\documentclass[a4paper,12pt]{amsart}
\usepackage{pstricks}

\newtheorem{thm}{Theorem}[section]
\newtheorem{prop}[thm]{Proposition}
\newtheorem{cor}[thm]{Corollary}
\theoremstyle{definition}
\newtheorem{defn}[thm]{Definition}
\newtheorem{exam}[thm]{Example}
\newtheorem{obs}[thm]{Observations}
\newtheorem{comment}[thm]{Comment}

\begin{document}

\title[On the Cyclic Deligne Conjecture]
{On the cyclic Deligne Conjecture}

\author[T.~Tradler]{Thomas~Tradler}
\address{Thomas Tradler, College of Technology of the City University
of New York, Department of Mathematics, 300 Jay Street, Brooklyn,
NY 11201, USA} \email{ttradler@citytech.cuny.edu}

\author[M.~Zeinalian]{Mahmoud~Zeinalian}
\address{Mahmoud Zeinalian, Department of Mathematics,C.W. Post Campus of Long Island University,
720 Northern Boulevard, Brookville, NY 11548, USA}
\email{mzeinalian@liu.edu}

\begin{abstract}
Let $A$ be a finite dimensional, unital, and associative algebra
which is endowed with a non-degenerate and invariant inner product.
We give an explicit description of an action of cyclic Sullivan chord
diagrams on the normalized Hochschild cochain complex of $A$. As a
corollary, the Hochschild cohomology of $A$ becomes a Frobenius
algebra which is endowed with a compatible $BV$ operator. If $A$ is
also commutative, then the discussion extends to an action of general
Sullivan chord diagrams. Some implications of this are discussed.
\end{abstract}

\maketitle

\section{Introduction}
This paper is concerned with algebraic structures on the
normalized Hochschild cochains, mirroring those of String
Topology. String Topology may be regarded as the study of the
algebraic topology of the free loop space of a manifold. Chas and
Sullivan, \cite{CS1}, \cite{CS2}, showed that the cohomology of
the free loop space of a manifold has the structure of a $BV$
algebra. Building on Sullivan's later work \cite{Su2}, Cohen and
Godin \cite{CG} showed that string topology operations give rise
to a two dimensional positive boundary TQFT. These were achieved
by looking at certain operations coming from what are known as
Sullivan chord diagrams.

The Hochschild cochain complex of the singular cochains on a
(simply connected) manifold gives a model for the chains on the
free loop space of that manifold. One expects analogs of the above
structures in a purely algebraic setting. Moreover, one is
interested in algebraic structures not only at the level of
Hochschild cohomology, but also, and more importantly so, at the
level of Hochschild cochains. The Deligne conjecture, which has
been proved in \cite{Ka1}, \cite{KS}, \cite{MS1}, \cite{Ta2},
\cite{Ta1}, and \cite{V}, partly addresses this issue. It states
that the chains on the little disc operad act on the Hochschild
cochain complex of an associative algebra. One relevant question
is whether the chains on the framed little disc operad, or
equivalently chains on cacti with marked points, act on the
Hochschild cochain complex of a unital and associative algebra
which has an invariant non-degenerate inner product. This question
has been affirmatively answered by Kaufmann and McClure and Smith;
see \cite{Ka2}, \cite{MS3}, and \cite{MS2}. The aim of the present
paper is to show that a much larger set of operations, with a
richer internal algebraic structure, act.

We give an explicit action of cyclic Sullivan chord diagrams (see
section \ref{sec1} for definition), which include the chains on
the cacti with marked points, on the normalized Hochschild cochain
complex. More precisely, we show the following:

\setcounter{section}{3} \setcounter{thm}{2}
\begin{thm}
Let $A$ be a finite dimensional, unital, and associative algebra
with a non-degenerate and invariant inner product. Then, the
normalized Hochschild cochain complex of $A$ is an algebra over
the PROP, $C_\ast \mathcal{S}^c$, of cyclic Sullivan chord
diagrams.
\end{thm}

\begin{cor}
Under the above assumptions, the Hochschild cohomology of $A$ is a
Frobenius algebra endowed with a compatible $BV$ operator.
\end{cor}

Cyclic Sullivan chord diagrams do not account for the operation
which reverses the orientation of a loop. As it turns out, the
concept of a Sullivan chord diagram is precisely the
generalization needed for labelling the orientation reversing
operations, in addition to the operations which are labelled by
the cyclic Sullivan chord diagrams. In section \ref{section-comm}
we show that this larger PROP still acts on the normalized
Hochschild complex, if the associative algebra $A$ happens to be
in addition commutative.

\setcounter{section}{4} \setcounter{thm}{2}
\begin{thm}
Let $A$ be a finite dimensional, unital, associative, and
commutative algebra endowed with a non-degenerate and invariant
inner product. Then, the normalized Hochschild cochain complex of
$A$ is an algebra over the PROP, $C_\ast \mathcal{S}$, of Sullivan
chord diagrams.
\end{thm}

\begin{cor}
Under the assumptions of Theorem \ref{theorem2}, the Hochschild
cohomology of $A$ is a Frobenius algebra, which is endowed with a
compatible $BV$ operator, $\Delta$, and an involution $\sim$. The
operator $\Delta$ maps each eigenspace of $\sim$ into the other, {\it
i.e.} $\Delta\left(HH^* (A;A)_\pm\right)\subset HH^*(A;A)_\mp$, where
$HH^*(A;A)_{\pm}$ are the $\pm 1$ eigenspaces of $\sim$. The map
$\sim$ is both an anti-algebra and an anti-coalgebra map. That is to
say $\widetilde{f\smile g}=\widetilde{ g}\smile \widetilde{f}$, and
$\vee_0\left(\widetilde{f}\right)= \sum_{(f)} \widetilde{f}''\otimes
\widetilde{f'}$, where $\vee_0(f)= \sum_{(f)} f'\otimes f''$.
\end{cor}

\setcounter{section}{1}

{\bf Acknowledgments.} We would like to thank James Stasheff and
Dennis Sullivan for their constructive comments.

\section{Cyclic Sullivan Chord Diagrams}\label{sec1}

In this section, we introduce a special kind of Sullivan chord
diagrams, called cyclic Sullivan chord diagrams. For them, one
defines a boundary operator, a composition, and an operation which
corresponds to a certain relabelling. These diagrams yield a PROP,
which is called the PROP of cyclic Sullivan chord diagrams. We
will discuss how the Frobenius PROP, as well as the $BV$ operad,
sit inside the homology of this PROP.

A \emph{cyclic Sullivan chord diagram} consists of a finite
collection of disjointly embedded planar circles which may be
connected using a finite number of immersed planar trees. Such a
tree is called a chord. An endpoint of a chord lies on a circle.
Different endpoints may lie on the same circle, and even on the
same point. However, there may exist circles to which no chords
are attached. The chords are not allowed to enter the circles. A
chord has two types of vertices, the inner vertices and the
endpoints, where it meets with the circles. The circles and the
chords together form a graph (with possibly a collection of
disjoint circles). At a vertex of this graph, there is a natural
cyclic ordering of the edges, which is induced by the orientation
of the plane. The cyclic ordering of the edges at each vertex
gives rise to a well defined thickening of the diagram to an
oriented surface with boundary. A diagram of type $(g; n, m)$ is
one for which this surface is of genus $g$, and has $n+m$ boundary
components, precisely $n$ of which are inside the original
circles. As part of the structure, these boundaries, which are
referred to as the inputs, are enumerated. Each input circle is
decorated with a marked point, called the input marked point, and
is oriented in a clockwise fashion. The remaining $m$ boundary
components, which are known as the outputs, are also enumerated
and decorated with output marked points. We reserve the term
\emph{special point} for collectively refereing to the input and
output marked points, as well as to the chord endpoints.

The thickened surface is merely an auxiliary tool for better
picturing the input and output circles. Mathematically, all that
matters is the combinatorial structure of the cyclic Sullivan
chord diagram. In fact, the input and output marked points are all
physically placed on the original circles. The chords are to be
thought of as objects of length zero. Consequently, an output
marked point or a chord endpoint which is located at an endpoint
of a chord, may slide from that endpoint to an adjacent one along
the perimeter of the output circle.

The output circles are oriented as follows. The induced
orientation of the surface from the plane induces an orientation
on its boundary components. This induced orientation, which
opposes the orientation of the input circles, should coincide with
the orientation of the output circles. Note that since the circles
are oriented and the chords do not enter the circles, at a vertex
on a circle, there is a natural linear ordering on the set of
chord endpoints union the set of output marked points at that
point. We would like to emphasize that in this linear ordering
output marked points may very well be positioned in between chord
endpoints, and are considered as part of the linear ordering.

The combinatorial dimension of a cyclic Sullivan chord diagram of
type $(g; n, m)$ is by definition the number of connected
components obtained by removing the special points (chord
endpoints, input and output marked points) from the input circles,
minus $n$. We only consider the cyclic Sullivan chord diagrams up
to abstract combinatorial isomorphism of graphs, sending input
circles to input circles, mapping chords to chords, respecting
special points, and cyclic orderings at the vertices.
Consequently, the orientations of the input and output circles
should match as well.

\begin{defn} [Cyclic Sullivan Chord Diagrams]\ Let $C_k \mathcal{S}^c(g; n, m)$ denote the vector space
generated by the cyclic Sullivan chord diagrams of type $(g; n,
m)$ and of combinatorial dimension $k$, up to isomorphism. Let
$C_k \mathcal{S}^c(n, m) = \bigoplus_{g=0}^\infty C_k
\mathcal{S}^c(g; n, m)$, and $C_\ast \mathcal{S}^c(n, m) =
\bigoplus_{k=0}^\infty C_k \mathcal{S}^c(n, m)$.
\end{defn}

\label{bdry}There is a natural boundary operator, $\partial$, on
$C_\ast \mathcal{S}^c(n, m)$. By linearity, it suffices to
describe $\partial$ on the basis elements. Consider a basis
element $s \in C_k \mathcal{S}^c(g; n, m)$. By removing the
special points (chord endpoints, the input, and the output marked
points) from the input circles of $s$, one obtains $k+n$ connected
components. Since these oriented circles are enumerated, there is
a natural numbering of these connected components from $1$ to
$n+k$. Let $\partial (s)\in C_{k-1} \mathcal{S}^c(g; n, m)$ denote
the alternating sum of all cyclic chord diagrams obtained by one
at a time collapsing of each of the connected components to a
point. Let us emphasize that in this paper, in defining the
boundary of a diagram, we collapse no chords or any segments
thereof. It is easy to verify that $\partial^2 =0$. See the
following example of the boundary of a diagram.

\ifx\setlinejoinmode\undefined
 \newcommand{\setlinejoinmode}[1]{}
\fi \ifx\setlinecaps\undefined
 \newcommand{\setlinecaps}[1]{}
\fi \ifx\setfont\undefined
 \newcommand{\setfont}[2]{}
\fi \[\pspicture(2.904294,-5.618330)(11.942900,-2.565420)
\scalebox{1 -1}{
\newrgbcolor{dialinecolor}{0 0 0}
\psset{linecolor=dialinecolor}
\newrgbcolor{diafillcolor}{1 1 1}
\psset{fillcolor=diafillcolor}
\newrgbcolor{dialinecolor}{0 0 0}
\psset{linecolor=dialinecolor}
\rput[r](3.746670,5.518330){\scalebox{1 -1}{}}
\newrgbcolor{dialinecolor}{1 1 1}
\psset{linecolor=dialinecolor}
\psellipse*(7.086580,2.994710)(0.387500,0.384290)
\psset{linewidth=0.040} \psset{linestyle=solid}
\psset{linestyle=solid}
\newrgbcolor{dialinecolor}{0 0 0}
\psset{linecolor=dialinecolor}
\psellipse(7.086580,2.994710)(0.387500,0.384290)
\newrgbcolor{dialinecolor}{1 1 1}
\psset{linecolor=dialinecolor}
\psellipse*(9.336580,2.969710)(0.387500,0.384290)
\psset{linewidth=0.040} \psset{linestyle=solid}
\psset{linestyle=solid}
\newrgbcolor{dialinecolor}{0 0 0}
\psset{linecolor=dialinecolor}
\psellipse(9.336580,2.969710)(0.387500,0.384290)

\newrgbcolor{dialinecolor}{0 0 0}
\psset{linecolor=dialinecolor}
\rput[r](6.089400,3.909300){\scalebox{1 -1}{=}}
\setfont{Helvetica}{0.8}
\newrgbcolor{dialinecolor}{0 0 0}
\psset{linecolor=dialinecolor} \rput[r](3.404294,3.7){\scalebox{1
-1}{$\partial$}} \psset{linewidth=0.040} \psset{linestyle=solid}
\psset{linestyle=solid}
\newrgbcolor{dialinecolor}{0 0 0}
\psset{linecolor=dialinecolor}
\psellipse(4.515420,4.488990)(0.381250,0.375000)
\newrgbcolor{dialinecolor}{1 1 1}
\psset{linecolor=dialinecolor}
\psellipse*(4.510330,3.042200)(0.387500,0.384290)
\psset{linewidth=0.040} \psset{linestyle=solid}
\psset{linestyle=solid}
\newrgbcolor{dialinecolor}{0 0 0}
\psset{linecolor=dialinecolor}
\psellipse(4.510330,3.042200)(0.387500,0.384290)
\psset{linewidth=0.050} \psset{linestyle=solid}
\psset{linestyle=solid} \setlinecaps{0}
\newrgbcolor{dialinecolor}{0 0 0}
\psset{linecolor=dialinecolor}
\psline(4.515330,3.240420)(4.510330,3.426490)
\psset{linewidth=0.050} \psset{linestyle=solid}
\psset{linestyle=solid} \setlinecaps{0}
\newrgbcolor{dialinecolor}{0 0 0}
\psset{linecolor=dialinecolor}
\psline(4.865330,4.315420)(5.027830,4.240420)
\psset{linewidth=0.050} \psset{linestyle=solid}
\psset{linestyle=solid} \setlinecaps{0}
\newrgbcolor{dialinecolor}{0 0 0}
\psset{linecolor=dialinecolor}
\psline(4.515420,4.863990)(4.515330,4.690420)
\psset{linewidth=0.010} \psset{linestyle=dashed,dash=1 1}
\psset{linestyle=dashed,dash=0.1 0.1} \setlinecaps{0}
\newrgbcolor{dialinecolor}{0 0 0}
\psset{linecolor=dialinecolor}
\psline(4.496670,4.201490)(4.510330,3.426490)
\setfont{Helvetica}{0.3}
\newrgbcolor{dialinecolor}{0 0 0}
\psset{linecolor=dialinecolor}
\rput[r](4.579170,4.478990){\scalebox{1 -1}{1}}
\setfont{Helvetica}{0.3}
\newrgbcolor{dialinecolor}{0 0 0}
\psset{linecolor=dialinecolor}
\rput[r](4.559170,3.026490){\scalebox{1 -1}{2}}
\setfont{Helvetica}{0.3}
\newrgbcolor{dialinecolor}{0 0 0}
\psset{linecolor=dialinecolor}
\rput[r](5.037830,4.115420){\scalebox{1 -1}{1'}}
\psset{linewidth=0.040} \psset{linestyle=solid}
\psset{linestyle=solid}
\newrgbcolor{dialinecolor}{0 0 0}
\psset{linecolor=dialinecolor}
\psellipse(7.091670,4.441500)(0.381250,0.375000)
\psset{linewidth=0.050} \psset{linestyle=solid}
\psset{linestyle=solid} \setlinecaps{0}
\newrgbcolor{dialinecolor}{0 0 0}
\psset{linecolor=dialinecolor}
\psline(7.091580,3.192940)(7.086580,3.379000)
\psset{linewidth=0.050} \psset{linestyle=solid}
\psset{linestyle=solid} \setlinecaps{0}
\newrgbcolor{dialinecolor}{0 0 0}
\psset{linecolor=dialinecolor}
\psline(7.429080,4.242940)(7.591580,4.167940)
\psset{linewidth=0.050} \psset{linestyle=solid}
\psset{linestyle=solid} \setlinecaps{0}
\newrgbcolor{dialinecolor}{0 0 0}
\psset{linecolor=dialinecolor}
\psline(7.091670,4.816500)(7.091580,4.642940)
\setfont{Helvetica}{0.3}
\newrgbcolor{dialinecolor}{0 0 0}
\psset{linecolor=dialinecolor}
\rput[r](7.155420,4.431500){\scalebox{1 -1}{1}}
\setfont{Helvetica}{0.3}
\newrgbcolor{dialinecolor}{0 0 0}
\psset{linecolor=dialinecolor}
\rput[r](7.135420,2.979000){\scalebox{1 -1}{2}}
\setfont{Helvetica}{0.3}
\newrgbcolor{dialinecolor}{0 0 0}
\psset{linecolor=dialinecolor}
\rput[r](7.530420,4.031500){\scalebox{1 -1}{1'}}
\psset{linewidth=0.040} \psset{linestyle=solid}
\psset{linestyle=solid}
\newrgbcolor{dialinecolor}{0 0 0}
\psset{linecolor=dialinecolor}
\psellipse(9.341670,4.416500)(0.381250,0.375000)
\psset{linewidth=0.050} \psset{linestyle=solid}
\psset{linestyle=solid} \setlinecaps{0}
\newrgbcolor{dialinecolor}{0 0 0}
\psset{linecolor=dialinecolor}
\psline(9.341580,3.167940)(9.336580,3.354000)
\psset{linewidth=0.050} \psset{linestyle=solid}
\psset{linestyle=solid} \setlinecaps{0}
\newrgbcolor{dialinecolor}{0 0 0}
\psset{linecolor=dialinecolor}
\psline(9.341670,4.041500)(9.465330,3.927920)
\psset{linewidth=0.050} \psset{linestyle=solid}
\psset{linestyle=solid} \setlinecaps{0}
\newrgbcolor{dialinecolor}{0 0 0}
\psset{linecolor=dialinecolor}
\psline(9.341670,4.791500)(9.341580,4.617940)
\psset{linewidth=0.010} \psset{linestyle=dashed,dash=1 1}
\psset{linestyle=dashed,dash=0.1 0.1} \setlinecaps{0}
\newrgbcolor{dialinecolor}{0 0 0}
\psset{linecolor=dialinecolor}
\psline(9.322920,4.129000)(9.336580,3.354000)
\setfont{Helvetica}{0.3}
\newrgbcolor{dialinecolor}{0 0 0}
\psset{linecolor=dialinecolor}
\rput[r](9.405420,4.406500){\scalebox{1 -1}{1}}
\setfont{Helvetica}{0.3}
\newrgbcolor{dialinecolor}{0 0 0}
\psset{linecolor=dialinecolor}
\rput[r](9.385420,2.954000){\scalebox{1 -1}{2}}
\setfont{Helvetica}{0.3}
\newrgbcolor{dialinecolor}{0 0 0}
\psset{linecolor=dialinecolor}
\rput[r](9.780420,4.006500){\scalebox{1 -1}{1'}}
\psset{linewidth=0.040} \psset{linestyle=solid}
\psset{linestyle=solid}
\newrgbcolor{dialinecolor}{0 0 0}
\psset{linecolor=dialinecolor}
\psellipse(11.504150,4.441500)(0.381250,0.375000)
\newrgbcolor{dialinecolor}{1 1 1}
\psset{linecolor=dialinecolor}
\psellipse*(11.499100,2.994710)(0.387500,0.384290)
\psset{linewidth=0.040} \psset{linestyle=solid}
\psset{linestyle=solid}
\newrgbcolor{dialinecolor}{0 0 0}
\psset{linecolor=dialinecolor}
\psellipse(11.499100,2.994710)(0.387500,0.384290)
\psset{linewidth=0.050} \psset{linestyle=solid}
\psset{linestyle=solid} \setlinecaps{0}
\newrgbcolor{dialinecolor}{0 0 0}
\psset{linecolor=dialinecolor}
\psline(11.504100,3.192940)(11.499100,3.379000)
\psset{linewidth=0.050} \psset{linestyle=solid}
\psset{linestyle=solid} \setlinecaps{0}
\newrgbcolor{dialinecolor}{0 0 0}
\psset{linecolor=dialinecolor}
\psline(11.502800,4.977920)(11.504100,4.642940)
\psset{linewidth=0.010} \psset{linestyle=dashed,dash=1 1}
\psset{linestyle=dashed,dash=0.1 0.1} \setlinecaps{0}
\newrgbcolor{dialinecolor}{0 0 0}
\psset{linecolor=dialinecolor}
\psline(11.485400,4.154000)(11.499100,3.379000)
\setfont{Helvetica}{0.3}
\newrgbcolor{dialinecolor}{0 0 0}
\psset{linecolor=dialinecolor}
\rput[r](11.567900,4.431500){\scalebox{1 -1}{1}}
\setfont{Helvetica}{0.3}
\newrgbcolor{dialinecolor}{0 0 0}
\psset{linecolor=dialinecolor}
\rput[r](11.547900,2.979000){\scalebox{1 -1}{2}}
\setfont{Helvetica}{0.3}
\newrgbcolor{dialinecolor}{0 0 0}
\psset{linecolor=dialinecolor}
\rput[r](11.942900,4.031500){\scalebox{1 -1}{1'}}

\newrgbcolor{dialinecolor}{0 0 0}
\psset{linecolor=dialinecolor}
\rput[r](8.241580,3.890420){\scalebox{1 -1}{$-$}}

\newrgbcolor{dialinecolor}{0 0 0}
\psset{linecolor=dialinecolor}
\rput[r](10.641600,3.902920){\scalebox{1 -1}{+}}
\psset{linewidth=0.010} \psset{linestyle=dashed,dash=0.1 0.1}
\psset{linestyle=dashed,dash=0.1 0.1} \setlinejoinmode{0}
\setlinecaps{0}
\newrgbcolor{dialinecolor}{0 0 0}
\psset{linecolor=dialinecolor} \pscustom{
\newpath
\moveto(7.086580,3.379000)
\curveto(6.540330,3.901120)(5.977830,5.201120)(7.091670,4.816500)
\stroke}
\newrgbcolor{dialinecolor}{0 0 0}
\psset{linecolor=dialinecolor}
\rput[r](4.862830,3.965420){\scalebox{1 -1}{}}
\newrgbcolor{dialinecolor}{0 0 0}
\psset{linecolor=dialinecolor}
\rput[r](5.412830,3.915420){\scalebox{1 -1}{}}
\psset{linewidth=0.050} \psset{linestyle=solid}
\psset{linestyle=solid} \setlinecaps{0}
\newrgbcolor{dialinecolor}{0 0 0}
\psset{linecolor=dialinecolor}
\psclip{\pswedge[linestyle=none,fillstyle=none](5.368945,3.752178){2.610188}{143.087664}{216.450350}}
\psellipse(5.368945,3.752178)(1.845682,1.845682)
\endpsclip
\psset{linewidth=0.050} \psset{linestyle=solid}
\psset{linestyle=solid} \setlinecaps{0}
\newrgbcolor{dialinecolor}{0 0 0}
\psset{linecolor=dialinecolor}
\psclip{\pswedge[linestyle=none,fillstyle=none](3.714017,3.777175){2.676227}{323.986426}{35.364731}}
\psellipse(3.714017,3.777175)(1.892378,1.892378)
\endpsclip
}\endpspicture \] Here the input circles are labelled by 1 and 2,
and the output circle is labelled by 1'. Throughout this paper, we
label the input circles using numbers 1, 2, 3, ..., and the output
circles with 1', 2', 3', .... For the inputs, these numbers are
written inside the input circles. In case of an output, these
numbers are written somewhere close to the diagram along the
perimeters of the appropriate output circles. In order to better
see the output circles one may, merely as a device, slightly
thicken the cyclic Sullivan chord diagram to obtain an auxiliary
orientable surface with boundary.

\label{comp}There is also a naturally defined composition. By
linearity, it suffices to define the composition $\circ : C_\ast
\mathcal{S}^c(k, l) \otimes C_\ast\mathcal{S}^c (m, k) \rightarrow
C_\ast \mathcal{S}^c (m, l)$ on the basis elements. For two such
elements $s \in C_\ast\mathcal{S}^c (k, l)$ and $s' \in C_\ast
\mathcal{S}^c(m, k)$, we want to define $s\circ s'$. For every
$1\leq i\leq k$ consider the $i^{th}$ output of $s'$ and the
$i^{th}$ input of $s$. Each of these is a circle with a certain
number of vertices and a particular marked point on it. Starting
from the marked point of the $i^{th}$ input circle of $s$, one can
read off the linear ordering of chords and output marked points
arranged around this input circle. After aligning the input of $s$
with that of $s'$, this linear ordering should be shuffled in
between the previously existing chords of the $i^{th}$ output
circle of $s'$ in all possible ways. The $i^{th}$ output circle of
$s'$ and the $i^{th}$ input circle of $s$ are now dissolved. With
respect to the total ordering of vertices on chord diagrams, for
each chord endpoint or marked point of $s$ which moves past an
output marked point or chord endpoint of $s'$, a sign factor of
$(-1)$ accrues. Summing over all possibilities gives rise to the
desired composition. See the following example:
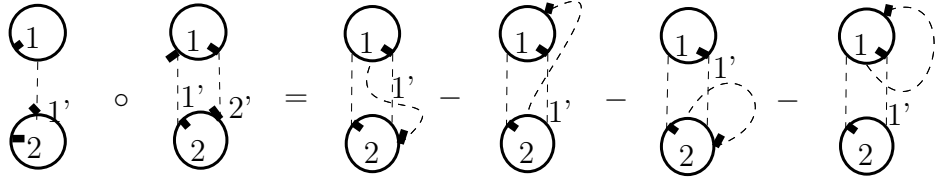
\ifx\setlinejoinmode\undefined
 \newcommand{\setlinejoinmode}[1]{}
\fi \ifx\setlinecaps\undefined
 \newcommand{\setlinecaps}[1]{}
\fi \ifx\setfont\undefined
 \newcommand{\setfont}[2]{}
\fi \[\pspicture(2.002830,-5.618330)(14.331767,-2.515424)
\scalebox{1 -1}{
\newrgbcolor{dialinecolor}{0 0 0}
\psset{linecolor=dialinecolor}
\newrgbcolor{diafillcolor}{1 1 1}
\psset{fillcolor=diafillcolor}
\newrgbcolor{dialinecolor}{0 0 0}
\psset{linecolor=dialinecolor}
\rput[r](3.746670,5.518330){\scalebox{1 -1}{}}

\newrgbcolor{dialinecolor}{0 0 0}
\psset{linecolor=dialinecolor}
\rput[r](6.008326,3.865079){\scalebox{1 -1}{=}}

\newrgbcolor{dialinecolor}{0 0 0}
\psset{linecolor=dialinecolor}
\rput[r](3.625330,3.882910){\scalebox{1 -1}{$\circ$}}
\psset{linewidth=0.040} \psset{linestyle=solid}
\psset{linestyle=solid}
\newrgbcolor{dialinecolor}{0 0 0}
\psset{linecolor=dialinecolor}
\psellipse(2.415420,4.413990)(0.381250,0.375000)
\newrgbcolor{dialinecolor}{1 1 1}
\psset{linecolor=dialinecolor}
\psellipse*(2.410330,2.967200)(0.387500,0.384290)
\psset{linewidth=0.040} \psset{linestyle=solid}
\psset{linestyle=solid}
\newrgbcolor{dialinecolor}{0 0 0}
\psset{linecolor=dialinecolor}
\psellipse(2.410330,2.967200)(0.387500,0.384290)
\psset{linewidth=0.1} \psset{linestyle=solid}
\psset{linestyle=solid} \setlinecaps{0}
\newrgbcolor{dialinecolor}{0 0 0}
\psset{linecolor=dialinecolor}
\psline(2.196670,3.101490)(2.096670,3.176490)
\psset{linewidth=0.1} \psset{linestyle=solid}
\psset{linestyle=solid} \setlinecaps{0}
\newrgbcolor{dialinecolor}{0 0 0}
\psset{linecolor=dialinecolor}
\psline(2.415420,4.038990)(2.321670,3.938990)
\psset{linewidth=0.1} \psset{linestyle=solid}
\psset{linestyle=solid} \setlinecaps{0}
\newrgbcolor{dialinecolor}{0 0 0}
\psset{linecolor=dialinecolor}
\psline(2.071670,4.376490)(2.234170,4.376490)
\psset{linewidth=0.010} \psset{linestyle=dashed,dash=1 1}
\psset{linestyle=dashed,dash=0.1 0.1} \setlinecaps{0}
\newrgbcolor{dialinecolor}{0 0 0}
\psset{linecolor=dialinecolor}
\psline(2.396670,4.126490)(2.410330,3.351490)
\setfont{Helvetica}{0.3}
\newrgbcolor{dialinecolor}{0 0 0}
\psset{linecolor=dialinecolor}
\rput[r](2.466670,4.503990){\scalebox{1 -1}{2}}
\setfont{Helvetica}{0.3}
\newrgbcolor{dialinecolor}{0 0 0}
\psset{linecolor=dialinecolor}
\rput[r](2.446670,3.038990){\scalebox{1 -1}{1}}
\setfont{Helvetica}{0.3}
\newrgbcolor{dialinecolor}{0 0 0}
\psset{linecolor=dialinecolor}
\rput[r](2.854170,4.003990){\scalebox{1 -1}{1'}}
\newrgbcolor{dialinecolor}{1 1 1}
\psset{linecolor=dialinecolor}
\psellipse*(4.576580,4.396660)(0.368750,0.393750)
\psset{linewidth=0.040} \psset{linestyle=solid}
\psset{linestyle=solid}
\newrgbcolor{dialinecolor}{0 0 0}
\psset{linecolor=dialinecolor}
\psellipse(4.576580,4.396660)(0.368750,0.393750)
\newrgbcolor{dialinecolor}{1 1 1}
\psset{linecolor=dialinecolor}
\psellipse*(4.545330,2.959160)(0.387500,0.381250)
\psset{linewidth=0.040} \psset{linestyle=solid}
\psset{linestyle=solid}
\newrgbcolor{dialinecolor}{0 0 0}
\psset{linecolor=dialinecolor}
\psellipse(4.545330,2.959160)(0.387500,0.381250)
\psset{linewidth=0.010} \psset{linestyle=dashed,dash=1 1}
\psset{linestyle=dashed,dash=0.1 0.1} \setlinecaps{0}
\newrgbcolor{dialinecolor}{0 0 0}
\psset{linecolor=dialinecolor}
\psline(4.271670,4.213990)(4.271326,3.228744)
\psset{linewidth=0.1} \psset{linestyle=solid}
\psset{linestyle=solid} \setlinecaps{0}
\newrgbcolor{dialinecolor}{0 0 0}
\psset{linecolor=dialinecolor}
\psline(4.819334,3.228744)(4.684170,3.138990)
\psset{linewidth=0.1} \psset{linestyle=solid}
\psset{linestyle=solid} \setlinecaps{0}
\newrgbcolor{dialinecolor}{0 0 0}
\psset{linecolor=dialinecolor}
\psline(4.134170,3.326490)(4.271326,3.228744)
\psset{linewidth=0.010} \psset{linestyle=dashed,dash=1 1}
\psset{linestyle=dashed,dash=0.1 0.1} \setlinecaps{0}
\newrgbcolor{dialinecolor}{0 0 0}
\psset{linecolor=dialinecolor}
\psline(4.837326,4.118237)(4.819334,3.228744)
\psset{linewidth=0.1} \psset{linestyle=solid}
\psset{linestyle=solid} \setlinecaps{0}
\newrgbcolor{dialinecolor}{0 0 0}
\psset{linecolor=dialinecolor}
\psline(4.837326,4.118237)(4.721670,3.963990)
\psset{linewidth=0.1} \psset{linestyle=solid}
\psset{linestyle=solid} \setlinecaps{0}
\newrgbcolor{dialinecolor}{0 0 0}
\psset{linecolor=dialinecolor}
\psline(4.409170,4.226490)(4.315834,4.118237)
\setfont{Helvetica}{0.3}
\newrgbcolor{dialinecolor}{0 0 0}
\psset{linecolor=dialinecolor}
\rput[r](4.641670,4.541490){\scalebox{1 -1}{2}}
\setfont{Helvetica}{0.3}
\newrgbcolor{dialinecolor}{0 0 0}
\psset{linecolor=dialinecolor}
\rput[r](4.579170,3.078990){\scalebox{1 -1}{1}}
\setfont{Helvetica}{0.3}
\newrgbcolor{dialinecolor}{0 0 0}
\psset{linecolor=dialinecolor}
\rput[r](4.629170,3.803990){\scalebox{1 -1}{1'}}
\setfont{Helvetica}{0.3}
\newrgbcolor{dialinecolor}{0 0 0}
\psset{linecolor=dialinecolor}
\rput[r](5.266670,3.953990){\scalebox{1 -1}{2'}}

\newrgbcolor{dialinecolor}{0 0 0}
\psset{linecolor=dialinecolor}
\rput[r](12.543935,3.883980){\scalebox{1 -1}{$-$}}

\newrgbcolor{dialinecolor}{0 0 0}
\psset{linecolor=dialinecolor}
\rput[r](8.059179,3.836517){\scalebox{1 -1}{$-$}}
\psset{linewidth=0.010} \psset{linestyle=dashed,dash=1 1}
\psset{linestyle=dashed,dash=0.1 0.1} \setlinejoinmode{0}
\setlinecaps{0}
\newrgbcolor{dialinecolor}{0 0 0}
\psset{linecolor=dialinecolor} \pscustom{
\newpath
\moveto(6.824443,3.412721)
\curveto(6.499443,4.250221)(8.186943,3.612721)(7.228543,4.442541)
\stroke}
\newrgbcolor{dialinecolor}{1 1 1}
\psset{linecolor=dialinecolor}
\psellipse*(6.859793,4.442541)(0.368750,0.393750)
\psset{linewidth=0.040} \psset{linestyle=solid}
\psset{linestyle=solid}
\newrgbcolor{dialinecolor}{0 0 0}
\psset{linecolor=dialinecolor}
\psellipse(6.859793,4.442541)(0.368750,0.393750)
\newrgbcolor{dialinecolor}{1 1 1}
\psset{linecolor=dialinecolor}
\psellipse*(6.859018,2.980041)(0.382075,0.381250)
\psset{linewidth=0.040} \psset{linestyle=solid}
\psset{linestyle=solid}
\newrgbcolor{dialinecolor}{0 0 0}
\psset{linecolor=dialinecolor}
\psellipse(6.859018,2.980041)(0.382075,0.381250)
\psset{linewidth=0.010} \psset{linestyle=dashed,dash=1 1}
\psset{linestyle=dashed,dash=0.1 0.1} \setlinecaps{0}
\newrgbcolor{dialinecolor}{0 0 0}
\psset{linecolor=dialinecolor}
\psline(6.579943,4.234861)(6.588851,3.249626)
\psset{linewidth=0.1} \psset{linestyle=solid}
\psset{linestyle=solid} \setlinecaps{0}
\newrgbcolor{dialinecolor}{0 0 0}
\psset{linecolor=dialinecolor}
\psline(7.129186,3.249626)(6.992443,3.159861)
\psset{linewidth=0.010} \psset{linestyle=dashed,dash=1 1}
\psset{linestyle=dashed,dash=0.1 0.1} \setlinecaps{0}
\newrgbcolor{dialinecolor}{0 0 0}
\psset{linecolor=dialinecolor}
\psline(7.120539,4.164118)(7.129186,3.249626)
\psset{linewidth=0.1} \psset{linestyle=solid}
\psset{linestyle=solid} \setlinecaps{0}
\newrgbcolor{dialinecolor}{0 0 0}
\psset{linecolor=dialinecolor}
\psline(7.228543,4.442541)(7.286943,4.262721)
\psset{linewidth=0.1} \psset{linestyle=solid}
\psset{linestyle=solid} \setlinecaps{0}
\newrgbcolor{dialinecolor}{0 0 0}
\psset{linecolor=dialinecolor}
\psline(6.717443,4.247361)(6.599048,4.164118)
\setfont{Helvetica}{0.3}
\newrgbcolor{dialinecolor}{0 0 0}
\psset{linecolor=dialinecolor}
\rput[r](6.949943,4.562361){\scalebox{1 -1}{2}}
\setfont{Helvetica}{0.3}
\newrgbcolor{dialinecolor}{0 0 0}
\psset{linecolor=dialinecolor}
\rput[r](6.887443,3.099861){\scalebox{1 -1}{1}}
\setfont{Helvetica}{0.3}
\newrgbcolor{dialinecolor}{0 0 0}
\psset{linecolor=dialinecolor}
\rput[r](7.424943,3.687361){\scalebox{1 -1}{1'}}
\psset{linewidth=0.010} \psset{linestyle=dashed,dash=1 1}
\psset{linestyle=dashed,dash=0.1 0.1} \setlinejoinmode{0}
\setlinecaps{0}
\newrgbcolor{dialinecolor}{0 0 0}
\psset{linecolor=dialinecolor} \pscustom{
\newpath
\moveto(13.430300,3.394159)
\curveto(14.270800,4.608089)(14.801846,2.293397)(13.716846,2.669457)
\stroke}
\newrgbcolor{dialinecolor}{1 1 1}
\psset{linecolor=dialinecolor}
\psellipse*(13.431075,4.475409)(0.368750,0.393750)
\psset{linewidth=0.040} \psset{linestyle=solid}
\psset{linestyle=solid}
\newrgbcolor{dialinecolor}{0 0 0}
\psset{linecolor=dialinecolor}
\psellipse(13.431075,4.475409)(0.368750,0.393750)
\newrgbcolor{dialinecolor}{1 1 1}
\psset{linecolor=dialinecolor}
\psellipse*(13.430300,3.012909)(0.382075,0.381250)
\psset{linewidth=0.040} \psset{linestyle=solid}
\psset{linestyle=solid}
\newrgbcolor{dialinecolor}{0 0 0}
\psset{linecolor=dialinecolor}
\psellipse(13.430300,3.012909)(0.382075,0.381250)
\psset{linewidth=0.010} \psset{linestyle=dashed,dash=1 1}
\psset{linestyle=dashed,dash=0.1 0.1} \setlinecaps{0}
\newrgbcolor{dialinecolor}{0 0 0}
\psset{linecolor=dialinecolor}
\psline(13.151225,4.267729)(13.160132,3.282494)
\psset{linewidth=0.1} \psset{linestyle=solid}
\psset{linestyle=solid} \setlinecaps{0}
\newrgbcolor{dialinecolor}{0 0 0}
\psset{linecolor=dialinecolor}
\psline(13.700468,3.282494)(13.563725,3.192729)
\psset{linewidth=0.010} \psset{linestyle=dashed,dash=1 1}
\psset{linestyle=dashed,dash=0.1 0.1} \setlinecaps{0}
\newrgbcolor{dialinecolor}{0 0 0}
\psset{linecolor=dialinecolor}
\psline(13.691821,4.196986)(13.700468,3.282494)
\psset{linewidth=0.1} \psset{linestyle=solid}
\psset{linestyle=solid} \setlinecaps{0}
\newrgbcolor{dialinecolor}{0 0 0}
\psset{linecolor=dialinecolor}
\psline(13.733225,2.595589)(13.700468,2.743325)
\psset{linewidth=0.1} \psset{linestyle=solid}
\psset{linestyle=solid} \setlinecaps{0}
\newrgbcolor{dialinecolor}{0 0 0}
\psset{linecolor=dialinecolor}
\psline(13.288725,4.280229)(13.170329,4.196986)
\setfont{Helvetica}{0.3}
\newrgbcolor{dialinecolor}{0 0 0}
\psset{linecolor=dialinecolor}
\rput[r](13.521225,4.595229){\scalebox{1 -1}{2}}
\setfont{Helvetica}{0.3}
\newrgbcolor{dialinecolor}{0 0 0}
\psset{linecolor=dialinecolor}
\rput[r](13.458725,3.132729){\scalebox{1 -1}{1}}
\setfont{Helvetica}{0.3}
\newrgbcolor{dialinecolor}{0 0 0}
\psset{linecolor=dialinecolor}
\rput[r](14.021225,4.045229){\scalebox{1 -1}{1'}}

\newrgbcolor{dialinecolor}{0 0 0}
\psset{linecolor=dialinecolor}
\rput[r](10.262610,3.865079){\scalebox{1 -1}{$-$}}
\psset{linewidth=0.010} \psset{linestyle=dashed,dash=1 1}
\psset{linestyle=dashed,dash=0.1 0.1} \setlinejoinmode{0}
\setlinecaps{0}
\newrgbcolor{dialinecolor}{0 0 0}
\psset{linecolor=dialinecolor} \pscustom{
\newpath
\moveto(8.938857,3.988990)
\curveto(8.988857,3.688990)(10.263857,2.138990)(9.238857,2.651490)
\stroke}
\newrgbcolor{dialinecolor}{1 1 1}
\psset{linecolor=dialinecolor}
\psellipse*(8.921757,4.442740)(0.368750,0.393750)
\psset{linewidth=0.040} \psset{linestyle=solid}
\psset{linestyle=solid}
\newrgbcolor{dialinecolor}{0 0 0}
\psset{linecolor=dialinecolor}
\psellipse(8.921757,4.442740)(0.368750,0.393750)
\newrgbcolor{dialinecolor}{1 1 1}
\psset{linecolor=dialinecolor}
\psellipse*(8.920932,2.980240)(0.382075,0.381250)
\psset{linewidth=0.040} \psset{linestyle=solid}
\psset{linestyle=solid}
\newrgbcolor{dialinecolor}{0 0 0}
\psset{linecolor=dialinecolor}
\psellipse(8.920932,2.980240)(0.382075,0.381250)
\psset{linewidth=0.010} \psset{linestyle=dashed,dash=1 1}
\psset{linestyle=dashed,dash=0.1 0.1} \setlinecaps{0}
\newrgbcolor{dialinecolor}{0 0 0}
\psset{linecolor=dialinecolor}
\psline(8.641847,4.235070)(8.650764,3.249824)
\psset{linewidth=0.1} \psset{linestyle=solid}
\psset{linestyle=solid} \setlinecaps{0}
\newrgbcolor{dialinecolor}{0 0 0}
\psset{linecolor=dialinecolor}
\psline(9.191100,3.249824)(9.054347,3.160070)
\psset{linewidth=0.010} \psset{linestyle=dashed,dash=1 1}
\psset{linestyle=dashed,dash=0.1 0.1} \setlinecaps{0}
\newrgbcolor{dialinecolor}{0 0 0}
\psset{linecolor=dialinecolor}
\psline(9.182503,4.164317)(9.191100,3.249824)
\psset{linewidth=0.1} \psset{linestyle=solid}
\psset{linestyle=solid} \setlinecaps{0}
\newrgbcolor{dialinecolor}{0 0 0}
\psset{linecolor=dialinecolor}
\psline(9.226357,2.576490)(9.191100,2.710656)
\psset{linewidth=0.1} \psset{linestyle=solid}
\psset{linestyle=solid} \setlinecaps{0}
\newrgbcolor{dialinecolor}{0 0 0}
\psset{linecolor=dialinecolor}
\psline(8.779347,4.247570)(8.661011,4.164317)
\setfont{Helvetica}{0.3}
\newrgbcolor{dialinecolor}{0 0 0}
\psset{linecolor=dialinecolor}
\rput[r](9.011847,4.562570){\scalebox{1 -1}{2}}
\setfont{Helvetica}{0.3}
\newrgbcolor{dialinecolor}{0 0 0}
\psset{linecolor=dialinecolor}
\rput[r](8.949347,3.100070){\scalebox{1 -1}{1}}
\setfont{Helvetica}{0.3}
\newrgbcolor{dialinecolor}{0 0 0}
\psset{linecolor=dialinecolor}
\rput[r](9.511847,4.012570){\scalebox{1 -1}{1'}}
\psset{linewidth=0.010} \psset{linestyle=dashed,dash=1 1}
\psset{linestyle=dashed,dash=0.1 0.1} \setlinejoinmode{0}
\setlinecaps{0}
\newrgbcolor{dialinecolor}{0 0 0}
\psset{linecolor=dialinecolor} \pscustom{
\newpath
\moveto(11.067241,4.027688)
\curveto(11.639741,3.089118)(12.531371,4.093938)(11.418891,4.481438)
\stroke}
\newrgbcolor{dialinecolor}{1 1 1}
\psset{linecolor=dialinecolor}
\psellipse*(11.050141,4.481438)(0.368750,0.393750)
\psset{linewidth=0.040} \psset{linestyle=solid}
\psset{linestyle=solid}
\newrgbcolor{dialinecolor}{0 0 0}
\psset{linecolor=dialinecolor}
\psellipse(11.050141,4.481438)(0.368750,0.393750)
\newrgbcolor{dialinecolor}{1 1 1}
\psset{linecolor=dialinecolor}
\psellipse*(11.061816,3.006438)(0.382075,0.381250)
\psset{linewidth=0.040} \psset{linestyle=solid}
\psset{linestyle=solid}
\newrgbcolor{dialinecolor}{0 0 0}
\psset{linecolor=dialinecolor}
\psellipse(11.061816,3.006438)(0.382075,0.381250)
\psset{linewidth=0.010} \psset{linestyle=dashed,dash=1 1}
\psset{linestyle=dashed,dash=0.1 0.1} \setlinecaps{0}
\newrgbcolor{dialinecolor}{0 0 0}
\psset{linecolor=dialinecolor}
\psline(10.770231,4.273758)(10.791648,3.276023)
\psset{linewidth=0.1} \psset{linestyle=solid}
\psset{linestyle=solid} \setlinecaps{0}
\newrgbcolor{dialinecolor}{0 0 0}
\psset{linecolor=dialinecolor}
\psline(11.331984,3.276023)(11.182731,3.198758)
\psset{linewidth=0.010} \psset{linestyle=dashed,dash=1 1}
\psset{linestyle=dashed,dash=0.1 0.1} \setlinecaps{0}
\newrgbcolor{dialinecolor}{0 0 0}
\psset{linecolor=dialinecolor}
\psline(11.310887,4.203015)(11.331984,3.276023)
\psset{linewidth=0.1} \psset{linestyle=solid}
\psset{linestyle=solid} \setlinecaps{0}
\newrgbcolor{dialinecolor}{0 0 0}
\psset{linecolor=dialinecolor}
\psline(11.502241,4.326618)(11.418891,4.481438)
\psset{linewidth=0.1} \psset{linestyle=solid}
\psset{linestyle=solid} \setlinecaps{0}
\newrgbcolor{dialinecolor}{0 0 0}
\psset{linecolor=dialinecolor}
\psline(10.907731,4.286258)(10.789395,4.203015)
\setfont{Helvetica}{0.3}
\newrgbcolor{dialinecolor}{0 0 0}
\psset{linecolor=dialinecolor}
\rput[r](11.140231,4.601258){\scalebox{1 -1}{2}}
\setfont{Helvetica}{0.3}
\newrgbcolor{dialinecolor}{0 0 0}
\psset{linecolor=dialinecolor}
\rput[r](11.077731,3.138758){\scalebox{1 -1}{1}}
\setfont{Helvetica}{0.3}
\newrgbcolor{dialinecolor}{0 0 0}
\psset{linecolor=dialinecolor}
\rput[r](11.702731,3.426258){\scalebox{1 -1}{1'}} }\endpspicture\]
The input circles are labelled by 1 and 2, and the output circles by
1' and 2'.

In the following composition example we consider a situation in which
two output points and three chord endpoints of $s$ coincide with one
of its input points. In the process of identifying the corresponding
circles, the marked point of the input circle of $s$, labelled by 1,
is to be identified with the marked point of the output circle of
$s'$, labelled by 1'. At this point, it is important to keep track of
the combinatorics of the chords and output points of $s$, together
with of those of $s'$. The chords and the output points of $s$ should
be inserted all together at one place in between those of $s'$,
respecting the linear ordering.
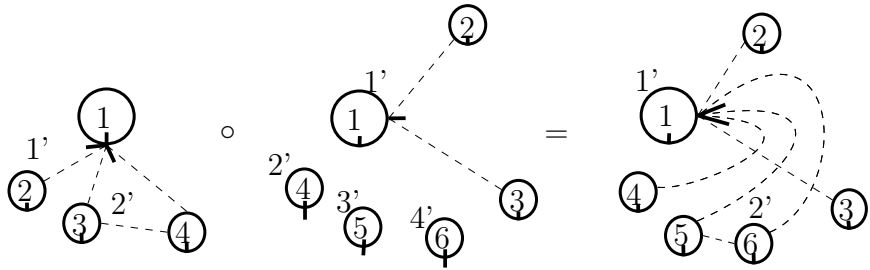
\ifx\setlinejoinmode\undefined
 \newcommand{\setlinejoinmode}[1]{}
\fi \ifx\setlinecaps\undefined
 \newcommand{\setlinecaps}[1]{}
\fi
\ifx\setfont\undefined
 \newcommand{\setfont}[2]{}
\fi \[\pspicture(2.761580,-5.618330)(14.239280,-1.529350)
\scalebox{1 -1}{
\newrgbcolor{dialinecolor}{0 0 0}
\psset{linecolor=dialinecolor}
\newrgbcolor{diafillcolor}{1 1 1}
\psset{fillcolor=diafillcolor}
\newrgbcolor{dialinecolor}{0 0 0}
\psset{linecolor=dialinecolor}
\rput[r](3.746670,5.518330){\scalebox{1 -1}{}}

\newrgbcolor{dialinecolor}{0 0 0}
\psset{linecolor=dialinecolor}
\rput[r](10.237800,3.307910){\scalebox{1 -1}{=}}

\newrgbcolor{dialinecolor}{0 0 0}
\psset{linecolor=dialinecolor}
\rput[r](5.812830,3.270410){\scalebox{1 -1}{$\circ$}}
\psset{linewidth=0.010} \psset{linestyle=dashed,dash=1 1}
\psset{linestyle=dashed,dash=0.1 0.1} \setlinejoinmode{0}
\setlinecaps{0}
\newrgbcolor{dialinecolor}{0 0 0}
\psset{linecolor=dialinecolor} \pscustom{
\newpath
\moveto(11.960300,3.031460)
\curveto(14.086600,2.673100)(13.251900,3.871510)(11.951900,4.434010)
\stroke}
\newrgbcolor{dialinecolor}{1 1 1}
\psset{linecolor=dialinecolor}
\psellipse*(7.460330,3.068960)(0.387500,0.384290)
\psset{linewidth=0.040} \psset{linestyle=solid}
\psset{linestyle=solid}
\newrgbcolor{dialinecolor}{0 0 0}
\psset{linecolor=dialinecolor}
\psellipse(7.460330,3.068960)(0.387500,0.384290)
\setfont{Helvetica}{0.3}
\newrgbcolor{dialinecolor}{0 0 0}
\psset{linecolor=dialinecolor}
\rput[r](7.496670,3.090750){\scalebox{1 -1}{1}}
\psset{linewidth=0.050} \psset{linestyle=solid}
\psset{linestyle=solid} \setlinecaps{0}
\newrgbcolor{dialinecolor}{0 0 0}
\psset{linecolor=dialinecolor}
\psline(7.465330,3.304680)(7.460330,3.453250)
\psset{linewidth=0.050} \psset{linestyle=solid}
\psset{linestyle=solid} \setlinecaps{0}
\newrgbcolor{dialinecolor}{0 0 0}
\psset{linecolor=dialinecolor}
\psline(7.847830,3.068960)(8.074080,3.066850)
\psset{linewidth=0.010} \psset{linestyle=dashed,dash=1 1}
\psset{linestyle=dashed,dash=0.1 0.1} \setlinecaps{0}
\newrgbcolor{dialinecolor}{0 0 0}
\psset{linecolor=dialinecolor}
\psline(7.847830,3.068960)(8.716290,2.020020)
\psset{linewidth=0.010} \psset{linestyle=dashed,dash=0.1 0.1}
\psset{linestyle=dashed,dash=0.1 0.1} \setlinecaps{0}
\newrgbcolor{dialinecolor}{0 0 0}
\psset{linecolor=dialinecolor}
\psline(9.349080,4.016850)(7.847830,3.068960)
\newrgbcolor{dialinecolor}{1 1 1}
\psset{linecolor=dialinecolor}
\psellipse*(11.572800,3.031460)(0.387500,0.384290)
\psset{linewidth=0.040} \psset{linestyle=solid}
\psset{linestyle=solid}
\newrgbcolor{dialinecolor}{0 0 0}
\psset{linecolor=dialinecolor}
\psellipse(11.572800,3.031460)(0.387500,0.384290)
\setfont{Helvetica}{0.3}
\newrgbcolor{dialinecolor}{0 0 0}
\psset{linecolor=dialinecolor}
\rput[r](11.634200,3.040750){\scalebox{1 -1}{1}}
\psset{linewidth=0.010} \psset{linestyle=dashed,dash=1 1}
\psset{linestyle=dashed,dash=0.1 0.1} \setlinecaps{0}
\newrgbcolor{dialinecolor}{0 0 0}
\psset{linecolor=dialinecolor}
\psline(12.027900,4.628960)(12.444100,4.735420)
\psset{linewidth=0.050} \psset{linestyle=solid}
\psset{linestyle=solid} \setlinecaps{0}
\newrgbcolor{dialinecolor}{0 0 0}
\psset{linecolor=dialinecolor}
\psline(11.960300,3.031460)(12.374100,3.148100)
\psset{linewidth=0.050} \psset{linestyle=solid}
\psset{linestyle=solid} \setlinecaps{0}
\newrgbcolor{dialinecolor}{0 0 0}
\psset{linecolor=dialinecolor}
\psline(11.577800,3.267180)(11.572800,3.415750)
\psset{linewidth=0.050} \psset{linestyle=solid}
\psset{linestyle=solid} \setlinecaps{0}
\newrgbcolor{dialinecolor}{0 0 0}
\psset{linecolor=dialinecolor}
\psline(11.960300,3.031460)(12.324100,2.885600)
\psset{linewidth=0.010} \psset{linestyle=dashed,dash=1 1}
\psset{linestyle=dashed,dash=0.1 0.1} \setlinecaps{0}
\newrgbcolor{dialinecolor}{0 0 0}
\psset{linecolor=dialinecolor}
\psline(11.960300,3.031460)(12.586600,2.048100)
\psset{linewidth=0.010} \psset{linestyle=dashed,dash=0.1 0.1}
\psset{linestyle=dashed,dash=0.1 0.1} \setlinecaps{0}
\newrgbcolor{dialinecolor}{0 0 0}
\psset{linecolor=dialinecolor}
\psline(13.749100,4.154350)(11.960300,3.031460)
\psset{linewidth=0.010} \psset{linestyle=dashed,dash=0.1 0.1}
\psset{linestyle=dashed,dash=0.1 0.1} \setlinejoinmode{0}
\setlinecaps{0}
\newrgbcolor{dialinecolor}{0 0 0}
\psset{linecolor=dialinecolor} \pscustom{
\newpath
\moveto(11.999100,3.004350)
\curveto(13.974100,1.348100)(14.014400,4.129350)(12.911600,4.585600)
\stroke} \psset{linewidth=0.010} \psset{linestyle=dashed,dash=0.1
0.1} \psset{linestyle=dashed,dash=0.1 0.1} \setlinejoinmode{0}
\setlinecaps{0}
\newrgbcolor{dialinecolor}{0 0 0}
\psset{linecolor=dialinecolor} \pscustom{
\newpath
\moveto(11.960300,3.031460)
\curveto(13.961600,3.123100)(11.997900,4.048220)(11.424100,3.985600)
\stroke} \psset{linewidth=0.040} \psset{linestyle=solid}
\psset{linestyle=solid}
\newrgbcolor{dialinecolor}{0 0 0}
\psset{linecolor=dialinecolor}
\psellipse(12.703430,4.735425)(0.259330,0.275715)
\psset{linewidth=0.050} \psset{linestyle=solid}
\psset{linestyle=solid} \setlinecaps{0}
\newrgbcolor{dialinecolor}{0 0 0}
\psset{linecolor=dialinecolor}
\psline(12.700200,4.861140)(12.703400,5.011140)
\setfont{Helvetica}{0.3}
\newrgbcolor{dialinecolor}{0 0 0}
\psset{linecolor=dialinecolor}
\rput[r](12.761800,4.742890){\scalebox{1 -1}{6}}
\newrgbcolor{dialinecolor}{1 1 1}
\psset{linecolor=dialinecolor}
\psellipse*(4.097830,3.042200)(0.387500,0.384290)
\psset{linewidth=0.040} \psset{linestyle=solid}
\psset{linestyle=solid}
\newrgbcolor{dialinecolor}{0 0 0}
\psset{linecolor=dialinecolor}
\psellipse(4.097830,3.042200)(0.387500,0.384290)
\psset{linewidth=0.050} \psset{linestyle=solid}
\psset{linestyle=solid} \setlinecaps{0}
\newrgbcolor{dialinecolor}{0 0 0}
\psset{linecolor=dialinecolor}
\psline(4.097830,3.426490)(3.824080,3.454350)
\setfont{Helvetica}{0.3}
\newrgbcolor{dialinecolor}{0 0 0}
\psset{linecolor=dialinecolor}
\rput[r](4.159170,3.063990){\scalebox{1 -1}{1}}
\psset{linewidth=0.010} \psset{linestyle=dashed,dash=1 1}
\psset{linestyle=dashed,dash=0.1 0.1} \setlinecaps{0}
\newrgbcolor{dialinecolor}{0 0 0}
\psset{linecolor=dialinecolor}
\psline(4.027830,4.464700)(4.906580,4.583660)
\psset{linewidth=0.050} \psset{linestyle=solid}
\psset{linestyle=solid} \setlinecaps{0}
\newrgbcolor{dialinecolor}{0 0 0}
\psset{linecolor=dialinecolor}
\psline(4.097830,3.426490)(4.202830,3.652170)
\psset{linewidth=0.050} \psset{linestyle=solid}
\psset{linestyle=solid} \setlinecaps{0}
\newrgbcolor{dialinecolor}{0 0 0}
\psset{linecolor=dialinecolor}
\psline(4.099080,3.254350)(4.097830,3.426490)
\psset{linewidth=0.010} \psset{linestyle=dashed,dash=1 1}
\psset{linestyle=dashed,dash=0.1 0.1} \setlinecaps{0}
\newrgbcolor{dialinecolor}{0 0 0}
\psset{linecolor=dialinecolor}
\psline(3.852200,4.198100)(4.097830,3.426490)
\psset{linewidth=0.010} \psset{linestyle=dashed,dash=0.1 0.1}
\psset{linestyle=dashed,dash=0.1 0.1} \setlinecaps{0}
\newrgbcolor{dialinecolor}{0 0 0}
\psset{linecolor=dialinecolor}
\psline(3.252200,3.910600)(4.097830,3.426490)
\psset{linewidth=0.010} \psset{linestyle=dashed,dash=0.1 0.1}
\psset{linestyle=dashed,dash=0.1 0.1} \setlinecaps{0}
\newrgbcolor{dialinecolor}{0 0 0}
\psset{linecolor=dialinecolor}
\psline(5.165910,4.307950)(4.097830,3.426490)
\setfont{Helvetica}{0.3}
\newrgbcolor{dialinecolor}{0 0 0}
\psset{linecolor=dialinecolor}
\rput[r](3.336580,3.466850){\scalebox{1 -1}{1'}}
\setfont{Helvetica}{0.3}
\newrgbcolor{dialinecolor}{0 0 0}
\psset{linecolor=dialinecolor}
\rput[r](4.474080,4.191850){\scalebox{1 -1}{2'}}
\setfont{Helvetica}{0.3}
\newrgbcolor{dialinecolor}{0 0 0}
\psset{linecolor=dialinecolor}
\rput[r](7.832830,2.566850){\scalebox{1 -1}{1'}}
\psset{linewidth=0.040} \psset{linestyle=solid}
\psset{linestyle=solid}
\newrgbcolor{dialinecolor}{0 0 0}
\psset{linecolor=dialinecolor}
\psellipse(13.971690,4.275065)(0.247590,0.275715)
\psset{linewidth=0.050} \psset{linestyle=solid}
\psset{linestyle=solid} \setlinecaps{0}
\newrgbcolor{dialinecolor}{0 0 0}
\psset{linecolor=dialinecolor}
\psline(13.971700,4.550780)(13.959300,4.431760)
\setfont{Helvetica}{0.3}
\newrgbcolor{dialinecolor}{0 0 0}
\psset{linecolor=dialinecolor}
\rput[r](14.036600,4.285600){\scalebox{1 -1}{3}}
\psset{linewidth=0.040} \psset{linestyle=solid}
\psset{linestyle=solid}
\newrgbcolor{dialinecolor}{0 0 0}
\psset{linecolor=dialinecolor}
\psellipse(11.165930,4.047985)(0.259330,0.275715)
\psset{linewidth=0.050} \psset{linestyle=solid}
\psset{linestyle=solid} \setlinecaps{0}
\newrgbcolor{dialinecolor}{0 0 0}
\psset{linecolor=dialinecolor}
\psline(11.165900,4.323700)(11.165300,4.204680)
\setfont{Helvetica}{0.3}
\newrgbcolor{dialinecolor}{0 0 0}
\psset{linecolor=dialinecolor}
\rput[r](11.211800,4.080450){\scalebox{1 -1}{4}}
\psset{linewidth=0.040} \psset{linestyle=solid}
\psset{linestyle=solid}
\newrgbcolor{dialinecolor}{0 0 0}
\psset{linecolor=dialinecolor}
\psellipse(11.768530,4.628965)(0.259330,0.275715)
\psset{linewidth=0.050} \psset{linestyle=solid}
\psset{linestyle=solid} \setlinecaps{0}
\newrgbcolor{dialinecolor}{0 0 0}
\psset{linecolor=dialinecolor}
\psline(11.765300,4.754680)(11.768500,4.904680)
\setfont{Helvetica}{0.3}
\newrgbcolor{dialinecolor}{0 0 0}
\psset{linecolor=dialinecolor}
\rput[r](11.839400,4.623930){\scalebox{1 -1}{5}}
\psset{linewidth=0.040} \psset{linestyle=solid}
\psset{linestyle=solid}
\newrgbcolor{dialinecolor}{0 0 0}
\psset{linecolor=dialinecolor}
\psellipse(12.809930,1.937565)(0.259330,0.275715)
\psset{linewidth=0.050} \psset{linestyle=solid}
\psset{linestyle=solid} \setlinecaps{0}
\newrgbcolor{dialinecolor}{0 0 0}
\psset{linecolor=dialinecolor}
\psline(12.809900,2.213280)(12.809300,2.094260)
\setfont{Helvetica}{0.3}
\newrgbcolor{dialinecolor}{0 0 0}
\psset{linecolor=dialinecolor}
\rput[r](12.880800,1.970030){\scalebox{1 -1}{2}}
\psset{linewidth=0.040} \psset{linestyle=solid}
\psset{linestyle=solid}
\newrgbcolor{dialinecolor}{0 0 0}
\psset{linecolor=dialinecolor}
\psellipse(6.728410,3.997985)(0.259330,0.275715)
\psset{linewidth=0.050} \psset{linestyle=solid}
\psset{linestyle=solid} \setlinecaps{0}
\newrgbcolor{dialinecolor}{0 0 0}
\psset{linecolor=dialinecolor}
\psline(6.736580,4.416850)(6.736580,4.166850)
\setfont{Helvetica}{0.3}
\newrgbcolor{dialinecolor}{0 0 0}
\psset{linecolor=dialinecolor}
\rput[r](6.811780,4.005450){\scalebox{1 -1}{4}}
\setfont{Helvetica}{0.3}
\newrgbcolor{dialinecolor}{0 0 0}
\psset{linecolor=dialinecolor}
\rput[r](6.557830,3.679350){\scalebox{1 -1}{2'}}
\psset{linewidth=0.040} \psset{linestyle=solid}
\psset{linestyle=solid}
\newrgbcolor{dialinecolor}{0 0 0}
\psset{linecolor=dialinecolor}
\psellipse(8.590910,4.660425)(0.259330,0.275715)
\psset{linewidth=0.050} \psset{linestyle=solid}
\psset{linestyle=solid} \setlinecaps{0}
\newrgbcolor{dialinecolor}{0 0 0}
\psset{linecolor=dialinecolor}
\psline(8.600240,4.811140)(8.599080,5.054350)
\setfont{Helvetica}{0.3}
\newrgbcolor{dialinecolor}{0 0 0}
\psset{linecolor=dialinecolor}
\rput[r](8.661580,4.654350){\scalebox{1 -1}{6}}
\setfont{Helvetica}{0.3}
\newrgbcolor{dialinecolor}{0 0 0}
\psset{linecolor=dialinecolor}
\rput[r](8.449080,4.379350){\scalebox{1 -1}{4'}}
\psset{linewidth=0.040} \psset{linestyle=solid}
\psset{linestyle=solid}
\newrgbcolor{dialinecolor}{0 0 0}
\psset{linecolor=dialinecolor}
\psellipse(8.899660,1.825065)(0.259330,0.275715)
\psset{linewidth=0.050} \psset{linestyle=solid}
\psset{linestyle=solid} \setlinecaps{0}
\newrgbcolor{dialinecolor}{0 0 0}
\psset{linecolor=dialinecolor}
\psline(8.899660,2.100780)(8.899080,1.981760)
\setfont{Helvetica}{0.3}
\newrgbcolor{dialinecolor}{0 0 0}
\psset{linecolor=dialinecolor}
\rput[r](8.995530,1.870030){\scalebox{1 -1}{2}}
\psset{linewidth=0.040} \psset{linestyle=solid}
\psset{linestyle=solid}
\newrgbcolor{dialinecolor}{0 0 0}
\psset{linecolor=dialinecolor}
\psellipse(3.040910,4.033725)(0.259330,0.275715)
\psset{linewidth=0.050} \psset{linestyle=solid}
\psset{linestyle=solid} \setlinecaps{0}
\newrgbcolor{dialinecolor}{0 0 0}
\psset{linecolor=dialinecolor}
\psline(3.040910,4.309440)(3.040330,4.190420)
\setfont{Helvetica}{0.3}
\newrgbcolor{dialinecolor}{0 0 0}
\psset{linecolor=dialinecolor}
\rput[r](3.111780,4.041190){\scalebox{1 -1}{2}}
\psset{linewidth=0.040} \psset{linestyle=solid}
\psset{linestyle=solid}
\newrgbcolor{dialinecolor}{0 0 0}
\psset{linecolor=dialinecolor}
\psellipse(3.768500,4.464705)(0.259330,0.275715)
\psset{linewidth=0.050} \psset{linestyle=solid}
\psset{linestyle=solid} \setlinecaps{0}
\newrgbcolor{dialinecolor}{0 0 0}
\psset{linecolor=dialinecolor}
\psline(3.765330,4.590420)(3.768500,4.740420)
\setfont{Helvetica}{0.3}
\newrgbcolor{dialinecolor}{0 0 0}
\psset{linecolor=dialinecolor}
\rput[r](3.839370,4.484670){\scalebox{1 -1}{3}}
\psset{linewidth=0.040} \psset{linestyle=solid}
\psset{linestyle=solid}
\newrgbcolor{dialinecolor}{0 0 0}
\psset{linecolor=dialinecolor}
\psellipse(5.165910,4.583665)(0.259330,0.275715)
\psset{linewidth=0.050} \psset{linestyle=solid}
\psset{linestyle=solid} \setlinecaps{0}
\newrgbcolor{dialinecolor}{0 0 0}
\psset{linecolor=dialinecolor}
\psline(5.162740,4.709380)(5.165910,4.859380)
\setfont{Helvetica}{0.3}
\newrgbcolor{dialinecolor}{0 0 0}
\psset{linecolor=dialinecolor}
\rput[r](5.211780,4.591130){\scalebox{1 -1}{4}}
\psset{linewidth=0.040} \psset{linestyle=solid}
\psset{linestyle=solid}
\newrgbcolor{dialinecolor}{0 0 0}
\psset{linecolor=dialinecolor}
\psellipse(9.574660,4.150065)(0.259330,0.275715)
\psset{linewidth=0.050} \psset{linestyle=solid}
\psset{linestyle=solid} \setlinecaps{0}
\newrgbcolor{dialinecolor}{0 0 0}
\psset{linecolor=dialinecolor}
\psline(9.574660,4.425780)(9.574080,4.306760)
\setfont{Helvetica}{0.3}
\newrgbcolor{dialinecolor}{0 0 0}
\psset{linecolor=dialinecolor}
\rput[r](9.658030,4.170030){\scalebox{1 -1}{3}}
\setfont{Helvetica}{0.3}
\newrgbcolor{dialinecolor}{0 0 0}
\psset{linecolor=dialinecolor}
\rput[r](12.951600,4.246850){\scalebox{1 -1}{2'}}
\setfont{Helvetica}{0.3}
\newrgbcolor{dialinecolor}{0 0 0}
\psset{linecolor=dialinecolor}
\rput[r](11.439100,2.559350){\scalebox{1 -1}{1'}}
\psset{linewidth=0.040} \psset{linestyle=solid}
\psset{linestyle=solid}
\newrgbcolor{dialinecolor}{0 0 0}
\psset{linecolor=dialinecolor}
\psellipse(7.506000,4.516465)(0.259330,0.275715)
\psset{linewidth=0.050} \psset{linestyle=solid}
\psset{linestyle=solid} \setlinecaps{0}
\newrgbcolor{dialinecolor}{0 0 0}
\psset{linecolor=dialinecolor}
\psline(7.524080,4.679350)(7.511580,4.891850)
\setfont{Helvetica}{0.3}
\newrgbcolor{dialinecolor}{0 0 0}
\psset{linecolor=dialinecolor}
\rput[r](7.576870,4.523930){\scalebox{1 -1}{5}}
\setfont{Helvetica}{0.3}
\newrgbcolor{dialinecolor}{0 0 0}
\psset{linecolor=dialinecolor}
\rput[r](7.470330,4.166850){\scalebox{1 -1}{3'}} }\endpspicture\]
If the combinatorial dimension of the composed object, $s \circ
s'$, is less than the sum of those of $s$ and $s'$, then the
composition is zero. This is precisely the situation in the
following composition. \label{comp-0-by-dimension}
\ifx\setlinejoinmode\undefined
 \newcommand{\setlinejoinmode}[1]{}
\fi \ifx\setlinecaps\undefined
 \newcommand{\setlinecaps}[1]{}
\fi
\ifx\setfont\undefined
 \newcommand{\setfont}[2]{}
\fi 
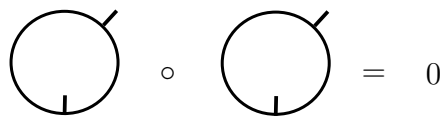
\[\pspicture(2.5,-3.458738)(8.262500,-1.999997) \scalebox{1
-1}{
\newrgbcolor{dialinecolor}{0 0 0}
\psset{linecolor=dialinecolor}
\newrgbcolor{diafillcolor}{1 1 1}
\psset{fillcolor=diafillcolor}
\newrgbcolor{dialinecolor}{0 0 0}
\psset{linecolor=dialinecolor}
\rput[r](4.728240,2.874890){\scalebox{1 -1}{$\circ$}}
\psset{linewidth=0.040} \psset{linestyle=solid}
\psset{linestyle=solid}
\newrgbcolor{dialinecolor}{0 0 0}
\psset{linecolor=dialinecolor}
\psellipse(3.251250,2.722820)(0.731250,0.7)
\psset{linewidth=0.050} \psset{linestyle=solid}
\psset{linestyle=solid} \setlinecaps{0}
\newrgbcolor{dialinecolor}{0 0 0}
\psset{linecolor=dialinecolor}
\psline(3.768320,2.227850)(3.945000,2.035320)
\psset{linewidth=0.050} \psset{linestyle=solid}
\psset{linestyle=solid} \setlinecaps{0}
\newrgbcolor{dialinecolor}{0 0 0}
\psset{linecolor=dialinecolor}
\psline(3.251250,3.422820)(3.257500,3.160320)
\psset{linewidth=0.040} \psset{linestyle=solid}
\psset{linestyle=solid}
\newrgbcolor{dialinecolor}{0 0 0}
\psset{linecolor=dialinecolor}
\psellipse(6.058750,2.733150)(0.731250,0.7)
\psset{linewidth=0.050} \psset{linestyle=solid}
\psset{linestyle=solid} \setlinecaps{0}
\newrgbcolor{dialinecolor}{0 0 0}
\psset{linecolor=dialinecolor}
\psline(6.575820,2.238170)(6.752500,2.045650)
\psset{linewidth=0.050} \psset{linestyle=solid}
\psset{linestyle=solid} \setlinecaps{0}
\newrgbcolor{dialinecolor}{0 0 0}
\psset{linecolor=dialinecolor}
\psline(6.058750,3.433150)(6.065000,3.170650)

\newrgbcolor{dialinecolor}{0 0 0}
\psset{linecolor=dialinecolor}
\rput[r](7.525000,2.887500){\scalebox{1 -1}{=}}

\newrgbcolor{dialinecolor}{0 0 0}
\psset{linecolor=dialinecolor}
\rput[r](8.262500,2.875000){\scalebox{1 -1}{0}} }\endpspicture\]
One can check that the above differential is a derivation of the
composition, {\it i.e.} the composition is a chain map between the
corresponding chain complexes.

Relabelling naturally gives rise to a map from the permutation
group $S_n$ to $C_0 \mathcal{S}^c(0, n,n) \subset C_0
\mathcal{S}^c(n,n)$, which induces an $S_n$-action on suitable
chord diagrams. More precisely, such a permutation corresponds to
$n$ disjoint circles in the plane, without chords, whose inputs
and outputs are numbered as prescribed by the permutation. The
following proposition organizes all of the above structures into a
single mathematical object.

\begin{prop}[PROP of Cyclic Sullivan Chord Diagrams]
The collection $C_\ast \mathcal{S}^c$, of chain complexes
$C_\ast\mathcal{S}^c (m, n)$, for $m, n \geq 1$, together with the
above composition rule is a PROP in the category of chain
complexes. The tensor product is the disjoint union.
\end{prop}

This PROP is referred to as the \emph{PROP of cyclic Sullivan
chord diagrams}.

\begin{exam}\label{brace-smile}
Let's look at the following four chord diagrams.
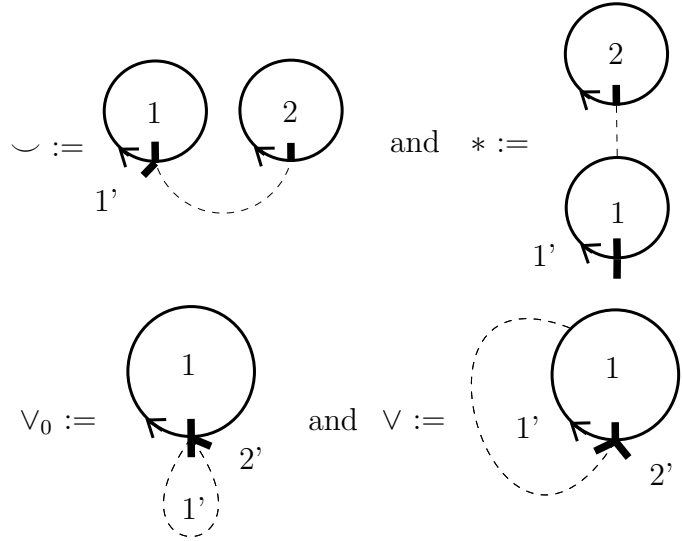
\ifx\setlinejoinmode\undefined
 \newcommand{\setlinejoinmode}[1]{}
\fi \ifx\setlinecaps\undefined
 \newcommand{\setlinecaps}[1]{}
\fi
\ifx\setfont\undefined
 \newcommand{\setfont}[2]{}
\fi \[\pspicture(1.621670,-6.330820)(10.275990,-2.582400)
\scalebox{1 -1}{
\newrgbcolor{dialinecolor}{0 0 0}
\psset{linecolor=dialinecolor}
\newrgbcolor{diafillcolor}{1 1 1}
\psset{fillcolor=diafillcolor}
\newrgbcolor{dialinecolor}{0 0 0}
\psset{linecolor=dialinecolor}
\rput[r](2.471670,4.555830){\scalebox{1 -1}{$\smile$ :=}}
\newrgbcolor{dialinecolor}{0 0 0}
\psset{linecolor=dialinecolor} \rput[r](7.2,4.45){\scalebox{1
-1}{and}}
\newrgbcolor{dialinecolor}{1 1 1}
\psset{linecolor=dialinecolor}
\psellipse*(3.416950,4.050830)(0.699400,0.687500)
\psset{linewidth=0.040} \psset{linestyle=solid}
\psset{linestyle=solid}
\newrgbcolor{dialinecolor}{0 0 0}
\psset{linecolor=dialinecolor}
\psellipse(3.416950,4.050830)(0.699400,0.687500)
\psset{linewidth=0.040} \psset{linestyle=solid}
\psset{linestyle=solid} \setlinecaps{0}
\newrgbcolor{dialinecolor}{0 0 0}
\psset{linecolor=dialinecolor}
\psline(3.065420,4.641770)(2.976175,4.576374)
\psset{linewidth=0.040} \psset{linestyle=solid}
\setlinejoinmode{0} \setlinecaps{0}
\newrgbcolor{dialinecolor}{0 0 0}
\psset{linecolor=dialinecolor}
\psline(3.199272,4.553892)(2.949287,4.556672)(3.021952,4.795879)
\newrgbcolor{dialinecolor}{1 1 1}
\psset{linecolor=dialinecolor}
\psellipse*(5.221070,4.039900)(0.699400,0.687500)
\psset{linewidth=0.040} \psset{linestyle=solid}
\psset{linestyle=solid}
\newrgbcolor{dialinecolor}{0 0 0}
\psset{linecolor=dialinecolor}
\psellipse(5.221070,4.039900)(0.699400,0.687500)
\psset{linewidth=0.040} \psset{linestyle=solid}
\psset{linestyle=solid} \setlinecaps{0}
\newrgbcolor{dialinecolor}{0 0 0}
\psset{linecolor=dialinecolor}
\psline(4.869540,4.630840)(4.780295,4.565444)
\psset{linewidth=0.040} \psset{linestyle=solid}
\setlinejoinmode{0} \setlinecaps{0}
\newrgbcolor{dialinecolor}{0 0 0}
\psset{linecolor=dialinecolor}
\psline(5.003392,4.542962)(4.753407,4.545742)(4.826072,4.784949)
\psset{linewidth=0.1} \psset{linestyle=solid}
\psset{linestyle=solid} \setlinecaps{0}
\newrgbcolor{dialinecolor}{0 0 0}
\psset{linecolor=dialinecolor}
\psline(5.219080,4.473410)(5.221070,4.727400)
\psset{linewidth=0.010} \psset{linestyle=dashed,dash=1 1}
\psset{linestyle=dashed,dash=0.1 0.1} \setlinecaps{0}
\newrgbcolor{dialinecolor}{0 0 0}
\psset{linecolor=dialinecolor}
\psclip{\pswedge[linestyle=none,fillstyle=none](4.317445,4.474475){1.327034}{15.636962}{163.668810}}
\psellipse(4.317445,4.474475)(0.938355,0.938355)
\endpsclip
\newrgbcolor{dialinecolor}{0 0 0}
\psset{linecolor=dialinecolor}
\rput(3.401670,4.069900){\scalebox{1 -1}{1}}
\newrgbcolor{dialinecolor}{0 0 0}
\psset{linecolor=dialinecolor}
\rput(5.214170,4.057400){\scalebox{1 -1}{2}}
\newrgbcolor{dialinecolor}{0 0 0}
\psset{linecolor=dialinecolor}
\rput(2.764170,5.232400){\scalebox{1 -1}{1'}}
\psset{linewidth=0.1} \psset{linestyle=solid}
\psset{linestyle=solid} \setlinecaps{0}
\newrgbcolor{dialinecolor}{0 0 0}
\psset{linecolor=dialinecolor}
\psline(3.416950,4.738330)(3.256850,4.906510)
\psset{linewidth=0.1} \psset{linestyle=solid}
\psset{linestyle=solid} \setlinecaps{0}
\newrgbcolor{dialinecolor}{0 0 0}
\psset{linecolor=dialinecolor}
\psline(3.413310,4.456620)(3.416950,4.738330)
\newrgbcolor{dialinecolor}{1 1 1}
\psset{linecolor=dialinecolor}
\psellipse*(9.557655,5.334910)(0.698335,0.687500)
\psset{linewidth=0.040} \psset{linestyle=solid}
\psset{linestyle=solid}
\newrgbcolor{dialinecolor}{0 0 0}
\psset{linecolor=dialinecolor}
\psellipse(9.557655,5.334910)(0.698335,0.687500)
\newrgbcolor{dialinecolor}{1 1 1}
\psset{linecolor=dialinecolor}
\psellipse*(9.547100,3.289900)(0.699400,0.687500)
\psset{linewidth=0.040} \psset{linestyle=solid}
\psset{linestyle=solid}
\newrgbcolor{dialinecolor}{0 0 0}
\psset{linecolor=dialinecolor}
\psellipse(9.547100,3.289900)(0.699400,0.687500)
\psset{linewidth=0.010} \psset{linestyle=dashed,dash=1 1}
\psset{linestyle=dashed,dash=0.1 0.1} \setlinecaps{0}
\newrgbcolor{dialinecolor}{0 0 0}
\psset{linecolor=dialinecolor}
\psline(9.557655,4.647410)(9.547100,3.977400)
\newrgbcolor{dialinecolor}{0 0 0}
\psset{linecolor=dialinecolor}
\rput[r](8.389320,4.494900){\scalebox{1 -1}{$*$ :=}}
\psset{linewidth=0.040} \psset{linestyle=solid}
\psset{linestyle=solid} \setlinecaps{0}
\newrgbcolor{dialinecolor}{0 0 0}
\psset{linecolor=dialinecolor}
\psline(9.226820,5.930840)(9.119146,5.858296)
\psset{linewidth=0.040} \psset{linestyle=solid}
\setlinejoinmode{0} \setlinecaps{0}
\newrgbcolor{dialinecolor}{0 0 0}
\psset{linecolor=dialinecolor}
\psline(9.341182,5.827022)(9.091502,5.839671)(9.173556,6.075822)
\psset{linewidth=0.040} \psset{linestyle=solid}
\psset{linestyle=solid} \setlinecaps{0}
\newrgbcolor{dialinecolor}{0 0 0}
\psset{linecolor=dialinecolor}
\psline(9.195570,3.880840)(9.106324,3.815441)
\psset{linewidth=0.040} \psset{linestyle=solid}
\setlinejoinmode{0} \setlinecaps{0}
\newrgbcolor{dialinecolor}{0 0 0}
\psset{linecolor=dialinecolor}
\psline(9.329421,3.792962)(9.079437,3.795739)(9.152098,4.034946)
\newrgbcolor{dialinecolor}{0 0 0}
\psset{linecolor=dialinecolor}
\rput(9.539300,3.294900){\scalebox{1 -1}{2}}
\newrgbcolor{dialinecolor}{0 0 0}
\psset{linecolor=dialinecolor}
\rput(9.564300,5.407400){\scalebox{1 -1}{1}}
\newrgbcolor{dialinecolor}{0 0 0}
\psset{linecolor=dialinecolor}
\rput(8.589320,5.969900){\scalebox{1 -1}{1'}}
\psset{linewidth=0.1} \psset{linestyle=solid}
\psset{linestyle=solid} \setlinecaps{0}
\newrgbcolor{dialinecolor}{0 0 0}
\psset{linecolor=dialinecolor}
\psline(9.544000,3.701660)(9.547100,3.977400)
\psset{linewidth=0.1} \psset{linestyle=solid}
\psset{linestyle=solid} \setlinecaps{0}
\newrgbcolor{dialinecolor}{0 0 0}
\psset{linecolor=dialinecolor}
\psline(9.552800,5.743430)(9.557655,6.022410)
\psset{linewidth=0.1} \psset{linestyle=solid}
\psset{linestyle=solid} \setlinecaps{0}
\newrgbcolor{dialinecolor}{0 0 0}
\psset{linecolor=dialinecolor}
\psline(9.557655,6.022410)(9.554800,6.280270) }\endpspicture\]
\ifx\setlinejoinmode\undefined
 \newcommand{\setlinejoinmode}[1]{}
\fi \ifx\setlinecaps\undefined
 \newcommand{\setlinecaps}[1]{}
\fi
\ifx\setfont\undefined
 \newcommand{\setfont}[2]{}
\fi \[ \pspicture(1.238340,-6.493967)(10.164583,-3.389170)
\scalebox{1 -1}{
\newrgbcolor{dialinecolor}{0 0 0}
\psset{linecolor=dialinecolor}
\newrgbcolor{diafillcolor}{1 1 1}
\psset{fillcolor=diafillcolor}
\newrgbcolor{dialinecolor}{0 0 0}
\psset{linecolor=dialinecolor}
\rput[r](5.865529,4.930065){\scalebox{1 -1}{and}}

\newrgbcolor{dialinecolor}{0 0 0}
\psset{linecolor=dialinecolor}
\rput[r](2.388340,4.972500){\scalebox{1 -1}{$\vee_0$ :=}}

\newrgbcolor{dialinecolor}{0 0 0}
\psset{linecolor=dialinecolor}
\rput[r](7.028980,4.945885){\scalebox{1 -1}{$\vee$ :=}}
\newrgbcolor{dialinecolor}{1 1 1}
\psset{linecolor=dialinecolor}
\psellipse*(3.643565,4.293695)(0.861435,0.884525)
\psset{linewidth=0.040} \psset{linestyle=solid}
\psset{linestyle=solid}
\newrgbcolor{dialinecolor}{0 0 0}
\psset{linecolor=dialinecolor}
\psellipse(3.643565,4.293695)(0.861435,0.884525)
\psset{linewidth=0.1} \psset{linestyle=solid}
\psset{linestyle=solid} \setlinecaps{0}
\newrgbcolor{dialinecolor}{0 0 0}
\psset{linecolor=dialinecolor}
\psline(3.641465,4.923545)(3.643565,5.178220)
\psset{linewidth=0.040} \psset{linestyle=solid}
\psset{linestyle=solid} \setlinecaps{0}
\newrgbcolor{dialinecolor}{0 0 0}
\psset{linecolor=dialinecolor}
\psline(3.260215,5.079795)(3.088758,4.957799)
\psset{linewidth=0.040} \psset{linestyle=solid}
\setlinejoinmode{0} \setlinecaps{0}
\newrgbcolor{dialinecolor}{0 0 0}
\psset{linecolor=dialinecolor}
\psline(3.311520,4.932204)(3.061598,4.938474)(3.137595,5.176643)
\psset{linewidth=0.010} \psset{linestyle=dashed,dash=1 1}
\psset{linestyle=dashed,dash=0.1 0.1} \setlinejoinmode{0}
\setlinecaps{0}
\newrgbcolor{dialinecolor}{0 0 0}
\psset{linecolor=dialinecolor} \pscustom{
\newpath
\moveto(3.643565,5.178220)
\curveto(2.372710,6.923545)(4.904995,6.928220)(3.643565,5.178220)
\stroke} \psset{linewidth=0.1} \psset{linestyle=solid}
\psset{linestyle=solid} \setlinecaps{0}
\newrgbcolor{dialinecolor}{0 0 0}
\psset{linecolor=dialinecolor}
\psline(3.906282,5.287653)(3.643565,5.178220)
\psset{linewidth=0.1} \psset{linestyle=solid}
\psset{linestyle=solid} \setlinecaps{0}
\newrgbcolor{dialinecolor}{0 0 0}
\psset{linecolor=dialinecolor}
\psline(3.643565,5.178220)(3.645536,5.433494)
\newrgbcolor{dialinecolor}{0 0 0}
\psset{linecolor=dialinecolor}
\rput(3.610181,4.204896){\scalebox{1 -1}{1}}
\newrgbcolor{dialinecolor}{0 0 0}
\psset{linecolor=dialinecolor}
\rput(3.571411,6.108382){\scalebox{1 -1}{\,\,\,1'}}
\newrgbcolor{dialinecolor}{0 0 0}
\psset{linecolor=dialinecolor}
\rput(4.443423,5.445469){\scalebox{1 -1}{2'}}
\newrgbcolor{dialinecolor}{1 1 1}
\psset{linecolor=dialinecolor}
\psellipse*(9.283148,4.337791)(0.861435,0.884525)
\psset{linewidth=0.040} \psset{linestyle=solid}
\psset{linestyle=solid}
\newrgbcolor{dialinecolor}{0 0 0}
\psset{linecolor=dialinecolor}
\psellipse(9.283148,4.337791)(0.861435,0.884525)
\psset{linewidth=0.1} \psset{linestyle=solid}
\psset{linestyle=solid} \setlinecaps{0}
\newrgbcolor{dialinecolor}{0 0 0}
\psset{linecolor=dialinecolor}
\psline(9.281048,4.967641)(9.283148,5.222316)
\psset{linewidth=0.040} \psset{linestyle=solid}
\psset{linestyle=solid} \setlinecaps{0}
\newrgbcolor{dialinecolor}{0 0 0}
\psset{linecolor=dialinecolor}
\psline(8.899798,5.123891)(8.728341,5.001895)
\psset{linewidth=0.040} \psset{linestyle=solid}
\setlinejoinmode{0} \setlinecaps{0}
\newrgbcolor{dialinecolor}{0 0 0}
\psset{linecolor=dialinecolor}
\psline(8.951103,4.976300)(8.701181,4.982570)(8.777178,5.220739)
\psset{linewidth=0.010} \psset{linestyle=dashed,dash=1 1}
\psset{linestyle=dashed,dash=0.1 0.1} \setlinejoinmode{0}
\setlinecaps{0}
\newrgbcolor{dialinecolor}{0 0 0}
\psset{linecolor=dialinecolor} \pscustom{
\newpath
\moveto(8.674022,3.712337)
\curveto(6.222056,2.759747)(7.689303,7.727172)(9.283148,5.222316)
\stroke} \psset{linewidth=0.1} \psset{linestyle=solid}
\psset{linestyle=solid} \setlinecaps{0}
\newrgbcolor{dialinecolor}{0 0 0}
\psset{linecolor=dialinecolor}
\psline(9.010709,5.353945)(9.283148,5.222316)
\psset{linewidth=0.1} \psset{linestyle=solid}
\psset{linestyle=solid} \setlinecaps{0}
\newrgbcolor{dialinecolor}{0 0 0}
\psset{linecolor=dialinecolor}
\psline(9.283148,5.222316)(9.458456,5.437914)
\newrgbcolor{dialinecolor}{0 0 0}
\psset{linecolor=dialinecolor}
\rput(9.249764,4.248992){\scalebox{1 -1}{1}}
\newrgbcolor{dialinecolor}{0 0 0}
\psset{linecolor=dialinecolor}
\rput(8.114979,5.029946){\scalebox{1 -1}{1'}}
\newrgbcolor{dialinecolor}{0 0 0}
\psset{linecolor=dialinecolor}
\rput(9.897390,5.657503){\scalebox{1
-1}{2'}} }\endpspicture \] %
Note that $\smile\in C_0\mathcal S^c(2,1)$ and $*\in C_1\mathcal
S^c(2,1)$. In each case the input circles are labelled by 1 and 2,
and the output circle by 1'. Similarly, $\vee_0\in C_0\mathcal
S^c(1,2)$ and $\vee\in C_1\mathcal S^c(1,2)$. Each input is labelled
by 1, and the output circles are labelled by 1' and 2'. The notation
$\smile$ was taken from \cite{G}, and $\vee_0$ and $\vee$ are
borrowed from \cite{Su}.

Observe that $\smile$ and $\vee_0$ are both closed elements. It is
straightforward to check that these elements satisfy associativity
and coassociativity,
$$\smile\circ (\smile\otimes id)=\smile\circ (id\otimes \smile)$$
and $$(id\otimes \vee_0)\circ \vee_0=(\vee_0\otimes id)\circ
\vee_0$$ Here $id\in C_0 \mathcal S^c(1,1)$ is the cyclic Sullivan
chord diagram consisting of one circle without chords, whose input
and output marked points coincide.

Let $\tau_2\in C_0 \mathcal S^c(2,2)$ denote the element that
switches the labelling, as defined by the map $S_2\to C_0 \mathcal
S^c(2,2)$. It is easy to see that $$\partial(*)=\smile -
(\smile\circ\tau_2)$$ and $$\partial(\vee)=\vee_0 - (\tau_2\circ
\vee_0) $$ This implies, that after passing to homology, $\smile$ and
$\vee_0$ are commutative and cocommutative, respectively.\\ Finally,
a word about the Frobenius compatibility conditions. We have, $$
\vee_0\circ \smile=(id\otimes \smile)\circ(\vee_0\otimes
id)=\tau_2\circ (\smile \otimes id) \circ (id \otimes (\tau_2\circ
\vee_0)) $$ Since $\vee_0$ is cocommutative on homology, this
equation implies the Frobenius compatibility condition at the level
of homology.

If we define
\[ \pspicture(3.940753,-3.458738)(6.810000,-2.010327)
\scalebox{1 -1}{
\newrgbcolor{dialinecolor}{0 0 0}
\psset{linecolor=dialinecolor}
\newrgbcolor{diafillcolor}{1 1 1}
\psset{fillcolor=diafillcolor} \psset{linewidth=0.04}
\psset{linestyle=solid} \psset{linestyle=solid}
\newrgbcolor{dialinecolor}{0 0 0}
\psset{linecolor=dialinecolor}
\psellipse(6.058750,2.733150)(0.731250,0.700000)
\psset{linewidth=0.050000} \psset{linestyle=solid}
\psset{linestyle=solid} \setlinecaps{0}
\newrgbcolor{dialinecolor}{0 0 0}
\psset{linecolor=dialinecolor}
\psline(6.575820,2.238170)(6.752500,2.045650)
\psset{linewidth=0.050000} \psset{linestyle=solid}
\psset{linestyle=solid} \setlinecaps{0}
\newrgbcolor{dialinecolor}{0 0 0}
\psset{linecolor=dialinecolor}
\psline(6.058750,3.433150)(6.065000,3.170650)
\setfont{Helvetica}{0.400000}
\newrgbcolor{dialinecolor}{0 0 0}
\psset{linecolor=dialinecolor}
\rput[r](4.872509,2.840134){\scalebox{1 -1}{$:=$}}
\setfont{Helvetica}{0.400000}
\newrgbcolor{dialinecolor}{0 0 0}
\psset{linecolor=dialinecolor}
\rput[r](4.253253,2.806798){\scalebox{1 -1}{$\Delta$}} }\endpspicture
\] then $\Delta^2=0$ (see page \pageref{comp-0-by-dimension}), and one may ascertain that the $BV$ relation is satisfied
on homology; see \cite[sec. 5]{CS1}.
\begin{multline*}
 \Delta\circ \smile\circ (\smile\otimes id) \,\,\simeq\,\,\, %
 \smile\circ (\Delta\circ \smile\otimes id)%
+ \smile\circ (id\otimes \Delta\circ \smile)\\%
+ \smile\circ (\Delta\circ \smile\otimes id)\circ(id \otimes\tau_2)%
+ \smile\circ (\smile\otimes id)\circ(\Delta\otimes id\otimes id)\\%
+ \smile\circ (\smile\otimes id)\circ(id\otimes \Delta\otimes id)%
+ \smile\circ (\smile\otimes id)\circ(id\otimes id\otimes \Delta).%
\end{multline*}
Similarly, the dual $coBV$ relations are satisfied. That is to
say,
\begin{multline*}
 (\vee_0\otimes id)\circ \vee_0\circ\Delta \,\, \simeq \,\, %
 (\vee_0\circ \Delta\otimes id)\circ \vee_0%
+ (id\otimes \vee_0\circ\Delta)\circ \vee_0\\%
+ (id\otimes \tau_2)\circ (\vee_0\circ \Delta\otimes id)\circ\vee_0%
+ (\Delta\otimes id\otimes id)\circ (\vee_0\otimes id)\circ \vee_0\\%
+ (id\otimes \Delta\otimes id)\circ (\vee_0\otimes id)\circ \vee_0%
+ (id\otimes id\otimes \Delta)\circ (\vee_0\otimes id)\circ \vee_0.%
\end{multline*}
This, however, turns out to be true for more trivial reasons. Each
individual term of the above equation is in fact homologous to
zero.
\end{exam}

\begin{comment}
Vector spaces generated by diagrams of type $(0; n, 1)$ whose
chords in the plane do not cross are closely related to the operad
of cacti. In a cyclic Sullivan chord diagram, collapsing each
chord to a point gives rise to a cactus. In fact, the operation of
collapsing chords establishes an isomorphism of operads.
\end{comment}

\section{The Associative Case}\label{sec2}

Now that the PROP of cyclic Sullivan chord diagrams is built, we
want to make it act. Let us first recall a few relevant notations
and definitions.

Let $(A,\cdot,1)$ be a finite dimensional, unital, and associative
algebra over a ground field $k$. $A$ and $A^*$ are both examples
of $A$-bimodules. More precisely, the left and right
multiplications give $A$ an $A$-bimodule structure. The
$A$-bimodule structure of $A^*:=Hom(A,k)$ is given by
$(a_1.a^*.a_2)(a_3):=a^*(a_2\cdot a_3\cdot a_1)$, for any
$a_1,a_2,a_3\in A$ and $a^*\in A^*$. Let $\beta:A\to A^*$ be an
isomorphism of $A$-bimodules whose inverse we denote by
$\gamma:=:A^*\to A$. Define an inner product $<-,->:A\otimes A\to
k$ by $<a_1,a_2>:= (\beta (a_1))(a_2)$. It is easy to verify that
$\beta$ is an $A$-bimodule isomorphism, if and only if, $<-,->$ is
a non-degenerate bilinear map, satisfying
\begin{equation*}\label{invar}
\begin{array}{c}
 <a\cdot b, c> = <a,b\cdot c> \\
 <a\cdot b, c> = <b,c\cdot a>
\end{array}
\end{equation*}
This implies that the map $(a_1, ..., a_r)\mapsto
<a_1\cdot...\cdot a_r,1>$ is invariant under a cyclic rotation of
$a_1,...,a_r$, {\it i.e.}

\begin{equation}\label{cyclic-inv}
<a_1\cdot ...\cdot a_r,1>=<a_r\cdot a_1\cdot...\cdot a_{r-1},1>
\end{equation}

Let us recall the definition of the normalized Hochschild cochain
complex and that of the endomorphism PROP.

\begin{defn}[Normalized Hochschild Cochain Complex]\label{DefHoch} Let $M$ be
an $A$-bimodule. The \emph{Hochschild cochain complex of $A$ with
values in $M$} is the graded vector space $ HC^*(A;M):=
\prod_{n\geq 0} Hom(A^{\otimes n},M)$, endowed with the
differential
\begin{eqnarray*}
\left(\delta(f)\right)(a_1, ...,a_n)&:=&a_1. f(a_2,...,a_n) \\%
&&+ \sum_{j=1}^{n-1}(-1)^{j}\cdot f(a_1,...,a_j\cdot a_{j+1},...,a_n)\\%
&&+ (-1)^{n}\cdot f(a_1,...,a_{n-1}). a_n,
\end{eqnarray*}
where ``$.$'' denotes the left and right module structures. A
straightforward check shows that $\delta^2=0$; see e.g. \cite[1.5.1]{L}.\\%
\
The \emph{normalized Hochschild cochain complex of $A$ with values
in $M$} is the subcomplex $$ \overline{HC^*}(A;M):= \{ f\in
HC^*(A;M) | f(a_1,...,a_n)=0 \text{ if one of the } a_j=1\} $$ It
is a well known fact that the inclusion $\overline{HC^*}(A;M)
\hookrightarrow HC^*(A;M)$ is a quasi-isomorphism; see e.g.
\cite[1.5.7]{L}.
\end{defn}

The bimodule isomorphism $\beta:A\stackrel{\cong}{\to} A^*$
induces an isomorphism of chain complexes $\beta_\sharp:\overline{
HC^*}(A;A) \stackrel{ \cong} {\to}\overline{ HC^*}(A;A^*)$, where
$\beta_\sharp(f):=\beta\circ f$.

\begin{defn}[Endomorphism PROP] Let $V$ be a differential graded
vector space over $k$. The \emph{endomorphism PROP of $V$} is
collection of differential graded vector spaces $\mathcal
End_V(k,l):=Hom(V^{\otimes k},V^{\otimes l})$. The map $S_k\to
\mathcal End_V(k,k)$ is given by $\sigma(v_1\otimes ...\otimes
v_k):= (-1)^{| \sigma|} v_{\sigma(1)}\otimes ... \otimes
v_{\sigma(k) }$, and the composition $\circ:\mathcal
End_V(k,l)\otimes \mathcal End_V(m,k)\to \mathcal End_V(m,l)$ is
defined by
\begin{equation}\label{End-comp}
(F\circ G) (v_1 ,... , v_{m}):= F(G(v_1,...,v_{m})).
\end{equation}
\end{defn}

\begin{thm}\label{theorem1} Let $A$ be a finite dimensional, unital, and
associative algebra with a non-degenerate and invariant inner
product. Then, the normalized Hochschild cochain complex of $A$ is
an algebra over the PROP, $C_\ast \mathcal{S}^c$, of cyclic
Sullivan chord diagrams.
\end{thm}

\begin{cor} Under the above assumptions, the Hochschild cohomology of
$A$ is a Frobenius algebra endowed with a compatible $BV$ operator.
\end{cor}

\begin{proof}[Proof of the Theorem \ref{theorem1}]

The objective is to establish a map $ \alpha:C_*\mathcal S^c\to
\mathcal End_{\overline{ HC^*}(A;A)}$ which respects the
differentials, composition, and symmetric group action. We will
achieve this in four steps.

{\it Step I: Construction of $\alpha$}

We need to define maps $$C_*\mathcal S^c (k,l)\otimes \overline{
HC^*}(A;A)^{\otimes k}\to \overline{ HC^*} (A; A)^{\otimes l}$$ Let
$s\in C_*\mathcal S^c (k,l)$, and $f_1,...,f_k\in \overline{
HC^*}(A;A)\cong \overline{ HC^*}(A;A^*)$. Each $f_i:A^{\otimes
n_i}\to A^*$ may be regarded as an element of $(A^*)^{\otimes
n_i}\otimes A^*$ and therefore be written as
$f_i=(c^i_1,...,c^i_{n_i};c^i_{n_i +1})$. In order to define
$\big(\alpha (s)\big)(f_1,...,f_k)\in \overline{ HC^*} (A;
A)^{\otimes l}$ proceed as follows:

\begin{itemize}
\item[(a)] Consider the cyclic Sullivan chord diagram $s$ and
place $f_i=(c^i_1,...,c^i_{n_i};c^i_{n_i +1})$, for $i=1,...,k$,
inside and around the $i^{th}$ input circle of $s$. To be more
precise, at the $i^{th}$ circle of $s$, start at the input marked
point with the last element $c^i_{n_i +1}$ and proceed with
$c^i_1$, $c^i_2$, ..., $c^i_{n_i}$ in the clockwise direction
along the input circle; see figure below.


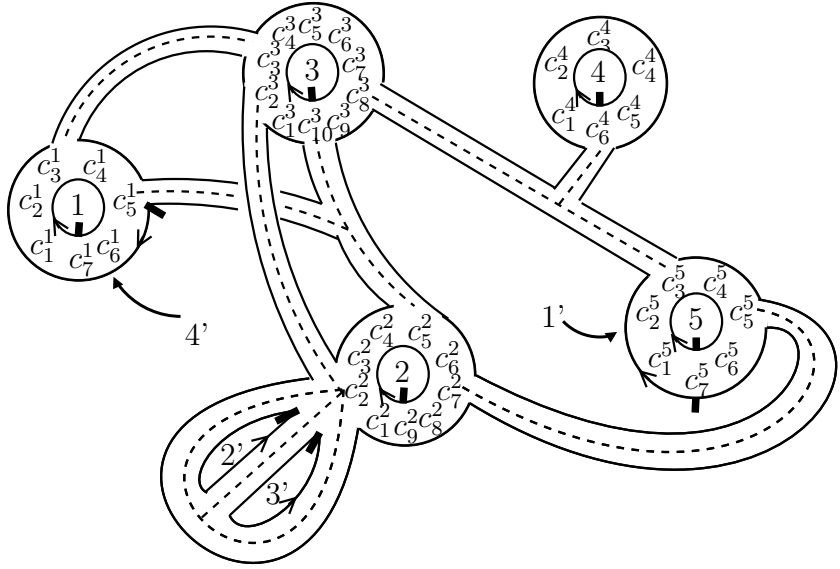
\ifx\setlinejoinmode\undefined
 \newcommand{\setlinejoinmode}[1]{}
\fi \ifx\setlinecaps\undefined
 \newcommand{\setlinecaps}[1]{}
\fi \ifx\setfont\undefined
 \newcommand{\setfont}[2]{}
\fi \[ \pspicture(0.272500,-9.967765)(11.320550,-2.380830)
\scalebox{1 -1}{
\newrgbcolor{dialinecolor}{0 0 0}
\psset{linecolor=dialinecolor}
\newrgbcolor{diafillcolor}{1 1 1}
\psset{fillcolor=diafillcolor} \psset{linewidth=0.03}
\psset{linestyle=solid} \psset{linestyle=solid} \setlinecaps{0}
\setlinejoinmode{0} \setlinecaps{0} \setlinejoinmode{0}
\psset{linestyle=solid}
\newrgbcolor{dialinecolor}{1 1 1}
\psset{linecolor=dialinecolor}
\psellipse*(4.356665,3.440005)(0.945835,0.945835)
\newrgbcolor{dialinecolor}{0 0 0}
\psset{linecolor=dialinecolor}
\psellipse(4.356665,3.440005)(0.945835,0.945835)
\psset{linewidth=0.003000} \setlinecaps{0} \setlinejoinmode{0}
\psset{linestyle=solid}
\newrgbcolor{dialinecolor}{0 0 0}
\psset{linecolor=dialinecolor}
\psellipse(4.356665,3.440005)(0.945835,0.945835)
\newrgbcolor{dialinecolor}{1 1 1}
\psset{linecolor=dialinecolor}
\pspolygon*(3.410630,2.916670)(3.410630,3.193340)(3.793960,3.193340)(3.793960,2.916670)
\psset{linewidth=0.04} \psset{linestyle=solid}
\psset{linestyle=solid} \setlinejoinmode{0}
\newrgbcolor{dialinecolor}{1 1 1}
\psset{linecolor=dialinecolor}
\pspolygon(3.410630,2.916670)(3.410630,3.193340)(3.793960,3.193340)(3.793960,2.916670)
\psset{linewidth=0.03} \psset{linestyle=solid}
\psset{linestyle=solid} \setlinecaps{0} \setlinejoinmode{0}
\setlinecaps{0} \setlinejoinmode{0} \psset{linestyle=solid}
\newrgbcolor{dialinecolor}{1 1 1}
\psset{linecolor=dialinecolor}
\psellipse*(1.233335,5.260005)(0.945835,0.945835)
\newrgbcolor{dialinecolor}{0 0 0}
\psset{linecolor=dialinecolor}
\psellipse(1.233335,5.260005)(0.945835,0.945835)
\psset{linewidth=0.003000} \setlinecaps{0} \setlinejoinmode{0}
\psset{linestyle=solid}
\newrgbcolor{dialinecolor}{0 0 0}
\psset{linecolor=dialinecolor}
\psellipse(1.233335,5.260005)(0.945835,0.945835)
\newrgbcolor{dialinecolor}{1 1 1}
\psset{linecolor=dialinecolor}
\pspolygon*(0.870630,4.2)(0.870630,4.533330)(1.193960,4.533330)(1.193960,4.2)
\psset{linewidth=0.04} \psset{linestyle=solid}
\psset{linestyle=solid} \setlinejoinmode{0}
\newrgbcolor{dialinecolor}{1 1 1}
\psset{linecolor=dialinecolor}
\pspolygon(0.870630,4.2)(0.870630,4.533330)(1.193960,4.533330)(1.193960,4.2)
\newrgbcolor{dialinecolor}{1 1 1}
\psset{linecolor=dialinecolor}
\pspolygon*(1.973970,4.88)(1.973970,5.135000)(2.190640,5.135000)(2.190640,4.88)
\psset{linewidth=0.04} \psset{linestyle=solid}
\psset{linestyle=solid} \setlinejoinmode{0}
\newrgbcolor{dialinecolor}{1 1 1}
\psset{linecolor=dialinecolor}
\pspolygon(1.973970,4.88)(1.973970,5.135000)(2.190640,5.135000)(2.190640,4.88)
\psset{linewidth=0.03} \psset{linestyle=solid}
\psset{linestyle=solid} \setlinecaps{0} \setlinejoinmode{0}
\setlinecaps{0} \setlinejoinmode{0} \psset{linestyle=solid}
\newrgbcolor{dialinecolor}{1 1 1}
\psset{linecolor=dialinecolor}
\psellipse*(5.573335,7.456665)(0.945835,0.945835)
\newrgbcolor{dialinecolor}{0 0 0}
\psset{linecolor=dialinecolor}
\psellipse(5.573335,7.456665)(0.945835,0.945835)
\psset{linewidth=0.003000} \setlinecaps{0} \setlinejoinmode{0}
\psset{linestyle=solid}
\newrgbcolor{dialinecolor}{0 0 0}
\psset{linecolor=dialinecolor}
\psellipse(5.573335,7.456665)(0.945835,0.945835)
\newrgbcolor{dialinecolor}{1 1 1}
\psset{linecolor=dialinecolor}
\pspolygon*(5.443970,6.33)(5.443970,6.818330)(6.043970,6.818330)(6.043970,6.33)
\psset{linewidth=0.04} \psset{linestyle=solid}
\psset{linestyle=solid} \setlinejoinmode{0}
\newrgbcolor{dialinecolor}{1 1 1}
\psset{linecolor=dialinecolor}
\pspolygon(5.443970,6.33)(5.443970,6.818330)(6.043970,6.818330)(6.043970,6.33)
\newrgbcolor{dialinecolor}{1 1 1}
\psset{linecolor=dialinecolor}
\pspolygon*(4.257300,4.213330)(4.257300,4.585000)(4.627300,4.585000)(4.627300,4.213330)
\psset{linewidth=0.04} \psset{linestyle=solid}
\psset{linestyle=solid} \setlinejoinmode{0}
\newrgbcolor{dialinecolor}{1 1 1}
\psset{linecolor=dialinecolor}
\pspolygon(4.257300,4.213330)(4.257300,4.585000)(4.627300,4.585000)(4.627300,4.213330)
\newrgbcolor{dialinecolor}{1 1 1}
\psset{linecolor=dialinecolor}
\pspolygon*(3.457300,3.58)(3.457300,4.085000)(3.877300,4.085000)(3.877300,3.58)
\psset{linewidth=0.04} \psset{linestyle=solid}
\psset{linestyle=solid} \setlinejoinmode{0}
\newrgbcolor{dialinecolor}{1 1 1}
\psset{linecolor=dialinecolor}
\pspolygon(3.457300,3.58)(3.457300,4.085000)(3.877300,4.085000)(3.877300,3.58)
\newrgbcolor{dialinecolor}{1 1 1}
\psset{linecolor=dialinecolor}
\pspolygon*(4.540640,7.301670)(4.540640,8.135000)(5.110640,8.135000)(5.110640,7.301670)
\psset{linewidth=0.04} \psset{linestyle=solid}
\psset{linestyle=solid} \setlinejoinmode{0}
\newrgbcolor{dialinecolor}{1 1 1}
\psset{linecolor=dialinecolor}
\pspolygon(4.540640,7.301670)(4.540640,8.135000)(5.110640,8.135000)(5.110640,7.301670)
\psset{linewidth=0.03} \psset{linestyle=solid}
\psset{linestyle=solid} \setlinecaps{0} \setlinejoinmode{0}
\setlinecaps{0} \setlinejoinmode{0} \psset{linestyle=solid}
\newrgbcolor{dialinecolor}{1 1 1}
\psset{linecolor=dialinecolor}
\psellipse*(9.440005,6.823335)(0.945835,0.945835)
\newrgbcolor{dialinecolor}{0 0 0}
\psset{linecolor=dialinecolor}
\psellipse(9.440005,6.823335)(0.945835,0.945835)
\psset{linewidth=0.003000} \setlinecaps{0} \setlinejoinmode{0}
\psset{linestyle=solid}
\newrgbcolor{dialinecolor}{0 0 0}
\psset{linecolor=dialinecolor}
\psellipse(9.440005,6.823335)(0.945835,0.945835)
\newrgbcolor{dialinecolor}{1 1 1}
\psset{linecolor=dialinecolor}
\pspolygon*(10.219200,6.464170)(10.219200,6.797500)(10.372530,6.797500)(10.372530,6.464170)
\psset{linewidth=0.04} \psset{linestyle=solid}
\psset{linestyle=solid} \setlinejoinmode{0}
\newrgbcolor{dialinecolor}{1 1 1}
\psset{linecolor=dialinecolor}
\pspolygon(10.219200,6.464170)(10.219200,6.797500)(10.372530,6.797500)(10.372530,6.464170)
\psset{linewidth=0.03} \psset{linestyle=solid}
\psset{linestyle=solid} \setlinecaps{0} \setlinejoinmode{0}
\setlinecaps{0} \setlinejoinmode{0} \psset{linestyle=solid}
\newrgbcolor{dialinecolor}{1 1 1}
\psset{linecolor=dialinecolor}
\psellipse*(8.172500,3.572500)(0.895000,0.895000)
\newrgbcolor{dialinecolor}{0 0 0}
\psset{linecolor=dialinecolor}
\psellipse(8.172500,3.572500)(0.895000,0.895000)
\psset{linewidth=0.003000} \setlinecaps{0} \setlinejoinmode{0}
\psset{linestyle=solid}
\newrgbcolor{dialinecolor}{0 0 0}
\psset{linecolor=dialinecolor}
\psellipse(8.172500,3.572500)(0.895000,0.895000)

\newrgbcolor{dialinecolor}{0 0 0}
\psset{linecolor=dialinecolor}
\rput(9.637500,6.214170){\scalebox{1 -1}{\,\,$c^5_4$}}
\newrgbcolor{dialinecolor}{1 1 1}
\psset{linecolor=dialinecolor}
\psellipse*(9.446670,6.741255)(0.357500,0.389585)
\psset{linewidth=0.03} \psset{linestyle=solid}
\psset{linestyle=solid}
\newrgbcolor{dialinecolor}{0 0 0}
\psset{linecolor=dialinecolor}
\psellipse(9.446670,6.741255)(0.357500,0.389585)
\newrgbcolor{dialinecolor}{0 0 0}
\psset{linecolor=dialinecolor}
\rput(9.441670,6.697500){\scalebox{1 -1}{5}}
\newrgbcolor{dialinecolor}{0 0 0}
\psset{linecolor=dialinecolor}
\rput(9.458330,7.514170){\scalebox{1 -1}{$c^5_7$}}
\newrgbcolor{dialinecolor}{0 0 0}
\psset{linecolor=dialinecolor}
\rput(8.975000,7.230830){\scalebox{1 -1}{$c^5_1$}}
\newrgbcolor{dialinecolor}{0 0 0}
\psset{linecolor=dialinecolor}
\rput(8.875000,6.614170){\scalebox{1 -1}{$c^5_2$\,\,}}
\newrgbcolor{dialinecolor}{0 0 0}
\psset{linecolor=dialinecolor}
\rput(9.920830,6.608330){\scalebox{1 -1}{\,\,\,\,$c^5_5$}}
\newrgbcolor{dialinecolor}{0 0 0}
\psset{linecolor=dialinecolor}
\rput(9.854170,7.230830){\scalebox{1 -1}{$c^5_6$}}
\newrgbcolor{dialinecolor}{0 0 0}
\psset{linecolor=dialinecolor}
\rput(8.154170,2.947500){\scalebox{1 -1}{$c^4_3$}}
\newrgbcolor{dialinecolor}{1 1 1}
\psset{linecolor=dialinecolor}
\psellipse*(8.163330,3.491255)(0.357500,0.389585)
\psset{linewidth=0.03} \psset{linestyle=solid}
\psset{linestyle=solid}
\newrgbcolor{dialinecolor}{0 0 0}
\psset{linecolor=dialinecolor}
\psellipse(8.163330,3.491255)(0.357500,0.389585)
\newrgbcolor{dialinecolor}{0 0 0}
\psset{linecolor=dialinecolor}
\rput(8.140250,3.418330){\scalebox{1 -1}{4}}
\newrgbcolor{dialinecolor}{0 0 0}
\psset{linecolor=dialinecolor}
\rput(7.758330,3.980830){\scalebox{1 -1}{$c^4_1$\,\,}}
\newrgbcolor{dialinecolor}{0 0 0}
\psset{linecolor=dialinecolor}
\rput(7.658330,3.297500){\scalebox{1 -1}{$c^4_2$\,\,}}
\newrgbcolor{dialinecolor}{0 0 0}
\psset{linecolor=dialinecolor}
\rput(8.620830,3.325000){\scalebox{1 -1}{\,\,\,\,$c^4_4$}}
\newrgbcolor{dialinecolor}{0 0 0}
\psset{linecolor=dialinecolor}
\rput(8.558330,3.947500){\scalebox{1 -1}{$c^4_5$}}
\newrgbcolor{dialinecolor}{0 0 0}
\psset{linecolor=dialinecolor}
\rput(5.304170,6.858330){\scalebox{1 -1}{$c^2_4$\,}}
\newrgbcolor{dialinecolor}{0 0 0}
\psset{linecolor=dialinecolor}
\rput(5.775000,6.864170){\scalebox{1 -1}{$c^2_5$}}
\newrgbcolor{dialinecolor}{1 1 1}
\psset{linecolor=dialinecolor}
\psellipse*(5.546670,7.441255)(0.357500,0.389585)
\psset{linewidth=0.03} \psset{linestyle=solid}
\psset{linestyle=solid}
\newrgbcolor{dialinecolor}{0 0 0}
\psset{linecolor=dialinecolor}
\psellipse(5.546670,7.441255)(0.357500,0.389585)
\newrgbcolor{dialinecolor}{0 0 0}
\psset{linecolor=dialinecolor}
\rput(5.540250,7.401670){\scalebox{1 -1}{2}}
\newrgbcolor{dialinecolor}{0 0 0}
\psset{linecolor=dialinecolor} \rput(5.591670,8.13){\scalebox{1
-1}{$c^2_9$}}
\newrgbcolor{dialinecolor}{0 0 0}
\psset{linecolor=dialinecolor}
\rput(5.220830,8.047500){\scalebox{1 -1}{$c^2_1$}}
\newrgbcolor{dialinecolor}{0 0 0}
\psset{linecolor=dialinecolor}
\rput(5.041670,7.197500){\scalebox{1 -1}{$c^2_3$\,\,}}
\newrgbcolor{dialinecolor}{0 0 0}
\psset{linecolor=dialinecolor}
\rput(6.070830,7.225000){\scalebox{1 -1}{\,\,$c^2_6$}}
\newrgbcolor{dialinecolor}{0 0 0}
\psset{linecolor=dialinecolor}
\rput(5.920830,8.030830){\scalebox{1 -1}{$c^2_8$}}
\newrgbcolor{dialinecolor}{0 0 0}
\psset{linecolor=dialinecolor}
\rput(4.937500,7.647500){\scalebox{1 -1}{$c^2_2$}}
\newrgbcolor{dialinecolor}{0 0 0}
\psset{linecolor=dialinecolor}
\rput(3.970830,2.897500){\scalebox{1 -1}{$c^3_4$}}
\newrgbcolor{dialinecolor}{0 0 0}
\psset{linecolor=dialinecolor}
\rput(4.320830,2.780830){\scalebox{1 -1}{$c^3_5$}}
\newrgbcolor{dialinecolor}{1 1 1}
\psset{linecolor=dialinecolor}
\psellipse*(4.342500,3.424585)(0.357500,0.389585)
\psset{linewidth=0.03} \psset{linestyle=solid}
\psset{linestyle=solid}
\newrgbcolor{dialinecolor}{0 0 0}
\psset{linecolor=dialinecolor}
\psellipse(4.342500,3.424585)(0.357500,0.389585)
\newrgbcolor{dialinecolor}{0 0 0}
\psset{linecolor=dialinecolor}
\rput(4.340250,3.385000){\scalebox{1 -1}{3}}
\newrgbcolor{dialinecolor}{0 0 0}
\psset{linecolor=dialinecolor} \rput(4.387500,4.1170){\scalebox{1
-1}{$c^3_{10}$}}
\newrgbcolor{dialinecolor}{0 0 0}
\psset{linecolor=dialinecolor}
\rput(4.037500,4.064170){\scalebox{1 -1}{$c^3_1$\,\,}}
\newrgbcolor{dialinecolor}{0 0 0}
\psset{linecolor=dialinecolor}
\rput(3.770830,3.230830){\scalebox{1 -1}{$c^3_3$}}
\newrgbcolor{dialinecolor}{0 0 0}
\psset{linecolor=dialinecolor} \rput(4.9,3.308330){\scalebox{1
-1}{$c^3_7$}}
\newrgbcolor{dialinecolor}{0 0 0}
\psset{linecolor=dialinecolor}
\rput(4.954170,3.714170){\scalebox{1 -1}{$c^3_8$}}
\newrgbcolor{dialinecolor}{0 0 0}
\psset{linecolor=dialinecolor}
\rput(4.712500,2.912500){\scalebox{1 -1}{$c^3_6$}}
\newrgbcolor{dialinecolor}{0 0 0}
\psset{linecolor=dialinecolor}
\rput(4.695830,4.062500){\scalebox{1 -1}{$c^3_9$}}
\newrgbcolor{dialinecolor}{0 0 0}
\psset{linecolor=dialinecolor}
\rput(3.745830,3.679170){\scalebox{1 -1}{$c^3_2$}}
\newrgbcolor{dialinecolor}{1 1 1}
\psset{linecolor=dialinecolor}
\pspolygon*(5.137500,3.580830)(5.137500,3.914160)(5.354170,3.914160)(5.354170,3.580830)
\psset{linewidth=0.04} \psset{linestyle=solid}
\psset{linestyle=solid} \setlinejoinmode{0}
\newrgbcolor{dialinecolor}{1 1 1}
\psset{linecolor=dialinecolor}
\pspolygon(5.137500,3.580830)(5.137500,3.914160)(5.354170,3.914160)(5.354170,3.580830)
\newrgbcolor{dialinecolor}{1 1 1}
\psset{linecolor=dialinecolor}
\pspolygon*(8.832500,5.897500)(8.832500,6.130830)(9.187500,6.130830)(9.187500,5.897500)
\psset{linewidth=0.04} \psset{linestyle=solid}
\psset{linestyle=solid} \setlinejoinmode{0}
\newrgbcolor{dialinecolor}{1 1 1}
\psset{linecolor=dialinecolor}
\pspolygon(8.832500,5.897500)(8.832500,6.130830)(9.187500,6.130830)(9.187500,5.897500)
\psset{linewidth=0.03} \psset{linestyle=dashed,dash=1 1}
\psset{linestyle=dashed,dash=0.1 0.1} \setlinecaps{0}
\newrgbcolor{dialinecolor}{0 0 0}
\psset{linecolor=dialinecolor}
\psline(5.137500,3.747500)(9.037500,6.014170)
\psset{linewidth=0.03} \psset{linestyle=solid}
\psset{linestyle=solid} \setlinecaps{0}
\newrgbcolor{dialinecolor}{0 0 0}
\psset{linecolor=dialinecolor}
\psline(5.137500,3.914170)(8.832500,6.130830)
\psset{linewidth=0.03} \psset{linestyle=solid}
\psset{linestyle=solid} \setlinecaps{0}
\newrgbcolor{dialinecolor}{0 0 0}
\psset{linecolor=dialinecolor}
\psline(5.245830,3.580830)(7.570830,4.964170)
\newrgbcolor{dialinecolor}{1 1 1}
\psset{linecolor=dialinecolor}
\pspolygon*(7.987500,4.330830)(7.987500,4.530830)(8.354170,4.530830)(8.354170,4.330830)
\psset{linewidth=0.04} \psset{linestyle=solid}
\psset{linestyle=solid} \setlinejoinmode{0}
\newrgbcolor{dialinecolor}{1 1 1}
\psset{linecolor=dialinecolor}
\pspolygon(7.987500,4.330830)(7.987500,4.530830)(8.354170,4.530830)(8.354170,4.330830)
\psset{linewidth=0.03} \psset{linestyle=dashed,dash=1 1}
\psset{linestyle=dashed,dash=0.1 0.1} \setlinecaps{0}
\newrgbcolor{dialinecolor}{0 0 0}
\psset{linecolor=dialinecolor}
\psline(7.620830,5.180830)(8.172500,4.467500)
\psset{linewidth=0.03} \psset{linestyle=solid}
\psset{linestyle=solid} \setlinecaps{0}
\newrgbcolor{dialinecolor}{0 0 0}
\psset{linecolor=dialinecolor}
\psline(7.554170,4.980830)(7.987500,4.430830)
\psset{linewidth=0.03} \psset{linestyle=solid}
\psset{linestyle=solid} \setlinecaps{0}
\newrgbcolor{dialinecolor}{0 0 0}
\psset{linecolor=dialinecolor}
\psline(7.887500,5.130830)(8.354170,4.430830)
\psset{linewidth=0.03} \psset{linestyle=solid}
\psset{linestyle=solid} \setlinecaps{0}
\newrgbcolor{dialinecolor}{0 0 0}
\psset{linecolor=dialinecolor}
\psline(9.187500,5.897500)(7.920830,5.164170)

\newrgbcolor{dialinecolor}{0 0 0}
\psset{linecolor=dialinecolor}
\rput(9.220830,6.214170){\scalebox{1 -1}{$c^5_3$\,\,}}
\newrgbcolor{dialinecolor}{0 0 0}
\psset{linecolor=dialinecolor}
\rput(8.208330,4.164170){\scalebox{1 -1}{$c^4_6$\,\,}}
\newrgbcolor{dialinecolor}{1 1 1}
\psset{linecolor=dialinecolor}
\pspolygon*(6.352500,7.492500)(6.352500,7.797500)(6.655830,7.797500)(6.655830,7.492500)
\psset{linewidth=0.04} \psset{linestyle=solid}
\psset{linestyle=solid} \setlinejoinmode{0}
\newrgbcolor{dialinecolor}{1 1 1}
\psset{linecolor=dialinecolor}
\pspolygon(6.352500,7.492500)(6.352500,7.797500)(6.655830,7.797500)(6.655830,7.492500)

\newrgbcolor{dialinecolor}{0 0 0}
\psset{linecolor=dialinecolor}
\rput(6.170830,7.664170){\scalebox{1 -1}{$c^2_7$}}
\psset{linewidth=0.03} \psset{linestyle=solid}
\psset{linestyle=solid} \setlinejoinmode{0} \setlinecaps{0}
\newrgbcolor{dialinecolor}{0 0 0}
\psset{linecolor=dialinecolor} \pscustom{
\newpath
\moveto(6.519170,7.456670)
\curveto(10.452500,9.54)(11.502500,6.756670)(10.385800,6.823340)
\stroke} \psset{linewidth=0.03} \psset{linestyle=dashed,dash=1 1}
\psset{linestyle=dashed,dash=0.1 0.1} \setlinejoinmode{0}
\setlinecaps{0}
\newrgbcolor{dialinecolor}{0 0 0}
\psset{linecolor=dialinecolor} \pscustom{
\newpath
\moveto(6.339170,7.614170)
\curveto(10.305800,9.914170)(12.172500,6.930830)(10.255800,6.580830)
\stroke} \psset{linewidth=0.03} \psset{linestyle=solid}
\psset{linestyle=solid} \setlinejoinmode{0} \setlinecaps{0}
\newrgbcolor{dialinecolor}{0 0 0}
\psset{linecolor=dialinecolor} \pscustom{
\newpath
\moveto(6.439170,7.830830)
\curveto(10.439200,10.297500)(12.772500,6.847500)(10.272500,6.447500)
\stroke} \psset{linewidth=0.03} \psset{linestyle=dashed,dash=1 1}
\psset{linestyle=dashed,dash=0.1 0.1} \setlinecaps{0}
\newrgbcolor{dialinecolor}{0 0 0}
\psset{linecolor=dialinecolor}
\psclip{\pswedge[linestyle=none,fillstyle=none](7.925752,3.728433){5.089093}{128.427664}{170.901133}}
\psellipse(7.925752,3.728433)(3.598532,3.598532)
\endpsclip
\psset{linewidth=0.03} \psset{linestyle=dashed,dash=0.1 0.1}
\psset{linestyle=dashed,dash=0.1 0.1} \setlinejoinmode{0}
\setlinecaps{0}
\newrgbcolor{dialinecolor}{0 0 0}
\psset{linecolor=dialinecolor} \pscustom{
\newpath
\moveto(4.755830,7.647500)
\curveto(4.522500,12.047500)(0.122500,8.014170)(4.739170,7.647500)
\stroke} \psset{linewidth=0.03} \psset{linestyle=dashed,dash=0.1
0.1} \psset{linestyle=dashed,dash=0.1 0.1} \setlinecaps{0}
\newrgbcolor{dialinecolor}{0 0 0}
\psset{linecolor=dialinecolor}
\psline(2.872500,9.347500)(4.739170,7.680830)

\newrgbcolor{dialinecolor}{0 0 0}
\psset{linecolor=dialinecolor}
\rput(0.987500,4.658330){\scalebox{1 -1}{$c^1_3$\,\,\,\,}}
\newrgbcolor{dialinecolor}{0 0 0}
\psset{linecolor=dialinecolor}
\rput(1.458330,4.664170){\scalebox{1 -1}{$c^1_4$}}
\newrgbcolor{dialinecolor}{1 1 1}
\psset{linecolor=dialinecolor}
\psellipse*(1.23,5.224585)(0.357500,0.389585)
\psset{linewidth=0.03} \psset{linestyle=solid}
\psset{linestyle=solid}
\newrgbcolor{dialinecolor}{0 0 0}
\psset{linecolor=dialinecolor}
\psellipse(1.23,5.224585)(0.357500,0.389585)

\newrgbcolor{dialinecolor}{0 0 0}
\psset{linecolor=dialinecolor}
\rput(1.223580,5.185000){\scalebox{1 -1}{1}}
\newrgbcolor{dialinecolor}{0 0 0}
\psset{linecolor=dialinecolor}
\rput(1.275000,5.914170){\scalebox{1 -1}{$c^1_7$}}
\newrgbcolor{dialinecolor}{0 0 0}
\psset{linecolor=dialinecolor}
\rput(0.825000,5.730830){\scalebox{1 -1}{$c^1_1$\,\,}}
\newrgbcolor{dialinecolor}{0 0 0}
\psset{linecolor=dialinecolor}
\rput(0.675000,5.130830){\scalebox{1 -1}{$c^1_2$\,\,}}
\newrgbcolor{dialinecolor}{0 0 0}
\psset{linecolor=dialinecolor}
\rput(1.720830,5.125000){\scalebox{1 -1}{\,\,\,\,$c^1_5$}}
\newrgbcolor{dialinecolor}{0 0 0}
\psset{linecolor=dialinecolor}
\rput(1.641670,5.730830){\scalebox{1 -1}{$c^1_6$}}
\psset{linewidth=0.03} \psset{linestyle=dashed,dash=1 1}
\psset{linestyle=dashed,dash=0.1 0.1} \setlinecaps{0}
\newrgbcolor{dialinecolor}{0 0 0}
\psset{linecolor=dialinecolor}
\psclip{\pswedge[linestyle=none,fillstyle=none](2.453166,10.709842){8.085545}{265.009691}{294.665967}}
\psellipse(2.453166,10.709842)(5.717344,5.717344)
\endpsclip
\psset{linewidth=0.03} \psset{linestyle=solid}
\psset{linestyle=solid} \setlinecaps{0}
\newrgbcolor{dialinecolor}{0 0 0}
\psset{linecolor=dialinecolor}
\psclip{\pswedge[linestyle=none,fillstyle=none](8.151628,3.775899){5.036551}{124.407721}{170.492272}}
\psellipse(8.151628,3.775899)(3.561379,3.561379)
\endpsclip
\psset{linewidth=0.03} \psset{linestyle=solid}
\psset{linestyle=solid} \setlinecaps{0}
\newrgbcolor{dialinecolor}{0 0 0}
\psset{linecolor=dialinecolor}
\psclip{\pswedge[linestyle=none,fillstyle=none](6.949547,3.953401){3.930701}{153.793184}{170.805819}}
\psellipse(6.949547,3.953401)(2.779425,2.779425)
\endpsclip
\psset{linewidth=0.03} \psset{linestyle=solid}
\psset{linestyle=solid} \setlinecaps{0}
\newrgbcolor{dialinecolor}{0 0 0}
\psset{linecolor=dialinecolor}
\psclip{\pswedge[linestyle=none,fillstyle=none](6.639493,4.868748){2.988645}{127.403456}{159.827841}}
\psellipse(6.639493,4.868748)(2.113291,2.113291)
\endpsclip
\psset{linewidth=0.03} \psset{linestyle=solid}
\psset{linestyle=solid} \setlinecaps{0}
\newrgbcolor{dialinecolor}{0 0 0}
\psset{linecolor=dialinecolor}
\psclip{\pswedge[linestyle=none,fillstyle=none](2.527483,10.652219){8.234291}{265.518228}{289.514811}}
\psellipse(2.527483,10.652219)(5.822523,5.822523)
\endpsclip
\psset{linewidth=0.03} \psset{linestyle=solid}
\psset{linestyle=solid} \setlinecaps{0}
\newrgbcolor{dialinecolor}{0 0 0}
\psset{linecolor=dialinecolor}
\psclip{\pswedge[linestyle=none,fillstyle=none](2.429666,10.796743){7.973960}{267.385861}{293.623946}}
\psellipse(2.429666,10.796743)(5.638441,5.638441)
\endpsclip
\psset{linewidth=0.03} \psset{linestyle=solid}
\psset{linestyle=solid} \setlinejoinmode{0} \setlinecaps{0}
\newrgbcolor{dialinecolor}{0 0 0}
\psset{linecolor=dialinecolor} \pscustom{
\newpath
\moveto(4.950490,8.164180)
\curveto(4.267160,12.297500)(-0.166170,8.047510)(4.300490,7.497510)
\stroke} \psset{linewidth=0.03} \psset{linestyle=solid}
\psset{linestyle=solid} \setlinejoinmode{0} \setlinecaps{0}
\newrgbcolor{dialinecolor}{0 0 0}
\psset{linecolor=dialinecolor} \pscustom{
\newpath
\moveto(2.943970,9.018330)
\curveto(2.627300,8.885000)(3.027300,8.051670)(4.138080,7.926720)
\stroke} \psset{linewidth=0.03} \psset{linestyle=solid}
\psset{linestyle=solid} \setlinejoinmode{0} \setlinecaps{0}
\newrgbcolor{dialinecolor}{0 0 0}
\psset{linecolor=dialinecolor} \pscustom{
\newpath
\moveto(3.240640,9.318190)
\curveto(3.443970,9.768330)(4.477300,9.118330)(4.434750,8.226580)
\stroke} \psset{linewidth=0.03} \psset{linestyle=solid}
\psset{linestyle=solid} \setlinecaps{0}
\newrgbcolor{dialinecolor}{0 0 0}
\psset{linecolor=dialinecolor}
\psline(2.927300,9.001670)(4.093970,7.968330)
\psset{linewidth=0.03} \psset{linestyle=solid}
\psset{linestyle=solid} \setlinecaps{0}
\newrgbcolor{dialinecolor}{0 0 0}
\psset{linecolor=dialinecolor}
\psline(3.243970,9.318330)(4.410640,8.235000)

\newrgbcolor{dialinecolor}{0 0 0}
\psset{linecolor=dialinecolor}
\rput(3.277300,8.501670){\scalebox{1 -1}{2'}}
\newrgbcolor{dialinecolor}{0 0 0}
\psset{linecolor=dialinecolor}
\rput(3.877300,9.018330){\scalebox{1 -1}{3'}}
\newrgbcolor{dialinecolor}{1 1 1}
\psset{linecolor=dialinecolor}
\pspolygon*(3.373970,4.868330)(3.373970,5.435000)(3.893970,5.435000)(3.893970,4.868330)
\psset{linewidth=0.04} \psset{linestyle=solid}
\psset{linestyle=solid} \setlinejoinmode{0}
\newrgbcolor{dialinecolor}{1 1 1}
\psset{linecolor=dialinecolor}
\pspolygon(3.373970,4.868330)(3.373970,5.435000)(3.893970,5.435000)(3.893970,4.868330)
\psset{linewidth=0.03} \psset{linestyle=solid}
\psset{linestyle=solid} \setlinecaps{0}
\newrgbcolor{dialinecolor}{0 0 0}
\psset{linecolor=dialinecolor}
\psclip{\pswedge[linestyle=none,fillstyle=none](9.545231,4.135415){8.234685}{147.108348}{179.717054}}
\psellipse(9.545231,4.135415)(5.822802,5.822802)
\endpsclip
\psset{linewidth=0.03} \psset{linestyle=solid}
\psset{linestyle=solid} \setlinecaps{0}
\newrgbcolor{dialinecolor}{0 0 0}
\psset{linecolor=dialinecolor}
\psclip{\pswedge[linestyle=none,fillstyle=none](9.909761,4.161467){9.206963}{148.832422}{184.822335}}
\psellipse(9.909761,4.161467)(6.510306,6.510306)
\endpsclip
\psset{linewidth=0.03} \psset{linestyle=dashed,dash=1 1}
\psset{linestyle=dashed,dash=0.1 0.1} \setlinecaps{0}
\newrgbcolor{dialinecolor}{0 0 0}
\psset{linecolor=dialinecolor}
\psclip{\pswedge[linestyle=none,fillstyle=none](9.181647,4.309297){7.915651}{142.532155}{185.546255}}
\psellipse(9.181647,4.309297)(5.597211,5.597211)
\endpsclip

\newrgbcolor{dialinecolor}{0 0 0}
\psset{linecolor=dialinecolor}
\rput(7.543970,6.658330){\scalebox{1 -1}{1'}}
\newrgbcolor{dialinecolor}{0 0 0}
\psset{linecolor=dialinecolor} \rput(2.806920,6.91){\scalebox{1
-1}{4'}} \psset{linewidth=0.04} \psset{linestyle=solid}
\psset{linestyle=solid} \setlinecaps{0}
\newrgbcolor{dialinecolor}{0 0 0}
\psset{linecolor=dialinecolor}
\psclip{\pswedge[linestyle=none,fillstyle=none](2.598009,5.636098){1.494395}{90.479476}{142.428963}}
\psellipse(2.598009,5.636098)(1.056697,1.056697)
\endpsclip
\setlinejoinmode{0}
\newrgbcolor{dialinecolor}{0 0 0}
\psset{linecolor=dialinecolor}
\pspolygon*(1.828557,6.248954)(1.697556,6.144246)(1.692391,6.311872)
\psset{linewidth=0.04} \psset{linestyle=solid}
\psset{linestyle=solid} \setlinecaps{0}
\newrgbcolor{dialinecolor}{0 0 0}
\psset{linecolor=dialinecolor}
\psclip{\pswedge[linestyle=none,fillstyle=none](8.163237,6.272151){0.982397}{81.421668}{136.346671}}
\psellipse(8.163237,6.272151)(0.694659,0.694659)
\endpsclip
\setlinejoinmode{0}
\newrgbcolor{dialinecolor}{0 0 0}
\psset{linecolor=dialinecolor}
\pspolygon*(8.306872,7.022470)(8.393715,6.879002)(8.226834,6.895608)
\psset{linewidth=0.03} \psset{linestyle=solid}
\psset{linestyle=solid} \setlinecaps{0}
\newrgbcolor{dialinecolor}{0 0 0}
\psset{linecolor=dialinecolor}
\psline(8.827300,7.518330)(8.727438,7.429566)
\psset{linewidth=0.03} \psset{linestyle=solid} \setlinejoinmode{0}
\setlinecaps{0}
\newrgbcolor{dialinecolor}{0 0 0}
\psset{linecolor=dialinecolor}
\psline(8.918288,7.465412)(8.702369,7.407283)(8.785417,7.614896)
\psset{linewidth=0.03} \psset{linestyle=solid}
\psset{linestyle=solid} \setlinecaps{0}
\newrgbcolor{dialinecolor}{0 0 0}
\psset{linecolor=dialinecolor}
\psline(2.093010,5.605350)(2.038252,5.709015)
\psset{linewidth=0.03} \psset{linestyle=solid} \setlinejoinmode{0}
\setlinecaps{0}
\newrgbcolor{dialinecolor}{0 0 0}
\psset{linecolor=dialinecolor}
\psline(2.027577,5.515121)(2.022586,5.738672)(2.204422,5.608535)
\psset{linewidth=0.03} \psset{linestyle=solid}
\psset{linestyle=solid} \setlinecaps{0}
\newrgbcolor{dialinecolor}{0 0 0}
\psset{linecolor=dialinecolor}
\psline(4.026350,9.172020)(4.090645,9.113489)
\psset{linewidth=0.03} \psset{linestyle=solid} \setlinejoinmode{0}
\setlinecaps{0}
\newrgbcolor{dialinecolor}{0 0 0}
\psset{linecolor=dialinecolor}
\psline(4.034871,9.299494)(4.115447,9.090909)(3.900234,9.151599)
\psset{linewidth=0.03} \psset{linestyle=solid}
\psset{linestyle=solid} \setlinecaps{0}
\newrgbcolor{dialinecolor}{0 0 0}
\psset{linecolor=dialinecolor}
\psline(3.656920,8.351670)(3.722649,8.295328)
\psset{linewidth=0.03} \psset{linestyle=solid} \setlinejoinmode{0}
\setlinecaps{0}
\newrgbcolor{dialinecolor}{0 0 0}
\psset{linecolor=dialinecolor}
\psline(3.661348,8.479585)(3.748114,8.273499)(3.531185,8.327738)
\psset{linewidth=0.03} \psset{linestyle=solid}
\psset{linestyle=solid} \setlinecaps{0}
\newrgbcolor{dialinecolor}{0 0 0}
\psset{linecolor=dialinecolor}
\psline(0.977210,5.500060)(0.924167,5.409546)
\psset{linewidth=0.03} \psset{linestyle=solid} \setlinejoinmode{0}
\setlinecaps{0}
\newrgbcolor{dialinecolor}{0 0 0}
\psset{linecolor=dialinecolor}
\psline(1.094605,5.502601)(0.907208,5.380608)(0.922052,5.603722)
\psset{linewidth=0.03} \psset{linestyle=solid}
\psset{linestyle=solid} \setlinecaps{0}
\newrgbcolor{dialinecolor}{0 0 0}
\psset{linecolor=dialinecolor}
\psline(4.089710,3.700060)(4.042616,3.625126)
\psset{linewidth=0.03} \psset{linestyle=solid} \setlinejoinmode{0}
\setlinecaps{0}
\newrgbcolor{dialinecolor}{0 0 0}
\psset{linecolor=dialinecolor}
\psline(4.215858,3.712851)(4.024768,3.596728)(4.046524,3.819274)
\psset{linewidth=0.03} \psset{linestyle=solid}
\psset{linestyle=solid} \setlinecaps{0}
\newrgbcolor{dialinecolor}{0 0 0}
\psset{linecolor=dialinecolor}
\psline(7.910540,3.766730)(7.871830,3.694184)
\psset{linewidth=0.03} \psset{linestyle=solid} \setlinejoinmode{0}
\setlinecaps{0}
\newrgbcolor{dialinecolor}{0 0 0}
\psset{linecolor=dialinecolor}
\psline(8.038419,3.793967)(7.856040,3.664592)(7.861967,3.888120)
\psset{linewidth=0.03} \psset{linestyle=solid}
\psset{linestyle=solid} \setlinecaps{0}
\newrgbcolor{dialinecolor}{0 0 0}
\psset{linecolor=dialinecolor}
\psline(5.293880,7.716730)(5.262576,7.672913)
\psset{linewidth=0.03} \psset{linestyle=solid} \setlinejoinmode{0}
\setlinecaps{0}
\newrgbcolor{dialinecolor}{0 0 0}
\psset{linecolor=dialinecolor}
\psline(5.440709,7.750225)(5.243078,7.645622)(5.277973,7.866489)
\psset{linewidth=0.03} \psset{linestyle=solid}
\psset{linestyle=solid} \setlinecaps{0}
\newrgbcolor{dialinecolor}{0 0 0}
\psset{linecolor=dialinecolor}
\psline(9.193880,7.016730)(9.147367,6.955187)
\psset{linewidth=0.03} \psset{linestyle=solid} \setlinejoinmode{0}
\setlinecaps{0}
\newrgbcolor{dialinecolor}{0 0 0}
\psset{linecolor=dialinecolor}
\psline(9.327511,7.027690)(9.127143,6.928428)(9.167955,7.148279)
\psset{linewidth=0.1} \psset{linestyle=solid}
\psset{linestyle=solid} \setlinecaps{0}
\newrgbcolor{dialinecolor}{0 0 0}
\psset{linecolor=dialinecolor}
\psline(5.561240,7.619320)(5.546670,7.830840)
\psset{linewidth=0.1} \psset{linestyle=solid}
\psset{linestyle=solid} \setlinecaps{0}
\newrgbcolor{dialinecolor}{0 0 0}
\psset{linecolor=dialinecolor}
\psline(9.444570,6.952660)(9.446670,7.130840)
\psset{linewidth=0.1} \psset{linestyle=solid}
\psset{linestyle=solid} \setlinecaps{0}
\newrgbcolor{dialinecolor}{0 0 0}
\psset{linecolor=dialinecolor}
\psline(1.244570,5.419320)(1.23,5.614170) \psset{linewidth=0.1}
\psset{linestyle=solid} \psset{linestyle=solid} \setlinecaps{0}
\newrgbcolor{dialinecolor}{0 0 0}
\psset{linecolor=dialinecolor}
\psline(8.161240,3.685990)(8.163330,3.880840)
\psset{linewidth=0.1} \psset{linestyle=solid}
\psset{linestyle=solid} \setlinecaps{0}
\newrgbcolor{dialinecolor}{0 0 0}
\psset{linecolor=dialinecolor}
\psline(4.327910,3.619320)(4.342500,3.814170)
\psset{linewidth=0.1} \psset{linestyle=solid}
\psset{linestyle=solid} \setlinecaps{0}
\newrgbcolor{dialinecolor}{0 0 0}
\psset{linecolor=dialinecolor}
\psline(2.123580,5.185000)(2.373580,5.335000)
\psset{linewidth=0.1} \psset{linestyle=solid}
\psset{linestyle=solid} \setlinecaps{0}
\newrgbcolor{dialinecolor}{0 0 0}
\psset{linecolor=dialinecolor}
\psline(9.440010,7.769170)(9.431240,7.952660)
\psset{linewidth=0.1} \psset{linestyle=solid}
\psset{linestyle=solid} \setlinecaps{0}
\newrgbcolor{dialinecolor}{0 0 0}
\psset{linecolor=dialinecolor}
\psline(4.156920,7.918330)(3.856920,8.068330)
\psset{linewidth=0.1} \psset{linestyle=solid}
\psset{linestyle=solid} \setlinecaps{0}
\newrgbcolor{dialinecolor}{0 0 0}
\psset{linecolor=dialinecolor}
\psline(4.290250,8.485000)(4.440250,8.218330)
\psset{linewidth=0.03} \psset{linestyle=dashed,dash=1 1}
\psset{linestyle=dashed,dash=0.1 0.1} \setlinecaps{0}
\newrgbcolor{dialinecolor}{0 0 0}
\psset{linecolor=dialinecolor}
\psclip{\pswedge[linestyle=none,fillstyle=none](2.991440,5.017404){2.944473}{199.375376}{285.389345}}
\psellipse(2.991440,5.017404)(2.082057,2.082057)
\endpsclip
\psset{linewidth=0.03} \psset{linestyle=solid}
\psset{linestyle=solid} \setlinecaps{0}
\newrgbcolor{dialinecolor}{0 0 0}
\psset{linecolor=dialinecolor}
\psclip{\pswedge[linestyle=none,fillstyle=none](2.993429,5.022918){3.142272}{197.178620}{286.575610}}
\psellipse(2.993429,5.022918)(2.221922,2.221922)
\endpsclip
\psset{linewidth=0.03} \psset{linestyle=solid}
\psset{linestyle=solid} \setlinecaps{0}
\newrgbcolor{dialinecolor}{0 0 0}
\psset{linecolor=dialinecolor}
\psclip{\pswedge[linestyle=none,fillstyle=none](2.932894,4.940503){2.561566}{200.230129}{285.292884}}
\psellipse(2.932894,4.940503)(1.811300,1.811300)
\endpsclip
}\endpspicture
\]\label{main-example}


Next, take the sum over all possibilities of placing the output
marked points and the chord endpoints on different $c^i_j$'s,
while respecting the cyclic ordering of the chord endpoints and
output marked points.

\item[(b)]\label{def-alpha-d} If none of the $c^i_j$s has more
than one special point attached to it, then go to the next step.
Otherwise, if some $c^i_j$s have several special points attached
to them, then use the dual of the product $-\cdot -:A\otimes A\to
A$ to pull things apart. More precisely, if there are $r$ such
things coming together at a $c^i_j$, we replace $c^i_j$ by $$
\quad\quad (\Delta\otimes id^{\otimes (r-2)} )\circ ...\circ
(\Delta\otimes id )\circ\Delta(c^i_j )= \sum_{ (c^i_j )} (c^i_j)
'\otimes (c^i_j)'' \otimes ... \otimes (c^i_j)^{(r)} $$ Here we
have used Sweedler's notation $\Delta:A^*\to A^*\otimes A^*$,
$\Delta(c)=\sum_{(c)} c'\otimes c''$. In the example on page
\pageref{main-example}, $c^2_2$ is replaced
by\ifx\setlinejoinmode\undefined
 \newcommand{\setlinejoinmode}[1]{}
\fi \ifx\setlinecaps\undefined
 \newcommand{\setlinecaps}[1]{}
\fi
\ifx\setfont\undefined
 \newcommand{\setfont}[2]{}
\fi 
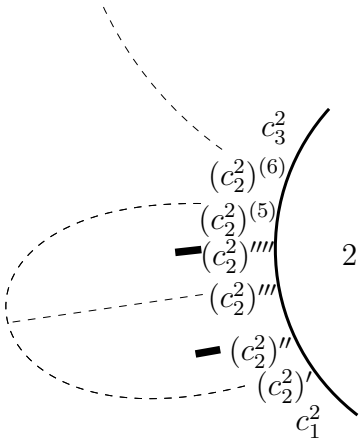
\[\pspicture(2.817621,-9.886767)(7.615442,-4.183467)
\scalebox{1 -1}{
\newrgbcolor{dialinecolor}{0 0 0}
\psset{linecolor=dialinecolor}
\newrgbcolor{diafillcolor}{1 1 1}
\psset{fillcolor=diafillcolor} \psset{linewidth=0.040}
\psset{linestyle=solid} \psset{linestyle=solid} \setlinecaps{0}
\newrgbcolor{dialinecolor}{0 0 0}
\psset{linecolor=dialinecolor}
\psclip{\pswedge[linestyle=none,fillstyle=none](9.192117,7.462849){3.968633}{127.457295}{221.962542}}
\psellipse(9.192117,7.462849)(2.806247,2.806247)
\endpsclip
\newrgbcolor{dialinecolor}{0 0 0}
\psset{linecolor=dialinecolor}
\rput(6.041442,6.479467){\scalebox{1 -1}{$(c^2_2)^{(6)}$}}
\psset{linewidth=0.010} \psset{linestyle=dashed,dash=1 1}
\psset{linestyle=dashed,dash=0.1 0.1} \setlinejoinmode{0}
\setlinecaps{0}
\newrgbcolor{dialinecolor}{0 0 0}
\psset{linecolor=dialinecolor} \pscustom{
\newpath
\moveto(6.001070,9.285710)
\curveto(2.218270,10.250700)(1.554442,6.801567)(5.414442,6.859467)
\stroke} \psset{linewidth=0.010} \psset{linestyle=dashed,dash=0.1
0.1} \psset{linestyle=dashed,dash=0.1 0.1} \setlinecaps{0}
\newrgbcolor{dialinecolor}{0 0 0}
\psset{linecolor=dialinecolor}
\psline(2.913070,8.436510)(5.441370,8.069810)
\psset{linewidth=0.1} \psset{linestyle=solid}
\psset{linestyle=solid} \setlinecaps{0}
\newrgbcolor{dialinecolor}{0 0 0}
\psset{linecolor=dialinecolor}
\psline(5.074670,7.510110)(5.422070,7.471510)
\psset{linewidth=0.1} \psset{linestyle=solid}
\psset{linestyle=solid} \setlinecaps{0}
\newrgbcolor{dialinecolor}{0 0 0}
\psset{linecolor=dialinecolor}
\psline(5.344870,8.861110)(5.672970,8.803210)
\psset{linewidth=0.010} \psset{linestyle=dashed,dash=1 1}
\psset{linestyle=dashed,dash=0.1 0.1} \setlinecaps{0}
\newrgbcolor{dialinecolor}{0 0 0}
\psset{linecolor=dialinecolor}
\psclip{\pswedge[linestyle=none,fillstyle=none](8.545899,2.169205){6.919163}{124.574446}{155.597997}}
\psellipse(8.545899,2.169205)(4.892587,4.892587)
\endpsclip
\newrgbcolor{dialinecolor}{0 0 0}
\psset{linecolor=dialinecolor}
\rput(5.913942,7.063767){\scalebox{1 -1}{$(c^2_2)^{(5)}$}}
\newrgbcolor{dialinecolor}{0 0 0}
\psset{linecolor=dialinecolor}
\rput(5.906942,7.557767){\scalebox{1 -1}{$(c^2_2)''''$}}
\newrgbcolor{dialinecolor}{0 0 0}
\psset{linecolor=dialinecolor}
\rput(5.970942,8.127767){\scalebox{1 -1}{$(c^2_2)'''$}}
\newrgbcolor{dialinecolor}{0 0 0}
\psset{linecolor=dialinecolor} \rput(6.20942,8.830767){\scalebox{1
-1}{$(c^2_2)''$}}
\newrgbcolor{dialinecolor}{0 0 0}
\psset{linecolor=dialinecolor}
\rput(6.526942,9.286767){\scalebox{1 -1}{$(c^2_2)'$}}
\newrgbcolor{dialinecolor}{0 0 0}
\psset{linecolor=dialinecolor}
\rput(6.388942,5.847767){\scalebox{1 -1}{$c^2_3$}}
\newrgbcolor{dialinecolor}{0 0 0}
\psset{linecolor=dialinecolor}
\rput(6.844942,9.761767){\scalebox{1 -1}{$c^2_1$}}
\setfont{Helvetica}{0.8}
\newrgbcolor{dialinecolor}{0 0 0}
\psset{linecolor=dialinecolor}
\rput(7.390442,7.543467){\scalebox{1 -1}{2}} }\endpspicture \]

If an input marked point is also involved, it can be placed
anywhere in between the above linear ordering of chord endpoints
and output marked points. The independence from this choice is
argued in the next step. Now that some $c_j^i$s are expanded,
relabel them so that the last one is again placed at the input
marked point. Let us use the same notation $c^i_{n_i+1}$ for the
new last element.

\item[(c)]\label{def-alpha-a} Note that now input marked points do
not coincide with any output marked points or chord endpoints. We
evaluate $c^i_{n_i+1}$ on the unit, to obtain $c^i_{n_i+1}(1) \in
k$. In the above picture we obtain $c^1_7(1)$ and $c^2_9(1)$. The
ambiguity of where to place the input marked point, which came up
in our previous step, is of no issue because the output marked
points and chord endpoints are linearly ordered, and the input
marked point is evaluated on the unit, $1$, which is in the center
of $A$.

\item[(d)]\label{def-alpha-c} We now deal with the $c^i_j$s which
are placed at the chord endpoints. The cyclic ordering at each vertex
of the chord induces a cyclic ordering of the chord endpoints. In the
example on page \pageref{main-example}, two of the chords have the
endpoints $(c^2_5, c^1_5, c^3_{10})$, and $(c^3_8 , c^4_6, c^5_3)$ up
to cyclic permutation.

We will multiply these elements in this cyclic order and evaluate
it on the unit. To be more precise, if
$c^{i_1}_{j_1},...,c^{i_r}_{j_r}$ are the endpoints of the chord
arranged in the cyclic order, then we obtain the term, $$
<\gamma(c^{i_1}_{j_1}) \cdot ...\cdot \gamma(c^{i_r}_{j_r}) ,1>$$
(see equation \eqref{cyclic-inv} on page \pageref{cyclic-inv}).
Here $\gamma:A^*\to A$ is the inverse of the $A$-bimodule
isomorphism $\beta:A\to A^*$.

\item[(e)]\label{def-alpha-b} For each of the $l$ output circles
of $s$, we look at its marked point. Following the orientation of
the output circle, we linearly read off the leftover $c^i_j$s (the
ones which did not correspond to input marked points or chord
endpoints) so that we end with the element at the output marked
point. For instance, in the example on page
\pageref{main-example}, the $1^{'th}$ output circle gives rise to
the term,
$$(c^5_1, c^5_2, c^3_9, c^2_6, c^5_6; c^5_7)$$
\item[(f)]\label{def-alpha-e} There is an overall sign factor
which is obtained in the following way. Note that, using the
ordering $f_1,...,f_k$, the $c^i_j$s can be linearly ordered,
\begin{equation}\label{f1-fk}
(c^1_1,...,c^1_{n_1};c^1_{n_1+1}),...,(c^k_1,...,c^k_{n_k};c^k_{n_k+1})
\end{equation}
The $c^i_j$'s, for $1\leq j\leq n_i$, are considered to be of degree
$1$, whereas $c^i_{n_i+1}$ is regarded as of degree $0$. Thus,
$f_1\otimes ...\otimes f_k$ has a total degree of $n_1+...+n_k$.
Having this in mind, the operation $\alpha(s)$ can be obtained by the
following two steps.

First, the $c^i_j$'s which correspond to chord endpoints are to be
evaluated using the inner product. These do not contribute to the
output total degree. Similarly, the $c^i_j$'s which correspond to
the output marked points change their degree from $1$ to $0$,
because they are to be positioned as the last entry of a
Hochschild element. The input marked points $c^i_{n_i+1}$ of
degree $0$ are evaluated on the unit, and do not change the total
degree. We see that $\alpha(s)$ changes the degree by the number
of special points on the input circles, minus $k$. This change of
degrees is obtained by applying a tensor product of shift and
identity maps to expression \eqref{f1-fk}, where the shift and
identity maps have degrees $1$ and $0$, respectively. In doing so,
the usual sign rule applies. That is, whenever something of degree
$r$ moves past something of degree $s$, a sign of $(-1)^{r\cdot
s}$ is introduced.

The second step is to rearrange expression \eqref{f1-fk} according to
the combinatorics of the output circles of $s$. This means that
blocks of $c^i_j$'s have to move past other blocks of $c^i_j$'s. We
introduce a sign of $(-1)^{r\cdot s}$ for each block of degree $r$
moving past a block of degree $s$.
\end{itemize}

{\it Step II: $\alpha$ is well-defined}

To ensure that the above procedure yields a well-defined map, the
following checks are in order. We first deal with the fact that
output marked points and chord endpoints may slide along chords.
For example in the chord diagram from page \pageref{main-example},
the $4'^{th}$ marked point may be put at $c^2_5$ instead of
$c^1_5$, while respecting the cyclic ordering of the chord. But we
can check that
\begin{multline*}
\quad\sum_{(c^1_5),(c^3_{10})} <\gamma(c^2_5)\cdot \gamma((c^1_5)'
) \cdot \gamma((c^3_{10})'),1>\cdot (c^3_{10})''(1)\cdot
(c^1_5)''(a)\\=<\gamma(c^2_5)\cdot a\cdot \gamma(c^1_5) \cdot
\gamma(c^3_{10}),1>\\=\sum_{(c^2_5),(c^3_{10})} <\gamma
((c^2_5)'')\cdot \gamma(c^1_5) \cdot \gamma ((c^3_{10})') ,1>\cdot
(c^3_{10})''(1)\cdot (c^2_5)'(a)
\end{multline*}
for all $a\in A$. Also, a chord may slide along another chord,
which does not change the outcome because:
\begin{multline*}
\quad\quad\sum_{(c^k_l)} <...\cdot \gamma(c^i_j)\cdot
\gamma((c^k_l)')\cdot ... ,1>\cdot <... \cdot
\gamma((c^k_l)'')\cdot ... ,1>\\= \sum_{(c^i_j)} <...\cdot
\gamma((c^i_j)'') \cdot \gamma (c^k_l)\cdot ... ,1>\cdot <...
\cdot \gamma((c^i_j)')\cdot ... ,1>
\end{multline*}
\ifx\setlinejoinmode\undefined
 \newcommand{\setlinejoinmode}[1]{}
\fi \ifx\setlinecaps\undefined
 \newcommand{\setlinecaps}[1]{}
\fi \ifx\setfont\undefined
 \newcommand{\setfont}[2]{}
\fi 
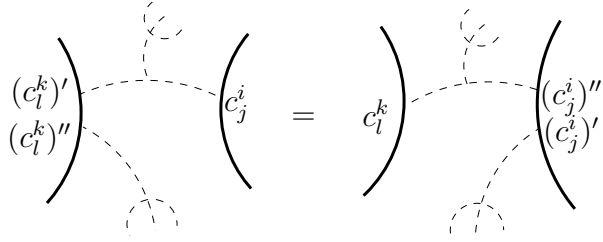
\[\pspicture(5.806030,-6.899360)(12.835434,-3.769099)
\scalebox{1 -1}{
\newrgbcolor{dialinecolor}{0 0 0}
\psset{linecolor=dialinecolor}
\newrgbcolor{diafillcolor}{1 1 1}
\psset{fillcolor=diafillcolor}
\newrgbcolor{dialinecolor}{0 0 0}
\psset{linecolor=dialinecolor}
\rput[r](9.385570,5.323900){\scalebox{1 -1}{=}}
\psset{linewidth=0.040} \psset{linestyle=solid}
\psset{linestyle=solid} \setlinecaps{0}
\newrgbcolor{dialinecolor}{0 0 0}
\psset{linecolor=dialinecolor}
\psclip{\pswedge[linestyle=none,fillstyle=none](14.506146,5.144420){3.106571}{140.383913}{210.885571}}
\psellipse(14.506146,5.144420)(2.196678,2.196678)
\endpsclip

\newrgbcolor{dialinecolor}{0 0 0}
\psset{linecolor=dialinecolor} \rput[r](13.2,5.039700){\scalebox{1
-1}{$(c^i_j)''$}}
\newrgbcolor{dialinecolor}{0 0 0}
\psset{linecolor=dialinecolor}
\rput[r](10.369870,5.304600){\scalebox{1 -1}{$c^k_l$}}

\newrgbcolor{dialinecolor}{0 0 0}
\psset{linecolor=dialinecolor}
\rput[r](13.16,5.536200){\scalebox{1 -1}{$(c^i_j)'$}}
\psset{linewidth=0.040} \psset{linestyle=solid}
\psset{linestyle=solid} \setlinecaps{0}
\newrgbcolor{dialinecolor}{0 0 0}
\psset{linecolor=dialinecolor}
\psclip{\pswedge[linestyle=none,fillstyle=none](8.813547,5.096742){2.496463}{325.386735}{46.226718}}
\psellipse(8.813547,5.096742)(1.765266,1.765266)
\endpsclip
\psset{linewidth=0.010} \psset{linestyle=dashed,dash=1 1}
\psset{linestyle=dashed,dash=0.1 0.1} \setlinecaps{0}
\newrgbcolor{dialinecolor}{0 0 0}
\psset{linecolor=dialinecolor}
\psclip{\pswedge[linestyle=none,fillstyle=none](11.678593,6.909514){2.930051}{237.233607}{288.367855}}
\psellipse(11.678593,6.909514)(2.071859,2.071859)
\endpsclip
\psset{linewidth=0.010} \psset{linestyle=dashed,dash=0.1 0.1}
\psset{linestyle=dashed,dash=0.1 0.1} \setlinecaps{0}
\newrgbcolor{dialinecolor}{0 0 0}
\psset{linecolor=dialinecolor}
\psclip{\pswedge[linestyle=none,fillstyle=none](13.315753,6.997535){2.569938}{183.313211}{237.977622}}
\psellipse(13.315753,6.997535)(1.817221,1.817221)
\endpsclip
\psset{linewidth=0.010} \psset{linestyle=dashed,dash=0.1 0.1}
\psset{linestyle=dashed,dash=0.1 0.1} \setlinecaps{0}
\newrgbcolor{dialinecolor}{0 0 0}
\psset{linecolor=dialinecolor}
\psclip{\pswedge[linestyle=none,fillstyle=none](12.041464,4.629594){0.980268}{158.114959}{237.306627}}
\psellipse(12.041464,4.629594)(0.693154,0.693154)
\endpsclip
\psset{linewidth=0.010} \psset{linestyle=dashed,dash=0.1 0.1}
\psset{linestyle=dashed,dash=0.1 0.1} \setlinecaps{0}
\newrgbcolor{dialinecolor}{0 0 0}
\psset{linecolor=dialinecolor}
\psclip{\pswedge[linestyle=none,fillstyle=none](11.527954,6.812613){0.508588}{170.299433}{0.429640}}
\psellipse(11.527954,6.812613)(0.359626,0.359626)
\endpsclip
\psset{linewidth=0.010} \psset{linestyle=dashed,dash=0.1 0.1}
\psset{linestyle=dashed,dash=0.1 0.1} \setlinecaps{0}
\newrgbcolor{dialinecolor}{0 0 0}
\psset{linecolor=dialinecolor}
\psclip{\pswedge[linestyle=none,fillstyle=none](11.586503,4.093761){0.391445}{16.393563}{232.821078}}
\psellipse(11.586503,4.093761)(0.276793,0.276793)
\endpsclip
\newrgbcolor{dialinecolor}{0 0 0}
\psset{linecolor=dialinecolor} \rput[r](8.5,5.136210){\scalebox{1
-1}{$c^i_j$}}
\newrgbcolor{dialinecolor}{0 0 0}
\psset{linecolor=dialinecolor}
\rput[r](6.104570,4.976500){\scalebox{1 -1}{$(c^k_l)'$}}
\newrgbcolor{dialinecolor}{0 0 0}
\psset{linecolor=dialinecolor}
\rput[r](6.143170,5.588310){\scalebox{1 -1}{$(c^k_l)''$}}
\psset{linewidth=0.010} \psset{linestyle=dashed,dash=0.1 0.1}
\psset{linestyle=dashed,dash=0.1 0.1} \setlinecaps{0}
\newrgbcolor{dialinecolor}{0 0 0}
\psset{linecolor=dialinecolor}
\psclip{\pswedge[linestyle=none,fillstyle=none](7.136408,7.034582){3.135997}{246.493854}{295.893125}}
\psellipse(7.136408,7.034582)(2.217485,2.217485)
\endpsclip
\psset{linewidth=0.040} \psset{linestyle=solid}
\psset{linestyle=solid} \setlinecaps{0}
\newrgbcolor{dialinecolor}{0 0 0}
\psset{linecolor=dialinecolor}
\psclip{\pswedge[linestyle=none,fillstyle=none](9.700854,5.210301){2.258920}{138.202000}{219.676177}}
\psellipse(9.700854,5.210301)(1.597297,1.597297)
\endpsclip
\psset{linewidth=0.040} \psset{linestyle=solid}
\psset{linestyle=solid} \setlinecaps{0}
\newrgbcolor{dialinecolor}{0 0 0}
\psset{linecolor=dialinecolor}
\psclip{\pswedge[linestyle=none,fillstyle=none](4.454388,5.257278){2.588711}{326.554412}{41.404353}}
\psellipse(4.454388,5.257278)(1.830495,1.830495)
\endpsclip
\psset{linewidth=0.010} \psset{linestyle=dashed,dash=1 1}
\psset{linestyle=dashed,dash=0.1 0.1} \setlinecaps{0}
\newrgbcolor{dialinecolor}{0 0 0}
\psset{linecolor=dialinecolor}
\psclip{\pswedge[linestyle=none,fillstyle=none](5.429606,7.066497){2.606622}{297.896327}{352.167311}}
\psellipse(5.429606,7.066497)(1.843160,1.843160)
\endpsclip
\psset{linewidth=0.010} \psset{linestyle=dashed,dash=0.1 0.1}
\psset{linestyle=dashed,dash=0.1 0.1} \setlinecaps{0}
\newrgbcolor{dialinecolor}{0 0 0}
\psset{linecolor=dialinecolor}
\psclip{\pswedge[linestyle=none,fillstyle=none](7.784372,4.522950){0.980268}{158.114959}{237.306627}}
\psellipse(7.784372,4.522950)(0.693154,0.693154)
\endpsclip
\psset{linewidth=0.010} \psset{linestyle=dashed,dash=0.1 0.1}
\psset{linestyle=dashed,dash=0.1 0.1} \setlinecaps{0}
\newrgbcolor{dialinecolor}{0 0 0}
\psset{linecolor=dialinecolor}
\psclip{\pswedge[linestyle=none,fillstyle=none](7.207804,6.736254){0.477904}{152.458352}{3.589383}}
\psellipse(7.207804,6.736254)(0.337929,0.337929)
\endpsclip
\psset{linewidth=0.010} \psset{linestyle=dashed,dash=0.1 0.1}
\psset{linestyle=dashed,dash=0.1 0.1} \setlinecaps{0}
\newrgbcolor{dialinecolor}{0 0 0}
\psset{linecolor=dialinecolor}
\psclip{\pswedge[linestyle=none,fillstyle=none](7.376030,3.996590){0.391445}{16.393563}{232.821078}}
\psellipse(7.376030,3.996590)(0.276793,0.276793)
\endpsclip
}\endpspicture\]

Thus, the example from page \pageref{main-example} yields the
following term in expression for $\alpha(s)(f_1,...,f_k)$:
\begin{multline}\label{exa-x}
\sum_{(c^1_5),(c^2_2),(c^3_{10}),(c^4_6),(c^5_7)}(-1)^\epsilon
\Big(c^5_1,c^5_2,c^3_9, c^2_6,c^5_6; (c^5_7)' \Big)\otimes
\Big(1_{TA^*};(c^2_2)''''\Big )\otimes \Big(1_{TA^*};(c^2_2)''\Big)\\
\otimes \Big(c^1_6, c^1_1, c^1_2, c^3_4, c^3_5, c^3_6, c^3_7,
c^4_1, c^4_2, c^4_3, c^4_4, c^4_5, c^5_4, c^2_8, c^2_1, c^1_4,
c^3_1, c^2_3, c^2_4 ; (c^1_5)'' \Big) \cdot c^1_7(1)\cdot
c^2_9(1)\\ \cdot (c^3_{10})'(1)\cdot (c^4_6)' (1)\cdot
(c^5_7)''(1)\cdot <\gamma(c^5_5)\cdot \gamma (c^2_7),1>\cdot
<\gamma(c^1_3)\cdot \gamma (c^3_3),1>\\ \cdot <\gamma(c^3_8)\cdot
\gamma((c^4_6)'')\cdot \gamma(c^5_3),1> \cdot <\gamma(c^2_5)\cdot
\gamma ((c^1_5)')\cdot \gamma((c^3_{10})''),1>\\ \cdot <\gamma
((c^2_2)')\cdot \gamma ((c^2_2)''')\cdot \gamma ((c^2_2)^{
(5)}),1> \cdot <\gamma((c^2_2)^{(6)})\cdot \gamma(c^3_2),1>
\end{multline}
Here, $\epsilon=(2+4+7+10+12+15+16+21+30+32)+(2\cdot 20+ 1\cdot
14+ 1\cdot 7+ 1\cdot 19+ 1\cdot 3+ 4\cdot 6+ 5\cdot 6+ 1\cdot 6+
1\cdot 4+ 1\cdot 2+ 2\cdot 2)\equiv 0\,\, (mod\,\, 2)$, and
$1_{TA^*}$ denotes the unit in the tensor algebra $TA^*$. Note
that we could use

\[ \sum_{(c^5_7)} (c^5_7)'\cdot ((c^5_7)''(1))=c^5_7 \]
\begin{multline*}
\sum_{(c^4_6)} (c^4_6)'(1)\cdot <\gamma (c^3_8)\cdot \gamma
((c^4_6)'')\cdot \gamma(c^5_3),1>\\= <\gamma(c^3_8)\cdot \gamma
(c^4_6)\cdot \gamma(c^5_3),1>
\end{multline*}

\begin{multline*}
\sum_{(c^1_5),(c^3_{10})} (c^3_{10})'(1)\cdot <\gamma(c^2_5)\cdot
\gamma ((c^1_5)')\cdot \gamma((c^3_{10})''),1>\cdot (c^1_5)''\\=
<\gamma(c^2_5)\cdot - \cdot\gamma (c^1_5)\cdot
\gamma(c^3_{10}),1>\in A^*
\end{multline*}where
$$<\gamma(c^2_5)\cdot - \cdot\gamma (c^1_5)\cdot
\gamma(c^3_{10}),1>: c \mapsto <\gamma(c^2_5)\cdot c \cdot\gamma
(c^1_5)\cdot \gamma(c^3_{10}),1>$$

and other identities to simplify this expression.

It remains to show that this well-defined map, $\alpha:C_*\mathcal
S^c\to \mathcal End_{\overline{ HC^*}(A;A)}$, respects the
differentials, composition, and symmetric group action. It must be
clear that the symmetric group action is respected, since it
simply has to do with relabelling.

{\it Step III: $\alpha$ respects differentials}

The boundary operator, $D$, on $\mathcal End_{\overline{
HC^*}(A;A)}$ is given by pre and post compositions with the
boundary operators on the domain and the range. More precisely,
for $\alpha(s) \in \mathcal End_{\overline{ HC^*}(A;A)}(k,l)$, we
have the following.
\begin{multline}\label{dif(a(s))}
D (\alpha(s))(f_1,...,f_k)= \quad \sum_{j=1}^{k} \alpha(s)\circ
\big( id^{\otimes (j-1)
}\otimes \delta\otimes id^{\otimes (k-j)} \big) (f_1,...,f_k)\\
-(-1)^{ |s|} \sum_{j=1}^l \big( id^{\otimes (j-1) }\otimes
\delta\otimes id^{\otimes (l-j)} \big) \circ
\alpha(s)(f_1,...,f_k)
\end{multline}
where $f_i=(c^i_1,...,c^i_{n_i};c^i_{n_i+1})\in
\overline{HC^*}(A;A)$, and $|s|$ denotes the degree of $s$. Here,
$\delta=\delta_1+\delta_2$ is the boundary operator on $\overline{
HC^*}(A;A)$, which is given by applying the comultiplication
$\Delta:A^*\to A^*\otimes A^*$ to all $c^i_j$'s as follows.
\begin{equation*}\delta_1(c^i_1,...,c^i_{n_i};c^i_{n_i+1})=\sum_{j=1}^{n_i}
\sum_{(c^i_j)} (-1)^{j-1}\cdot (c^i_1,..., (c^i_j)',(c^i_j)'',...,
c^i_{n_i} ;c^i_{ n_i+1})
\end{equation*}
\begin{eqnarray*}
\delta_2(c^i_1,...,c^i_{n_i};c^i_{n_i+1})&=&\sum_{(c^i_{n_i+1})}
(-1)^{n_i}( c^i_1, ...,c^i_{n_i}, (c^i_{n_i+1} )';(c^i_{n_i+1}
)'')\\ &+& (-1)^{n_i+ 1}((c^i_{n_i+1} )'', c^i_1, ..., c^i_{n_i} ;
(c^i_{n_i+1} )')
\end{eqnarray*}
For an $s\in C_* \mathcal S^c$, let us compare $D(\alpha(s))$ with
$\alpha (\partial (s))$. The right hand side of the equation
(\ref{dif(a(s))}) is of the form $S+T_1+T_2$, where
$$S=\sum_{j=1}^{k} \alpha(s)\circ \big( id^{\otimes (j-1) }\otimes
\delta\otimes id^{\otimes (k-j)} \big) (f_1,...,f_k)$$
$$T_1=-(-1)^{ |s|} \sum_{j=1}^l \big( id^{\otimes (j-1) }\otimes
\delta_1 \otimes id^{\otimes (l-j)} \big) \circ
\alpha(s)(f_1,...,f_k)$$ $$T_2=-(-1)^{ |s|} \sum_{j=1}^l \big(
id^{\otimes (j-1) }\otimes \delta_2 \otimes id^{\otimes (l-j)}
\big) \circ \alpha(s)(f_1,...,f_k)$$

Each term in $S$ is obtained by first applying $\Delta$ to a
$c^i_j$, and then placing special points of $\alpha$ on the
outcome in all possible ways as described in {\it Step I}. In
doing so, there are terms in which none, one, or both of the
tensor factors of $\Delta(c^i_j)=\sum_{(c^i_j)} (c^i_j)'\otimes
(c^i_j)''$ come in contact with the special points. For the
purposes of this proof, refer to these terms as $S_0$, $S_1$ and
$S_2$, respectively. It is easy to see that $\alpha (\partial (s))
=S_2$. Note that the alternating signs which appears in the
definition of the boundary operator corresponds to those arising
from the fact that if an element $c^i_j$ is at the $r^{th}$ one of
the $|s|+k$ special points on the input circles of $s$, then
$|s|+k-r$ shift maps move over an additional factor of $(c^i_j)'$,
giving rise to the appropriate sign. It is also easy to see that
$S_0= -T_1$. It remains to understand what happens with $S_1$ and
$T_2$. The claim is that $S_1=-T_2$. Seeing this is a bit less
straightforward, since some of the terms in $S_1$ cancel amongst
themselves, whereas other terms cancel with $T_2$. The following
helps to better understand the situation. Consider a simple
diagram in which the special points do not coalesce, and let us
concentrate on an element in $S_1$. If an output marked point is
placed on one of the factor in $\Delta$, then it corresponds
exactly to a term in $T_2$. If an input marked point is placed on
one of the factor in $\Delta$, then it cancels out with a similar
term of $S_1$. This is because for a single input marked point, we
are dealing with evaluation on the unit, and the two terms
$\sum_{(c^i_j)} ((c^i_j)'(1))\otimes (c^i_j)''=c^i_j$ and
$\sum_{(c^i_j)} (c^i_j)'\otimes ((c^i_j)''(1))=c^i_j$ cancel. Note
that in one of the two expressions the shift map has moved past
$(c^i_j)'$, giving rise to a desired negative sign. Now consider
the case in which a single chord endpoint is attached to one of
the tensor factors, for example, $(c^i_j)''$. Follow the cyclic
ordering of the chord's endpoints to go the next chord endpoint
which is attached to, let's say, $c^p_q$. In this case the term
involving $(c^i_j)''$ cancels out with a term in the sum where
$c^p_q$ is split and the chord is attached to $(c^p_q)'$.
\ifx\setlinejoinmode\undefined
 \newcommand{\setlinejoinmode}[1]{}
\fi \ifx\setlinecaps\undefined
 \newcommand{\setlinecaps}[1]{}
\fi \ifx\setfont\undefined
 \newcommand{\setfont}[2]{}
\fi \[\pspicture(4.976079,-6.715574)(12.476431,-3.902295)
\scalebox{1 -1}{
\newrgbcolor{dialinecolor}{0 0 0}
\psset{linecolor=dialinecolor}
\newrgbcolor{diafillcolor}{1 1 1}
\psset{fillcolor=diafillcolor}
\newrgbcolor{dialinecolor}{0 0 0}
\psset{linecolor=dialinecolor}
\rput[r](8.843170,5.323900){\scalebox{1 -1}{=}}

\newrgbcolor{dialinecolor}{0 0 0}
\psset{linecolor=dialinecolor}
\rput[r](12.455700,5.337210){\scalebox{1 -1}{$c^i_j$}}

\newrgbcolor{dialinecolor}{0 0 0}
\psset{linecolor=dialinecolor}
\rput[r](9.971000,5.109700){\scalebox{1 -1}{$(c^p_q)'$}}

\newrgbcolor{dialinecolor}{0 0 0}
\psset{linecolor=dialinecolor}
\rput[r](9.984400,5.721510){\scalebox{1 -1}{$(c^p_q)''$}}
\psset{linewidth=0.010} \psset{linestyle=dashed,dash=1 1}
\psset{linestyle=dashed,dash=0.1 0.1} \setlinecaps{0}
\newrgbcolor{dialinecolor}{0 0 0}
\psset{linecolor=dialinecolor}
\psclip{\pswedge[linestyle=none,fillstyle=none](11.070639,7.167736){3.135938}{246.493361}{295.893618}}
\psellipse(11.070639,7.167736)(2.217443,2.217443)
\endpsclip
\psset{linewidth=0.040} \psset{linestyle=solid}
\psset{linestyle=solid} \setlinecaps{0}
\newrgbcolor{dialinecolor}{0 0 0}
\psset{linecolor=dialinecolor}
\psclip{\pswedge[linestyle=none,fillstyle=none](13.753968,5.449889){2.484393}{137.633871}{219.871283}}
\psellipse(13.753968,5.449889)(1.756731,1.756731)
\endpsclip
\psset{linewidth=0.040} \psset{linestyle=solid}
\psset{linestyle=solid} \setlinecaps{0}
\newrgbcolor{dialinecolor}{0 0 0}
\psset{linecolor=dialinecolor}
\psclip{\pswedge[linestyle=none,fillstyle=none](8.388599,5.390479){2.588698}{326.554195}{41.404570}}
\psellipse(8.388599,5.390479)(1.830486,1.830486)
\endpsclip
\psset{linewidth=0.010} \psset{linestyle=dashed,dash=1 1}
\psset{linestyle=dashed,dash=0.1 0.1} \setlinecaps{0}
\newrgbcolor{dialinecolor}{0 0 0}
\psset{linecolor=dialinecolor}
\psclip{\pswedge[linestyle=none,fillstyle=none](11.718592,4.656148){0.980257}{158.114498}{237.307256}}
\psellipse(11.718592,4.656148)(0.693146,0.693146)
\endpsclip
\psset{linewidth=0.010} \psset{linestyle=dashed,dash=0.1 0.1}
\psset{linestyle=dashed,dash=0.1 0.1} \setlinecaps{0}
\newrgbcolor{dialinecolor}{0 0 0}
\psset{linecolor=dialinecolor}
\psclip{\pswedge[linestyle=none,fillstyle=none](11.310261,4.129789){0.391443}{16.393748}{232.822626}}
\psellipse(11.310261,4.129789)(0.276792,0.276792)
\endpsclip
\psset{linewidth=0.040} \psset{linestyle=solid}
\psset{linestyle=solid} \setlinecaps{0}
\newrgbcolor{dialinecolor}{0 0 0}
\psset{linecolor=dialinecolor}
\psclip{\pswedge[linestyle=none,fillstyle=none](9.469364,5.293421){3.106557}{140.383727}{210.885758}}
\psellipse(9.469364,5.293421)(2.196667,2.196667)
\endpsclip

\newrgbcolor{dialinecolor}{0 0 0}
\psset{linecolor=dialinecolor} \rput[l](7.5,5.163695){\scalebox{1
-1}{$(c^i_j)''$}}
\newrgbcolor{dialinecolor}{0 0 0}
\psset{linecolor=dialinecolor}
\rput[r](5.333100,5.453600){\scalebox{1 -1}{$c^p_q$}}

\newrgbcolor{dialinecolor}{0 0 0}
\psset{linecolor=dialinecolor}
\rput[l](7.536100,5.685200){\scalebox{1 -1}{$(c^i_j)'$}}
\psset{linewidth=0.040} \psset{linestyle=solid}
\psset{linestyle=solid} \setlinecaps{0}
\newrgbcolor{dialinecolor}{0 0 0}
\psset{linecolor=dialinecolor}
\psclip{\pswedge[linestyle=none,fillstyle=none](3.776788,5.245743){2.496452}{325.386523}{46.226931}}
\psellipse(3.776788,5.245743)(1.765258,1.765258)
\endpsclip
\psset{linewidth=0.010} \psset{linestyle=dashed,dash=1 1}
\psset{linestyle=dashed,dash=0.1 0.1} \setlinecaps{0}
\newrgbcolor{dialinecolor}{0 0 0}
\psset{linecolor=dialinecolor}
\psclip{\pswedge[linestyle=none,fillstyle=none](6.641816,7.058474){2.930001}{237.233132}{288.368367}}
\psellipse(6.641816,7.058474)(2.071823,2.071823)
\endpsclip
\psset{linewidth=0.010} \psset{linestyle=dashed,dash=0.1 0.1}
\psset{linestyle=dashed,dash=0.1 0.1} \setlinecaps{0}
\newrgbcolor{dialinecolor}{0 0 0}
\psset{linecolor=dialinecolor}
\psclip{\pswedge[linestyle=none,fillstyle=none](7.004692,4.778588){0.980257}{158.114498}{237.307256}}
\psellipse(7.004692,4.778588)(0.693146,0.693146)
\endpsclip
\psset{linewidth=0.010} \psset{linestyle=dashed,dash=0.1 0.1}
\psset{linestyle=dashed,dash=0.1 0.1} \setlinecaps{0}
\newrgbcolor{dialinecolor}{0 0 0}
\psset{linecolor=dialinecolor}
\psclip{\pswedge[linestyle=none,fillstyle=none](6.549661,4.242767){0.391440}{16.392263}{232.822317}}
\psellipse(6.549661,4.242767)(0.276790,0.276790)
\endpsclip
}\endpspicture\] In other words, we use the algebraic fact that
\begin{multline*}
\sum_{(c^i_j)} (c^i_j)'\cdot <...\cdot \gamma((c^i_j)'') \cdot
\gamma(c^p_q) \cdot..., 1>\\=\sum_{(c^p_q)} <...\cdot
\gamma(c^i_j) \cdot \gamma((c^p_q)') \cdot..., 1>\cdot (c^p_q)''
\end{multline*}
Note that the signs become opposite when applying $\Delta$ to
$c^i_j$ on the left side of this equation and moving $(c^i_j)'$ to
the spot of $c^p_q$, instead of applying $\Delta$ to $c^p_q$.

Now, if several special points coalesce, then we can do the same
steps as above in the cyclic order specified at this point. One
can slide the tensor factors of $\sum_{(c^i_j)}(c^i_j)'\otimes
(c^i_j)'' \otimes ... \otimes (c^i_j)^{(r)}$ which are not
attached to anything from one side to the other in order for them
to cancel out.

{\it Step IV: $\alpha$ respects compositions}

Let us argue why $\alpha$ respects the composition. Let $s\in C_*
\mathcal S^c(k,l)$ and $s'\in C_* \mathcal S^c(m,k)$, and consider
the composition $s\circ s'\in C_* \mathcal S^c(m,l)$. Recall that
we have to identify the $j^{th}$ output circle of $s'$ with the
$j^{th}$ input circle of $s$ starting at the respective marked
points; see page \pageref{comp}. We need to show that
$\alpha(s\circ s')=\alpha(s)\circ \alpha(s')$, where the
composition in $\mathcal End_{\overline{HC^*}(A,A)}$ is given by
\eqref{End-comp} on page \pageref{End-comp}. Thus, for
$\alpha(s)\circ \alpha(s')$, we need to apply $\alpha(s)$ to the
output $\alpha(s')(f_1,...,f_m)\in \overline {HC^*}(A,A)^{\otimes
k}$. Assuming again that each $f_i$ is of the form
$f_i=(c^i_1,...,c^i_{n_i};c^i_{n_i+1})\in (A^*)^{\otimes n_i}
\otimes A^*$, we see that the $k$ tensor factors of $\alpha(s')
(f_1, ...,f_m)$ consist of $c^i_j$'s following the direction of
the output circles of $s'$, together with some coefficients,
compare \eqref{exa-x}. We need to apply $\alpha(s)$ to this, which
means that output marked points and chords need to be added to the
$k$ tensor factors of $\alpha(s') (f_1, ...,f_m)$ in all possible
ways. Thus one sums over all possibilities of placing output
marked points (see item (b) on page \pageref{def-alpha-b}), and
chord endpoints (see item (c) on page \pageref{def-alpha-c}) on
the $c^i_j$'s
according to the combinatorics given by $s$. But this means
exactly that we apply chords and output marked points at the
points specified by the composition $s\circ s'$. Since everything
is graded, and we have to move the same number of elements past
each other to obtain the same expression, we also obtain the same
overall sign. Notice that this argument also works if several
special points coincide at some point, since this only means that
the coproduct $\Delta$ has to be applied to $c^i_j$; see item (d)
on page \pageref{def-alpha-d}. The above arguments can be applied
to all $c^i_j$'s of the $k$ factors of $\alpha(s') (f_1,
...,f_m)$, which are not output marked points. Now, let
$\tilde{c}\in A^*$ represent one of the output marked points. The
definition of $\alpha(s)$ in item (a) on page
\pageref{def-alpha-a} requires to apply the unit $1$ to this
element $\tilde{c}$. Note that $\tilde{c}$ might either be a
factor of some coproduct -such as $(c^5_7)'$, $(c^2_2)''''$,
$(c^2_2)''$ and $(c^1_5)''$ in equation \eqref{exa-x}- or not. In
the first case, we can completely eliminate this marked point by
using the algebraic fact
\begin{multline*}
\sum_{(c^i_j)} (c^i_j)'\otimes (c^i_j)''\otimes...\otimes
(c^i_j)^{(p)}(1)\otimes ...\otimes (c^i_j)^{(r)}\\=\sum_{(c^i_j)}
(c^i_j)'\otimes (c^i_j)''\otimes...\otimes (c^i_j)^{(r-1)}
\end{multline*}
\ifx\setlinejoinmode\undefined
 \newcommand{\setlinejoinmode}[1]{}
\fi \ifx\setlinecaps\undefined
 \newcommand{\setlinecaps}[1]{}
\fi \ifx\setfont\undefined
 \newcommand{\setfont}[2]{}
\fi \[\pspicture(4.485320,-9.386770)(10.133398,-6.079470)
\scalebox{1 -1}{
\newrgbcolor{dialinecolor}{0 0 0}
\psset{linecolor=dialinecolor}
\newrgbcolor{diafillcolor}{1 1 1}
\psset{fillcolor=diafillcolor} \psset{linewidth=0.040}
\psset{linestyle=solid} \psset{linestyle=solid} \setlinecaps{0}
\newrgbcolor{dialinecolor}{0 0 0}
\psset{linecolor=dialinecolor}
\psclip{\pswedge[linestyle=none,fillstyle=none](5.358357,9.465992){2.557318}{252.271558}{268.313872}}
\psellipse(5.358357,9.465992)(1.808297,1.808297)
\endpsclip
\psset{linewidth=0.040} \psset{linestyle=solid}
\setlinejoinmode{0} \setlinecaps{0}
\newrgbcolor{dialinecolor}{0 0 0}
\psset{linecolor=dialinecolor}
\psline(5.247907,7.706217)(5.349835,7.660274)(5.251924,7.606298)

\newrgbcolor{dialinecolor}{0 0 0}
\psset{linecolor=dialinecolor}
\rput(4.585320,7.759570){\scalebox{1 -1}{$1$}}
\psset{linewidth=0.040} \psset{linestyle=solid}
\psset{linestyle=solid} \setlinecaps{0}
\newrgbcolor{dialinecolor}{0 0 0}
\psset{linecolor=dialinecolor}
\psclip{\pswedge[linestyle=none,fillstyle=none](9.064675,7.466525){3.825024}{135.984331}{209.227689}}
\psellipse(9.064675,7.466525)(2.704701,2.704701)
\endpsclip
\newrgbcolor{dialinecolor}{0 0 0}
\psset{linecolor=dialinecolor}
\rput(6.041440,6.479470){\scalebox{1 -1}{$(c^i_j)^{(r)}$}}
\newrgbcolor{dialinecolor}{0 0 0}
\psset{linecolor=dialinecolor}
\rput(5.913940,7.063770){\scalebox{1 -1}{...}}
\newrgbcolor{dialinecolor}{0 0 0}
\psset{linecolor=dialinecolor}
\rput(5.873240,7.622970){\scalebox{1 -1}{$(c^i_j)^{(p)}$}}
\newrgbcolor{dialinecolor}{0 0 0}
\psset{linecolor=dialinecolor}
\rput(5.889440,8.241870){\scalebox{1 -1}{...}}
\newrgbcolor{dialinecolor}{0 0 0}
\psset{linecolor=dialinecolor}
\rput(6.084940,8.830770){\scalebox{1 -1}{$(c^i_j)''$}}
\newrgbcolor{dialinecolor}{0 0 0}
\psset{linecolor=dialinecolor}
\rput(6.426940,9.286770){\scalebox{1 -1}{$(c^i_j)'$}}
\newrgbcolor{dialinecolor}{0 0 0}
\psset{linecolor=dialinecolor}
\rput(7.794070,7.684570){\scalebox{1 -1}{=}}
\psset{linewidth=0.040} \psset{linestyle=solid}
\psset{linestyle=solid} \setlinecaps{0}
\newrgbcolor{dialinecolor}{0 0 0}
\psset{linecolor=dialinecolor}
\psclip{\pswedge[linestyle=none,fillstyle=none](11.790802,7.484180){3.443157}{133.629229}{210.980161}}
\psellipse(11.790802,7.484180)(2.434680,2.434680)
\endpsclip
\newrgbcolor{dialinecolor}{0 0 0}
\psset{linecolor=dialinecolor} \rput(8.8,6.693800){\scalebox{1
-1}{$(c^i_j)^{(r-1)}$}}
\newrgbcolor{dialinecolor}{0 0 0}
\psset{linecolor=dialinecolor}
\rput(8.872360,7.330700){\scalebox{1 -1}{...}}
\newrgbcolor{dialinecolor}{0 0 0}
\psset{linecolor=dialinecolor}
\rput(8.810870,8.016190){\scalebox{1 -1}{$(c^i_j)''$}}
\newrgbcolor{dialinecolor}{0 0 0}
\psset{linecolor=dialinecolor}
\rput(8.995270,8.582890){\scalebox{1 -1}{$(c^i_j)'$}}
\newrgbcolor{dialinecolor}{0 0 0}
\psset{linecolor=dialinecolor} }\endpspicture
 \] In the second case, we apply the unit $1$
to some $\tilde{c}=c^i_j$, where $1\leq j\leq n_i$. But in the
normalized Hochschild complex, we have $f_i( ...,1,...)=0$, or
$c^i_j(1)=0$. Thus, the composition vanishes. This is consistent with
the fact that $s\circ s'=0$, since the dimension of $s\circ s'$ is
less than the sum of the dimensions of $s$ and $s'$; see page
\pageref{comp-0-by-dimension}.

We have shown that the action of $\alpha(s) \circ \alpha(s')$ is the
same as the action of $\alpha(s\circ s')$, and this completes the
proof of the theorem.
\end{proof}

\section{The commutative Case}\label{section-comm}

Throughout this section $A$ denotes an associative, commutative,
and unital algebra, which is endowed with a non-degenerate
invariant inner product. We describe how the PROP of cyclic
Sullivan chord diagrams, $C_\ast \mathcal{S}^c$, can be enlarged
to include orientation-reversing chords in the action. Thickening
such a chord gives rise to a non-orientable surface with boundary.
This enlarged PROP, denoted by $C_\ast \mathcal{S}$, will then act
on the Hochschild complex of the algebra $A$.

\begin{obs}\
\begin{itemize}
\item[(i)]\label{observat-i} Define the
\emph{orientation-reversing} operation $$\sim : HC^*(A;A)\to
HC^*(A;A)$$ $$\sim : f\mapsto \tilde{f}$$ $$\tilde{f}( a_1,a_2,\dots
,a_{n-1},a_n):= (-1)^{\frac{n(n+1)}{2}}\cdot
f(a_n,a_{n-1},\dots,a_2,a_1)$$ One can check that $\left(
\delta\Big(\tilde{f}\Big)-\widetilde{\delta(f)} \right)(a_1, \dots,
a_{n+1})$ is equal to
\begin{eqnarray*}
&&\sum_{j=1}^n \pm f(a_{n+1}, \dots,a_ja_{j+1}- a_{j+1}a_j
,\dots,a_1)\\ &&\pm \big(a_1\cdot f(a_2,\dots, a_{n+1})- f(a_2,\dots,
a_{n+1})\cdot a_1\big) \\ &&\pm \big(a_{n+1}\cdot f(a_1,\dots,
a_{n})- f(a_1,\dots, a_{n})\cdot a_{n+1}\big)
\end{eqnarray*}
Since $A$ is commutative, this expression vanishes. This means
that for commutative $A$ the map $\sim$ is a chain map of the
Hochschild complex into itself. The operation $\sim$ can be
obtained from the following chord diagram, where we insert a
string of elements in one direction, and read them off in the
opposite direction. We refer to this diagram as the
orientation-reversing chord diagram.
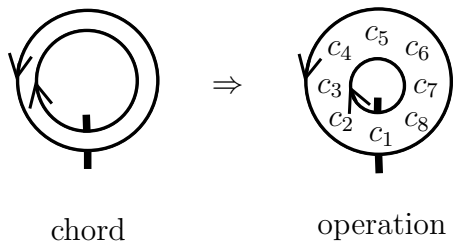
\[ \pspicture(1.747413,-11.748770)(7.868260,-8.659467)
\scalebox{1 -1}{
\newrgbcolor{dialinecolor}{0 0 0}
\psset{linecolor=dialinecolor}
\newrgbcolor{diafillcolor}{1 1 1}
\psset{fillcolor=diafillcolor} \setfont{Helvetica}{0.40}
\newrgbcolor{dialinecolor}{0 0 0}
\psset{linecolor=dialinecolor}
\rput(4.858620,9.742467){\scalebox{1 -1}{$\Rightarrow$}}
\setfont{Helvetica}{0.40}
\newrgbcolor{dialinecolor}{0 0 0}
\psset{linecolor=dialinecolor}
\rput(3.006940,11.623770){\scalebox{1 -1}{chord}}
\setfont{Helvetica}{0.40}
\newrgbcolor{dialinecolor}{0 0 0}
\psset{linecolor=dialinecolor}
\rput(6.910620,11.585467){\scalebox{1 -1}{operation}}
\psset{linewidth=0.040} \psset{linestyle=solid}
\psset{linestyle=solid} \setlinecaps{0} \setlinejoinmode{0}
\setlinecaps{0} \setlinejoinmode{0} \psset{linestyle=solid}
\newrgbcolor{dialinecolor}{1 1 1}
\psset{linecolor=dialinecolor}
\psellipse*(6.864080,9.663647)(0.984180,0.984180)
\newrgbcolor{dialinecolor}{0 0 0}
\psset{linecolor=dialinecolor}
\psellipse(6.864080,9.663647)(0.984180,0.984180)
\psset{linewidth=0.004000} \setlinecaps{0} \setlinejoinmode{0}
\psset{linestyle=solid}
\newrgbcolor{dialinecolor}{0 0 0}
\psset{linecolor=dialinecolor}
\psellipse(6.864080,9.663647)(0.984180,0.984180)
\psset{linewidth=0.040} \psset{linestyle=solid}
\psset{linestyle=solid} \setlinecaps{0} \setlinejoinmode{0}
\setlinecaps{0} \setlinejoinmode{0} \psset{linestyle=solid}
\newrgbcolor{dialinecolor}{1 1 1}
\psset{linecolor=dialinecolor}
\psellipse*(6.849599,9.716690)(0.375899,0.375899)
\newrgbcolor{dialinecolor}{0 0 0}
\psset{linecolor=dialinecolor}
\psellipse(6.849599,9.716690)(0.375899,0.375899)
\psset{linewidth=0.004000} \setlinecaps{0} \setlinejoinmode{0}
\psset{linestyle=solid}
\newrgbcolor{dialinecolor}{0 0 0}
\psset{linecolor=dialinecolor}
\psellipse(6.849599,9.716690)(0.375899,0.375899)
\psset{linewidth=0.10} \psset{linestyle=solid}
\psset{linestyle=solid} \setlinecaps{0}
\newrgbcolor{dialinecolor}{0 0 0}
\psset{linecolor=dialinecolor}
\psline(6.849599,10.092589)(6.853620,9.875467)
\psset{linewidth=0.10} \psset{linestyle=solid}
\psset{linestyle=solid} \setlinecaps{0}
\newrgbcolor{dialinecolor}{0 0 0}
\psset{linecolor=dialinecolor}
\psline(6.872620,10.882467)(6.864080,10.647826)
\psset{linewidth=0.040} \psset{linestyle=solid}
\psset{linestyle=solid} \setlinecaps{0}
\newrgbcolor{dialinecolor}{0 0 0}
\psset{linecolor=dialinecolor}
\psline(6.509715,9.798562)(6.568620,9.932467)
\psset{linewidth=0.040} \psset{linestyle=solid}
\setlinejoinmode{0} \setlinecaps{0}
\newrgbcolor{dialinecolor}{0 0 0}
\psset{linecolor=dialinecolor}
\psline(6.749808,9.971832)(6.491708,9.757626)(6.475203,10.092630)
\psset{linewidth=0.040} \psset{linestyle=solid}
\psset{linestyle=solid} \setlinecaps{0}
\newrgbcolor{dialinecolor}{0 0 0}
\psset{linecolor=dialinecolor}
\psline(5.899012,9.576270)(5.941620,9.381467)
\psset{linewidth=0.040} \psset{linestyle=solid}
\setlinejoinmode{0} \setlinecaps{0}
\newrgbcolor{dialinecolor}{0 0 0}
\psset{linecolor=dialinecolor}
\psline(5.807022,9.294836)(5.889456,9.619958)(6.100093,9.358937)
\setfont{Helvetica}{0.40}
\newrgbcolor{dialinecolor}{0 0 0}
\psset{linecolor=dialinecolor} \rput(6.910620,10.4){\scalebox{1
-1}{$c_1$}} \setfont{Helvetica}{0.40}
\newrgbcolor{dialinecolor}{0 0 0}
\psset{linecolor=dialinecolor} \rput(6.373420,10.2){\scalebox{1
-1}{$c_2$}} \setfont{Helvetica}{0.40}
\newrgbcolor{dialinecolor}{0 0 0}
\psset{linecolor=dialinecolor} \rput(6.221420,9.75){\scalebox{1
-1}{$c_3$}} \setfont{Helvetica}{0.40}
\newrgbcolor{dialinecolor}{0 0 0}
\psset{linecolor=dialinecolor} \rput(6.354420,9.25){\scalebox{1
-1}{$c_4$}} \setfont{Helvetica}{0.40}
\newrgbcolor{dialinecolor}{0 0 0}
\psset{linecolor=dialinecolor} \rput(6.848420,9.05){\scalebox{1
-1}{$c_5$}} \setfont{Helvetica}{0.40}
\newrgbcolor{dialinecolor}{0 0 0}
\psset{linecolor=dialinecolor} \rput(7.385620,9.25){\scalebox{1
-1}{$c_6$}} \setfont{Helvetica}{0.40}
\newrgbcolor{dialinecolor}{0 0 0}
\psset{linecolor=dialinecolor} \rput(7.499620,9.75){\scalebox{1
-1}{$c_7$}} \setfont{Helvetica}{0.40}
\newrgbcolor{dialinecolor}{0 0 0}
\psset{linecolor=dialinecolor} \rput(7.380420,10.2){\scalebox{1
-1}{$c_8$}} \psset{linewidth=0.040} \psset{linestyle=solid}
\psset{linestyle=solid} \setlinecaps{0} \setlinejoinmode{0}
\setlinecaps{0} \setlinejoinmode{0} \psset{linestyle=solid}
\newrgbcolor{dialinecolor}{1 1 1}
\psset{linecolor=dialinecolor}
\psellipse*(2.996620,9.647467)(0.950,0.950)
\newrgbcolor{dialinecolor}{0 0 0}
\psset{linecolor=dialinecolor}
\psellipse(2.996620,9.647467)(0.950,0.950)
\psset{linewidth=0.004000} \setlinecaps{0} \setlinejoinmode{0}
\psset{linestyle=solid}
\newrgbcolor{dialinecolor}{0 0 0}
\psset{linecolor=dialinecolor}
\psellipse(2.996620,9.647467)(0.950,0.950) \psset{linewidth=0.040}
\psset{linestyle=solid} \psset{linestyle=solid} \setlinecaps{0}
\setlinejoinmode{0} \setlinecaps{0} \setlinejoinmode{0}
\psset{linestyle=solid}
\newrgbcolor{dialinecolor}{1 1 1}
\psset{linecolor=dialinecolor}
\psellipse*(2.985720,9.651767)(0.687300,0.687300)
\newrgbcolor{dialinecolor}{0 0 0}
\psset{linecolor=dialinecolor}
\psellipse(2.985720,9.651767)(0.687300,0.687300)
\psset{linewidth=0.004000} \setlinecaps{0} \setlinejoinmode{0}
\psset{linestyle=solid}
\newrgbcolor{dialinecolor}{0 0 0}
\psset{linecolor=dialinecolor}
\psellipse(2.985720,9.651767)(0.687300,0.687300)
\psset{linewidth=0.10} \psset{linestyle=solid}
\psset{linestyle=solid} \setlinecaps{0}
\newrgbcolor{dialinecolor}{0 0 0}
\psset{linecolor=dialinecolor}
\psline(2.985720,10.339067)(2.977620,10.103467)
\psset{linewidth=0.10} \psset{linestyle=solid}
\psset{linestyle=solid} \setlinecaps{0}
\newrgbcolor{dialinecolor}{0 0 0}
\psset{linecolor=dialinecolor}
\psline(2.996620,10.825467)(2.996620,10.597467)
\psset{linewidth=0.040} \psset{linestyle=solid}
\psset{linestyle=solid} \setlinecaps{0}
\newrgbcolor{dialinecolor}{0 0 0}
\psset{linecolor=dialinecolor}
\psline(2.318035,9.739032)(2.331620,9.799467)
\psset{linewidth=0.040} \psset{linestyle=solid}
\setlinejoinmode{0} \setlinecaps{0}
\newrgbcolor{dialinecolor}{0 0 0}
\psset{linecolor=dialinecolor}
\psline(2.520368,9.955200)(2.308227,9.695399)(2.227671,10.020992)
\psset{linewidth=0.040} \psset{linestyle=solid}
\psset{linestyle=solid} \setlinecaps{0}
\newrgbcolor{dialinecolor}{0 0 0}
\psset{linecolor=dialinecolor}
\psline(2.061324,9.559241)(2.084620,9.419467)
\psset{linewidth=0.040} \psset{linestyle=solid}
\setlinejoinmode{0} \setlinecaps{0}
\newrgbcolor{dialinecolor}{0 0 0}
\psset{linecolor=dialinecolor}
\psline(1.955332,9.282776)(2.053972,9.603354)(2.251251,9.332096)
}\endpspicture \] Note that this chord diagram is a closed element in
the complex of chord diagrams, as defined in section \ref{sec1}.
\item[(ii)]\label{observat-ii} Let's look at the brace operation $*$
from page \pageref{brace-smile}. For $f,g\in HC^*(A;A)$ we have:
\begin{multline*}
\quad\quad\quad (f*g)(a_1,\dots, a_n)\\=\sum_k \pm f(a_1,
\dots,a_k,g(a_{k+1},\dots, a_{k+l}), a_{k+l+1},\dots, a_n)
\end{multline*}
We now want to allow reversing of orientations, as described in (i).
For example the operation $\tilde{f}*g$, defined below, is also
legitimate. In this case $(\tilde{f}*g)(a_1,\dots, a_n)$ is equal to
\begin{equation*}
\quad\quad\quad \sum_k \pm f(a_n, \dots,a_{k+l+1},g(a_{k+1},\dots,
a_{k+l}), a_k,\dots, a_1)
\end{equation*}
Note that the elements plugged into $f$ are reversed, while those
plugged into $g$ have preserved their linear ordering. The
following figure shows that this phenomenon can be expressed by
considering diagrams with twisted chords.
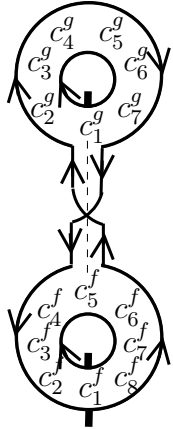
\[ \pspicture(5.581056,-14.390381)(8.149742,-8.659470)
\scalebox{1 -1}{
\newrgbcolor{dialinecolor}{0 0 0}
\psset{linecolor=dialinecolor}
\newrgbcolor{diafillcolor}{1 1 1}
\psset{fillcolor=diafillcolor} \psset{linewidth=0.040}
\psset{linestyle=solid} \psset{linestyle=solid} \setlinecaps{0}
\setlinejoinmode{0} \setlinecaps{0} \setlinejoinmode{0}
\psset{linestyle=solid}
\newrgbcolor{dialinecolor}{1 1 1}
\psset{linecolor=dialinecolor}
\psellipse*(6.864080,9.663650)(0.984180,0.984180)
\newrgbcolor{dialinecolor}{0 0 0}
\psset{linecolor=dialinecolor}
\psellipse(6.864080,9.663650)(0.984180,0.984180)
\psset{linewidth=0.004000} \setlinecaps{0} \setlinejoinmode{0}
\psset{linestyle=solid}
\newrgbcolor{dialinecolor}{0 0 0}
\psset{linecolor=dialinecolor}
\psellipse(6.864080,9.663650)(0.984180,0.984180)
\psset{linewidth=0.040} \psset{linestyle=solid}
\psset{linestyle=solid} \setlinecaps{0} \setlinejoinmode{0}
\setlinecaps{0} \setlinejoinmode{0} \psset{linestyle=solid}
\newrgbcolor{dialinecolor}{1 1 1}
\psset{linecolor=dialinecolor}
\psellipse*(6.849600,9.716690)(0.375900,0.375900)
\newrgbcolor{dialinecolor}{0 0 0}
\psset{linecolor=dialinecolor}
\psellipse(6.849600,9.716690)(0.375900,0.375900)
\psset{linewidth=0.004000} \setlinecaps{0} \setlinejoinmode{0}
\psset{linestyle=solid}
\newrgbcolor{dialinecolor}{0 0 0}
\psset{linecolor=dialinecolor}
\psellipse(6.849600,9.716690)(0.375900,0.375900)
\psset{linewidth=0.10} \psset{linestyle=solid}
\psset{linestyle=solid} \setlinecaps{0}
\newrgbcolor{dialinecolor}{0 0 0}
\psset{linecolor=dialinecolor}
\psline(6.849600,10.092600)(6.853620,9.875470)
\psset{linewidth=0.040} \psset{linestyle=solid}
\psset{linestyle=solid} \setlinecaps{0}
\newrgbcolor{dialinecolor}{0 0 0}
\psset{linecolor=dialinecolor}
\psline(6.509715,9.798562)(6.568620,9.932470)
\psset{linewidth=0.040} \psset{linestyle=solid}
\setlinejoinmode{0} \setlinecaps{0}
\newrgbcolor{dialinecolor}{0 0 0}
\psset{linecolor=dialinecolor}
\psline(6.749807,9.971833)(6.491707,9.757626)(6.475201,10.092630)
\psset{linewidth=0.040} \psset{linestyle=solid}
\psset{linestyle=solid} \setlinecaps{0}
\newrgbcolor{dialinecolor}{0 0 0}
\psset{linecolor=dialinecolor}
\psline(5.895624,9.751700)(5.927014,9.927470)
\psset{linewidth=0.040} \psset{linestyle=solid}
\setlinejoinmode{0} \setlinecaps{0}
\newrgbcolor{dialinecolor}{0 0 0}
\psset{linecolor=dialinecolor}
\psline(6.088167,9.976632)(5.887762,9.707675)(5.792839,10.029373)
\setfont{Helvetica}{0.40}
\newrgbcolor{dialinecolor}{0 0 0}
\psset{linecolor=dialinecolor} \rput(6.910620,10.35){\scalebox{1
-1}{$c^g_1$}} \setfont{Helvetica}{0.40}
\newrgbcolor{dialinecolor}{0 0 0}
\psset{linecolor=dialinecolor} \rput(6.258620,10.1){\scalebox{1
-1}{$c^g_2$}} \setfont{Helvetica}{0.40}
\newrgbcolor{dialinecolor}{0 0 0}
\psset{linecolor=dialinecolor} \rput(6.221420,9.5){\scalebox{1
-1}{$c^g_3$}} \setfont{Helvetica}{0.40}
\newrgbcolor{dialinecolor}{0 0 0}
\psset{linecolor=dialinecolor} \rput(6.518420,9.1){\scalebox{1
-1}{$c^g_4$}} \setfont{Helvetica}{0.40}
\newrgbcolor{dialinecolor}{0 0 0}
\psset{linecolor=dialinecolor} \rput(7.172420,9.1){\scalebox{1
-1}{$c^g_5$}} \setfont{Helvetica}{0.40}
\newrgbcolor{dialinecolor}{0 0 0}
\psset{linecolor=dialinecolor} \rput(7.499620,9.5){\scalebox{1
-1}{$c^g_6$}} \setfont{Helvetica}{0.40}
\newrgbcolor{dialinecolor}{0 0 0}
\psset{linecolor=dialinecolor} \rput(7.413220,10.1){\scalebox{1
-1}{$c^g_7$}} \psset{linewidth=0.040} \psset{linestyle=solid}
\psset{linestyle=solid} \setlinecaps{0} \setlinejoinmode{0}
\setlinecaps{0} \setlinejoinmode{0} \psset{linestyle=solid}
\newrgbcolor{dialinecolor}{1 1 1}
\psset{linecolor=dialinecolor}
\psellipse*(6.867419,13.144010)(0.984180,0.984180)
\newrgbcolor{dialinecolor}{0 0 0}
\psset{linecolor=dialinecolor}
\psellipse(6.867419,13.144010)(0.984180,0.984180)
\psset{linewidth=0.004000} \setlinecaps{0} \setlinejoinmode{0}
\psset{linestyle=solid}
\newrgbcolor{dialinecolor}{0 0 0}
\psset{linecolor=dialinecolor}
\psellipse(6.867419,13.144010)(0.984180,0.984180)
\psset{linewidth=0.040} \psset{linestyle=solid}
\psset{linestyle=solid} \setlinecaps{0} \setlinejoinmode{0}
\setlinecaps{0} \setlinejoinmode{0} \psset{linestyle=solid}
\newrgbcolor{dialinecolor}{1 1 1}
\psset{linecolor=dialinecolor}
\psellipse*(6.852939,13.197050)(0.375900,0.375900)
\newrgbcolor{dialinecolor}{0 0 0}
\psset{linecolor=dialinecolor}
\psellipse(6.852939,13.197050)(0.375900,0.375900)
\psset{linewidth=0.004000} \setlinecaps{0} \setlinejoinmode{0}
\psset{linestyle=solid}
\newrgbcolor{dialinecolor}{0 0 0}
\psset{linecolor=dialinecolor}
\psellipse(6.852939,13.197050)(0.375900,0.375900)
\psset{linewidth=0.10} \psset{linestyle=solid}
\psset{linestyle=solid} \setlinecaps{0}
\newrgbcolor{dialinecolor}{0 0 0}
\psset{linecolor=dialinecolor}
\psline(6.852939,13.572950)(6.856959,13.355830)
\psset{linewidth=0.10} \psset{linestyle=solid}
\psset{linestyle=solid} \setlinecaps{0}
\newrgbcolor{dialinecolor}{0 0 0}
\psset{linecolor=dialinecolor}
\psline(6.861814,14.339070)(6.867419,14.128190)
\psset{linewidth=0.040} \psset{linestyle=solid}
\psset{linestyle=solid} \setlinecaps{0}
\newrgbcolor{dialinecolor}{0 0 0}
\psset{linecolor=dialinecolor}
\psline(6.513053,13.278922)(6.571959,13.412830)
\psset{linewidth=0.040} \psset{linestyle=solid}
\setlinejoinmode{0} \setlinecaps{0}
\newrgbcolor{dialinecolor}{0 0 0}
\psset{linecolor=dialinecolor}
\psline(6.753146,13.452193)(6.495046,13.237986)(6.478540,13.572990)
\psset{linewidth=0.040} \psset{linestyle=solid}
\psset{linestyle=solid} \setlinecaps{0}
\newrgbcolor{dialinecolor}{0 0 0}
\psset{linecolor=dialinecolor}
\psline(5.899841,13.056122)(5.927014,12.912270)
\psset{linewidth=0.040} \psset{linestyle=solid}
\setlinejoinmode{0} \setlinecaps{0}
\newrgbcolor{dialinecolor}{0 0 0}
\psset{linecolor=dialinecolor}
\psline(5.799831,12.777437)(5.891540,13.100066)(6.094618,12.833121)
\setfont{Helvetica}{0.40}
\newrgbcolor{dialinecolor}{0 0 0}
\psset{linecolor=dialinecolor} \rput(6.913959,13.8){\scalebox{1
-1}{$c^f_1$}} \setfont{Helvetica}{0.40}
\newrgbcolor{dialinecolor}{0 0 0}
\psset{linecolor=dialinecolor} \rput(6.376759,13.65){\scalebox{1
-1}{$c^f_2$}} \setfont{Helvetica}{0.40}
\newrgbcolor{dialinecolor}{0 0 0}
\psset{linecolor=dialinecolor} \rput(6.224759,13.2){\scalebox{1
-1}{$c^f_3$}} \setfont{Helvetica}{0.40}
\newrgbcolor{dialinecolor}{0 0 0}
\psset{linecolor=dialinecolor} \rput(6.357759,12.75){\scalebox{1
-1}{$c^f_4$}} \setfont{Helvetica}{0.40}
\newrgbcolor{dialinecolor}{0 0 0}
\psset{linecolor=dialinecolor} \rput(6.851759,12.5){\scalebox{1
-1}{$c^f_5$}} \setfont{Helvetica}{0.40}
\newrgbcolor{dialinecolor}{0 0 0}
\psset{linecolor=dialinecolor} \rput(7.388959,12.75){\scalebox{1
-1}{$c^f_6$}} \setfont{Helvetica}{0.40}
\newrgbcolor{dialinecolor}{0 0 0}
\psset{linecolor=dialinecolor} \rput(7.502959,13.2){\scalebox{1
-1}{$c^f_7$}} \setfont{Helvetica}{0.40}
\newrgbcolor{dialinecolor}{0 0 0}
\psset{linecolor=dialinecolor} \rput(7.383759,13.65){\scalebox{1
-1}{$c^f_8$}}
\newrgbcolor{dialinecolor}{1 1 1}
\psset{linecolor=dialinecolor}
\pspolygon*(6.697814,10.616270)(6.697814,12.190670)(6.976614,12.190670)(6.976614,10.616270)
\psset{linewidth=0.10} \psset{linestyle=solid}
\psset{linestyle=solid} \setlinejoinmode{0}
\newrgbcolor{dialinecolor}{1 1 1}
\psset{linecolor=dialinecolor}
\pspolygon(6.697814,10.616270)(6.697814,12.190670)(6.976614,12.190670)(6.976614,10.616270)
\psset{linewidth=0.040} \psset{linestyle=solid}
\psset{linestyle=solid} \setlinecaps{0}
\newrgbcolor{dialinecolor}{0 0 0}
\psset{linecolor=dialinecolor}
\psline(6.648614,10.599870)(6.648614,11.157470)
\psset{linewidth=0.040} \psset{linestyle=solid}
\psset{linestyle=solid} \setlinecaps{0}
\newrgbcolor{dialinecolor}{0 0 0}
\psset{linecolor=dialinecolor}
\psline(6.630894,11.618630)(6.630894,12.176230)
\psset{linewidth=0.040} \psset{linestyle=solid}
\psset{linestyle=solid} \setlinecaps{0}
\newrgbcolor{dialinecolor}{0 0 0}
\psset{linecolor=dialinecolor}
\psline(7.042214,10.632670)(7.025814,11.255870)
\psset{linewidth=0.040} \psset{linestyle=solid}
\psset{linestyle=solid} \setlinecaps{0}
\newrgbcolor{dialinecolor}{0 0 0}
\psset{linecolor=dialinecolor}
\psline(7.057294,11.635030)(7.057294,12.192630)
\psset{linewidth=0.010} \psset{linestyle=dashed,dash=1 1}
\psset{linestyle=dashed,dash=0.10 0.10} \setlinecaps{0}
\newrgbcolor{dialinecolor}{0 0 0}
\psset{linecolor=dialinecolor}
\psline(6.845414,10.534270)(6.845414,12.289070)
\psset{linewidth=0.040} \psset{linestyle=solid}
\psset{linestyle=solid} \setlinecaps{0}
\newrgbcolor{dialinecolor}{0 0 0}
\psset{linecolor=dialinecolor}
\psclip{\pswedge[linestyle=none,fillstyle=none](7.207884,11.110978){0.818965}{106.623674}{173.765184}}
\psellipse(7.207884,11.110978)(0.579095,0.579095)
\endpsclip
\psset{linewidth=0.040} \psset{linestyle=solid}
\psset{linestyle=solid} \setlinecaps{0}
\newrgbcolor{dialinecolor}{0 0 0}
\psset{linecolor=dialinecolor}
\psclip{\pswedge[linestyle=none,fillstyle=none](6.518172,11.092628){0.786247}{13.569452}{76.430548}}
\psellipse(6.518172,11.092628)(0.555961,0.555961)
\endpsclip
\psset{linewidth=0.040} \psset{linestyle=solid}
\psset{linestyle=solid} \setlinecaps{0}
\newrgbcolor{dialinecolor}{0 0 0}
\psset{linecolor=dialinecolor}
\psline(7.834163,13.231737)(7.796614,13.420670)
\psset{linewidth=0.040} \psset{linestyle=solid}
\setlinejoinmode{0} \setlinecaps{0}
\newrgbcolor{dialinecolor}{0 0 0}
\psset{linecolor=dialinecolor}
\psline(7.931524,13.511358)(7.842881,13.187873)(7.637279,13.452878)
\psset{linewidth=0.040} \psset{linestyle=solid}
\psset{linestyle=solid} \setlinecaps{0}
\newrgbcolor{dialinecolor}{0 0 0}
\psset{linecolor=dialinecolor}
\psline(7.832381,9.575628)(7.813014,9.468270)
\psset{linewidth=0.040} \psset{linestyle=solid}
\setlinejoinmode{0} \setlinecaps{0}
\newrgbcolor{dialinecolor}{0 0 0}
\psset{linecolor=dialinecolor}
\psline(7.639444,9.351034)(7.840321,9.619639)(7.934678,9.297775)
\psset{linewidth=0.040} \psset{linestyle=solid}
\psset{linestyle=solid} \setlinecaps{0}
\newrgbcolor{dialinecolor}{0 0 0}
\psset{linecolor=dialinecolor}
\psline(7.040353,10.978784)(7.042214,10.796670)
\psset{linewidth=0.040} \psset{linestyle=solid}
\setlinejoinmode{0} \setlinecaps{0}
\newrgbcolor{dialinecolor}{0 0 0}
\psset{linecolor=dialinecolor}
\psline(6.892969,10.721985)(7.039896,11.023503)(7.192954,10.725051)
\psset{linewidth=0.040} \psset{linestyle=solid}
\psset{linestyle=solid} \setlinecaps{0}
\newrgbcolor{dialinecolor}{0 0 0}
\psset{linecolor=dialinecolor}
\psline(6.625211,11.962879)(6.615814,11.764270)
\psset{linewidth=0.040} \psset{linestyle=solid}
\setlinejoinmode{0} \setlinecaps{0}
\newrgbcolor{dialinecolor}{0 0 0}
\psset{linecolor=dialinecolor}
\psline(6.463314,11.714975)(6.627325,12.007550)(6.762978,11.700796)
\psset{linewidth=0.040} \psset{linestyle=solid}
\psset{linestyle=solid} \setlinecaps{0}
\newrgbcolor{dialinecolor}{0 0 0}
\psset{linecolor=dialinecolor}
\psline(6.646819,10.846059)(6.648614,11.009870)
\psset{linewidth=0.040} \psset{linestyle=solid}
\setlinejoinmode{0} \setlinecaps{0}
\newrgbcolor{dialinecolor}{0 0 0}
\psset{linecolor=dialinecolor}
\psline(6.799607,11.099679)(6.646329,10.801340)(6.499625,11.102966)
\psset{linewidth=0.040} \psset{linestyle=solid}
\psset{linestyle=solid} \setlinecaps{0}
\newrgbcolor{dialinecolor}{0 0 0}
\psset{linecolor=dialinecolor}
\psline(7.063753,11.846061)(7.057294,11.913830)
\psset{linewidth=0.040} \psset{linestyle=solid}
\setlinejoinmode{0} \setlinecaps{0}
\newrgbcolor{dialinecolor}{0 0 0}
\psset{linecolor=dialinecolor}
\psline(7.188856,12.114419)(7.067996,11.801541)(6.890209,12.085956)
}\endpspicture \]
\end{itemize}
\end{obs}
These two observations demonstrate all the new features of diagrams
describing orientation-reversing operations. Chord diagrams with
possible twisted chords form a PROP which acts on the normalized
Hochschild cochain complex of $A$. The relevant definition and its
application to the normalized Hochschild cochain complex will occupy
the rest of the paper.

\begin{defn}[Sullivan Chord Diagram]A \emph{Sullivan chord diagram} is a
generalization of a cyclic Sullivan chord diagram, where chords
may have twists in them and the orientation of output circles may
be arbitrary.
\end{defn}

Note that, when moving along an output circle, one may alternate
between going in the direction compatible with those of the input
circles, and the opposite direction. This is shown in the below
figure, where the direction of $1^{'th}$ output circle is
compatible with that of the $1^{st}$ input circle, but in
opposition to that of the $2^{nd}$ input circle. A similar remark
applies to going along the input circles, as seen for example the
$1^{st}$ input circle.

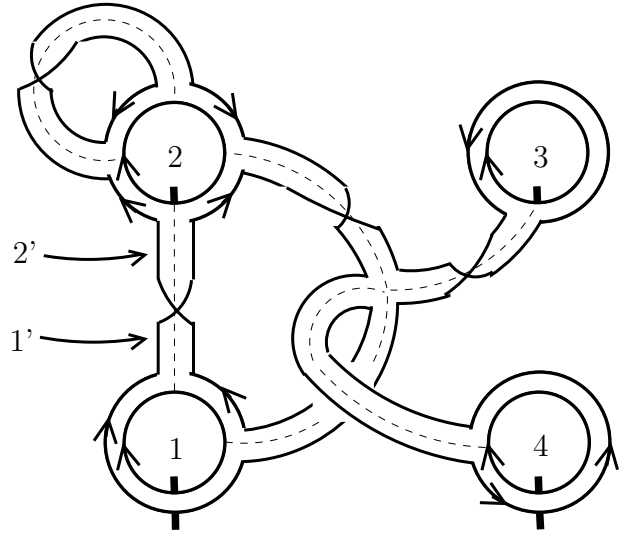
\[ \pspicture(1.080814,-15.086341)(9.613568,-8.004242)
\scalebox{1 -1}{
\newrgbcolor{dialinecolor}{0 0 0}
\psset{linecolor=dialinecolor}
\newrgbcolor{diafillcolor}{1 1 1}
\psset{fillcolor=diafillcolor} \psset{linewidth=0.040}
\psset{linestyle=solid} \psset{linestyle=solid} \setlinecaps{0}
\setlinejoinmode{0} \setlinecaps{0} \setlinejoinmode{0}
\psset{linestyle=solid}
\newrgbcolor{dialinecolor}{1 1 1}
\psset{linecolor=dialinecolor}
\psellipse*(3.517720,10.014170)(0.950,0.950)
\newrgbcolor{dialinecolor}{0 0 0}
\psset{linecolor=dialinecolor}
\psellipse(3.517720,10.014170)(0.950,0.950)
\psset{linewidth=0.004000} \setlinecaps{0} \setlinejoinmode{0}
\psset{linestyle=solid}
\newrgbcolor{dialinecolor}{0 0 0}
\psset{linecolor=dialinecolor}
\psellipse(3.517720,10.014170)(0.950,0.950)
\psset{linewidth=0.040} \psset{linestyle=solid}
\psset{linestyle=solid} \setlinecaps{0} \setlinejoinmode{0}
\setlinecaps{0} \setlinejoinmode{0} \psset{linestyle=solid}
\newrgbcolor{dialinecolor}{1 1 1}
\psset{linecolor=dialinecolor}
\psellipse*(3.506820,10.018470)(0.687300,0.687300)
\newrgbcolor{dialinecolor}{0 0 0}
\psset{linecolor=dialinecolor}
\psellipse(3.506820,10.018470)(0.687300,0.687300)
\psset{linewidth=0.004000} \setlinecaps{0} \setlinejoinmode{0}
\psset{linestyle=solid}
\newrgbcolor{dialinecolor}{0 0 0}
\psset{linecolor=dialinecolor}
\psellipse(3.506820,10.018470)(0.687300,0.687300)
\psset{linewidth=0.10} \psset{linestyle=solid}
\psset{linestyle=solid} \setlinecaps{0}
\newrgbcolor{dialinecolor}{0 0 0}
\psset{linecolor=dialinecolor}
\psline(3.506820,10.705770)(3.498720,10.470200)
\psset{linewidth=0.040} \psset{linestyle=solid}
\psset{linestyle=solid} \setlinecaps{0}
\newrgbcolor{dialinecolor}{0 0 0}
\psset{linecolor=dialinecolor}
\psline(2.839135,10.105735)(2.852720,10.166170)
\psset{linewidth=0.040} \psset{linestyle=solid}
\setlinejoinmode{0} \setlinecaps{0}
\newrgbcolor{dialinecolor}{0 0 0}
\psset{linecolor=dialinecolor}
\psline(3.041468,10.321903)(2.829328,10.062103)(2.748772,10.387696)
\psset{linewidth=0.040} \psset{linestyle=solid}
\psset{linestyle=solid} \setlinecaps{0}
\newrgbcolor{dialinecolor}{0 0 0}
\psset{linecolor=dialinecolor}
\psline(2.719190,9.496012)(2.752014,9.438570)
\psset{linewidth=0.040} \psset{linestyle=solid}
\setlinejoinmode{0} \setlinecaps{0}
\newrgbcolor{dialinecolor}{0 0 0}
\psset{linecolor=dialinecolor}
\psline(2.715607,9.199947)(2.697002,9.534841)(2.976080,9.348789)
\psset{linewidth=0.040} \psset{linestyle=solid}
\psset{linestyle=solid} \setlinecaps{0} \setlinejoinmode{0}
\setlinecaps{0} \setlinejoinmode{0} \psset{linestyle=solid}
\newrgbcolor{dialinecolor}{1 1 1}
\psset{linecolor=dialinecolor}
\psellipse*(8.345650,10.009060)(0.950,0.950)
\newrgbcolor{dialinecolor}{0 0 0}
\psset{linecolor=dialinecolor}
\psellipse(8.345650,10.009060)(0.950,0.950)
\psset{linewidth=0.004000} \setlinecaps{0} \setlinejoinmode{0}
\psset{linestyle=solid}
\newrgbcolor{dialinecolor}{0 0 0}
\psset{linecolor=dialinecolor}
\psellipse(8.345650,10.009060)(0.950,0.950)
\psset{linewidth=0.040} \psset{linestyle=solid}
\psset{linestyle=solid} \setlinecaps{0} \setlinejoinmode{0}
\setlinecaps{0} \setlinejoinmode{0} \psset{linestyle=solid}
\newrgbcolor{dialinecolor}{1 1 1}
\psset{linecolor=dialinecolor}
\psellipse*(8.334750,10.013360)(0.687300,0.687300)
\newrgbcolor{dialinecolor}{0 0 0}
\psset{linecolor=dialinecolor}
\psellipse(8.334750,10.013360)(0.687300,0.687300)
\psset{linewidth=0.004000} \setlinecaps{0} \setlinejoinmode{0}
\psset{linestyle=solid}
\newrgbcolor{dialinecolor}{0 0 0}
\psset{linecolor=dialinecolor}
\psellipse(8.334750,10.013360)(0.687300,0.687300)
\psset{linewidth=0.10} \psset{linestyle=solid}
\psset{linestyle=solid} \setlinecaps{0}
\newrgbcolor{dialinecolor}{0 0 0}
\psset{linecolor=dialinecolor}
\psline(8.334750,10.700660)(8.326650,10.465090)
\psset{linewidth=0.040} \psset{linestyle=solid}
\psset{linestyle=solid} \setlinecaps{0}
\newrgbcolor{dialinecolor}{0 0 0}
\psset{linecolor=dialinecolor}
\psline(7.667065,10.100625)(7.680650,10.161060)
\psset{linewidth=0.040} \psset{linestyle=solid}
\setlinejoinmode{0} \setlinecaps{0}
\newrgbcolor{dialinecolor}{0 0 0}
\psset{linecolor=dialinecolor}
\psline(7.869398,10.316793)(7.657258,10.056993)(7.576702,10.382586)
\psset{linewidth=0.040} \psset{linestyle=solid}
\psset{linestyle=solid} \setlinecaps{0}
\newrgbcolor{dialinecolor}{0 0 0}
\psset{linecolor=dialinecolor}
\psline(7.410354,9.920834)(7.433650,9.781060)
\psset{linewidth=0.040} \psset{linestyle=solid}
\setlinejoinmode{0} \setlinecaps{0}
\newrgbcolor{dialinecolor}{0 0 0}
\psset{linecolor=dialinecolor}
\psline(7.304363,9.644369)(7.403002,9.964947)(7.600281,9.693689)
\psset{linewidth=0.040} \psset{linestyle=solid}
\psset{linestyle=solid} \setlinecaps{0} \setlinejoinmode{0}
\setlinecaps{0} \setlinejoinmode{0} \psset{linestyle=solid}
\newrgbcolor{dialinecolor}{1 1 1}
\psset{linecolor=dialinecolor}
\psellipse*(3.520650,13.869060)(0.950,0.950)
\newrgbcolor{dialinecolor}{0 0 0}
\psset{linecolor=dialinecolor}
\psellipse(3.520650,13.869060)(0.950,0.950)
\psset{linewidth=0.004000} \setlinecaps{0} \setlinejoinmode{0}
\psset{linestyle=solid}
\newrgbcolor{dialinecolor}{0 0 0}
\psset{linecolor=dialinecolor}
\psellipse(3.520650,13.869060)(0.950,0.950)
\psset{linewidth=0.040} \psset{linestyle=solid}
\psset{linestyle=solid} \setlinecaps{0} \setlinejoinmode{0}
\setlinecaps{0} \setlinejoinmode{0} \psset{linestyle=solid}
\newrgbcolor{dialinecolor}{1 1 1}
\psset{linecolor=dialinecolor}
\psellipse*(3.509750,13.873360)(0.687300,0.687300)
\newrgbcolor{dialinecolor}{0 0 0}
\psset{linecolor=dialinecolor}
\psellipse(3.509750,13.873360)(0.687300,0.687300)
\psset{linewidth=0.004000} \setlinecaps{0} \setlinejoinmode{0}
\psset{linestyle=solid}
\newrgbcolor{dialinecolor}{0 0 0}
\psset{linecolor=dialinecolor}
\psellipse(3.509750,13.873360)(0.687300,0.687300)
\psset{linewidth=0.10} \psset{linestyle=solid}
\psset{linestyle=solid} \setlinecaps{0}
\newrgbcolor{dialinecolor}{0 0 0}
\psset{linecolor=dialinecolor}
\psline(3.509750,14.560660)(3.501650,14.325090)
\psset{linewidth=0.040} \psset{linestyle=solid}
\psset{linestyle=solid} \setlinecaps{0}
\newrgbcolor{dialinecolor}{0 0 0}
\psset{linecolor=dialinecolor}
\psline(2.842065,13.960625)(2.855650,14.021060)
\psset{linewidth=0.040} \psset{linestyle=solid}
\setlinejoinmode{0} \setlinecaps{0}
\newrgbcolor{dialinecolor}{0 0 0}
\psset{linecolor=dialinecolor}
\psline(3.044398,14.176793)(2.832258,13.916993)(2.751702,14.242586)
\psset{linewidth=0.040} \psset{linestyle=solid}
\psset{linestyle=solid} \setlinecaps{0}
\newrgbcolor{dialinecolor}{0 0 0}
\psset{linecolor=dialinecolor}
\psline(2.623565,13.618714)(2.597614,13.800370)
\psset{linewidth=0.040} \psset{linestyle=solid}
\setlinejoinmode{0} \setlinecaps{0}
\newrgbcolor{dialinecolor}{0 0 0}
\psset{linecolor=dialinecolor}
\psline(2.735955,13.892640)(2.629889,13.574442)(2.438970,13.850214)
\psset{linewidth=0.040} \psset{linestyle=solid}
\psset{linestyle=solid} \setlinecaps{0} \setlinejoinmode{0}
\setlinecaps{0} \setlinejoinmode{0} \psset{linestyle=solid}
\newrgbcolor{dialinecolor}{1 1 1}
\psset{linecolor=dialinecolor}
\psellipse*(8.364950,13.869060)(0.950,0.950)
\newrgbcolor{dialinecolor}{0 0 0}
\psset{linecolor=dialinecolor}
\psellipse(8.364950,13.869060)(0.950,0.950)
\psset{linewidth=0.004000} \setlinecaps{0} \setlinejoinmode{0}
\psset{linestyle=solid}
\newrgbcolor{dialinecolor}{0 0 0}
\psset{linecolor=dialinecolor}
\psellipse(8.364950,13.869060)(0.950,0.950)
\psset{linewidth=0.040} \psset{linestyle=solid}
\psset{linestyle=solid} \setlinecaps{0} \setlinejoinmode{0}
\setlinecaps{0} \setlinejoinmode{0} \psset{linestyle=solid}
\newrgbcolor{dialinecolor}{1 1 1}
\psset{linecolor=dialinecolor}
\psellipse*(8.354050,13.873360)(0.687300,0.687300)
\newrgbcolor{dialinecolor}{0 0 0}
\psset{linecolor=dialinecolor}
\psellipse(8.354050,13.873360)(0.687300,0.687300)
\psset{linewidth=0.004000} \setlinecaps{0} \setlinejoinmode{0}
\psset{linestyle=solid}
\newrgbcolor{dialinecolor}{0 0 0}
\psset{linecolor=dialinecolor}
\psellipse(8.354050,13.873360)(0.687300,0.687300)
\psset{linewidth=0.10} \psset{linestyle=solid}
\psset{linestyle=solid} \setlinecaps{0}
\newrgbcolor{dialinecolor}{0 0 0}
\psset{linecolor=dialinecolor}
\psline(8.354050,14.560660)(8.345950,14.325090)
\psset{linewidth=0.040} \psset{linestyle=solid}
\psset{linestyle=solid} \setlinecaps{0}
\newrgbcolor{dialinecolor}{0 0 0}
\psset{linecolor=dialinecolor}
\psline(7.686365,13.960625)(7.699950,14.021060)
\psset{linewidth=0.040} \psset{linestyle=solid}
\setlinejoinmode{0} \setlinecaps{0}
\newrgbcolor{dialinecolor}{0 0 0}
\psset{linecolor=dialinecolor}
\psline(7.888698,14.176793)(7.676558,13.916993)(7.596002,14.242586)
\psset{linewidth=0.040} \psset{linestyle=solid}
\psset{linestyle=solid} \setlinecaps{0}
\newrgbcolor{dialinecolor}{0 0 0}
\psset{linecolor=dialinecolor}
\psline(7.865192,14.664640)(7.731414,14.591670)
\psset{linewidth=0.040} \psset{linestyle=solid}
\setlinejoinmode{0} \setlinecaps{0}
\newrgbcolor{dialinecolor}{0 0 0}
\psset{linecolor=dialinecolor}
\psline(7.569257,14.674084)(7.904453,14.686055)(7.712912,14.410715)
\setfont{Helvetica}{0.40}
\newrgbcolor{dialinecolor}{0 0 0}
\psset{linecolor=dialinecolor}
\rput[l](8.273744,10.075960){\scalebox{1 -1}{3}}
\setfont{Helvetica}{0.40}
\newrgbcolor{dialinecolor}{0 0 0}
\psset{linecolor=dialinecolor}
\rput[l](3.427514,13.954770){\scalebox{1 -1}{1}}
\setfont{Helvetica}{0.40}
\newrgbcolor{dialinecolor}{0 0 0}
\psset{linecolor=dialinecolor}
\rput[l](8.273744,13.916660){\scalebox{1 -1}{4}}
\psset{linewidth=0.010} \psset{linestyle=dashed,dash=1 1}
\psset{linestyle=dashed,dash=0.10 0.10} \setlinecaps{0}
\newrgbcolor{dialinecolor}{0 0 0}
\psset{linecolor=dialinecolor}
\psclip{\pswedge[linestyle=none,fillstyle=none](6.023207,12.534937){1.019856}{153.681055}{292.847536}}
\psellipse(6.023207,12.534937)(0.721147,0.721147)
\endpsclip
\setfont{Helvetica}{0.40}
\newrgbcolor{dialinecolor}{0 0 0}
\psset{linecolor=dialinecolor}
\rput[l](1.356624,11.340110){\scalebox{1 -1}{2'}}
\setfont{Helvetica}{0.40}
\newrgbcolor{dialinecolor}{0 0 0}
\psset{linecolor=dialinecolor}
\rput[l](1.314164,12.507270){\scalebox{1 -1}{1'}}
\setfont{Helvetica}{0.40}
\newrgbcolor{dialinecolor}{0 0 0}
\psset{linecolor=dialinecolor}
\rput[l](3.414004,10.066310){\scalebox{1 -1}{2}}
\psset{linewidth=0.040} \psset{linestyle=solid}
\psset{linestyle=solid} \setlinecaps{0}
\newrgbcolor{dialinecolor}{0 0 0}
\psset{linecolor=dialinecolor}
\psclip{\pswedge[linestyle=none,fillstyle=none](2.368791,7.885821){5.051969}{78.677556}{99.686745}}
\psellipse(2.368791,7.885821)(3.572281,3.572281)
\endpsclip
\psset{linewidth=0.040} \psset{linestyle=solid}
\setlinejoinmode{0} \setlinecaps{0}
\newrgbcolor{dialinecolor}{0 0 0}
\psset{linecolor=dialinecolor}
\psline(2.941163,11.520778)(3.113726,11.378574)(2.896425,11.325847)
\psset{linewidth=0.040} \psset{linestyle=solid}
\psset{linestyle=solid} \setlinecaps{0}
\newrgbcolor{dialinecolor}{0 0 0}
\psset{linecolor=dialinecolor}
\psclip{\pswedge[linestyle=none,fillstyle=none](2.430507,8.731220){5.413174}{80.389397}{100.738158}}
\psellipse(2.430507,8.731220)(3.827692,3.827692)
\endpsclip
\psset{linewidth=0.040} \psset{linestyle=solid}
\setlinejoinmode{0} \setlinecaps{0}
\newrgbcolor{dialinecolor}{0 0 0}
\psset{linecolor=dialinecolor}
\psline(2.936424,12.633216)(3.113429,12.496580)(2.897918,12.436958)
\newrgbcolor{dialinecolor}{1 1 1}
\psset{linecolor=dialinecolor}
\pspolygon*(4.334614,9.863170)(4.334614,10.191270)(4.662714,10.191270)(4.662714,9.863170)
\psset{linewidth=0.10} \psset{linestyle=solid}
\psset{linestyle=solid} \setlinejoinmode{0}
\newrgbcolor{dialinecolor}{1 1 1}
\psset{linecolor=dialinecolor}
\pspolygon(4.334614,9.863170)(4.334614,10.191270)(4.662714,10.191270)(4.662714,9.863170)
\newrgbcolor{dialinecolor}{1 1 1}
\psset{linecolor=dialinecolor}
\pspolygon*(3.328904,8.925010)(3.328904,9.253110)(3.657004,9.253110)(3.657004,8.925010)
\psset{linewidth=0.10} \psset{linestyle=solid}
\psset{linestyle=solid} \setlinejoinmode{0}
\newrgbcolor{dialinecolor}{1 1 1}
\psset{linecolor=dialinecolor}
\pspolygon(3.328904,8.925010)(3.328904,9.253110)(3.657004,9.253110)(3.657004,8.925010)
\newrgbcolor{dialinecolor}{1 1 1}
\psset{linecolor=dialinecolor}
\pspolygon*(2.383204,9.928610)(2.383204,10.256710)(2.674814,10.256710)(2.674814,9.928610)
\psset{linewidth=0.10} \psset{linestyle=solid}
\psset{linestyle=solid} \setlinejoinmode{0}
\newrgbcolor{dialinecolor}{1 1 1}
\psset{linecolor=dialinecolor}
\pspolygon(2.383204,9.928610)(2.383204,10.256710)(2.674814,10.256710)(2.674814,9.928610)
\newrgbcolor{dialinecolor}{1 1 1}
\psset{linecolor=dialinecolor}
\pspolygon*(7.980204,10.797110)(7.980204,11.125210)(8.308304,11.125210)(8.308304,10.797110)
\psset{linewidth=0.10} \psset{linestyle=solid}
\psset{linestyle=solid} \setlinejoinmode{0}
\newrgbcolor{dialinecolor}{1 1 1}
\psset{linecolor=dialinecolor}
\pspolygon(7.980204,10.797110)(7.980204,11.125210)(8.308304,11.125210)(8.308304,10.797110)
\newrgbcolor{dialinecolor}{1 1 1}
\psset{linecolor=dialinecolor}
\pspolygon*(3.328904,12.719570)(3.328904,13.055210)(3.657004,13.055210)(3.657004,12.719570)
\psset{linewidth=0.10} \psset{linestyle=solid}
\psset{linestyle=solid} \setlinejoinmode{0}
\newrgbcolor{dialinecolor}{1 1 1}
\psset{linecolor=dialinecolor}
\pspolygon(3.328904,12.719570)(3.328904,13.055210)(3.657004,13.055210)(3.657004,12.719570)
\newrgbcolor{dialinecolor}{1 1 1}
\psset{linecolor=dialinecolor}
\pspolygon*(3.328904,10.855010)(3.328904,11.183110)(3.657004,11.183110)(3.657004,10.855010)
\psset{linewidth=0.10} \psset{linestyle=solid}
\psset{linestyle=solid} \setlinejoinmode{0}
\newrgbcolor{dialinecolor}{1 1 1}
\psset{linecolor=dialinecolor}
\pspolygon(3.328904,10.855010)(3.328904,11.183110)(3.657004,11.183110)(3.657004,10.855010)
\newrgbcolor{dialinecolor}{1 1 1}
\psset{linecolor=dialinecolor}
\pspolygon*(4.274604,13.730710)(4.274604,14.058810)(4.602704,14.058810)(4.602704,13.730710)
\psset{linewidth=0.10} \psset{linestyle=solid}
\psset{linestyle=solid} \setlinejoinmode{0}
\newrgbcolor{dialinecolor}{1 1 1}
\psset{linecolor=dialinecolor}
\pspolygon(4.274604,13.730710)(4.274604,14.058810)(4.602704,14.058810)(4.602704,13.730710)
\psset{linewidth=0.010} \psset{linestyle=dashed,dash=1 1}
\psset{linestyle=dashed,dash=0.10 0.10} \setlinecaps{0}
\newrgbcolor{dialinecolor}{0 0 0}
\psset{linecolor=dialinecolor}
\psclip{\pswedge[linestyle=none,fillstyle=none](4.433205,11.976449){2.703346}{356.818816}{97.096499}}
\psellipse(4.433205,11.976449)(1.911554,1.911554)
\endpsclip
\newrgbcolor{dialinecolor}{1 1 1}
\psset{linecolor=dialinecolor}
\pspolygon*(7.188904,13.750010)(7.188904,14.078110)(7.517004,14.078110)(7.517004,13.750010)
\psset{linewidth=0.10} \psset{linestyle=solid}
\psset{linestyle=solid} \setlinejoinmode{0}
\newrgbcolor{dialinecolor}{1 1 1}
\psset{linecolor=dialinecolor}
\pspolygon(7.188904,13.750010)(7.188904,14.078110)(7.517004,14.078110)(7.517004,13.750010)
\psset{linewidth=0.040} \psset{linestyle=solid}
\psset{linestyle=solid} \setlinecaps{0}
\newrgbcolor{dialinecolor}{0 0 0}
\psset{linecolor=dialinecolor}
\psclip{\pswedge[linestyle=none,fillstyle=none](4.352090,12.005012){3.032377}{0.527958}{87.887927}}
\psellipse(4.352090,12.005012)(2.144214,2.144214)
\endpsclip
\psset{linewidth=0.040} \psset{linestyle=solid}
\psset{linestyle=solid} \setlinecaps{0}
\newrgbcolor{dialinecolor}{0 0 0}
\psset{linecolor=dialinecolor}
\psclip{\pswedge[linestyle=none,fillstyle=none](4.249658,11.826151){2.706683}{7.119814}{84.332902}}
\psellipse(4.249658,11.826151)(1.913914,1.913914)
\endpsclip
\newrgbcolor{dialinecolor}{1 1 1}
\psset{linecolor=dialinecolor}
\pspolygon*(5.396114,12.854670)(5.396114,13.549470)(6.110214,13.549470)(6.110214,12.854670)
\psset{linewidth=0.10} \psset{linestyle=solid}
\psset{linestyle=solid} \setlinejoinmode{0}
\newrgbcolor{dialinecolor}{1 1 1}
\psset{linecolor=dialinecolor}
\pspolygon(5.396114,12.854670)(5.396114,13.549470)(6.110214,13.549470)(6.110214,12.854670)
\psset{linewidth=0.010} \psset{linestyle=dashed,dash=1 1}
\psset{linestyle=dashed,dash=0.10 0.10} \setlinecaps{0}
\newrgbcolor{dialinecolor}{0 0 0}
\psset{linecolor=dialinecolor}
\psclip{\pswedge[linestyle=none,fillstyle=none](7.793227,10.770868){4.494712}{91.981201}{139.490543}}
\psellipse(7.793227,10.770868)(3.178241,3.178241)
\endpsclip
\psset{linewidth=0.040} \psset{linestyle=solid}
\psset{linestyle=solid} \setlinecaps{0}
\newrgbcolor{dialinecolor}{0 0 0}
\psset{linecolor=dialinecolor}
\psclip{\pswedge[linestyle=none,fillstyle=none](7.855991,10.528819){4.497691}{97.130521}{136.461642}}
\psellipse(7.855991,10.528819)(3.180348,3.180348)
\endpsclip
\psset{linewidth=0.040} \psset{linestyle=solid}
\psset{linestyle=solid} \setlinecaps{0}
\newrgbcolor{dialinecolor}{0 0 0}
\psset{linecolor=dialinecolor}
\psclip{\pswedge[linestyle=none,fillstyle=none](7.945296,10.561055){5.141861}{97.344402}{138.495130}}
\psellipse(7.945296,10.561055)(3.635845,3.635845)
\endpsclip
\psset{linewidth=0.10} \psset{linestyle=solid}
\psset{linestyle=solid} \setlinecaps{0}
\newrgbcolor{dialinecolor}{0 0 0}
\psset{linecolor=dialinecolor}
\psline(3.524014,15.035570)(3.520650,14.819060)
\psset{linewidth=0.10} \psset{linestyle=solid}
\psset{linestyle=solid} \setlinecaps{0}
\newrgbcolor{dialinecolor}{0 0 0}
\psset{linecolor=dialinecolor}
\psline(8.375993,15.029469)(8.364950,14.819060)
\psset{linewidth=0.040} \psset{linestyle=solid}
\psset{linestyle=solid} \setlinecaps{0}
\newrgbcolor{dialinecolor}{0 0 0}
\psset{linecolor=dialinecolor}
\psclip{\pswedge[linestyle=none,fillstyle=none](5.989859,12.536498){0.680040}{162.515738}{289.303023}}
\psellipse(5.989859,12.536498)(0.480861,0.480861)
\endpsclip
\psset{linewidth=0.040} \psset{linestyle=solid}
\psset{linestyle=solid} \setlinecaps{0}
\newrgbcolor{dialinecolor}{0 0 0}
\psset{linecolor=dialinecolor}
\psclip{\pswedge[linestyle=none,fillstyle=none](5.919690,12.487699){1.215129}{144.244392}{278.879232}}
\psellipse(5.919690,12.487699)(0.859226,0.859226)
\endpsclip
\psset{linewidth=0.040} \psset{linestyle=solid}
\psset{linestyle=solid} \setlinecaps{0}
\newrgbcolor{dialinecolor}{0 0 0}
\psset{linecolor=dialinecolor}
\psclip{\pswedge[linestyle=none,fillstyle=none](6.500957,9.862612){3.085657}{29.127907}{91.138259}}
\psellipse(6.500957,9.862612)(2.181889,2.181889)
\endpsclip
\psset{linewidth=0.040} \psset{linestyle=solid}
\psset{linestyle=solid} \setlinecaps{0}
\newrgbcolor{dialinecolor}{0 0 0}
\psset{linecolor=dialinecolor}
\psclip{\pswedge[linestyle=none,fillstyle=none](6.656710,10.155400){2.055782}{29.298902}{96.338876}}
\psellipse(6.656710,10.155400)(1.453657,1.453657)
\endpsclip
\psset{linewidth=0.040} \psset{linestyle=solid}
\psset{linestyle=solid} \setlinecaps{0}
\newrgbcolor{dialinecolor}{0 0 0}
\psset{linecolor=dialinecolor}
\psclip{\pswedge[linestyle=none,fillstyle=none](4.099193,12.280008){3.531148}{277.639118}{343.739630}}
\psellipse(4.099193,12.280008)(2.496898,2.496898)
\endpsclip
\psset{linewidth=0.040} \psset{linestyle=solid}
\psset{linestyle=solid} \setlinecaps{0}
\newrgbcolor{dialinecolor}{0 0 0}
\psset{linecolor=dialinecolor}
\psclip{\pswedge[linestyle=none,fillstyle=none](3.935806,12.534124){3.279793}{282.331761}{338.319911}}
\psellipse(3.935806,12.534124)(2.319164,2.319164)
\endpsclip
\psset{linewidth=0.040} \psset{linestyle=solid}
\psset{linestyle=solid} \setlinecaps{0}
\newrgbcolor{dialinecolor}{0 0 0}
\psset{linecolor=dialinecolor}
\psclip{\pswedge[linestyle=none,fillstyle=none](2.597986,9.204957){1.669783}{87.207257}{355.409971}}
\psellipse(2.597986,9.204957)(1.180715,1.180715)
\endpsclip
\psset{linewidth=0.040} \psset{linestyle=solid}
\psset{linestyle=solid} \setlinecaps{0}
\newrgbcolor{dialinecolor}{0 0 0}
\psset{linecolor=dialinecolor}
\psclip{\pswedge[linestyle=none,fillstyle=none](2.570015,9.179267){1.022519}{87.812410}{357.606370}}
\psellipse(2.570015,9.179267)(0.723030,0.723030)
\endpsclip
\psset{linewidth=0.040} \psset{linestyle=solid}
\psset{linestyle=solid} \setlinecaps{0}
\newrgbcolor{dialinecolor}{0 0 0}
\psset{linecolor=dialinecolor}
\psline(3.292414,10.924670)(3.292414,12.970470)
\psset{linewidth=0.040} \psset{linestyle=solid}
\psset{linestyle=solid} \setlinecaps{0}
\newrgbcolor{dialinecolor}{0 0 0}
\psset{linecolor=dialinecolor}
\psline(3.762104,10.921510)(3.762104,12.967310)
\newrgbcolor{dialinecolor}{1 1 1}
\psset{linecolor=dialinecolor}
\pspolygon*(3.155204,11.754570)(3.155204,12.217770)(3.869304,12.217770)(3.869304,11.754570)
\psset{linewidth=0.10} \psset{linestyle=solid}
\psset{linestyle=solid} \setlinejoinmode{0}
\newrgbcolor{dialinecolor}{1 1 1}
\psset{linecolor=dialinecolor}
\pspolygon(3.155204,11.754570)(3.155204,12.217770)(3.869304,12.217770)(3.869304,11.754570)
\newrgbcolor{dialinecolor}{1 1 1}
\psset{linecolor=dialinecolor}
\pspolygon*(7.227504,11.136970)(7.227504,11.928270)(7.654214,11.928270)(7.654214,11.136970)
\psset{linewidth=0.10} \psset{linestyle=solid}
\psset{linestyle=solid} \setlinejoinmode{0}
\newrgbcolor{dialinecolor}{1 1 1}
\psset{linecolor=dialinecolor}
\pspolygon(7.227504,11.136970)(7.227504,11.928270)(7.654214,11.928270)(7.654214,11.136970)
\newrgbcolor{dialinecolor}{1 1 1}
\psset{linecolor=dialinecolor}
\pspolygon*(4.913614,10.519370)(4.913614,10.943970)(6.476914,10.943970)(6.476914,10.519370)
\psset{linewidth=0.10} \psset{linestyle=solid}
\psset{linestyle=solid} \setlinejoinmode{0}
\newrgbcolor{dialinecolor}{1 1 1}
\psset{linecolor=dialinecolor}
\pspolygon(4.913614,10.519370)(4.913614,10.943970)(6.476914,10.943970)(6.476914,10.519370)
\psset{linewidth=0.010} \psset{linestyle=dashed,dash=1 1}
\psset{linestyle=dashed,dash=0.10 0.10} \setlinecaps{0}
\newrgbcolor{dialinecolor}{0 0 0}
\psset{linecolor=dialinecolor}
\psclip{\pswedge[linestyle=none,fillstyle=none](6.424155,9.690144){3.056637}{27.874423}{91.159623}}
\psellipse(6.424155,9.690144)(2.161369,2.161369)
\endpsclip
\psset{linewidth=0.010} \psset{linestyle=dashed,dash=0.10 0.10}
\psset{linestyle=dashed,dash=0.10 0.10} \setlinecaps{0}
\newrgbcolor{dialinecolor}{0 0 0}
\psset{linecolor=dialinecolor}
\psclip{\pswedge[linestyle=none,fillstyle=none](4.055652,12.361000){3.318620}{273.382840}{346.966030}}
\psellipse(4.055652,12.361000)(2.346619,2.346619)
\endpsclip
\psset{linewidth=0.010} \psset{linestyle=dashed,dash=0.10 0.10}
\psset{linestyle=dashed,dash=0.10 0.10} \setlinecaps{0}
\newrgbcolor{dialinecolor}{0 0 0}
\psset{linecolor=dialinecolor}
\psline(3.506820,10.705770)(3.509750,13.186060)
\newrgbcolor{dialinecolor}{1 1 1}
\psset{linecolor=dialinecolor}
\pspolygon*(1.130814,8.589370)(1.130814,9.110470)(2.385314,9.110470)(2.385314,8.589370)
\psset{linewidth=0.10} \psset{linestyle=solid}
\psset{linestyle=solid} \setlinejoinmode{0}
\newrgbcolor{dialinecolor}{1 1 1}
\psset{linecolor=dialinecolor}
\pspolygon(1.130814,8.589370)(1.130814,9.110470)(2.385314,9.110470)(2.385314,8.589370)
\psset{linewidth=0.010} \psset{linestyle=dashed,dash=1 1}
\psset{linestyle=dashed,dash=0.10 0.10} \setlinecaps{0}
\newrgbcolor{dialinecolor}{0 0 0}
\psset{linecolor=dialinecolor}
\psclip{\pswedge[linestyle=none,fillstyle=none](2.575445,9.186547){1.332948}{73.944535}{8.826319}}
\psellipse(2.575445,9.186547)(0.942537,0.942537)
\endpsclip
\psset{linewidth=0.040} \psset{linestyle=solid}
\psset{linestyle=solid} \setlinecaps{0}
\newrgbcolor{dialinecolor}{0 0 0}
\psset{linecolor=dialinecolor}
\psclip{\pswedge[linestyle=none,fillstyle=none](6.180716,10.539611){2.366028}{20.213249}{56.101204}}
\psellipse(6.180716,10.539611)(1.673035,1.673035)
\endpsclip
\psset{linewidth=0.040} \psset{linestyle=solid}
\psset{linestyle=solid} \setlinecaps{0}
\newrgbcolor{dialinecolor}{0 0 0}
\psset{linecolor=dialinecolor}
\psclip{\pswedge[linestyle=none,fillstyle=none](7.511464,11.252084){0.620021}{56.926993}{144.981132}}
\psellipse(7.511464,11.252084)(0.438421,0.438421)
\endpsclip
\psset{linewidth=0.040} \psset{linestyle=solid}
\psset{linestyle=solid} \setlinecaps{0}
\newrgbcolor{dialinecolor}{0 0 0}
\psset{linecolor=dialinecolor}
\psclip{\pswedge[linestyle=none,fillstyle=none](3.065670,11.646675){0.980637}{5.735333}{69.217026}}
\psellipse(3.065670,11.646675)(0.693415,0.693415)
\endpsclip
\psset{linewidth=0.040} \psset{linestyle=solid}
\psset{linestyle=solid} \setlinecaps{0}
\newrgbcolor{dialinecolor}{0 0 0}
\psset{linecolor=dialinecolor}
\psclip{\pswedge[linestyle=none,fillstyle=none](4.428498,11.322594){1.717467}{123.646905}{162.059722}}
\psellipse(4.428498,11.322594)(1.214433,1.214433)
\endpsclip
\psset{linewidth=0.040} \psset{linestyle=solid}
\psset{linestyle=solid} \setlinecaps{0}
\newrgbcolor{dialinecolor}{0 0 0}
\psset{linecolor=dialinecolor}
\psclip{\pswedge[linestyle=none,fillstyle=none](0.920102,7.748236){2.173609}{31.476717}{69.478021}}
\psellipse(0.920102,7.748236)(1.536974,1.536974)
\endpsclip
\psset{linewidth=0.040} \psset{linestyle=solid}
\psset{linestyle=solid} \setlinecaps{0}
\newrgbcolor{dialinecolor}{0 0 0}
\psset{linecolor=dialinecolor}
\psclip{\pswedge[linestyle=none,fillstyle=none](2.163501,8.722485){0.815098}{119.063878}{197.333303}}
\psellipse(2.163501,8.722485)(0.576361,0.576361)
\endpsclip
\psset{linewidth=0.040} \psset{linestyle=solid}
\psset{linestyle=solid} \setlinecaps{0}
\newrgbcolor{dialinecolor}{0 0 0}
\psset{linecolor=dialinecolor}
\psclip{\pswedge[linestyle=none,fillstyle=none](5.397483,10.689027){0.642699}{332.725402}{50.653336}}
\psellipse(5.397483,10.689027)(0.454457,0.454457)
\endpsclip
\psset{linewidth=0.040} \psset{linestyle=solid}
\psset{linestyle=solid} \setlinecaps{0}
\newrgbcolor{dialinecolor}{0 0 0}
\psset{linecolor=dialinecolor}
\psclip{\pswedge[linestyle=none,fillstyle=none](8.385134,4.577867){9.653383}{108.269202}{119.819245}}
\psellipse(8.385134,4.577867)(6.825973,6.825973)
\endpsclip
\psset{linewidth=0.040} \psset{linestyle=solid}
\psset{linestyle=solid} \setlinecaps{0}
\newrgbcolor{dialinecolor}{0 0 0}
\psset{linecolor=dialinecolor}
\psline(4.230140,10.602945)(4.141614,10.693070)
\psset{linewidth=0.040} \psset{linestyle=solid}
\setlinejoinmode{0} \setlinecaps{0}
\newrgbcolor{dialinecolor}{0 0 0}
\psset{linecolor=dialinecolor}
\psline(4.158265,10.890175)(4.261479,10.571040)(3.944243,10.679950)
\psset{linewidth=0.040} \psset{linestyle=solid}
\psset{linestyle=solid} \setlinecaps{0}
\newrgbcolor{dialinecolor}{0 0 0}
\psset{linecolor=dialinecolor}
\psline(9.298643,13.957004)(9.256114,14.186370)
\psset{linewidth=0.040} \psset{linestyle=solid}
\setlinejoinmode{0} \setlinecaps{0}
\newrgbcolor{dialinecolor}{0 0 0}
\psset{linecolor=dialinecolor}
\psline(9.399588,14.235351)(9.306797,13.913032)(9.104616,14.180657)
\psset{linewidth=0.040} \psset{linestyle=solid}
\psset{linestyle=solid} \setlinecaps{0}
\newrgbcolor{dialinecolor}{0 0 0}
\psset{linecolor=dialinecolor}
\psline(4.181212,13.192183)(4.257414,13.279270)
\psset{linewidth=0.040} \psset{linestyle=solid}
\setlinejoinmode{0} \setlinecaps{0}
\newrgbcolor{dialinecolor}{0 0 0}
\psset{linecolor=dialinecolor}
\psline(4.462201,13.285524)(4.151763,13.158526)(4.236428,13.483075)
\psset{linewidth=0.040} \psset{linestyle=solid}
\psset{linestyle=solid} \setlinecaps{0}
\newrgbcolor{dialinecolor}{0 0 0}
\psset{linecolor=dialinecolor}
\psline(2.837280,10.673049)(3.002914,10.808870)
\psset{linewidth=0.040} \psset{linestyle=solid}
\setlinejoinmode{0} \setlinecaps{0}
\newrgbcolor{dialinecolor}{0 0 0}
\psset{linecolor=dialinecolor}
\psline(3.129791,10.718927)(2.802699,10.644692)(2.939566,10.950907)
\psset{linewidth=0.040} \psset{linestyle=solid}
\psset{linestyle=solid} \setlinecaps{0}
\newrgbcolor{dialinecolor}{0 0 0}
\psset{linecolor=dialinecolor}
\psline(4.302622,9.519696)(4.218814,9.399970)
\psset{linewidth=0.040} \psset{linestyle=solid}
\setlinejoinmode{0} \setlinecaps{0}
\newrgbcolor{dialinecolor}{0 0 0}
\psset{linecolor=dialinecolor}
\psline(4.033344,9.396583)(4.328268,9.556333)(4.279114,9.224544)
}\endpspicture \]

We also remark that Sullivan chord diagrams are considered up to
abstract isomorphism of the thickened surfaces respecting all
labelling and orientations (see figure below).

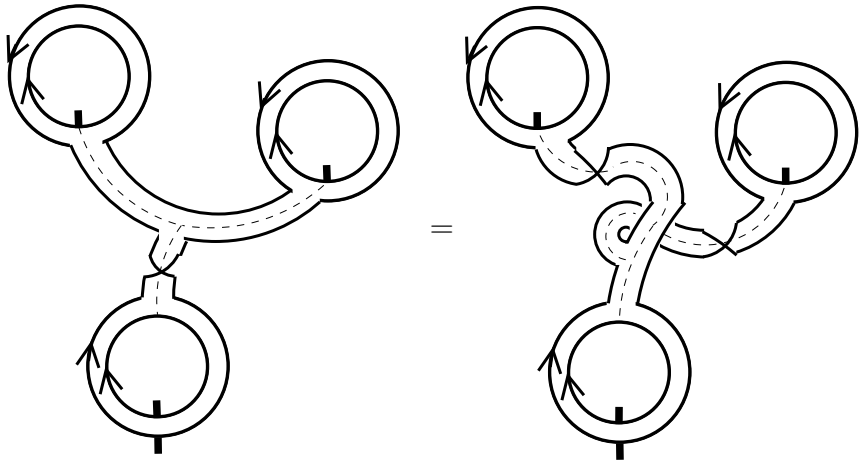
\[ \pspicture(1.602736,-15.575879)(13.211200,-9.449470)
\scalebox{1 -1}{
\newrgbcolor{dialinecolor}{0 0 0}
\psset{linecolor=dialinecolor}
\newrgbcolor{diafillcolor}{1 1 1}
\psset{fillcolor=diafillcolor} \psset{linewidth=0.040}
\psset{linestyle=solid} \psset{linestyle=solid} \setlinecaps{0}
\setlinejoinmode{0} \setlinecaps{0} \setlinejoinmode{0}
\psset{linestyle=solid}
\newrgbcolor{dialinecolor}{1 1 1}
\psset{linecolor=dialinecolor}
\psellipse*(2.842220,10.419470)(0.950,0.950)
\newrgbcolor{dialinecolor}{0 0 0}
\psset{linecolor=dialinecolor}
\psellipse(2.842220,10.419470)(0.950,0.950)
\psset{linewidth=0.004000} \setlinecaps{0} \setlinejoinmode{0}
\psset{linestyle=solid}
\newrgbcolor{dialinecolor}{0 0 0}
\psset{linecolor=dialinecolor}
\psellipse(2.842220,10.419470)(0.950,0.950)
\psset{linewidth=0.040} \psset{linestyle=solid}
\psset{linestyle=solid} \setlinecaps{0} \setlinejoinmode{0}
\setlinecaps{0} \setlinejoinmode{0} \psset{linestyle=solid}
\newrgbcolor{dialinecolor}{1 1 1}
\psset{linecolor=dialinecolor}
\psellipse*(2.831320,10.423770)(0.687300,0.687300)
\newrgbcolor{dialinecolor}{0 0 0}
\psset{linecolor=dialinecolor}
\psellipse(2.831320,10.423770)(0.687300,0.687300)
\psset{linewidth=0.004000} \setlinecaps{0} \setlinejoinmode{0}
\psset{linestyle=solid}
\newrgbcolor{dialinecolor}{0 0 0}
\psset{linecolor=dialinecolor}
\psellipse(2.831320,10.423770)(0.687300,0.687300)
\psset{linewidth=0.10} \psset{linestyle=solid}
\psset{linestyle=solid} \setlinecaps{0}
\newrgbcolor{dialinecolor}{0 0 0}
\psset{linecolor=dialinecolor}
\psline(2.831320,11.111100)(2.823220,10.875500)
\psset{linewidth=0.040} \psset{linestyle=solid}
\psset{linestyle=solid} \setlinecaps{0}
\newrgbcolor{dialinecolor}{0 0 0}
\psset{linecolor=dialinecolor}
\psline(2.163635,10.511065)(2.177220,10.571500)
\psset{linewidth=0.040} \psset{linestyle=solid}
\setlinejoinmode{0} \setlinecaps{0}
\newrgbcolor{dialinecolor}{0 0 0}
\psset{linecolor=dialinecolor}
\psline(2.365968,10.727233)(2.153828,10.467433)(2.073272,10.793026)
\psset{linewidth=0.040} \psset{linestyle=solid}
\psset{linestyle=solid} \setlinecaps{0}
\newrgbcolor{dialinecolor}{0 0 0}
\psset{linecolor=dialinecolor}
\psline(1.920988,10.167988)(1.999310,9.998270)
\psset{linewidth=0.040} \psset{linestyle=solid}
\setlinejoinmode{0} \setlinecaps{0}
\newrgbcolor{dialinecolor}{0 0 0}
\psset{linecolor=dialinecolor}
\psline(1.891757,9.873348)(1.902249,10.208594)(2.164151,9.999053)
\psset{linewidth=0.040} \psset{linestyle=solid}
\psset{linestyle=solid} \setlinecaps{0} \setlinejoinmode{0}
\setlinecaps{0} \setlinejoinmode{0} \psset{linestyle=solid}
\newrgbcolor{dialinecolor}{1 1 1}
\psset{linecolor=dialinecolor}
\psellipse*(6.145450,11.147800)(0.950,0.950)
\newrgbcolor{dialinecolor}{0 0 0}
\psset{linecolor=dialinecolor}
\psellipse(6.145450,11.147800)(0.950,0.950)
\psset{linewidth=0.004000} \setlinecaps{0} \setlinejoinmode{0}
\psset{linestyle=solid}
\newrgbcolor{dialinecolor}{0 0 0}
\psset{linecolor=dialinecolor}
\psellipse(6.145450,11.147800)(0.950,0.950)
\psset{linewidth=0.040} \psset{linestyle=solid}
\psset{linestyle=solid} \setlinecaps{0} \setlinejoinmode{0}
\setlinecaps{0} \setlinejoinmode{0} \psset{linestyle=solid}
\newrgbcolor{dialinecolor}{1 1 1}
\psset{linecolor=dialinecolor}
\psellipse*(6.134550,11.152100)(0.687300,0.687300)
\newrgbcolor{dialinecolor}{0 0 0}
\psset{linecolor=dialinecolor}
\psellipse(6.134550,11.152100)(0.687300,0.687300)
\psset{linewidth=0.004000} \setlinecaps{0} \setlinejoinmode{0}
\psset{linestyle=solid}
\newrgbcolor{dialinecolor}{0 0 0}
\psset{linecolor=dialinecolor}
\psellipse(6.134550,11.152100)(0.687300,0.687300)
\psset{linewidth=0.10} \psset{linestyle=solid}
\psset{linestyle=solid} \setlinecaps{0}
\newrgbcolor{dialinecolor}{0 0 0}
\psset{linecolor=dialinecolor}
\psline(6.134550,11.839400)(6.126450,11.603800)
\psset{linewidth=0.040} \psset{linestyle=solid}
\psset{linestyle=solid} \setlinecaps{0}
\newrgbcolor{dialinecolor}{0 0 0}
\psset{linecolor=dialinecolor}
\psline(5.466865,11.239365)(5.480450,11.299800)
\psset{linewidth=0.040} \psset{linestyle=solid}
\setlinejoinmode{0} \setlinecaps{0}
\newrgbcolor{dialinecolor}{0 0 0}
\psset{linecolor=dialinecolor}
\psline(5.669198,11.455533)(5.457058,11.195733)(5.376502,11.521326)
\psset{linewidth=0.040} \psset{linestyle=solid}
\psset{linestyle=solid} \setlinecaps{0}
\newrgbcolor{dialinecolor}{0 0 0}
\psset{linecolor=dialinecolor}
\psline(5.275220,10.783872)(5.302120,10.717300)
\psset{linewidth=0.040} \psset{linestyle=solid}
\setlinejoinmode{0} \setlinecaps{0}
\newrgbcolor{dialinecolor}{0 0 0}
\psset{linecolor=dialinecolor}
\psline(5.231785,10.490988)(5.258465,10.825336)(5.509935,10.603383)
\newrgbcolor{dialinecolor}{1 1 1}
\psset{linecolor=dialinecolor}
\pspolygon*(5.724210,11.897200)(5.724210,12.225300)(5.985490,12.225300)(5.985490,11.897200)
\psset{linewidth=0.10} \psset{linestyle=solid}
\psset{linestyle=solid} \setlinejoinmode{0}
\newrgbcolor{dialinecolor}{1 1 1}
\psset{linecolor=dialinecolor}
\pspolygon(5.724210,11.897200)(5.724210,12.225300)(5.985490,12.225300)(5.985490,11.897200)
\setfont{Helvetica}{0.40}
\newrgbcolor{dialinecolor}{0 0 0}
\psset{linecolor=dialinecolor}
\rput[l](7.497060,12.479200){\scalebox{1 -1}{=}}
\psset{linewidth=0.040} \psset{linestyle=solid}
\psset{linestyle=solid} \setlinecaps{0} \setlinejoinmode{0}
\setlinecaps{0} \setlinejoinmode{0} \psset{linestyle=solid}
\newrgbcolor{dialinecolor}{1 1 1}
\psset{linecolor=dialinecolor}
\psellipse*(3.884970,14.278500)(0.950,0.950)
\newrgbcolor{dialinecolor}{0 0 0}
\psset{linecolor=dialinecolor}
\psellipse(3.884970,14.278500)(0.950,0.950)
\psset{linewidth=0.004000} \setlinecaps{0} \setlinejoinmode{0}
\psset{linestyle=solid}
\newrgbcolor{dialinecolor}{0 0 0}
\psset{linecolor=dialinecolor}
\psellipse(3.884970,14.278500)(0.950,0.950)
\psset{linewidth=0.040} \psset{linestyle=solid}
\psset{linestyle=solid} \setlinecaps{0} \setlinejoinmode{0}
\setlinecaps{0} \setlinejoinmode{0} \psset{linestyle=solid}
\newrgbcolor{dialinecolor}{1 1 1}
\psset{linecolor=dialinecolor}
\psellipse*(3.874070,14.282800)(0.687300,0.687300)
\newrgbcolor{dialinecolor}{0 0 0}
\psset{linecolor=dialinecolor}
\psellipse(3.874070,14.282800)(0.687300,0.687300)
\psset{linewidth=0.004000} \setlinecaps{0} \setlinejoinmode{0}
\psset{linestyle=solid}
\newrgbcolor{dialinecolor}{0 0 0}
\psset{linecolor=dialinecolor}
\psellipse(3.874070,14.282800)(0.687300,0.687300)
\psset{linewidth=0.10} \psset{linestyle=solid}
\psset{linestyle=solid} \setlinecaps{0}
\newrgbcolor{dialinecolor}{0 0 0}
\psset{linecolor=dialinecolor}
\psline(3.874070,14.970100)(3.865960,14.734600)
\psset{linewidth=0.040} \psset{linestyle=solid}
\psset{linestyle=solid} \setlinecaps{0}
\newrgbcolor{dialinecolor}{0 0 0}
\psset{linecolor=dialinecolor}
\psline(3.206385,14.370065)(3.219970,14.430500)
\psset{linewidth=0.040} \psset{linestyle=solid}
\setlinejoinmode{0} \setlinecaps{0}
\newrgbcolor{dialinecolor}{0 0 0}
\psset{linecolor=dialinecolor}
\psline(3.408718,14.586233)(3.196578,14.326433)(3.116022,14.652026)
\psset{linewidth=0.040} \psset{linestyle=solid}
\psset{linestyle=solid} \setlinecaps{0}
\newrgbcolor{dialinecolor}{0 0 0}
\psset{linecolor=dialinecolor}
\psline(2.987881,14.028144)(2.961930,14.209800)
\psset{linewidth=0.040} \psset{linestyle=solid}
\setlinejoinmode{0} \setlinecaps{0}
\newrgbcolor{dialinecolor}{0 0 0}
\psset{linecolor=dialinecolor}
\psline(3.100271,14.302070)(2.994205,13.983872)(2.803287,14.259644)
\psset{linewidth=0.10} \psset{linestyle=solid}
\psset{linestyle=solid} \setlinecaps{0}
\newrgbcolor{dialinecolor}{0 0 0}
\psset{linecolor=dialinecolor}
\psline(3.888360,15.445000)(3.884970,15.228500)
\newrgbcolor{dialinecolor}{1 1 1}
\psset{linecolor=dialinecolor}
\pspolygon*(3.715460,13.224900)(3.715460,13.475800)(4.045670,13.475800)(4.045670,13.224900)
\psset{linewidth=0.10} \psset{linestyle=solid}
\psset{linestyle=solid} \setlinejoinmode{0}
\newrgbcolor{dialinecolor}{1 1 1}
\psset{linecolor=dialinecolor}
\pspolygon(3.715460,13.224900)(3.715460,13.475800)(4.045670,13.475800)(4.045670,13.224900)
\newrgbcolor{dialinecolor}{1 1 1}
\psset{linecolor=dialinecolor}
\pspolygon*(2.763880,11.225800)(2.763880,11.553900)(3.095800,11.553900)(3.095800,11.225800)
\psset{linewidth=0.10} \psset{linestyle=solid}
\psset{linestyle=solid} \setlinejoinmode{0}
\newrgbcolor{dialinecolor}{1 1 1}
\psset{linecolor=dialinecolor}
\pspolygon(2.763880,11.225800)(2.763880,11.553900)(3.095800,11.553900)(3.095800,11.225800)
\psset{linewidth=0.040} \psset{linestyle=solid}
\psset{linestyle=solid} \setlinecaps{0}
\newrgbcolor{dialinecolor}{0 0 0}
\psset{linecolor=dialinecolor}
\psclip{\pswedge[linestyle=none,fillstyle=none](5.709757,13.348339){2.314474}{178.744419}{206.801634}}
\psellipse(5.709757,13.348339)(1.636580,1.636580)
\endpsclip
\psset{linewidth=0.040} \psset{linestyle=solid}
\psset{linestyle=solid} \setlinecaps{0}
\newrgbcolor{dialinecolor}{0 0 0}
\psset{linecolor=dialinecolor}
\psclip{\pswedge[linestyle=none,fillstyle=none](5.497250,13.424088){2.604586}{181.888342}{209.628097}}
\psellipse(5.497250,13.424088)(1.841721,1.841721)
\endpsclip
\psset{linewidth=0.040} \psset{linestyle=solid}
\psset{linestyle=solid} \setlinecaps{0}
\newrgbcolor{dialinecolor}{0 0 0}
\psset{linecolor=dialinecolor}
\psclip{\pswedge[linestyle=none,fillstyle=none](4.651512,10.531583){2.982508}{48.397845}{157.265831}}
\psellipse(4.651512,10.531583)(2.108952,2.108952)
\endpsclip
\psset{linewidth=0.040} \psset{linestyle=solid}
\psset{linestyle=solid} \setlinecaps{0}
\newrgbcolor{dialinecolor}{0 0 0}
\psset{linecolor=dialinecolor}
\psclip{\pswedge[linestyle=none,fillstyle=none](4.690410,10.621595){2.338497}{54.037787}{155.827751}}
\psellipse(4.690410,10.621595)(1.653567,1.653567)
\endpsclip
\newrgbcolor{dialinecolor}{1 1 1}
\psset{linecolor=dialinecolor}
\pspolygon*(3.947730,12.350300)(3.947730,12.568550)(4.178440,12.568550)(4.178440,12.350300)
\psset{linewidth=0.10} \psset{linestyle=solid}
\psset{linestyle=solid} \setlinejoinmode{0}
\newrgbcolor{dialinecolor}{1 1 1}
\psset{linecolor=dialinecolor}
\pspolygon(3.947730,12.350300)(3.947730,12.568550)(4.178440,12.568550)(4.178440,12.350300)
\psset{linewidth=0.010} \psset{linestyle=dashed,dash=1 1}
\psset{linestyle=dashed,dash=0.10 0.10} \setlinecaps{0}
\newrgbcolor{dialinecolor}{0 0 0}
\psset{linecolor=dialinecolor}
\psclip{\pswedge[linestyle=none,fillstyle=none](4.587090,10.160077){3.229476}{47.340100}{99.567171}}
\psellipse(4.587090,10.160077)(2.283584,2.283584)
\endpsclip
\psset{linewidth=0.010} \psset{linestyle=dashed,dash=0.10 0.10}
\psset{linestyle=dashed,dash=0.10 0.10} \setlinecaps{0}
\newrgbcolor{dialinecolor}{0 0 0}
\psset{linecolor=dialinecolor}
\psclip{\pswedge[linestyle=none,fillstyle=none](4.663275,10.556108){2.707055}{106.003610}{163.145783}}
\psellipse(4.663275,10.556108)(1.914177,1.914177)
\endpsclip
\newrgbcolor{dialinecolor}{1 1 1}
\psset{linecolor=dialinecolor}
\pspolygon*(3.601320,12.862300)(3.601320,13.038310)(4.261430,13.038310)(4.261430,12.862300)
\psset{linewidth=0.10} \psset{linestyle=solid}
\psset{linestyle=solid} \setlinejoinmode{0}
\newrgbcolor{dialinecolor}{1 1 1}
\psset{linecolor=dialinecolor}
\pspolygon(3.601320,12.862300)(3.601320,13.038310)(4.261430,13.038310)(4.261430,12.862300)
\psset{linewidth=0.040} \psset{linestyle=solid}
\psset{linestyle=solid} \setlinecaps{0}
\newrgbcolor{dialinecolor}{0 0 0}
\psset{linecolor=dialinecolor}
\psclip{\pswedge[linestyle=none,fillstyle=none](3.675567,12.560578){0.780575}{25.281825}{88.353102}}
\psellipse(3.675567,12.560578)(0.551950,0.551950)
\endpsclip
\psset{linewidth=0.040} \psset{linestyle=solid}
\psset{linestyle=solid} \setlinecaps{0}
\newrgbcolor{dialinecolor}{0 0 0}
\psset{linecolor=dialinecolor}
\psclip{\pswedge[linestyle=none,fillstyle=none](4.128632,12.724799){0.536058}{93.836505}{167.924001}}
\psellipse(4.128632,12.724799)(0.379050,0.379050)
\endpsclip
\psset{linewidth=0.010} \psset{linestyle=dashed,dash=1 1}
\psset{linestyle=dashed,dash=0.10 0.10} \setlinecaps{0}
\newrgbcolor{dialinecolor}{0 0 0}
\psset{linecolor=dialinecolor}
\psclip{\pswedge[linestyle=none,fillstyle=none](5.549719,13.404622){2.385051}{173.501288}{216.726945}}
\psellipse(5.549719,13.404622)(1.686486,1.686486)
\endpsclip
\psset{linewidth=0.040} \psset{linestyle=solid}
\psset{linestyle=solid} \setlinecaps{0} \setlinejoinmode{0}
\setlinecaps{0} \setlinejoinmode{0} \psset{linestyle=solid}
\newrgbcolor{dialinecolor}{1 1 1}
\psset{linecolor=dialinecolor}
\psellipse*(8.937930,10.441210)(0.950,0.950)
\newrgbcolor{dialinecolor}{0 0 0}
\psset{linecolor=dialinecolor}
\psellipse(8.937930,10.441210)(0.950,0.950)
\psset{linewidth=0.004000} \setlinecaps{0} \setlinejoinmode{0}
\psset{linestyle=solid}
\newrgbcolor{dialinecolor}{0 0 0}
\psset{linecolor=dialinecolor}
\psellipse(8.937930,10.441210)(0.950,0.950)
\psset{linewidth=0.040} \psset{linestyle=solid}
\psset{linestyle=solid} \setlinecaps{0} \setlinejoinmode{0}
\setlinecaps{0} \setlinejoinmode{0} \psset{linestyle=solid}
\newrgbcolor{dialinecolor}{1 1 1}
\psset{linecolor=dialinecolor}
\psellipse*(8.927030,10.445510)(0.687300,0.687300)
\newrgbcolor{dialinecolor}{0 0 0}
\psset{linecolor=dialinecolor}
\psellipse(8.927030,10.445510)(0.687300,0.687300)
\psset{linewidth=0.004000} \setlinecaps{0} \setlinejoinmode{0}
\psset{linestyle=solid}
\newrgbcolor{dialinecolor}{0 0 0}
\psset{linecolor=dialinecolor}
\psellipse(8.927030,10.445510)(0.687300,0.687300)
\psset{linewidth=0.10} \psset{linestyle=solid}
\psset{linestyle=solid} \setlinecaps{0}
\newrgbcolor{dialinecolor}{0 0 0}
\psset{linecolor=dialinecolor}
\psline(8.927030,11.132800)(8.926350,10.899400)
\psset{linewidth=0.040} \psset{linestyle=solid}
\psset{linestyle=solid} \setlinecaps{0}
\newrgbcolor{dialinecolor}{0 0 0}
\psset{linecolor=dialinecolor}
\psline(8.259345,10.532765)(8.272930,10.593200)
\psset{linewidth=0.040} \psset{linestyle=solid}
\setlinejoinmode{0} \setlinecaps{0}
\newrgbcolor{dialinecolor}{0 0 0}
\psset{linecolor=dialinecolor}
\psline(8.461678,10.748933)(8.249538,10.489133)(8.168982,10.814726)
\psset{linewidth=0.040} \psset{linestyle=solid}
\psset{linestyle=solid} \setlinecaps{0}
\newrgbcolor{dialinecolor}{0 0 0}
\psset{linecolor=dialinecolor}
\psline(8.016702,10.189690)(8.095020,10.020)
\psset{linewidth=0.040} \psset{linestyle=solid}
\setlinejoinmode{0} \setlinecaps{0}
\newrgbcolor{dialinecolor}{0 0 0}
\psset{linecolor=dialinecolor}
\psline(7.987484,9.895048)(7.997961,10.230295)(8.259872,10.020766)
\psset{linewidth=0.040} \psset{linestyle=solid}
\psset{linestyle=solid} \setlinecaps{0} \setlinejoinmode{0}
\setlinecaps{0} \setlinejoinmode{0} \psset{linestyle=solid}
\newrgbcolor{dialinecolor}{1 1 1}
\psset{linecolor=dialinecolor}
\psellipse*(12.241200,11.169500)(0.950,0.950)
\newrgbcolor{dialinecolor}{0 0 0}
\psset{linecolor=dialinecolor}
\psellipse(12.241200,11.169500)(0.950,0.950)
\psset{linewidth=0.004000} \setlinecaps{0} \setlinejoinmode{0}
\psset{linestyle=solid}
\newrgbcolor{dialinecolor}{0 0 0}
\psset{linecolor=dialinecolor}
\psellipse(12.241200,11.169500)(0.950,0.950)
\psset{linewidth=0.040} \psset{linestyle=solid}
\psset{linestyle=solid} \setlinecaps{0} \setlinejoinmode{0}
\setlinecaps{0} \setlinejoinmode{0} \psset{linestyle=solid}
\newrgbcolor{dialinecolor}{1 1 1}
\psset{linecolor=dialinecolor}
\psellipse*(12.230300,11.173800)(0.687300,0.687300)
\newrgbcolor{dialinecolor}{0 0 0}
\psset{linecolor=dialinecolor}
\psellipse(12.230300,11.173800)(0.687300,0.687300)
\psset{linewidth=0.004000} \setlinecaps{0} \setlinejoinmode{0}
\psset{linestyle=solid}
\newrgbcolor{dialinecolor}{0 0 0}
\psset{linecolor=dialinecolor}
\psellipse(12.230300,11.173800)(0.687300,0.687300)
\psset{linewidth=0.10} \psset{linestyle=solid}
\psset{linestyle=solid} \setlinecaps{0}
\newrgbcolor{dialinecolor}{0 0 0}
\psset{linecolor=dialinecolor}
\psline(12.230300,11.861100)(12.229300,11.637400)
\psset{linewidth=0.040} \psset{linestyle=solid}
\psset{linestyle=solid} \setlinecaps{0}
\newrgbcolor{dialinecolor}{0 0 0}
\psset{linecolor=dialinecolor}
\psline(11.562615,11.261065)(11.576200,11.321500)
\psset{linewidth=0.040} \psset{linestyle=solid}
\setlinejoinmode{0} \setlinecaps{0}
\newrgbcolor{dialinecolor}{0 0 0}
\psset{linecolor=dialinecolor}
\psline(11.764948,11.477233)(11.552808,11.217433)(11.472252,11.543026)
\psset{linewidth=0.040} \psset{linestyle=solid}
\psset{linestyle=solid} \setlinecaps{0}
\newrgbcolor{dialinecolor}{0 0 0}
\psset{linecolor=dialinecolor}
\psline(11.376011,10.807773)(11.390500,10.777400)
\psset{linewidth=0.040} \psset{linestyle=solid}
\setlinejoinmode{0} \setlinecaps{0}
\newrgbcolor{dialinecolor}{0 0 0}
\psset{linecolor=dialinecolor}
\psline(11.350541,10.512784)(11.356755,10.848136)(11.621309,10.641953)
\newrgbcolor{dialinecolor}{1 1 1}
\psset{linecolor=dialinecolor}
\pspolygon*(11.995500,11.938300)(11.995500,12.266400)(12.283230,12.266400)(12.283230,11.938300)
\psset{linewidth=0.10} \psset{linestyle=solid}
\psset{linestyle=solid} \setlinejoinmode{0}
\newrgbcolor{dialinecolor}{1 1 1}
\psset{linecolor=dialinecolor}
\pspolygon(11.995500,11.938300)(11.995500,12.266400)(12.283230,12.266400)(12.283230,11.938300)
\psset{linewidth=0.040} \psset{linestyle=solid}
\psset{linestyle=solid} \setlinecaps{0}
\newrgbcolor{dialinecolor}{0 0 0}
\psset{linecolor=dialinecolor}
\psclip{\pswedge[linestyle=none,fillstyle=none](11.172087,11.511040){1.889292}{27.132940}{132.334144}}
\psellipse(11.172087,11.511040)(1.335931,1.335931)
\endpsclip
\psset{linewidth=0.040} \psset{linestyle=solid}
\psset{linestyle=solid} \setlinecaps{0} \setlinejoinmode{0}
\setlinecaps{0} \setlinejoinmode{0} \psset{linestyle=solid}
\newrgbcolor{dialinecolor}{1 1 1}
\psset{linecolor=dialinecolor}
\psellipse*(10.023750,14.358600)(0.950,0.950)
\newrgbcolor{dialinecolor}{0 0 0}
\psset{linecolor=dialinecolor}
\psellipse(10.023750,14.358600)(0.950,0.950)
\psset{linewidth=0.004000} \setlinecaps{0} \setlinejoinmode{0}
\psset{linestyle=solid}
\newrgbcolor{dialinecolor}{0 0 0}
\psset{linecolor=dialinecolor}
\psellipse(10.023750,14.358600)(0.950,0.950)
\psset{linewidth=0.040} \psset{linestyle=solid}
\psset{linestyle=solid} \setlinecaps{0} \setlinejoinmode{0}
\setlinecaps{0} \setlinejoinmode{0} \psset{linestyle=solid}
\newrgbcolor{dialinecolor}{1 1 1}
\psset{linecolor=dialinecolor}
\psellipse*(10.012850,14.362900)(0.687300,0.687300)
\newrgbcolor{dialinecolor}{0 0 0}
\psset{linecolor=dialinecolor}
\psellipse(10.012850,14.362900)(0.687300,0.687300)
\psset{linewidth=0.004000} \setlinecaps{0} \setlinejoinmode{0}
\psset{linestyle=solid}
\newrgbcolor{dialinecolor}{0 0 0}
\psset{linecolor=dialinecolor}
\psellipse(10.012850,14.362900)(0.687300,0.687300)
\psset{linewidth=0.10} \psset{linestyle=solid}
\psset{linestyle=solid} \setlinecaps{0}
\newrgbcolor{dialinecolor}{0 0 0}
\psset{linecolor=dialinecolor}
\psline(10.012900,15.050200)(10.014600,14.821300)
\psset{linewidth=0.040} \psset{linestyle=solid}
\psset{linestyle=solid} \setlinecaps{0}
\newrgbcolor{dialinecolor}{0 0 0}
\psset{linecolor=dialinecolor}
\psline(9.345165,14.450165)(9.358750,14.510600)
\psset{linewidth=0.040} \psset{linestyle=solid}
\setlinejoinmode{0} \setlinecaps{0}
\newrgbcolor{dialinecolor}{0 0 0}
\psset{linecolor=dialinecolor}
\psline(9.547498,14.666333)(9.335358,14.406533)(9.254802,14.732126)
\psset{linewidth=0.040} \psset{linestyle=solid}
\psset{linestyle=solid} \setlinecaps{0}
\newrgbcolor{dialinecolor}{0 0 0}
\psset{linecolor=dialinecolor}
\psline(9.126701,14.108244)(9.100750,14.289900)
\psset{linewidth=0.040} \psset{linestyle=solid}
\setlinejoinmode{0} \setlinecaps{0}
\newrgbcolor{dialinecolor}{0 0 0}
\psset{linecolor=dialinecolor}
\psline(9.239091,14.382170)(9.133025,14.063972)(8.942107,14.339744)
\psset{linewidth=0.10} \psset{linestyle=solid}
\psset{linestyle=solid} \setlinecaps{0}
\newrgbcolor{dialinecolor}{0 0 0}
\psset{linecolor=dialinecolor}
\psline(10.027200,15.525100)(10.023800,15.308600)
\newrgbcolor{dialinecolor}{1 1 1}
\psset{linecolor=dialinecolor}
\pspolygon*(9.900550,13.305000)(9.900550,13.555900)(10.184460,13.555900)(10.184460,13.305000)
\psset{linewidth=0.10} \psset{linestyle=solid}
\psset{linestyle=solid} \setlinejoinmode{0}
\newrgbcolor{dialinecolor}{1 1 1}
\psset{linecolor=dialinecolor}
\pspolygon(9.900550,13.305000)(9.900550,13.555900)(10.184460,13.555900)(10.184460,13.305000)
\newrgbcolor{dialinecolor}{1 1 1}
\psset{linecolor=dialinecolor}
\pspolygon*(8.954150,11.248200)(8.954150,11.576300)(9.282250,11.576300)(9.282250,11.248200)
\psset{linewidth=0.10} \psset{linestyle=solid}
\psset{linestyle=solid} \setlinejoinmode{0}
\newrgbcolor{dialinecolor}{1 1 1}
\psset{linecolor=dialinecolor}
\pspolygon(8.954150,11.248200)(8.954150,11.576300)(9.282250,11.576300)(9.282250,11.248200)
\psset{linewidth=0.040} \psset{linestyle=solid}
\psset{linestyle=solid} \setlinecaps{0}
\newrgbcolor{dialinecolor}{0 0 0}
\psset{linecolor=dialinecolor}
\psclip{\pswedge[linestyle=none,fillstyle=none](10.122678,12.526388){0.640192}{107.168914}{310.056643}}
\psellipse(10.122678,12.526388)(0.452684,0.452684)
\endpsclip
\psset{linewidth=0.040} \psset{linestyle=solid}
\psset{linestyle=solid} \setlinecaps{0}
\newrgbcolor{dialinecolor}{0 0 0}
\psset{linecolor=dialinecolor}
\psclip{\pswedge[linestyle=none,fillstyle=none](11.145698,11.518738){1.375271}{34.783057}{135.726453}}
\psellipse(11.145698,11.518738)(0.972463,0.972463)
\endpsclip
\psset{linewidth=0.040} \psset{linestyle=solid}
\psset{linestyle=solid} \setlinecaps{0}
\newrgbcolor{dialinecolor}{0 0 0}
\psset{linecolor=dialinecolor}
\psclip{\pswedge[linestyle=none,fillstyle=none](10.106609,12.526493){0.164405}{77.989362}{309.274288}}
\psellipse(10.106609,12.526493)(0.116252,0.116252)
\endpsclip
\psset{linewidth=0.040} \psset{linestyle=solid}
\psset{linestyle=solid} \setlinecaps{0}
\newrgbcolor{dialinecolor}{0 0 0}
\psset{linecolor=dialinecolor}
\psclip{\pswedge[linestyle=none,fillstyle=none](9.625444,11.123227){0.481195}{60.237366}{150.214588}}
\psellipse(9.625444,11.123227)(0.340256,0.340256)
\endpsclip
\psset{linewidth=0.040} \psset{linestyle=solid}
\psset{linestyle=solid} \setlinecaps{0}
\newrgbcolor{dialinecolor}{0 0 0}
\psset{linecolor=dialinecolor}
\psclip{\pswedge[linestyle=none,fillstyle=none](9.540981,11.209854){0.956800}{61.726911}{163.509460}}
\psellipse(9.540981,11.209854)(0.676560,0.676560)
\endpsclip
\psset{linewidth=0.040} \psset{linestyle=solid}
\psset{linestyle=solid} \setlinecaps{0}
\newrgbcolor{dialinecolor}{0 0 0}
\psset{linecolor=dialinecolor}
\psclip{\pswedge[linestyle=none,fillstyle=none](10.162155,12.028939){1.016590}{240.515351}{7.824409}}
\psellipse(10.162155,12.028939)(0.718838,0.718838)
\endpsclip
\newrgbcolor{dialinecolor}{1 1 1}
\psset{linecolor=dialinecolor}
\pspolygon*(11.195200,12.247600)(11.195200,12.976590)(11.523300,12.976590)(11.523300,12.247600)
\psset{linewidth=0.10} \psset{linestyle=solid}
\psset{linestyle=solid} \setlinejoinmode{0}
\newrgbcolor{dialinecolor}{1 1 1}
\psset{linecolor=dialinecolor}
\pspolygon(11.195200,12.247600)(11.195200,12.976590)(11.523300,12.976590)(11.523300,12.247600)
\psset{linewidth=0.010} \psset{linestyle=dashed,dash=1 1}
\psset{linestyle=dashed,dash=0.10 0.10} \setlinecaps{0}
\newrgbcolor{dialinecolor}{0 0 0}
\psset{linecolor=dialinecolor}
\psclip{\pswedge[linestyle=none,fillstyle=none](11.156113,11.555344){1.579472}{15.888450}{127.785980}}
\psellipse(11.156113,11.555344)(1.116855,1.116855)
\endpsclip
\psset{linewidth=0.040} \psset{linestyle=solid}
\psset{linestyle=solid} \setlinecaps{0}
\newrgbcolor{dialinecolor}{0 0 0}
\psset{linecolor=dialinecolor}
\psclip{\pswedge[linestyle=none,fillstyle=none](11.062320,12.313630){0.757637}{6.927913}{82.079977}}
\psellipse(11.062320,12.313630)(0.535730,0.535730)
\endpsclip
\psset{linewidth=0.040} \psset{linestyle=solid}
\psset{linestyle=solid} \setlinecaps{0}
\newrgbcolor{dialinecolor}{0 0 0}
\psset{linecolor=dialinecolor}
\psclip{\pswedge[linestyle=none,fillstyle=none](10.075718,14.501808){3.234114}{297.534195}{311.469849}}
\psellipse(10.075718,14.501808)(2.286864,2.286864)
\endpsclip
\newrgbcolor{dialinecolor}{1 1 1}
\psset{linecolor=dialinecolor}
\pspolygon*(9.495550,11.309300)(9.495550,12.038290)(9.823650,12.038290)(9.823650,11.309300)
\psset{linewidth=0.10} \psset{linestyle=solid}
\psset{linestyle=solid} \setlinejoinmode{0}
\newrgbcolor{dialinecolor}{1 1 1}
\psset{linecolor=dialinecolor}
\pspolygon(9.495550,11.309300)(9.495550,12.038290)(9.823650,12.038290)(9.823650,11.309300)
\psset{linewidth=0.010} \psset{linestyle=dashed,dash=1 1}
\psset{linestyle=dashed,dash=0.10 0.10} \setlinecaps{0}
\newrgbcolor{dialinecolor}{0 0 0}
\psset{linecolor=dialinecolor}
\psclip{\pswedge[linestyle=none,fillstyle=none](9.715897,10.871196){1.175371}{79.465339}{161.653467}}
\psellipse(9.715897,10.871196)(0.831113,0.831113)
\endpsclip
\psset{linewidth=0.040} \psset{linestyle=solid}
\psset{linestyle=solid} \setlinecaps{0}
\newrgbcolor{dialinecolor}{0 0 0}
\psset{linecolor=dialinecolor}
\psclip{\pswedge[linestyle=none,fillstyle=none](7.557838,13.595842){4.086542}{310.490651}{322.886354}}
\psellipse(7.557838,13.595842)(2.889622,2.889622)
\endpsclip
\psset{linewidth=0.040} \psset{linestyle=solid}
\psset{linestyle=solid} \setlinecaps{0}
\newrgbcolor{dialinecolor}{0 0 0}
\psset{linecolor=dialinecolor}
\psclip{\pswedge[linestyle=none,fillstyle=none](9.361094,11.353607){0.753990}{1.143928}{81.038984}}
\psellipse(9.361094,11.353607)(0.533151,0.533151)
\endpsclip
\newrgbcolor{dialinecolor}{1 1 1}
\psset{linecolor=dialinecolor}
\pspolygon*(10.367400,12.054000)(10.367400,12.381310)(10.707470,12.381310)(10.707470,12.054000)
\psset{linewidth=0.10} \psset{linestyle=solid}
\psset{linestyle=solid} \setlinejoinmode{0}
\newrgbcolor{dialinecolor}{1 1 1}
\psset{linecolor=dialinecolor}
\pspolygon(10.367400,12.054000)(10.367400,12.381310)(10.707470,12.381310)(10.707470,12.054000)
\newrgbcolor{dialinecolor}{1 1 1}
\psset{linecolor=dialinecolor}
\pspolygon*(10.205900,12.296300)(10.205900,12.708630)(10.507710,12.708630)(10.507710,12.296300)
\psset{linewidth=0.10} \psset{linestyle=solid}
\psset{linestyle=solid} \setlinejoinmode{0}
\newrgbcolor{dialinecolor}{1 1 1}
\psset{linecolor=dialinecolor}
\pspolygon(10.205900,12.296300)(10.205900,12.708630)(10.507710,12.708630)(10.507710,12.296300)
\psset{linewidth=0.040} \psset{linestyle=solid}
\psset{linestyle=solid} \setlinecaps{0}
\newrgbcolor{dialinecolor}{0 0 0}
\psset{linecolor=dialinecolor}
\psclip{\pswedge[linestyle=none,fillstyle=none](12.635733,13.857167){3.979488}{188.531757}{219.112957}}
\psellipse(12.635733,13.857167)(2.813923,2.813923)
\endpsclip
\psset{linewidth=0.010} \psset{linestyle=dashed,dash=1 1}
\psset{linestyle=dashed,dash=0.10 0.10} \setlinecaps{0}
\newrgbcolor{dialinecolor}{0 0 0}
\psset{linecolor=dialinecolor}
\psclip{\pswedge[linestyle=none,fillstyle=none](10.187345,12.025522){0.672006}{237.244973}{54.215701}}
\psellipse(10.187345,12.025522)(0.475180,0.475180)
\endpsclip
\psset{linewidth=0.040} \psset{linestyle=solid}
\psset{linestyle=solid} \setlinecaps{0}
\newrgbcolor{dialinecolor}{0 0 0}
\psset{linecolor=dialinecolor}
\psclip{\pswedge[linestyle=none,fillstyle=none](10.108205,12.070914){0.485037}{217.948851}{3.022834}}
\psellipse(10.108205,12.070914)(0.342973,0.342973)
\endpsclip
\psset{linewidth=0.040} \psset{linestyle=solid}
\psset{linestyle=solid} \setlinecaps{0}
\newrgbcolor{dialinecolor}{0 0 0}
\psset{linecolor=dialinecolor}
\psclip{\pswedge[linestyle=none,fillstyle=none](13.033461,13.928238){4.004798}{189.934220}{220.418053}}
\psellipse(13.033461,13.928238)(2.831820,2.831820)
\endpsclip
\newrgbcolor{dialinecolor}{1 1 1}
\psset{linecolor=dialinecolor}
\pspolygon*(10.018800,12.700100)(10.018800,12.882890)(10.269600,12.882890)(10.269600,12.700100)
\psset{linewidth=0.10} \psset{linestyle=solid}
\psset{linestyle=solid} \setlinejoinmode{0}
\newrgbcolor{dialinecolor}{1 1 1}
\psset{linecolor=dialinecolor}
\pspolygon(10.018800,12.700100)(10.018800,12.882890)(10.269600,12.882890)(10.269600,12.700100)
\psset{linewidth=0.010} \psset{linestyle=dashed,dash=1 1}
\psset{linestyle=dashed,dash=0.10 0.10} \setlinecaps{0}
\newrgbcolor{dialinecolor}{0 0 0}
\psset{linecolor=dialinecolor}
\psclip{\pswedge[linestyle=none,fillstyle=none](12.252680,13.758126){3.169677}{182.110133}{216.906517}}
\psellipse(12.252680,13.758126)(2.241300,2.241300)
\endpsclip
\psset{linewidth=0.010} \psset{linestyle=dashed,dash=0.10 0.10}
\psset{linestyle=dashed,dash=0.10 0.10} \setlinecaps{0}
\newrgbcolor{dialinecolor}{0 0 0}
\psset{linecolor=dialinecolor}
\psclip{\pswedge[linestyle=none,fillstyle=none](10.152940,12.553151){0.450789}{84.426153}{285.029823}}
\psellipse(10.152940,12.553151)(0.318756,0.318756)
\endpsclip
}\endpspicture \]

In particular, a chord with two adjacent twists is identified with
a chord without a twist. Also, note that the relations among
diagrams, such as sliding along a chord, now have to respect the
twists of that chord (see figure below).

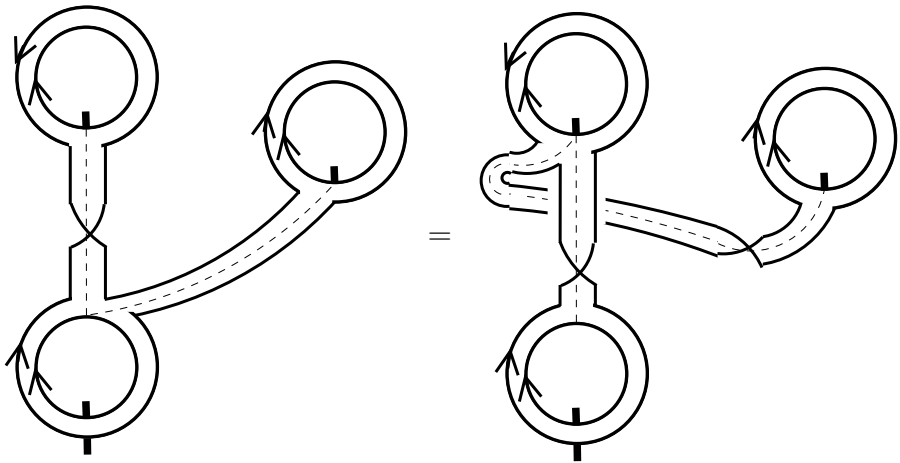
\[ \pspicture(1.602736,-15.578677)(13.628700,-9.449470)
\scalebox{1 -1}{
\newrgbcolor{dialinecolor}{0 0 0}
\psset{linecolor=dialinecolor}
\newrgbcolor{diafillcolor}{1 1 1}
\psset{fillcolor=diafillcolor} \psset{linewidth=0.01}
\psset{linestyle=dashed,dash=1 1} \psset{linestyle=dashed,dash=0.1
0.1} \setlinecaps{0}
\newrgbcolor{dialinecolor}{0 0 0}
\psset{linecolor=dialinecolor}
\psclip{\pswedge[linestyle=none,fillstyle=none](6.749514,25.744536){19.629216}{277.872548}{289.295737}}
\psellipse(6.749514,25.744536)(13.879952,13.879952)
\endpsclip
\psset{linewidth=0.04} \psset{linestyle=solid}
\psset{linestyle=solid} \setlinecaps{0} \setlinejoinmode{0}
\setlinecaps{0} \setlinejoinmode{0} \psset{linestyle=solid}
\newrgbcolor{dialinecolor}{1 1 1}
\psset{linecolor=dialinecolor}
\psellipse*(2.842220,10.419470)(0.950000,0.950000)
\newrgbcolor{dialinecolor}{0 0 0}
\psset{linecolor=dialinecolor}
\psellipse(2.842220,10.419470)(0.950000,0.950000)
\psset{linewidth=0.004000} \setlinecaps{0} \setlinejoinmode{0}
\psset{linestyle=solid}
\newrgbcolor{dialinecolor}{0 0 0}
\psset{linecolor=dialinecolor}
\psellipse(2.842220,10.419470)(0.950000,0.950000)
\psset{linewidth=0.04} \psset{linestyle=solid}
\psset{linestyle=solid} \setlinecaps{0} \setlinejoinmode{0}
\setlinecaps{0} \setlinejoinmode{0} \psset{linestyle=solid}
\newrgbcolor{dialinecolor}{1 1 1}
\psset{linecolor=dialinecolor}
\psellipse*(2.831320,10.423770)(0.687300,0.687300)
\newrgbcolor{dialinecolor}{0 0 0}
\psset{linecolor=dialinecolor}
\psellipse(2.831320,10.423770)(0.687300,0.687300)
\psset{linewidth=0.004000} \setlinecaps{0} \setlinejoinmode{0}
\psset{linestyle=solid}
\newrgbcolor{dialinecolor}{0 0 0}
\psset{linecolor=dialinecolor}
\psellipse(2.831320,10.423770)(0.687300,0.687300)
\psset{linewidth=0.1} \psset{linestyle=solid}
\psset{linestyle=solid} \setlinecaps{0}
\newrgbcolor{dialinecolor}{0 0 0}
\psset{linecolor=dialinecolor}
\psline(2.831320,11.111100)(2.823220,10.875500)
\psset{linewidth=0.04} \psset{linestyle=solid}
\psset{linestyle=solid} \setlinecaps{0}
\newrgbcolor{dialinecolor}{0 0 0}
\psset{linecolor=dialinecolor}
\psline(2.163635,10.511065)(2.177220,10.571500)
\psset{linewidth=0.04} \psset{linestyle=solid} \setlinejoinmode{0}
\setlinecaps{0}
\newrgbcolor{dialinecolor}{0 0 0}
\psset{linecolor=dialinecolor}
\psline(2.365968,10.727233)(2.153828,10.467433)(2.073272,10.793026)
\psset{linewidth=0.04} \psset{linestyle=solid}
\psset{linestyle=solid} \setlinecaps{0} \setlinejoinmode{0}
\setlinecaps{0} \setlinejoinmode{0} \psset{linestyle=solid}
\newrgbcolor{dialinecolor}{1 1 1}
\psset{linecolor=dialinecolor}
\psellipse*(2.845150,14.274400)(0.950000,0.950000)
\newrgbcolor{dialinecolor}{0 0 0}
\psset{linecolor=dialinecolor}
\psellipse(2.845150,14.274400)(0.950000,0.950000)
\psset{linewidth=0.004000} \setlinecaps{0} \setlinejoinmode{0}
\psset{linestyle=solid}
\newrgbcolor{dialinecolor}{0 0 0}
\psset{linecolor=dialinecolor}
\psellipse(2.845150,14.274400)(0.950000,0.950000)
\psset{linewidth=0.04} \psset{linestyle=solid}
\psset{linestyle=solid} \setlinecaps{0} \setlinejoinmode{0}
\setlinecaps{0} \setlinejoinmode{0} \psset{linestyle=solid}
\newrgbcolor{dialinecolor}{1 1 1}
\psset{linecolor=dialinecolor}
\psellipse*(2.834250,14.278700)(0.687300,0.687300)
\newrgbcolor{dialinecolor}{0 0 0}
\psset{linecolor=dialinecolor}
\psellipse(2.834250,14.278700)(0.687300,0.687300)
\psset{linewidth=0.004000} \setlinecaps{0} \setlinejoinmode{0}
\psset{linestyle=solid}
\newrgbcolor{dialinecolor}{0 0 0}
\psset{linecolor=dialinecolor}
\psellipse(2.834250,14.278700)(0.687300,0.687300)
\psset{linewidth=0.1} \psset{linestyle=solid}
\psset{linestyle=solid} \setlinecaps{0}
\newrgbcolor{dialinecolor}{0 0 0}
\psset{linecolor=dialinecolor}
\psline(2.834250,14.966000)(2.826150,14.730400)
\psset{linewidth=0.04} \psset{linestyle=solid}
\psset{linestyle=solid} \setlinecaps{0}
\newrgbcolor{dialinecolor}{0 0 0}
\psset{linecolor=dialinecolor}
\psline(2.166565,14.365965)(2.180150,14.426400)
\psset{linewidth=0.04} \psset{linestyle=solid} \setlinejoinmode{0}
\setlinecaps{0}
\newrgbcolor{dialinecolor}{0 0 0}
\psset{linecolor=dialinecolor}
\psline(2.368898,14.582133)(2.156758,14.322333)(2.076202,14.647926)
\psset{linewidth=0.04} \psset{linestyle=solid}
\psset{linestyle=solid} \setlinecaps{0}
\newrgbcolor{dialinecolor}{0 0 0}
\psset{linecolor=dialinecolor}
\psline(1.948061,14.024044)(1.922110,14.205700)
\psset{linewidth=0.04} \psset{linestyle=solid} \setlinejoinmode{0}
\setlinecaps{0}
\newrgbcolor{dialinecolor}{0 0 0}
\psset{linecolor=dialinecolor}
\psline(2.060451,14.297970)(1.954385,13.979772)(1.763467,14.255544)
\newrgbcolor{dialinecolor}{1 1 1}
\psset{linecolor=dialinecolor}
\pspolygon*(2.692000,13.317900)(2.692000,13.491600)(4.006510,13.491600)(4.006510,13.317900)
\psset{linewidth=0.1} \psset{linestyle=solid}
\psset{linestyle=solid} \setlinejoinmode{0}
\newrgbcolor{dialinecolor}{1 1 1}
\psset{linecolor=dialinecolor}
\pspolygon(2.692000,13.317900)(2.692000,13.491600)(4.006510,13.491600)(4.006510,13.317900)
\newrgbcolor{dialinecolor}{1 1 1}
\psset{linecolor=dialinecolor}
\pspolygon*(2.653400,11.260300)(2.653400,11.588400)(2.981500,11.588400)(2.981500,11.260300)
\psset{linewidth=0.1} \psset{linestyle=solid}
\psset{linestyle=solid} \setlinejoinmode{0}
\newrgbcolor{dialinecolor}{1 1 1}
\psset{linecolor=dialinecolor}
\pspolygon(2.653400,11.260300)(2.653400,11.588400)(2.981500,11.588400)(2.981500,11.260300)
\psset{linewidth=0.1} \psset{linestyle=solid}
\psset{linestyle=solid} \setlinecaps{0}
\newrgbcolor{dialinecolor}{0 0 0}
\psset{linecolor=dialinecolor}
\psline(2.848510,15.440900)(2.845150,15.224400)
\psset{linewidth=0.04} \psset{linestyle=solid}
\psset{linestyle=solid} \setlinecaps{0}
\newrgbcolor{dialinecolor}{0 0 0}
\psset{linecolor=dialinecolor}
\psline(2.618930,11.336000)(2.618930,13.381800)
\psset{linewidth=0.04} \psset{linestyle=solid}
\psset{linestyle=solid} \setlinecaps{0}
\newrgbcolor{dialinecolor}{0 0 0}
\psset{linecolor=dialinecolor}
\psline(3.086600,11.326800)(3.086600,13.372600)
\newrgbcolor{dialinecolor}{1 1 1}
\psset{linecolor=dialinecolor}
\pspolygon*(2.479700,12.159900)(2.479700,12.623100)(3.193800,12.623100)(3.193800,12.159900)
\psset{linewidth=0.1} \psset{linestyle=solid}
\psset{linestyle=solid} \setlinejoinmode{0}
\newrgbcolor{dialinecolor}{1 1 1}
\psset{linecolor=dialinecolor}
\pspolygon(2.479700,12.159900)(2.479700,12.623100)(3.193800,12.623100)(3.193800,12.159900)
\psset{linewidth=0.01} \psset{linestyle=dashed,dash=1 1}
\psset{linestyle=dashed,dash=0.1 0.1} \setlinecaps{0}
\newrgbcolor{dialinecolor}{0 0 0}
\psset{linecolor=dialinecolor}
\psline(2.831320,11.111100)(2.834250,13.591400)
\psset{linewidth=0.04} \psset{linestyle=solid}
\psset{linestyle=solid} \setlinecaps{0}
\newrgbcolor{dialinecolor}{0 0 0}
\psset{linecolor=dialinecolor}
\psclip{\pswedge[linestyle=none,fillstyle=none](2.377058,12.031056){1.014936}{6.459998}{68.832104}}
\psellipse(2.377058,12.031056)(0.717668,0.717668)
\endpsclip
\psset{linewidth=0.04} \psset{linestyle=solid}
\psset{linestyle=solid} \setlinecaps{0}
\newrgbcolor{dialinecolor}{0 0 0}
\psset{linecolor=dialinecolor}
\psclip{\pswedge[linestyle=none,fillstyle=none](3.663581,11.744830){1.573099}{120.909821}{161.064040}}
\psellipse(3.663581,11.744830)(1.112349,1.112349)
\endpsclip
\psset{linewidth=0.04} \psset{linestyle=solid}
\psset{linestyle=solid} \setlinecaps{0}
\newrgbcolor{dialinecolor}{0 0 0}
\psset{linecolor=dialinecolor}
\psline(1.920988,10.167988)(1.999310,9.998270)
\psset{linewidth=0.04} \psset{linestyle=solid} \setlinejoinmode{0}
\setlinecaps{0}
\newrgbcolor{dialinecolor}{0 0 0}
\psset{linecolor=dialinecolor}
\psline(1.891757,9.873348)(1.902249,10.208594)(2.164151,9.999053)
\psset{linewidth=0.04} \psset{linestyle=solid}
\psset{linestyle=solid} \setlinecaps{0} \setlinejoinmode{0}
\setlinecaps{0} \setlinejoinmode{0} \psset{linestyle=solid}
\newrgbcolor{dialinecolor}{1 1 1}
\psset{linecolor=dialinecolor}
\psellipse*(6.145450,11.147800)(0.950000,0.950000)
\newrgbcolor{dialinecolor}{0 0 0}
\psset{linecolor=dialinecolor}
\psellipse(6.145450,11.147800)(0.950000,0.950000)
\psset{linewidth=0.004000} \setlinecaps{0} \setlinejoinmode{0}
\psset{linestyle=solid}
\newrgbcolor{dialinecolor}{0 0 0}
\psset{linecolor=dialinecolor}
\psellipse(6.145450,11.147800)(0.950000,0.950000)
\psset{linewidth=0.04} \psset{linestyle=solid}
\psset{linestyle=solid} \setlinecaps{0} \setlinejoinmode{0}
\setlinecaps{0} \setlinejoinmode{0} \psset{linestyle=solid}
\newrgbcolor{dialinecolor}{1 1 1}
\psset{linecolor=dialinecolor}
\psellipse*(6.134550,11.152100)(0.687300,0.687300)
\newrgbcolor{dialinecolor}{0 0 0}
\psset{linecolor=dialinecolor}
\psellipse(6.134550,11.152100)(0.687300,0.687300)
\psset{linewidth=0.004000} \setlinecaps{0} \setlinejoinmode{0}
\psset{linestyle=solid}
\newrgbcolor{dialinecolor}{0 0 0}
\psset{linecolor=dialinecolor}
\psellipse(6.134550,11.152100)(0.687300,0.687300)
\psset{linewidth=0.1} \psset{linestyle=solid}
\psset{linestyle=solid} \setlinecaps{0}
\newrgbcolor{dialinecolor}{0 0 0}
\psset{linecolor=dialinecolor}
\psline(6.134550,11.839400)(6.126450,11.603800)
\psset{linewidth=0.04} \psset{linestyle=solid}
\psset{linestyle=solid} \setlinecaps{0}
\newrgbcolor{dialinecolor}{0 0 0}
\psset{linecolor=dialinecolor}
\psline(5.466865,11.239365)(5.480450,11.299800)
\psset{linewidth=0.04} \psset{linestyle=solid} \setlinejoinmode{0}
\setlinecaps{0}
\newrgbcolor{dialinecolor}{0 0 0}
\psset{linecolor=dialinecolor}
\psline(5.669198,11.455533)(5.457058,11.195733)(5.376502,11.521326)
\psset{linewidth=0.04} \psset{linestyle=solid}
\psset{linestyle=solid} \setlinecaps{0}
\newrgbcolor{dialinecolor}{0 0 0}
\psset{linecolor=dialinecolor}
\psline(5.227181,10.955055)(5.195450,11.147800)
\psset{linewidth=0.04} \psset{linestyle=solid} \setlinejoinmode{0}
\setlinecaps{0}
\newrgbcolor{dialinecolor}{0 0 0}
\psset{linecolor=dialinecolor}
\psline(5.333721,11.231309)(5.234445,10.910927)(5.037706,11.182577)
\newrgbcolor{dialinecolor}{1 1 1}
\psset{linecolor=dialinecolor}
\pspolygon*(5.724210,11.897200)(5.724210,12.225300)(6.052310,12.225300)(6.052310,11.897200)
\psset{linewidth=0.1} \psset{linestyle=solid}
\psset{linestyle=solid} \setlinejoinmode{0}
\newrgbcolor{dialinecolor}{1 1 1}
\psset{linecolor=dialinecolor}
\pspolygon(5.724210,11.897200)(5.724210,12.225300)(6.052310,12.225300)(6.052310,11.897200)
\psset{linewidth=0.01} \psset{linestyle=dashed,dash=1 1}
\psset{linestyle=dashed,dash=0.1 0.1} \setlinecaps{0}
\newrgbcolor{dialinecolor}{0 0 0}
\psset{linecolor=dialinecolor}
\psclip{\pswedge[linestyle=none,fillstyle=none](2.042204,8.114956){7.825442}{42.305331}{81.770506}}
\psellipse(2.042204,8.114956)(5.533423,5.533423)
\endpsclip
\psset{linewidth=0.04} \psset{linestyle=solid}
\psset{linestyle=solid} \setlinecaps{0}
\newrgbcolor{dialinecolor}{0 0 0}
\psset{linecolor=dialinecolor}
\psclip{\pswedge[linestyle=none,fillstyle=none](2.266871,8.866818){6.506411}{42.362913}{79.818719}}
\psellipse(2.266871,8.866818)(4.600728,4.600728)
\endpsclip
\psset{linewidth=0.04} \psset{linestyle=solid}
\psset{linestyle=solid} \setlinecaps{0}
\newrgbcolor{dialinecolor}{0 0 0}
\psset{linecolor=dialinecolor}
\psclip{\pswedge[linestyle=none,fillstyle=none](2.440798,8.668302){7.139454}{42.791276}{78.058014}}
\psellipse(2.440798,8.668302)(5.048356,5.048356)
\endpsclip
\setfont{Helvetica}{0.400000}
\newrgbcolor{dialinecolor}{0 0 0}
\psset{linecolor=dialinecolor}
\rput[l](7.374800,12.562900){\scalebox{1 -1}{=}}
\psset{linewidth=0.04} \psset{linestyle=solid}
\psset{linestyle=solid} \setlinecaps{0} \setlinejoinmode{0}
\setlinecaps{0} \setlinejoinmode{0} \psset{linestyle=solid}
\newrgbcolor{dialinecolor}{1 1 1}
\psset{linecolor=dialinecolor}
\psellipse*(9.355520,10.506560)(0.950000,0.950000)
\newrgbcolor{dialinecolor}{0 0 0}
\psset{linecolor=dialinecolor}
\psellipse(9.355520,10.506560)(0.950000,0.950000)
\psset{linewidth=0.004000} \setlinecaps{0} \setlinejoinmode{0}
\psset{linestyle=solid}
\newrgbcolor{dialinecolor}{0 0 0}
\psset{linecolor=dialinecolor}
\psellipse(9.355520,10.506560)(0.950000,0.950000)
\psset{linewidth=0.04} \psset{linestyle=solid}
\psset{linestyle=solid} \setlinecaps{0} \setlinejoinmode{0}
\setlinecaps{0} \setlinejoinmode{0} \psset{linestyle=solid}
\newrgbcolor{dialinecolor}{1 1 1}
\psset{linecolor=dialinecolor}
\psellipse*(9.344620,10.510860)(0.687300,0.687300)
\newrgbcolor{dialinecolor}{0 0 0}
\psset{linecolor=dialinecolor}
\psellipse(9.344620,10.510860)(0.687300,0.687300)
\psset{linewidth=0.004000} \setlinecaps{0} \setlinejoinmode{0}
\psset{linestyle=solid}
\newrgbcolor{dialinecolor}{0 0 0}
\psset{linecolor=dialinecolor}
\psellipse(9.344620,10.510860)(0.687300,0.687300)
\psset{linewidth=0.1} \psset{linestyle=solid}
\psset{linestyle=solid} \setlinecaps{0}
\newrgbcolor{dialinecolor}{0 0 0}
\psset{linecolor=dialinecolor}
\psline(9.344620,11.198200)(9.336540,10.962600)
\psset{linewidth=0.04} \psset{linestyle=solid}
\psset{linestyle=solid} \setlinecaps{0}
\newrgbcolor{dialinecolor}{0 0 0}
\psset{linecolor=dialinecolor}
\psline(8.676948,10.598162)(8.690520,10.658500)
\psset{linewidth=0.04} \psset{linestyle=solid} \setlinejoinmode{0}
\setlinecaps{0}
\newrgbcolor{dialinecolor}{0 0 0}
\psset{linecolor=dialinecolor}
\psline(8.879312,10.814301)(8.667134,10.554531)(8.586625,10.880136)
\psset{linewidth=0.04} \psset{linestyle=solid}
\psset{linestyle=solid} \setlinecaps{0} \setlinejoinmode{0}
\setlinecaps{0} \setlinejoinmode{0} \psset{linestyle=solid}
\newrgbcolor{dialinecolor}{1 1 1}
\psset{linecolor=dialinecolor}
\psellipse*(9.358450,14.361400)(0.950000,0.950000)
\newrgbcolor{dialinecolor}{0 0 0}
\psset{linecolor=dialinecolor}
\psellipse(9.358450,14.361400)(0.950000,0.950000)
\psset{linewidth=0.004000} \setlinecaps{0} \setlinejoinmode{0}
\psset{linestyle=solid}
\newrgbcolor{dialinecolor}{0 0 0}
\psset{linecolor=dialinecolor}
\psellipse(9.358450,14.361400)(0.950000,0.950000)
\psset{linewidth=0.04} \psset{linestyle=solid}
\psset{linestyle=solid} \setlinecaps{0} \setlinejoinmode{0}
\setlinecaps{0} \setlinejoinmode{0} \psset{linestyle=solid}
\newrgbcolor{dialinecolor}{1 1 1}
\psset{linecolor=dialinecolor}
\psellipse*(9.347550,14.365700)(0.687300,0.687300)
\newrgbcolor{dialinecolor}{0 0 0}
\psset{linecolor=dialinecolor}
\psellipse(9.347550,14.365700)(0.687300,0.687300)
\psset{linewidth=0.004000} \setlinecaps{0} \setlinejoinmode{0}
\psset{linestyle=solid}
\newrgbcolor{dialinecolor}{0 0 0}
\psset{linecolor=dialinecolor}
\psellipse(9.347550,14.365700)(0.687300,0.687300)
\psset{linewidth=0.1} \psset{linestyle=solid}
\psset{linestyle=solid} \setlinecaps{0}
\newrgbcolor{dialinecolor}{0 0 0}
\psset{linecolor=dialinecolor}
\psline(9.347550,15.053000)(9.339440,14.817500)
\psset{linewidth=0.04} \psset{linestyle=solid}
\psset{linestyle=solid} \setlinecaps{0}
\newrgbcolor{dialinecolor}{0 0 0}
\psset{linecolor=dialinecolor}
\psline(8.679865,14.452965)(8.693450,14.513400)
\psset{linewidth=0.04} \psset{linestyle=solid} \setlinejoinmode{0}
\setlinecaps{0}
\newrgbcolor{dialinecolor}{0 0 0}
\psset{linecolor=dialinecolor}
\psline(8.882198,14.669133)(8.670058,14.409333)(8.589502,14.734926)
\psset{linewidth=0.04} \psset{linestyle=solid}
\psset{linestyle=solid} \setlinecaps{0}
\newrgbcolor{dialinecolor}{0 0 0}
\psset{linecolor=dialinecolor}
\psline(8.461361,14.111044)(8.435410,14.292700)
\psset{linewidth=0.04} \psset{linestyle=solid} \setlinejoinmode{0}
\setlinecaps{0}
\newrgbcolor{dialinecolor}{0 0 0}
\psset{linecolor=dialinecolor}
\psline(8.573751,14.384970)(8.467685,14.066772)(8.276767,14.342544)
\newrgbcolor{dialinecolor}{1 1 1}
\psset{linecolor=dialinecolor}
\pspolygon*(8.862950,11.347400)(8.862950,11.675500)(9.494800,11.675500)(9.494800,11.347400)
\psset{linewidth=0.1} \psset{linestyle=solid}
\psset{linestyle=solid} \setlinejoinmode{0}
\newrgbcolor{dialinecolor}{1 1 1}
\psset{linecolor=dialinecolor}
\pspolygon(8.862950,11.347400)(8.862950,11.675500)(9.494800,11.675500)(9.494800,11.347400)
\psset{linewidth=0.1} \psset{linestyle=solid}
\psset{linestyle=solid} \setlinecaps{0}
\newrgbcolor{dialinecolor}{0 0 0}
\psset{linecolor=dialinecolor}
\psline(9.361840,15.527900)(9.358450,15.311400)
\psset{linewidth=0.04} \psset{linestyle=solid}
\psset{linestyle=solid} \setlinecaps{0}
\newrgbcolor{dialinecolor}{0 0 0}
\psset{linecolor=dialinecolor}
\psline(8.434304,10.254995)(8.512610,10.085400)
\psset{linewidth=0.04} \psset{linestyle=solid} \setlinejoinmode{0}
\setlinecaps{0}
\newrgbcolor{dialinecolor}{0 0 0}
\psset{linecolor=dialinecolor}
\psline(8.405131,9.960350)(8.415557,10.295598)(8.677500,10.086108)
\psset{linewidth=0.04} \psset{linestyle=solid}
\psset{linestyle=solid} \setlinecaps{0} \setlinejoinmode{0}
\setlinecaps{0} \setlinejoinmode{0} \psset{linestyle=solid}
\newrgbcolor{dialinecolor}{1 1 1}
\psset{linecolor=dialinecolor}
\psellipse*(12.658700,11.234800)(0.950000,0.950000)
\newrgbcolor{dialinecolor}{0 0 0}
\psset{linecolor=dialinecolor}
\psellipse(12.658700,11.234800)(0.950000,0.950000)
\psset{linewidth=0.004000} \setlinecaps{0} \setlinejoinmode{0}
\psset{linestyle=solid}
\newrgbcolor{dialinecolor}{0 0 0}
\psset{linecolor=dialinecolor}
\psellipse(12.658700,11.234800)(0.950000,0.950000)
\psset{linewidth=0.04} \psset{linestyle=solid}
\psset{linestyle=solid} \setlinecaps{0} \setlinejoinmode{0}
\setlinecaps{0} \setlinejoinmode{0} \psset{linestyle=solid}
\newrgbcolor{dialinecolor}{1 1 1}
\psset{linecolor=dialinecolor}
\psellipse*(12.647800,11.239100)(0.687300,0.687300)
\newrgbcolor{dialinecolor}{0 0 0}
\psset{linecolor=dialinecolor}
\psellipse(12.647800,11.239100)(0.687300,0.687300)
\psset{linewidth=0.004000} \setlinecaps{0} \setlinejoinmode{0}
\psset{linestyle=solid}
\newrgbcolor{dialinecolor}{0 0 0}
\psset{linecolor=dialinecolor}
\psellipse(12.647800,11.239100)(0.687300,0.687300)
\psset{linewidth=0.1} \psset{linestyle=solid}
\psset{linestyle=solid} \setlinecaps{0}
\newrgbcolor{dialinecolor}{0 0 0}
\psset{linecolor=dialinecolor}
\psline(12.647800,11.926400)(12.639700,11.690900)
\psset{linewidth=0.04} \psset{linestyle=solid}
\psset{linestyle=solid} \setlinecaps{0}
\newrgbcolor{dialinecolor}{0 0 0}
\psset{linecolor=dialinecolor}
\psline(11.980115,11.326365)(11.993700,11.386800)
\psset{linewidth=0.04} \psset{linestyle=solid} \setlinejoinmode{0}
\setlinecaps{0}
\newrgbcolor{dialinecolor}{0 0 0}
\psset{linecolor=dialinecolor}
\psline(12.182448,11.542533)(11.970308,11.282733)(11.889752,11.608326)
\psset{linewidth=0.04} \psset{linestyle=solid}
\psset{linestyle=solid} \setlinecaps{0}
\newrgbcolor{dialinecolor}{0 0 0}
\psset{linecolor=dialinecolor}
\psline(11.740459,11.042053)(11.708700,11.234800)
\psset{linewidth=0.04} \psset{linestyle=solid} \setlinejoinmode{0}
\setlinecaps{0}
\newrgbcolor{dialinecolor}{0 0 0}
\psset{linecolor=dialinecolor}
\psline(11.846961,11.318322)(11.747729,10.997926)(11.550952,11.269549)
\newrgbcolor{dialinecolor}{1 1 1}
\psset{linecolor=dialinecolor}
\pspolygon*(12.372600,12.003600)(12.372600,12.331700)(12.700700,12.331700)(12.700700,12.003600)
\psset{linewidth=0.1} \psset{linestyle=solid}
\psset{linestyle=solid} \setlinejoinmode{0}
\newrgbcolor{dialinecolor}{1 1 1}
\psset{linecolor=dialinecolor}
\pspolygon(12.372600,12.003600)(12.372600,12.331700)(12.700700,12.331700)(12.700700,12.003600)
\psset{linewidth=0.01} \psset{linestyle=dashed,dash=1 1}
\psset{linestyle=dashed,dash=0.1 0.1} \setlinecaps{0}
\newrgbcolor{dialinecolor}{0 0 0}
\psset{linecolor=dialinecolor}
\psclip{\pswedge[linestyle=none,fillstyle=none](8.401944,11.776562){0.301467}{69.398740}{285.148638}}
\psellipse(8.401944,11.776562)(0.213170,0.213170)
\endpsclip
\psset{linewidth=0.01} \psset{linestyle=dashed,dash=0.1 0.1}
\psset{linestyle=dashed,dash=0.1 0.1} \setlinecaps{0}
\newrgbcolor{dialinecolor}{0 0 0}
\psset{linecolor=dialinecolor}
\psclip{\pswedge[linestyle=none,fillstyle=none](8.629936,10.759181){1.186180}{31.561762}{97.838192}}
\psellipse(8.629936,10.759181)(0.838756,0.838756)
\endpsclip
\newrgbcolor{dialinecolor}{1 1 1}
\psset{linecolor=dialinecolor}
\pspolygon*(9.188940,13.307800)(9.188940,13.558700)(9.519150,13.558700)(9.519150,13.307800)
\psset{linewidth=0.1} \psset{linestyle=solid}
\psset{linestyle=solid} \setlinejoinmode{0}
\newrgbcolor{dialinecolor}{1 1 1}
\psset{linecolor=dialinecolor}
\pspolygon(9.188940,13.307800)(9.188940,13.558700)(9.519150,13.558700)(9.519150,13.307800)
\psset{linewidth=0.04} \psset{linestyle=solid}
\psset{linestyle=solid} \setlinecaps{0}
\newrgbcolor{dialinecolor}{0 0 0}
\psset{linecolor=dialinecolor}
\psclip{\pswedge[linestyle=none,fillstyle=none](11.624129,11.663429){1.810309}{23.267849}{93.517146}}
\psellipse(11.624129,11.663429)(1.280082,1.280082)
\endpsclip
\psset{linewidth=0.04} \psset{linestyle=solid}
\psset{linestyle=solid} \setlinecaps{0}
\newrgbcolor{dialinecolor}{0 0 0}
\psset{linecolor=dialinecolor}
\psclip{\pswedge[linestyle=none,fillstyle=none](11.738022,11.928847){0.904135}{14.776071}{120.816972}}
\psellipse(11.738022,11.928847)(0.639320,0.639320)
\endpsclip
\psset{linewidth=0.04} \psset{linestyle=solid}
\psset{linestyle=solid} \setlinecaps{0}
\newrgbcolor{dialinecolor}{0 0 0}
\psset{linecolor=dialinecolor}
\psclip{\pswedge[linestyle=none,fillstyle=none](6.841741,27.099960){21.704186}{276.012247}{288.151276}}
\psellipse(6.841741,27.099960)(15.347177,15.347177)
\endpsclip
\psset{linewidth=0.04} \psset{linestyle=solid}
\psset{linestyle=solid} \setlinecaps{0}
\newrgbcolor{dialinecolor}{0 0 0}
\psset{linecolor=dialinecolor}
\psclip{\pswedge[linestyle=none,fillstyle=none](7.582807,22.374965){14.513397}{274.890191}{293.182722}}
\psellipse(7.582807,22.374965)(10.262521,10.262521)
\endpsclip
\newrgbcolor{dialinecolor}{1 1 1}
\psset{linecolor=dialinecolor}
\pspolygon*(9.015240,11.864000)(9.015240,12.327200)(9.685690,12.327200)(9.685690,11.864000)
\psset{linewidth=0.1} \psset{linestyle=solid}
\psset{linestyle=solid} \setlinejoinmode{0}
\newrgbcolor{dialinecolor}{1 1 1}
\psset{linecolor=dialinecolor}
\pspolygon(9.015240,11.864000)(9.015240,12.327200)(9.685690,12.327200)(9.685690,11.864000)
\psset{linewidth=0.04} \psset{linestyle=solid}
\psset{linestyle=solid} \setlinecaps{0}
\newrgbcolor{dialinecolor}{0 0 0}
\psset{linecolor=dialinecolor}
\psline(9.599940,11.413900)(9.599940,13.459700)
\psset{linewidth=0.04} \psset{linestyle=solid}
\psset{linestyle=solid} \setlinecaps{0}
\newrgbcolor{dialinecolor}{0 0 0}
\psset{linecolor=dialinecolor}
\psline(9.133140,11.570800)(9.130240,13.462800)
\newrgbcolor{dialinecolor}{1 1 1}
\psset{linecolor=dialinecolor}
\pspolygon*(9.036640,12.671600)(9.036640,13.134120)(9.707090,13.134120)(9.707090,12.671600)
\psset{linewidth=0.1} \psset{linestyle=solid}
\psset{linestyle=solid} \setlinejoinmode{0}
\newrgbcolor{dialinecolor}{1 1 1}
\psset{linecolor=dialinecolor}
\pspolygon(9.036640,12.671600)(9.036640,13.134120)(9.707090,13.134120)(9.707090,12.671600)
\psset{linewidth=0.01} \psset{linestyle=dashed,dash=1 1}
\psset{linestyle=dashed,dash=0.1 0.1} \setlinecaps{0}
\newrgbcolor{dialinecolor}{0 0 0}
\psset{linecolor=dialinecolor}
\psline(9.344620,11.198200)(9.347550,13.678400)
\psset{linewidth=0.04} \psset{linestyle=solid}
\psset{linestyle=solid} \setlinecaps{0}
\newrgbcolor{dialinecolor}{0 0 0}
\psset{linecolor=dialinecolor}
\psclip{\pswedge[linestyle=none,fillstyle=none](10.181615,12.221776){1.612420}{120.238103}{159.894678}}
\psellipse(10.181615,12.221776)(1.140153,1.140153)
\endpsclip
\psset{linewidth=0.04} \psset{linestyle=solid}
\psset{linestyle=solid} \setlinecaps{0}
\newrgbcolor{dialinecolor}{0 0 0}
\psset{linecolor=dialinecolor}
\psclip{\pswedge[linestyle=none,fillstyle=none](8.903498,12.563606){0.980635}{5.735253}{69.217106}}
\psellipse(8.903498,12.563606)(0.693413,0.693413)
\endpsclip
\psset{linewidth=0.04} \psset{linestyle=solid}
\psset{linestyle=solid} \setlinecaps{0}
\newrgbcolor{dialinecolor}{0 0 0}
\psset{linecolor=dialinecolor}
\psclip{\pswedge[linestyle=none,fillstyle=none](8.424913,11.795777){0.513778}{84.337380}{274.293223}}
\psellipse(8.424913,11.795777)(0.363296,0.363296)
\endpsclip
\psset{linewidth=0.04} \psset{linestyle=solid}
\psset{linestyle=solid} \setlinecaps{0}
\newrgbcolor{dialinecolor}{0 0 0}
\psset{linecolor=dialinecolor}
\psclip{\pswedge[linestyle=none,fillstyle=none](8.554168,11.108965){0.485126}{34.002392}{107.177227}}
\psellipse(8.554168,11.108965)(0.343036,0.343036)
\endpsclip
\psset{linewidth=0.04} \psset{linestyle=solid}
\psset{linestyle=solid} \setlinecaps{0}
\newrgbcolor{dialinecolor}{0 0 0}
\psset{linecolor=dialinecolor}
\psclip{\pswedge[linestyle=none,fillstyle=none](8.689255,10.900224){1.160141}{57.241698}{104.086120}}
\psellipse(8.689255,10.900224)(0.820343,0.820343)
\endpsclip
\psset{linewidth=0.04} \psset{linestyle=solid}
\psset{linestyle=solid} \setlinecaps{0}
\newrgbcolor{dialinecolor}{0 0 0}
\psset{linecolor=dialinecolor}
\psclip{\pswedge[linestyle=none,fillstyle=none](8.431705,11.759303){0.131999}{68.092737}{309.806791}}
\psellipse(8.431705,11.759303)(0.093337,0.093337)
\endpsclip
\newrgbcolor{dialinecolor}{1 1 1}
\psset{linecolor=dialinecolor}
\pspolygon*(11.274900,12.151100)(11.274900,13.126100)(11.774900,13.126100)(11.774900,12.151100)
\psset{linewidth=0.1} \psset{linestyle=solid}
\psset{linestyle=solid} \setlinejoinmode{0}
\newrgbcolor{dialinecolor}{1 1 1}
\psset{linecolor=dialinecolor}
\pspolygon(11.274900,12.151100)(11.274900,13.126100)(11.774900,13.126100)(11.774900,12.151100)
\psset{linewidth=0.01} \psset{linestyle=dashed,dash=1 1}
\psset{linestyle=dashed,dash=0.1 0.1} \setlinecaps{0}
\newrgbcolor{dialinecolor}{0 0 0}
\psset{linecolor=dialinecolor}
\psclip{\pswedge[linestyle=none,fillstyle=none](11.661759,11.730468){1.421736}{11.238620}{116.225676}}
\psellipse(11.661759,11.730468)(1.005319,1.005319)
\endpsclip
\psset{linewidth=0.04} \psset{linestyle=solid}
\psset{linestyle=solid} \setlinecaps{0}
\newrgbcolor{dialinecolor}{0 0 0}
\psset{linecolor=dialinecolor}
\psclip{\pswedge[linestyle=none,fillstyle=none](10.901785,13.360747){1.442558}{288.497093}{335.898004}}
\psellipse(10.901785,13.360747)(1.020043,1.020043)
\endpsclip
\psset{linewidth=0.04} \psset{linestyle=solid}
\psset{linestyle=solid} \setlinecaps{0}
\newrgbcolor{dialinecolor}{0 0 0}
\psset{linecolor=dialinecolor}
\psclip{\pswedge[linestyle=none,fillstyle=none](11.303411,12.047637){1.049999}{43.754345}{96.380582}}
\psellipse(11.303411,12.047637)(0.742462,0.742462)
\endpsclip
}\endpspicture \]

These more general diagrams are made into a PROP, denoted by
$C_*\mathcal S$, similarly to the case of $C_*\mathcal S^c$
described in section \ref{sec1}. In fact, the tensor product, the
symmetric group action and the differential are exactly the same.
As for the composition, the following comments are in order.
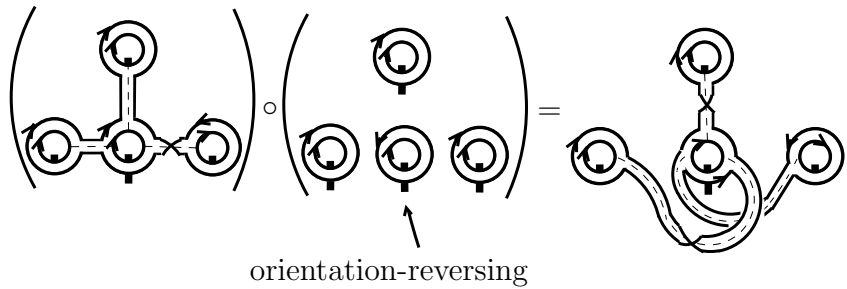
\[ \pspicture(2.029420,-13.302800)(13.162600,-9.633770)
\scalebox{1 -1}{
\newrgbcolor{dialinecolor}{0 0 0}
\psset{linecolor=dialinecolor}
\newrgbcolor{diafillcolor}{1 1 1}
\psset{fillcolor=diafillcolor} \psset{linewidth=0.04}
\psset{linestyle=solid} \psset{linestyle=solid} \setlinecaps{0}
\newrgbcolor{dialinecolor}{0 0 0}
\psset{linecolor=dialinecolor}
\psclip{\pswedge[linestyle=none,fillstyle=none](15.774738,10.233001){12.320986}{162.350483}{165.700255}}
\psellipse(15.774738,10.233001)(8.712253,8.712253)
\endpsclip
\psset{linewidth=0.04} \psset{linestyle=solid} \setlinejoinmode{0}
\setlinecaps{0}
\newrgbcolor{dialinecolor}{0 0 0}
\psset{linecolor=dialinecolor}
\psline(7.393922,12.426990)(7.321999,12.341391)(7.296674,12.450288)
\psset{linewidth=0.04} \psset{linestyle=solid}
\psset{linestyle=solid} \setlinecaps{0} \setlinejoinmode{0}
\setlinecaps{0} \setlinejoinmode{0} \psset{linestyle=solid}
\newrgbcolor{dialinecolor}{1 1 1}
\psset{linecolor=dialinecolor}
\psellipse*(11.290550,10.332470)(0.383350,0.383350)
\newrgbcolor{dialinecolor}{0 0 0}
\psset{linecolor=dialinecolor}
\psellipse(11.290550,10.332470)(0.383350,0.383350)
\psset{linewidth=0.004000} \setlinecaps{0} \setlinejoinmode{0}
\psset{linestyle=solid}
\newrgbcolor{dialinecolor}{0 0 0}
\psset{linecolor=dialinecolor}
\psellipse(11.290550,10.332470)(0.383350,0.383350)
\psset{linewidth=0.04} \psset{linestyle=solid}
\psset{linestyle=solid} \setlinecaps{0} \setlinejoinmode{0}
\setlinecaps{0} \setlinejoinmode{0} \psset{linestyle=solid}
\newrgbcolor{dialinecolor}{1 1 1}
\psset{linecolor=dialinecolor}
\psellipse*(11.288295,10.338495)(0.220595,0.220595)
\newrgbcolor{dialinecolor}{0 0 0}
\psset{linecolor=dialinecolor}
\psellipse(11.288295,10.338495)(0.220595,0.220595)
\psset{linewidth=0.004000} \setlinecaps{0} \setlinejoinmode{0}
\psset{linestyle=solid}
\newrgbcolor{dialinecolor}{0 0 0}
\psset{linecolor=dialinecolor}
\psellipse(11.288295,10.338495)(0.220595,0.220595)
\psset{linewidth=0.1} \psset{linestyle=solid} \psset{linestyle=solid}
\setlinecaps{0}
\newrgbcolor{dialinecolor}{0 0 0}
\psset{linecolor=dialinecolor}
\psline(11.288300,10.559100)(11.287700,10.442900)
\psset{linewidth=0.04} \psset{linestyle=solid}
\psset{linestyle=solid} \setlinecaps{0}
\newrgbcolor{dialinecolor}{0 0 0}
\psset{linecolor=dialinecolor}
\psline(11.083608,10.281365)(11.067700,10.338500)
\psset{linewidth=0.04} \psset{linestyle=solid} \setlinejoinmode{0}
\setlinecaps{0}
\newrgbcolor{dialinecolor}{0 0 0}
\psset{linecolor=dialinecolor}
\psline(11.127621,10.402903)(11.095604,10.238283)(10.983118,10.362668)
\psset{linewidth=0.04} \psset{linestyle=solid}
\psset{linestyle=solid} \setlinecaps{0}
\newrgbcolor{dialinecolor}{0 0 0}
\psset{linecolor=dialinecolor}
\psline(10.907214,10.332008)(10.907200,10.332500)
\psset{linewidth=0.04} \psset{linestyle=solid} \setlinejoinmode{0}
\setlinecaps{0}
\newrgbcolor{dialinecolor}{0 0 0}
\psset{linecolor=dialinecolor}
\psline(10.979258,10.439331)(10.908457,10.287304)(10.829316,10.435161)
\psset{linewidth=0.04} \psset{linestyle=solid}
\psset{linestyle=solid} \setlinecaps{0} \setlinejoinmode{0}
\setlinecaps{0} \setlinejoinmode{0} \psset{linestyle=solid}
\newrgbcolor{dialinecolor}{1 1 1}
\psset{linecolor=dialinecolor}
\psellipse*(11.338150,11.643150)(0.383350,0.383350)
\newrgbcolor{dialinecolor}{0 0 0}
\psset{linecolor=dialinecolor}
\psellipse(11.338150,11.643150)(0.383350,0.383350)
\psset{linewidth=0.004000} \setlinecaps{0} \setlinejoinmode{0}
\psset{linestyle=solid}
\newrgbcolor{dialinecolor}{0 0 0}
\psset{linecolor=dialinecolor}
\psellipse(11.338150,11.643150)(0.383350,0.383350)
\psset{linewidth=0.04} \psset{linestyle=solid}
\psset{linestyle=solid} \setlinecaps{0} \setlinejoinmode{0}
\setlinecaps{0} \setlinejoinmode{0} \psset{linestyle=solid}
\newrgbcolor{dialinecolor}{1 1 1}
\psset{linecolor=dialinecolor}
\psellipse*(11.325795,11.659295)(0.220595,0.220595)
\newrgbcolor{dialinecolor}{0 0 0}
\psset{linecolor=dialinecolor}
\psellipse(11.325795,11.659295)(0.220595,0.220595)
\psset{linewidth=0.004000} \setlinecaps{0} \setlinejoinmode{0}
\psset{linestyle=solid}
\newrgbcolor{dialinecolor}{0 0 0}
\psset{linecolor=dialinecolor}
\psellipse(11.325795,11.659295)(0.220595,0.220595)
\psset{linewidth=0.1} \psset{linestyle=solid} \psset{linestyle=solid}
\setlinecaps{0}
\newrgbcolor{dialinecolor}{0 0 0}
\psset{linecolor=dialinecolor}
\psline(11.325800,11.879900)(11.325200,11.763700)
\psset{linewidth=0.04} \psset{linestyle=solid}
\psset{linestyle=solid} \setlinecaps{0} \setlinejoinmode{0}
\setlinecaps{0} \setlinejoinmode{0} \psset{linestyle=solid}
\newrgbcolor{dialinecolor}{1 1 1}
\psset{linecolor=dialinecolor}
\psellipse*(9.896730,11.653350)(0.383350,0.383350)
\newrgbcolor{dialinecolor}{0 0 0}
\psset{linecolor=dialinecolor}
\psellipse(9.896730,11.653350)(0.383350,0.383350)
\psset{linewidth=0.004000} \setlinecaps{0} \setlinejoinmode{0}
\psset{linestyle=solid}
\newrgbcolor{dialinecolor}{0 0 0}
\psset{linecolor=dialinecolor}
\psellipse(9.896730,11.653350)(0.383350,0.383350)
\psset{linewidth=0.04} \psset{linestyle=solid}
\psset{linestyle=solid} \setlinecaps{0} \setlinejoinmode{0}
\setlinecaps{0} \setlinejoinmode{0} \psset{linestyle=solid}
\newrgbcolor{dialinecolor}{1 1 1}
\psset{linecolor=dialinecolor}
\psellipse*(9.894515,11.659295)(0.220595,0.220595)
\newrgbcolor{dialinecolor}{0 0 0}
\psset{linecolor=dialinecolor}
\psellipse(9.894515,11.659295)(0.220595,0.220595)
\psset{linewidth=0.004000} \setlinecaps{0} \setlinejoinmode{0}
\psset{linestyle=solid}
\newrgbcolor{dialinecolor}{0 0 0}
\psset{linecolor=dialinecolor}
\psellipse(9.894515,11.659295)(0.220595,0.220595)
\psset{linewidth=0.1} \psset{linestyle=solid} \psset{linestyle=solid}
\setlinecaps{0}
\newrgbcolor{dialinecolor}{0 0 0}
\psset{linecolor=dialinecolor}
\psline(9.894520,11.879900)(9.893920,11.763700)
\psset{linewidth=0.04} \psset{linestyle=solid}
\psset{linestyle=solid} \setlinecaps{0}
\newrgbcolor{dialinecolor}{0 0 0}
\psset{linecolor=dialinecolor}
\psline(9.694545,11.603448)(9.681270,11.663700)
\psset{linewidth=0.04} \psset{linestyle=solid} \setlinejoinmode{0}
\setlinecaps{0}
\newrgbcolor{dialinecolor}{0 0 0}
\psset{linecolor=dialinecolor}
\psline(9.745136,11.722398)(9.704168,11.559774)(9.598650,11.690123)
\psset{linewidth=0.04} \psset{linestyle=solid}
\psset{linestyle=solid} \setlinecaps{0}
\newrgbcolor{dialinecolor}{0 0 0}
\psset{linecolor=dialinecolor}
\psline(9.565468,11.472419)(9.579960,11.443400)
\psset{linewidth=0.04} \psset{linestyle=solid} \setlinejoinmode{0}
\setlinecaps{0}
\newrgbcolor{dialinecolor}{0 0 0}
\psset{linecolor=dialinecolor}
\psline(9.585528,11.600114)(9.585449,11.432409)(9.451332,11.533096)
\psset{linewidth=0.04} \psset{linestyle=solid}
\psset{linestyle=solid} \setlinecaps{0} \setlinejoinmode{0}
\setlinecaps{0} \setlinejoinmode{0} \psset{linestyle=solid}
\newrgbcolor{dialinecolor}{1 1 1}
\psset{linecolor=dialinecolor}
\psellipse*(12.759250,11.653350)(0.383350,0.383350)
\newrgbcolor{dialinecolor}{0 0 0}
\psset{linecolor=dialinecolor}
\psellipse(12.759250,11.653350)(0.383350,0.383350)
\psset{linewidth=0.004000} \setlinecaps{0} \setlinejoinmode{0}
\psset{linestyle=solid}
\newrgbcolor{dialinecolor}{0 0 0}
\psset{linecolor=dialinecolor}
\psellipse(12.759250,11.653350)(0.383350,0.383350)
\psset{linewidth=0.04} \psset{linestyle=solid}
\psset{linestyle=solid} \setlinecaps{0} \setlinejoinmode{0}
\setlinecaps{0} \setlinejoinmode{0} \psset{linestyle=solid}
\newrgbcolor{dialinecolor}{1 1 1}
\psset{linecolor=dialinecolor}
\psellipse*(12.756995,11.659295)(0.220595,0.220595)
\newrgbcolor{dialinecolor}{0 0 0}
\psset{linecolor=dialinecolor}
\psellipse(12.756995,11.659295)(0.220595,0.220595)
\psset{linewidth=0.004000} \setlinecaps{0} \setlinejoinmode{0}
\psset{linestyle=solid}
\newrgbcolor{dialinecolor}{0 0 0}
\psset{linecolor=dialinecolor}
\psellipse(12.756995,11.659295)(0.220595,0.220595)
\psset{linewidth=0.1} \psset{linestyle=solid} \psset{linestyle=solid}
\setlinecaps{0}
\newrgbcolor{dialinecolor}{0 0 0}
\psset{linecolor=dialinecolor}
\psline(12.757000,11.879900)(12.756400,11.763700)
\psset{linewidth=0.04} \psset{linestyle=solid}
\psset{linestyle=solid} \setlinecaps{0}
\newrgbcolor{dialinecolor}{0 0 0}
\psset{linecolor=dialinecolor}
\psline(12.860858,11.480818)(12.881100,11.494300)
\psset{linewidth=0.04} \psset{linestyle=solid} \setlinejoinmode{0}
\setlinecaps{0}
\newrgbcolor{dialinecolor}{0 0 0}
\psset{linecolor=dialinecolor}
\psline(12.731660,11.484879)(12.898079,11.505609)(12.814812,11.360036)
\psset{linewidth=0.04} \psset{linestyle=solid}
\psset{linestyle=solid} \setlinecaps{0}
\newrgbcolor{dialinecolor}{0 0 0}
\psset{linecolor=dialinecolor}
\psline(12.413984,11.481440)(12.435900,11.430900)
\psset{linewidth=0.04} \psset{linestyle=solid} \setlinejoinmode{0}
\setlinecaps{0}
\newrgbcolor{dialinecolor}{0 0 0}
\psset{linecolor=dialinecolor}
\psline(12.387059,11.355014)(12.396192,11.522470)(12.524677,11.414690)
\setfont{Helvetica}{0.400000}
\newrgbcolor{dialinecolor}{0 0 0}
\psset{linecolor=dialinecolor}
\rput[l](5.402880,11.015500){\scalebox{1 -1}{$\circ$}}
\setfont{Helvetica}{0.400000}
\newrgbcolor{dialinecolor}{0 0 0}
\psset{linecolor=dialinecolor}
\rput[l](9.067480,11.087500){\scalebox{1 -1}{=}}
\psset{linewidth=0.04} \psset{linestyle=solid}
\psset{linestyle=solid} \setlinecaps{0}
\newrgbcolor{dialinecolor}{0 0 0}
\psset{linecolor=dialinecolor}
\psclip{\pswedge[linestyle=none,fillstyle=none](8.141566,10.970715){3.490585}{150.447253}{209.552747}}
\psellipse(8.141566,10.970715)(2.468216,2.468216)
\endpsclip
\psset{linewidth=0.04} \psset{linestyle=solid}
\psset{linestyle=solid} \setlinecaps{0}
\newrgbcolor{dialinecolor}{0 0 0}
\psset{linecolor=dialinecolor}
\psclip{\pswedge[linestyle=none,fillstyle=none](6.058209,11.005115){4.075228}{335.009668}{24.990332}}
\psellipse(6.058209,11.005115)(2.881621,2.881621)
\endpsclip
\psset{linewidth=0.04} \psset{linestyle=solid}
\psset{linestyle=solid} \setlinecaps{0} \setlinejoinmode{0}
\setlinecaps{0} \setlinejoinmode{0} \psset{linestyle=solid}
\newrgbcolor{dialinecolor}{1 1 1}
\psset{linecolor=dialinecolor}
\psellipse*(7.265110,10.332470)(0.383350,0.383350)
\newrgbcolor{dialinecolor}{0 0 0}
\psset{linecolor=dialinecolor}
\psellipse(7.265110,10.332470)(0.383350,0.383350)
\psset{linewidth=0.004000} \setlinecaps{0} \setlinejoinmode{0}
\psset{linestyle=solid}
\newrgbcolor{dialinecolor}{0 0 0}
\psset{linecolor=dialinecolor}
\psellipse(7.265110,10.332470)(0.383350,0.383350)
\psset{linewidth=0.04} \psset{linestyle=solid}
\psset{linestyle=solid} \setlinecaps{0} \setlinejoinmode{0}
\setlinecaps{0} \setlinejoinmode{0} \psset{linestyle=solid}
\newrgbcolor{dialinecolor}{1 1 1}
\psset{linecolor=dialinecolor}
\psellipse*(7.262895,10.338495)(0.220595,0.220595)
\newrgbcolor{dialinecolor}{0 0 0}
\psset{linecolor=dialinecolor}
\psellipse(7.262895,10.338495)(0.220595,0.220595)
\psset{linewidth=0.004000} \setlinecaps{0} \setlinejoinmode{0}
\psset{linestyle=solid}
\newrgbcolor{dialinecolor}{0 0 0}
\psset{linecolor=dialinecolor}
\psellipse(7.262895,10.338495)(0.220595,0.220595)
\psset{linewidth=0.1} \psset{linestyle=solid} \psset{linestyle=solid}
\setlinecaps{0}
\newrgbcolor{dialinecolor}{0 0 0}
\psset{linecolor=dialinecolor}
\psline(7.262900,10.559100)(7.262300,10.442900)
\psset{linewidth=0.04} \psset{linestyle=solid}
\psset{linestyle=solid} \setlinecaps{0}
\newrgbcolor{dialinecolor}{0 0 0}
\psset{linecolor=dialinecolor}
\psline(7.058195,10.281370)(7.042300,10.338500)
\psset{linewidth=0.04} \psset{linestyle=solid} \setlinejoinmode{0}
\setlinecaps{0}
\newrgbcolor{dialinecolor}{0 0 0}
\psset{linecolor=dialinecolor}
\psline(7.102231,10.402899)(7.070183,10.238285)(6.957720,10.362692)
\psset{linewidth=0.04} \psset{linestyle=solid}
\psset{linestyle=solid} \setlinecaps{0}
\newrgbcolor{dialinecolor}{0 0 0}
\psset{linecolor=dialinecolor}
\psline(6.974294,10.088161)(6.938500,10.135900)
\psset{linewidth=0.04} \psset{linestyle=solid} \setlinejoinmode{0}
\setlinecaps{0}
\newrgbcolor{dialinecolor}{0 0 0}
\psset{linecolor=dialinecolor}
\psline(6.971144,10.217385)(7.001122,10.052381)(6.851132,10.127401)
\psset{linewidth=0.1} \psset{linestyle=solid} \psset{linestyle=solid}
\setlinecaps{0}
\newrgbcolor{dialinecolor}{0 0 0}
\psset{linecolor=dialinecolor}
\psline(7.263680,10.833800)(7.265110,10.715800)
\psset{linewidth=0.04} \psset{linestyle=solid}
\psset{linestyle=solid} \setlinecaps{0} \setlinejoinmode{0}
\setlinecaps{0} \setlinejoinmode{0} \psset{linestyle=solid}
\newrgbcolor{dialinecolor}{1 1 1}
\psset{linecolor=dialinecolor}
\psellipse*(7.299040,11.627550)(0.383350,0.383350)
\newrgbcolor{dialinecolor}{0 0 0}
\psset{linecolor=dialinecolor}
\psellipse(7.299040,11.627550)(0.383350,0.383350)
\psset{linewidth=0.004000} \setlinecaps{0} \setlinejoinmode{0}
\psset{linestyle=solid}
\newrgbcolor{dialinecolor}{0 0 0}
\psset{linecolor=dialinecolor}
\psellipse(7.299040,11.627550)(0.383350,0.383350)
\psset{linewidth=0.04} \psset{linestyle=solid}
\psset{linestyle=solid} \setlinecaps{0} \setlinejoinmode{0}
\setlinecaps{0} \setlinejoinmode{0} \psset{linestyle=solid}
\newrgbcolor{dialinecolor}{1 1 1}
\psset{linecolor=dialinecolor}
\psellipse*(7.296825,11.633595)(0.220595,0.220595)
\newrgbcolor{dialinecolor}{0 0 0}
\psset{linecolor=dialinecolor}
\psellipse(7.296825,11.633595)(0.220595,0.220595)
\psset{linewidth=0.004000} \setlinecaps{0} \setlinejoinmode{0}
\psset{linestyle=solid}
\newrgbcolor{dialinecolor}{0 0 0}
\psset{linecolor=dialinecolor}
\psellipse(7.296825,11.633595)(0.220595,0.220595)
\psset{linewidth=0.1} \psset{linestyle=solid} \psset{linestyle=solid}
\setlinecaps{0}
\newrgbcolor{dialinecolor}{0 0 0}
\psset{linecolor=dialinecolor}
\psline(7.296830,11.854200)(7.296230,11.738000)
\psset{linewidth=0.04} \psset{linestyle=solid}
\psset{linestyle=solid} \setlinecaps{0}
\newrgbcolor{dialinecolor}{0 0 0}
\psset{linecolor=dialinecolor}
\psline(7.096859,11.577647)(7.083580,11.637900)
\psset{linewidth=0.04} \psset{linestyle=solid} \setlinejoinmode{0}
\setlinecaps{0}
\newrgbcolor{dialinecolor}{0 0 0}
\psset{linecolor=dialinecolor}
\psline(7.147443,11.696600)(7.106485,11.533973)(7.000958,11.664316)
\psset{linewidth=0.04} \psset{linestyle=solid}
\psset{linestyle=solid} \setlinecaps{0}
\newrgbcolor{dialinecolor}{0 0 0}
\psset{linecolor=dialinecolor}
\psline(6.953825,11.455649)(6.975760,11.405100)
\psset{linewidth=0.04} \psset{linestyle=solid} \setlinejoinmode{0}
\setlinecaps{0}
\newrgbcolor{dialinecolor}{0 0 0}
\psset{linecolor=dialinecolor}
\psline(6.926932,11.329216)(6.936022,11.496675)(7.064535,11.388927)
\psset{linewidth=0.1} \psset{linestyle=solid} \psset{linestyle=solid}
\setlinecaps{0}
\newrgbcolor{dialinecolor}{0 0 0}
\psset{linecolor=dialinecolor}
\psline(7.301230,12.129300)(7.299040,12.010900)
\psset{linewidth=0.04} \psset{linestyle=solid}
\psset{linestyle=solid} \setlinecaps{0} \setlinejoinmode{0}
\setlinecaps{0} \setlinejoinmode{0} \psset{linestyle=solid}
\newrgbcolor{dialinecolor}{1 1 1}
\psset{linecolor=dialinecolor}
\psellipse*(8.297770,11.636150)(0.383350,0.383350)
\newrgbcolor{dialinecolor}{0 0 0}
\psset{linecolor=dialinecolor}
\psellipse(8.297770,11.636150)(0.383350,0.383350)
\psset{linewidth=0.004000} \setlinecaps{0} \setlinejoinmode{0}
\psset{linestyle=solid}
\newrgbcolor{dialinecolor}{0 0 0}
\psset{linecolor=dialinecolor}
\psellipse(8.297770,11.636150)(0.383350,0.383350)
\psset{linewidth=0.04} \psset{linestyle=solid}
\psset{linestyle=solid} \setlinecaps{0} \setlinejoinmode{0}
\setlinecaps{0} \setlinejoinmode{0} \psset{linestyle=solid}
\newrgbcolor{dialinecolor}{1 1 1}
\psset{linecolor=dialinecolor}
\psellipse*(8.295555,11.642195)(0.220595,0.220595)
\newrgbcolor{dialinecolor}{0 0 0}
\psset{linecolor=dialinecolor}
\psellipse(8.295555,11.642195)(0.220595,0.220595)
\psset{linewidth=0.004000} \setlinecaps{0} \setlinejoinmode{0}
\psset{linestyle=solid}
\newrgbcolor{dialinecolor}{0 0 0}
\psset{linecolor=dialinecolor}
\psellipse(8.295555,11.642195)(0.220595,0.220595)
\psset{linewidth=0.1} \psset{linestyle=solid} \psset{linestyle=solid}
\setlinecaps{0}
\newrgbcolor{dialinecolor}{0 0 0}
\psset{linecolor=dialinecolor}
\psline(8.295560,11.862800)(8.294960,11.746600)
\psset{linewidth=0.04} \psset{linestyle=solid}
\psset{linestyle=solid} \setlinecaps{0}
\newrgbcolor{dialinecolor}{0 0 0}
\psset{linecolor=dialinecolor}
\psline(8.095585,11.586248)(8.082310,11.646500)
\psset{linewidth=0.04} \psset{linestyle=solid} \setlinejoinmode{0}
\setlinecaps{0}
\newrgbcolor{dialinecolor}{0 0 0}
\psset{linecolor=dialinecolor}
\psline(8.146176,11.705198)(8.105208,11.542574)(7.999690,11.672923)
\psset{linewidth=0.04} \psset{linestyle=solid}
\psset{linestyle=solid} \setlinecaps{0}
\newrgbcolor{dialinecolor}{0 0 0}
\psset{linecolor=dialinecolor}
\psline(8.006950,11.391858)(7.971150,11.439600)
\psset{linewidth=0.04} \psset{linestyle=solid} \setlinejoinmode{0}
\setlinecaps{0}
\newrgbcolor{dialinecolor}{0 0 0}
\psset{linecolor=dialinecolor}
\psline(8.003793,11.521082)(8.033780,11.356079)(7.883786,11.431091)
\psset{linewidth=0.1} \psset{linestyle=solid} \psset{linestyle=solid}
\setlinecaps{0}
\newrgbcolor{dialinecolor}{0 0 0}
\psset{linecolor=dialinecolor}
\psline(8.296340,12.137500)(8.297770,12.019500)
\psset{linewidth=0.04} \psset{linestyle=solid}
\psset{linestyle=solid} \setlinecaps{0} \setlinejoinmode{0}
\setlinecaps{0} \setlinejoinmode{0} \psset{linestyle=solid}
\newrgbcolor{dialinecolor}{1 1 1}
\psset{linecolor=dialinecolor}
\psellipse*(6.314890,11.606350)(0.383350,0.383350)
\newrgbcolor{dialinecolor}{0 0 0}
\psset{linecolor=dialinecolor}
\psellipse(6.314890,11.606350)(0.383350,0.383350)
\psset{linewidth=0.004000} \setlinecaps{0} \setlinejoinmode{0}
\psset{linestyle=solid}
\newrgbcolor{dialinecolor}{0 0 0}
\psset{linecolor=dialinecolor}
\psellipse(6.314890,11.606350)(0.383350,0.383350)
\psset{linewidth=0.04} \psset{linestyle=solid}
\psset{linestyle=solid} \setlinecaps{0} \setlinejoinmode{0}
\setlinecaps{0} \setlinejoinmode{0} \psset{linestyle=solid}
\newrgbcolor{dialinecolor}{1 1 1}
\psset{linecolor=dialinecolor}
\psellipse*(6.312675,11.612395)(0.220595,0.220595)
\newrgbcolor{dialinecolor}{0 0 0}
\psset{linecolor=dialinecolor}
\psellipse(6.312675,11.612395)(0.220595,0.220595)
\psset{linewidth=0.004000} \setlinecaps{0} \setlinejoinmode{0}
\psset{linestyle=solid}
\newrgbcolor{dialinecolor}{0 0 0}
\psset{linecolor=dialinecolor}
\psellipse(6.312675,11.612395)(0.220595,0.220595)
\psset{linewidth=0.1} \psset{linestyle=solid} \psset{linestyle=solid}
\setlinecaps{0}
\newrgbcolor{dialinecolor}{0 0 0}
\psset{linecolor=dialinecolor}
\psline(6.312680,11.833000)(6.312080,11.716800)
\psset{linewidth=0.04} \psset{linestyle=solid}
\psset{linestyle=solid} \setlinecaps{0}
\newrgbcolor{dialinecolor}{0 0 0}
\psset{linecolor=dialinecolor}
\psline(6.112705,11.556448)(6.099430,11.616700)
\psset{linewidth=0.04} \psset{linestyle=solid} \setlinejoinmode{0}
\setlinecaps{0}
\newrgbcolor{dialinecolor}{0 0 0}
\psset{linecolor=dialinecolor}
\psline(6.163296,11.675398)(6.122328,11.512774)(6.016810,11.643123)
\psset{linewidth=0.04} \psset{linestyle=solid}
\psset{linestyle=solid} \setlinecaps{0}
\newrgbcolor{dialinecolor}{0 0 0}
\psset{linecolor=dialinecolor}
\psline(6.024041,11.362037)(5.988270,11.409700)
\psset{linewidth=0.04} \psset{linestyle=solid} \setlinejoinmode{0}
\setlinecaps{0}
\newrgbcolor{dialinecolor}{0 0 0}
\psset{linecolor=dialinecolor}
\psline(6.020833,11.491259)(6.050886,11.326268)(5.900862,11.401220)
\psset{linewidth=0.1} \psset{linestyle=solid} \psset{linestyle=solid}
\setlinecaps{0}
\newrgbcolor{dialinecolor}{0 0 0}
\psset{linecolor=dialinecolor}
\psline(6.313460,12.107600)(6.314890,11.989700)
\psset{linewidth=0.04} \psset{linestyle=solid}
\psset{linestyle=solid} \setlinecaps{0}
\newrgbcolor{dialinecolor}{0 0 0}
\psset{linecolor=dialinecolor}
\psclip{\pswedge[linestyle=none,fillstyle=none](4.517750,10.871185){3.490746}{150.447950}{209.552050}}
\psellipse(4.517750,10.871185)(2.468330,2.468330)
\endpsclip
\psset{linewidth=0.04} \psset{linestyle=solid}
\psset{linestyle=solid} \setlinecaps{0}
\newrgbcolor{dialinecolor}{0 0 0}
\psset{linecolor=dialinecolor}
\psclip{\pswedge[linestyle=none,fillstyle=none](2.434166,10.905590){4.075387}{335.010165}{24.989835}}
\psellipse(2.434166,10.905590)(2.881734,2.881734)
\endpsclip
\psset{linewidth=0.04} \psset{linestyle=solid}
\psset{linestyle=solid} \setlinecaps{0} \setlinejoinmode{0}
\setlinecaps{0} \setlinejoinmode{0} \psset{linestyle=solid}
\newrgbcolor{dialinecolor}{1 1 1}
\psset{linecolor=dialinecolor}
\psellipse*(3.626270,10.232910)(0.383350,0.383350)
\newrgbcolor{dialinecolor}{0 0 0}
\psset{linecolor=dialinecolor}
\psellipse(3.626270,10.232910)(0.383350,0.383350)
\psset{linewidth=0.004000} \setlinecaps{0} \setlinejoinmode{0}
\psset{linestyle=solid}
\newrgbcolor{dialinecolor}{0 0 0}
\psset{linecolor=dialinecolor}
\psellipse(3.626270,10.232910)(0.383350,0.383350)
\psset{linewidth=0.04} \psset{linestyle=solid}
\psset{linestyle=solid} \setlinecaps{0} \setlinejoinmode{0}
\setlinecaps{0} \setlinejoinmode{0} \psset{linestyle=solid}
\newrgbcolor{dialinecolor}{1 1 1}
\psset{linecolor=dialinecolor}
\psellipse*(3.624065,10.238895)(0.220595,0.220595)
\newrgbcolor{dialinecolor}{0 0 0}
\psset{linecolor=dialinecolor}
\psellipse(3.624065,10.238895)(0.220595,0.220595)
\psset{linewidth=0.004000} \setlinecaps{0} \setlinejoinmode{0}
\psset{linestyle=solid}
\newrgbcolor{dialinecolor}{0 0 0}
\psset{linecolor=dialinecolor}
\psellipse(3.624065,10.238895)(0.220595,0.220595)
\psset{linewidth=0.1} \psset{linestyle=solid} \psset{linestyle=solid}
\setlinecaps{0}
\newrgbcolor{dialinecolor}{0 0 0}
\psset{linecolor=dialinecolor}
\psline(3.624060,10.459500)(3.623470,10.343300)
\psset{linewidth=0.04} \psset{linestyle=solid}
\psset{linestyle=solid} \setlinecaps{0}
\newrgbcolor{dialinecolor}{0 0 0}
\psset{linecolor=dialinecolor}
\psline(3.419350,10.181865)(3.403470,10.238900)
\psset{linewidth=0.04} \psset{linestyle=solid} \setlinejoinmode{0}
\setlinecaps{0}
\newrgbcolor{dialinecolor}{0 0 0}
\psset{linecolor=dialinecolor}
\psline(3.463364,10.303403)(3.431345,10.138783)(3.318860,10.263170)
\psset{linewidth=0.04} \psset{linestyle=solid}
\psset{linestyle=solid} \setlinecaps{0}
\newrgbcolor{dialinecolor}{0 0 0}
\psset{linecolor=dialinecolor}
\psline(3.335442,9.988593)(3.299660,10.036300) \psset{linewidth=0.04}
\psset{linestyle=solid} \setlinejoinmode{0} \setlinecaps{0}
\newrgbcolor{dialinecolor}{0 0 0}
\psset{linecolor=dialinecolor}
\psline(3.332272,10.117816)(3.362276,9.952816)(3.212274,10.027812)
\psset{linewidth=0.04} \psset{linestyle=solid}
\psset{linestyle=solid} \setlinecaps{0} \setlinejoinmode{0}
\setlinecaps{0} \setlinejoinmode{0} \psset{linestyle=solid}
\newrgbcolor{dialinecolor}{1 1 1}
\psset{linecolor=dialinecolor}
\psellipse*(4.748380,11.536650)(0.383350,0.383350)
\newrgbcolor{dialinecolor}{0 0 0}
\psset{linecolor=dialinecolor}
\psellipse(4.748380,11.536650)(0.383350,0.383350)
\psset{linewidth=0.004000} \setlinecaps{0} \setlinejoinmode{0}
\psset{linestyle=solid}
\newrgbcolor{dialinecolor}{0 0 0}
\psset{linecolor=dialinecolor}
\psellipse(4.748380,11.536650)(0.383350,0.383350)
\psset{linewidth=0.04} \psset{linestyle=solid}
\psset{linestyle=solid} \setlinecaps{0} \setlinejoinmode{0}
\setlinecaps{0} \setlinejoinmode{0} \psset{linestyle=solid}
\newrgbcolor{dialinecolor}{1 1 1}
\psset{linecolor=dialinecolor}
\psellipse*(4.746175,11.542595)(0.220595,0.220595)
\newrgbcolor{dialinecolor}{0 0 0}
\psset{linecolor=dialinecolor}
\psellipse(4.746175,11.542595)(0.220595,0.220595)
\psset{linewidth=0.004000} \setlinecaps{0} \setlinejoinmode{0}
\psset{linestyle=solid}
\newrgbcolor{dialinecolor}{0 0 0}
\psset{linecolor=dialinecolor}
\psellipse(4.746175,11.542595)(0.220595,0.220595)
\psset{linewidth=0.1} \psset{linestyle=solid} \psset{linestyle=solid}
\setlinecaps{0}
\newrgbcolor{dialinecolor}{0 0 0}
\psset{linecolor=dialinecolor}
\psline(4.746170,11.763200)(4.745570,11.647000)
\psset{linewidth=0.04} \psset{linestyle=solid}
\psset{linestyle=solid} \setlinecaps{0}
\newrgbcolor{dialinecolor}{0 0 0}
\psset{linecolor=dialinecolor}
\psline(4.663525,11.356178)(4.693460,11.343800)
\psset{linewidth=0.04} \psset{linestyle=solid} \setlinejoinmode{0}
\setlinecaps{0}
\newrgbcolor{dialinecolor}{0 0 0}
\psset{linecolor=dialinecolor}
\psline(4.594895,11.465716)(4.704853,11.339089)(4.537577,11.327100)
\psset{linewidth=0.04} \psset{linestyle=solid}
\psset{linestyle=solid} \setlinecaps{0}
\newrgbcolor{dialinecolor}{0 0 0}
\psset{linecolor=dialinecolor}
\psline(4.595708,11.191345)(4.635270,11.177600)
\psset{linewidth=0.04} \psset{linestyle=solid} \setlinejoinmode{0}
\setlinecaps{0}
\newrgbcolor{dialinecolor}{0 0 0}
\psset{linecolor=dialinecolor}
\psline(4.670541,11.085948)(4.553464,11.206023)(4.719770,11.227639)
\psset{linewidth=0.04} \psset{linestyle=solid}
\psset{linestyle=solid} \setlinecaps{0} \setlinejoinmode{0}
\setlinecaps{0} \setlinejoinmode{0} \psset{linestyle=solid}
\newrgbcolor{dialinecolor}{1 1 1}
\psset{linecolor=dialinecolor}
\psellipse*(2.646230,11.521750)(0.383350,0.383350)
\newrgbcolor{dialinecolor}{0 0 0}
\psset{linecolor=dialinecolor}
\psellipse(2.646230,11.521750)(0.383350,0.383350)
\psset{linewidth=0.004000} \setlinecaps{0} \setlinejoinmode{0}
\psset{linestyle=solid}
\newrgbcolor{dialinecolor}{0 0 0}
\psset{linecolor=dialinecolor}
\psellipse(2.646230,11.521750)(0.383350,0.383350)
\psset{linewidth=0.04} \psset{linestyle=solid}
\psset{linestyle=solid} \setlinecaps{0} \setlinejoinmode{0}
\setlinecaps{0} \setlinejoinmode{0} \psset{linestyle=solid}
\newrgbcolor{dialinecolor}{1 1 1}
\psset{linecolor=dialinecolor}
\psellipse*(2.644025,11.527695)(0.220595,0.220595)
\newrgbcolor{dialinecolor}{0 0 0}
\psset{linecolor=dialinecolor}
\psellipse(2.644025,11.527695)(0.220595,0.220595)
\psset{linewidth=0.004000} \setlinecaps{0} \setlinejoinmode{0}
\psset{linestyle=solid}
\newrgbcolor{dialinecolor}{0 0 0}
\psset{linecolor=dialinecolor}
\psellipse(2.644025,11.527695)(0.220595,0.220595)
\psset{linewidth=0.1} \psset{linestyle=solid} \psset{linestyle=solid}
\setlinecaps{0}
\newrgbcolor{dialinecolor}{0 0 0}
\psset{linecolor=dialinecolor}
\psline(2.644030,11.748300)(2.643430,11.632100)
\psset{linewidth=0.04} \psset{linestyle=solid}
\psset{linestyle=solid} \setlinecaps{0}
\newrgbcolor{dialinecolor}{0 0 0}
\psset{linecolor=dialinecolor}
\psline(2.444055,11.471748)(2.430780,11.532000)
\psset{linewidth=0.04} \psset{linestyle=solid} \setlinejoinmode{0}
\setlinecaps{0}
\newrgbcolor{dialinecolor}{0 0 0}
\psset{linecolor=dialinecolor}
\psline(2.494646,11.590698)(2.453678,11.428074)(2.348160,11.558423)
\psset{linewidth=0.04} \psset{linestyle=solid}
\psset{linestyle=solid} \setlinecaps{0}
\newrgbcolor{dialinecolor}{0 0 0}
\psset{linecolor=dialinecolor}
\psline(2.375141,11.220373)(2.349440,11.260700)
\psset{linewidth=0.04} \psset{linestyle=solid} \setlinejoinmode{0}
\setlinecaps{0}
\newrgbcolor{dialinecolor}{0 0 0}
\psset{linecolor=dialinecolor}
\psline(2.333735,11.424890)(2.351106,11.258087)(2.207241,11.344272)
\psset{linewidth=0.04} \psset{linestyle=solid}
\psset{linestyle=solid} \setlinecaps{0} \setlinejoinmode{0}
\setlinecaps{0} \setlinejoinmode{0} \psset{linestyle=solid}
\newrgbcolor{dialinecolor}{1 1 1}
\psset{linecolor=dialinecolor}
\psellipse*(3.631290,11.516850)(0.383350,0.383350)
\newrgbcolor{dialinecolor}{0 0 0}
\psset{linecolor=dialinecolor}
\psellipse(3.631290,11.516850)(0.383350,0.383350)
\psset{linewidth=0.004000} \setlinecaps{0} \setlinejoinmode{0}
\psset{linestyle=solid}
\newrgbcolor{dialinecolor}{0 0 0}
\psset{linecolor=dialinecolor}
\psellipse(3.631290,11.516850)(0.383350,0.383350)
\psset{linewidth=0.04} \psset{linestyle=solid}
\psset{linestyle=solid} \setlinecaps{0} \setlinejoinmode{0}
\setlinecaps{0} \setlinejoinmode{0} \psset{linestyle=solid}
\newrgbcolor{dialinecolor}{1 1 1}
\psset{linecolor=dialinecolor}
\psellipse*(3.629085,11.522895)(0.220595,0.220595)
\newrgbcolor{dialinecolor}{0 0 0}
\psset{linecolor=dialinecolor}
\psellipse(3.629085,11.522895)(0.220595,0.220595)
\psset{linewidth=0.004000} \setlinecaps{0} \setlinejoinmode{0}
\psset{linestyle=solid}
\newrgbcolor{dialinecolor}{0 0 0}
\psset{linecolor=dialinecolor}
\psellipse(3.629085,11.522895)(0.220595,0.220595)
\psset{linewidth=0.1} \psset{linestyle=solid} \psset{linestyle=solid}
\setlinecaps{0}
\newrgbcolor{dialinecolor}{0 0 0}
\psset{linecolor=dialinecolor}
\psline(3.629090,11.743500)(3.628490,11.627300)
\psset{linewidth=0.04} \psset{linestyle=solid}
\psset{linestyle=solid} \setlinecaps{0}
\newrgbcolor{dialinecolor}{0 0 0}
\psset{linecolor=dialinecolor}
\psline(3.429115,11.466948)(3.415840,11.527200)
\psset{linewidth=0.04} \psset{linestyle=solid} \setlinejoinmode{0}
\setlinecaps{0}
\newrgbcolor{dialinecolor}{0 0 0}
\psset{linecolor=dialinecolor}
\psline(3.479706,11.585898)(3.438738,11.423274)(3.333220,11.553623)
\psset{linewidth=0.04} \psset{linestyle=solid}
\psset{linestyle=solid} \setlinecaps{0}
\newrgbcolor{dialinecolor}{0 0 0}
\psset{linecolor=dialinecolor}
\psline(3.340474,11.272561)(3.304680,11.320300)
\psset{linewidth=0.04} \psset{linestyle=solid} \setlinejoinmode{0}
\setlinecaps{0}
\newrgbcolor{dialinecolor}{0 0 0}
\psset{linecolor=dialinecolor}
\psline(3.337324,11.401785)(3.367302,11.236781)(3.217312,11.311801)
\psset{linewidth=0.1} \psset{linestyle=solid} \psset{linestyle=solid}
\setlinecaps{0}
\newrgbcolor{dialinecolor}{0 0 0}
\psset{linecolor=dialinecolor}
\psline(3.629870,12.018200)(3.631290,11.900200)
\psset{linewidth=0.04} \psset{linestyle=solid}
\psset{linestyle=solid} \setlinecaps{0}
\newrgbcolor{dialinecolor}{0 0 0}
\psset{linecolor=dialinecolor}
\psline(3.997540,11.418400)(4.376480,11.416000)
\psset{linewidth=0.04} \psset{linestyle=solid}
\psset{linestyle=solid} \setlinecaps{0}
\newrgbcolor{dialinecolor}{0 0 0}
\psset{linecolor=dialinecolor}
\psline(4.007770,11.619500)(4.371740,11.622000)
\psset{linewidth=0.04} \psset{linestyle=solid}
\psset{linestyle=solid} \setlinecaps{0}
\newrgbcolor{dialinecolor}{0 0 0}
\psset{linecolor=dialinecolor}
\psline(3.528610,11.160200)(3.526240,10.601300)
\psset{linewidth=0.04} \psset{linestyle=solid}
\psset{linestyle=solid} \setlinecaps{0}
\newrgbcolor{dialinecolor}{0 0 0}
\psset{linecolor=dialinecolor}
\psline(3.706240,11.138900)(3.710970,10.603600)
\newrgbcolor{dialinecolor}{1 1 1}
\psset{linecolor=dialinecolor}
\pspolygon*(3.566510,10.551500)(3.566510,11.183850)(3.670720,11.183850)(3.670720,10.551500)
\psset{linewidth=0.04} \psset{linestyle=solid}
\psset{linestyle=solid} \setlinejoinmode{0}
\newrgbcolor{dialinecolor}{1 1 1}
\psset{linecolor=dialinecolor}
\pspolygon(3.566510,10.551500)(3.566510,11.183850)(3.670720,11.183850)(3.670720,10.551500)
\psset{linewidth=0.01} \psset{linestyle=dashed,dash=1 1}
\psset{linestyle=dashed,dash=0.1 0.1} \setlinecaps{0}
\newrgbcolor{dialinecolor}{0 0 0}
\psset{linecolor=dialinecolor}
\psline(3.623760,10.401400)(3.629090,11.302300)
\newrgbcolor{dialinecolor}{1 1 1}
\psset{linecolor=dialinecolor}
\pspolygon*(3.960230,11.456200)(3.960230,11.586460)(4.443681,11.586460)(4.443681,11.456200)
\psset{linewidth=0.04} \psset{linestyle=solid}
\psset{linestyle=solid} \setlinejoinmode{0}
\newrgbcolor{dialinecolor}{1 1 1}
\psset{linecolor=dialinecolor}
\pspolygon(3.960230,11.456200)(3.960230,11.586460)(4.443681,11.586460)(4.443681,11.456200)
\newrgbcolor{dialinecolor}{1 1 1}
\psset{linecolor=dialinecolor}
\pspolygon*(4.081630,11.394000)(4.081630,11.678280)(4.284000,11.678280)(4.284000,11.394000)
\psset{linewidth=0.04} \psset{linestyle=solid}
\psset{linestyle=solid} \setlinejoinmode{0}
\newrgbcolor{dialinecolor}{1 1 1}
\psset{linecolor=dialinecolor}
\pspolygon(4.081630,11.394000)(4.081630,11.678280)(4.284000,11.678280)(4.284000,11.394000)
\psset{linewidth=0.04} \psset{linestyle=solid}
\psset{linestyle=solid} \setlinecaps{0}
\newrgbcolor{dialinecolor}{0 0 0}
\psset{linecolor=dialinecolor}
\psclip{\pswedge[linestyle=none,fillstyle=none](3.867599,11.923293){0.771714}{290.389987}{328.157784}}
\psellipse(3.867599,11.923293)(0.545684,0.545684)
\endpsclip
\psset{linewidth=0.04} \psset{linestyle=solid}
\psset{linestyle=solid} \setlinecaps{0}
\newrgbcolor{dialinecolor}{0 0 0}
\psset{linecolor=dialinecolor}
\psclip{\pswedge[linestyle=none,fillstyle=none](4.496874,11.885800){0.708495}{210.544160}{249.977082}}
\psellipse(4.496874,11.885800)(0.500981,0.500981)
\endpsclip
\setfont{Helvetica}{0.400000}
\newrgbcolor{dialinecolor}{0 0 0}
\psset{linecolor=dialinecolor}
\rput[l](5.239740,13.202800){\scalebox{1
-1}{orientation-reversing}} \psset{linewidth=0.01}
\psset{linestyle=dashed,dash=1 1} \psset{linestyle=dashed,dash=0.1
0.1} \setlinecaps{0}
\newrgbcolor{dialinecolor}{0 0 0}
\psset{linecolor=dialinecolor}
\psclip{\pswedge[linestyle=none,fillstyle=none](9.744554,10.379347){4.343428}{24.629618}{38.092197}}
\psellipse(9.744554,10.379347)(3.071267,3.071267)
\endpsclip
\psset{linewidth=0.1} \psset{linestyle=solid} \psset{linestyle=solid}
\setlinecaps{0}
\newrgbcolor{dialinecolor}{0 0 0}
\psset{linecolor=dialinecolor}
\psline(11.337000,12.015200)(11.337500,12.166800)
\psset{linewidth=0.04} \psset{linestyle=solid}
\psset{linestyle=solid} \setlinecaps{0}
\newrgbcolor{dialinecolor}{0 0 0}
\psset{linecolor=dialinecolor}
\psline(11.218700,11.288300)(11.220300,10.693600)
\psset{linewidth=0.04} \psset{linestyle=solid}
\psset{linestyle=solid} \setlinecaps{0}
\newrgbcolor{dialinecolor}{0 0 0}
\psset{linecolor=dialinecolor}
\psline(11.398400,11.268400)(11.396000,10.709500)
\newrgbcolor{dialinecolor}{1 1 1}
\psset{linecolor=dialinecolor}
\pspolygon*(11.156400,10.869100)(11.156400,11.053830)(11.447710,11.053830)(11.447710,10.869100)
\psset{linewidth=0.04} \psset{linestyle=solid}
\psset{linestyle=solid} \setlinejoinmode{0}
\newrgbcolor{dialinecolor}{1 1 1}
\psset{linecolor=dialinecolor}
\pspolygon(11.156400,10.869100)(11.156400,11.053830)(11.447710,11.053830)(11.447710,10.869100)
\newrgbcolor{dialinecolor}{1 1 1}
\psset{linecolor=dialinecolor}
\pspolygon*(11.254100,10.600600)(11.254100,11.318270)(11.360120,11.318270)(11.360120,10.600600)
\psset{linewidth=0.04} \psset{linestyle=solid}
\psset{linestyle=solid} \setlinejoinmode{0}
\newrgbcolor{dialinecolor}{1 1 1}
\psset{linecolor=dialinecolor}
\pspolygon(11.254100,10.600600)(11.254100,11.318270)(11.360120,11.318270)(11.360120,10.600600)
\psset{linewidth=0.01} \psset{linestyle=dashed,dash=1 1}
\psset{linestyle=dashed,dash=0.1 0.1} \setlinecaps{0}
\newrgbcolor{dialinecolor}{0 0 0}
\psset{linecolor=dialinecolor}
\psline(11.288300,10.559100)(11.325800,11.438700)
\psset{linewidth=0.04} \psset{linestyle=solid}
\psset{linestyle=solid} \setlinecaps{0}
\newrgbcolor{dialinecolor}{0 0 0}
\psset{linecolor=dialinecolor}
\psclip{\pswedge[linestyle=none,fillstyle=none](11.895749,10.511369){1.072073}{130.115578}{153.948252}}
\psellipse(11.895749,10.511369)(0.758070,0.758070)
\endpsclip
\psset{linewidth=0.04} \psset{linestyle=solid}
\psset{linestyle=solid} \setlinecaps{0}
\newrgbcolor{dialinecolor}{0 0 0}
\psset{linecolor=dialinecolor}
\psclip{\pswedge[linestyle=none,fillstyle=none](11.062562,10.767772){0.496257}{12.078589}{64.577989}}
\psellipse(11.062562,10.767772)(0.350906,0.350906)
\endpsclip
\newrgbcolor{dialinecolor}{1 1 1}
\psset{linecolor=dialinecolor}
\pspolygon*(11.597000,11.686900)(11.597000,11.815490)(11.775880,11.815490)(11.775880,11.686900)
\psset{linewidth=0.04} \psset{linestyle=solid}
\psset{linestyle=solid} \setlinejoinmode{0}
\newrgbcolor{dialinecolor}{1 1 1}
\psset{linecolor=dialinecolor}
\pspolygon(11.597000,11.686900)(11.597000,11.815490)(11.775880,11.815490)(11.775880,11.686900)
\psset{linewidth=0.04} \psset{linestyle=solid}
\psset{linestyle=solid} \setlinecaps{0}
\newrgbcolor{dialinecolor}{0 0 0}
\psset{linecolor=dialinecolor}
\psclip{\pswedge[linestyle=none,fillstyle=none](9.684967,12.691672){1.470013}{302.221585}{347.127527}}
\psellipse(9.684967,12.691672)(1.039456,1.039456)
\endpsclip
\psset{linewidth=0.04} \psset{linestyle=solid}
\psset{linestyle=solid} \setlinecaps{0}
\newrgbcolor{dialinecolor}{0 0 0}
\psset{linecolor=dialinecolor}
\psclip{\pswedge[linestyle=none,fillstyle=none](11.310680,12.194306){0.775267}{65.064321}{154.986558}}
\psellipse(11.310680,12.194306)(0.548196,0.548196)
\endpsclip
\psset{linewidth=0.04} \psset{linestyle=solid}
\psset{linestyle=solid} \setlinecaps{0}
\newrgbcolor{dialinecolor}{0 0 0}
\psset{linecolor=dialinecolor}
\psclip{\pswedge[linestyle=none,fillstyle=none](11.365040,12.185267){1.070003}{80.593325}{159.252127}}
\psellipse(11.365040,12.185267)(0.756607,0.756607)
\endpsclip
\psset{linewidth=0.04} \psset{linestyle=solid}
\psset{linestyle=solid} \setlinecaps{0}
\newrgbcolor{dialinecolor}{0 0 0}
\psset{linecolor=dialinecolor}
\psclip{\pswedge[linestyle=none,fillstyle=none](11.631587,11.799286){0.913500}{40.454248}{175.897303}}
\psellipse(11.631587,11.799286)(0.645942,0.645942)
\endpsclip
\psset{linewidth=0.04} \psset{linestyle=solid}
\psset{linestyle=solid} \setlinecaps{0}
\newrgbcolor{dialinecolor}{0 0 0}
\psset{linecolor=dialinecolor}
\psclip{\pswedge[linestyle=none,fillstyle=none](11.627687,11.805151){1.123182}{43.864232}{154.625176}}
\psellipse(11.627687,11.805151)(0.794210,0.794210)
\endpsclip
\psset{linewidth=0.04} \psset{linestyle=solid}
\psset{linestyle=solid} \setlinecaps{0}
\newrgbcolor{dialinecolor}{0 0 0}
\psset{linecolor=dialinecolor}
\psline(12.183100,12.364100)(12.473100,11.899400)
\newrgbcolor{dialinecolor}{1 1 1}
\psset{linecolor=dialinecolor}
\pspolygon*(10.199700,11.691200)(10.199700,11.788344)(10.352870,11.788344)(10.352870,11.691200)
\psset{linewidth=0.04} \psset{linestyle=solid}
\psset{linestyle=solid} \setlinejoinmode{0}
\newrgbcolor{dialinecolor}{1 1 1}
\psset{linecolor=dialinecolor}
\pspolygon(10.199700,11.691200)(10.199700,11.788344)(10.352870,11.788344)(10.352870,11.691200)
\psset{linewidth=0.04} \psset{linestyle=solid}
\psset{linestyle=solid} \setlinecaps{0}
\newrgbcolor{dialinecolor}{0 0 0}
\psset{linecolor=dialinecolor}
\psclip{\pswedge[linestyle=none,fillstyle=none](9.807164,12.601982){1.499111}{296.497099}{351.193837}}
\psellipse(9.807164,12.601982)(1.060031,1.060031)
\endpsclip
\psset{linewidth=0.01} \psset{linestyle=dashed,dash=1 1}
\psset{linestyle=dashed,dash=0.1 0.1} \setlinecaps{0}
\newrgbcolor{dialinecolor}{0 0 0}
\psset{linecolor=dialinecolor}
\psclip{\pswedge[linestyle=none,fillstyle=none](9.691332,12.686362){1.571264}{292.421163}{350.239229}}
\psellipse(9.691332,12.686362)(1.111052,1.111052)
\endpsclip
\newrgbcolor{dialinecolor}{1 1 1}
\psset{linecolor=dialinecolor}
\pspolygon*(10.872700,11.592600)(10.872700,11.815490)(11.013010,11.815490)(11.013010,11.592600)
\psset{linewidth=0.04} \psset{linestyle=solid}
\psset{linestyle=solid} \setlinejoinmode{0}
\newrgbcolor{dialinecolor}{1 1 1}
\psset{linecolor=dialinecolor}
\pspolygon(10.872700,11.592600)(10.872700,11.815490)(11.013010,11.815490)(11.013010,11.592600)
\newrgbcolor{dialinecolor}{1 1 1}
\psset{linecolor=dialinecolor}
\pspolygon*(12.333100,11.776900)(12.333100,11.875485)(12.538840,11.875485)(12.538840,11.776900)
\psset{linewidth=0.04} \psset{linestyle=solid}
\psset{linestyle=solid} \setlinejoinmode{0}
\newrgbcolor{dialinecolor}{1 1 1}
\psset{linecolor=dialinecolor}
\pspolygon(12.333100,11.776900)(12.333100,11.875485)(12.538840,11.875485)(12.538840,11.776900)
\psset{linewidth=0.04} \psset{linestyle=solid}
\psset{linestyle=solid} \setlinecaps{0}
\newrgbcolor{dialinecolor}{0 0 0}
\psset{linecolor=dialinecolor}
\psline(11.242628,11.471600)(11.240100,11.472600)
\psset{linewidth=0.04} \psset{linestyle=solid} \setlinejoinmode{0}
\setlinecaps{0}
\newrgbcolor{dialinecolor}{0 0 0}
\psset{linecolor=dialinecolor}
\psline(11.172318,11.580067)(11.284214,11.455150)(11.117143,11.440583)
\psset{linewidth=0.04} \psset{linestyle=solid}
\psset{linestyle=solid} \setlinecaps{0}
\newrgbcolor{dialinecolor}{0 0 0}
\psset{linecolor=dialinecolor}
\psclip{\pswedge[linestyle=none,fillstyle=none](11.378316,11.894027){0.755131}{150.322393}{217.199540}}
\psellipse(11.378316,11.894027)(0.533958,0.533958)
\endpsclip
\psset{linewidth=0.01} \psset{linestyle=dashed,dash=1 1}
\psset{linestyle=dashed,dash=0.1 0.1} \setlinecaps{0}
\newrgbcolor{dialinecolor}{0 0 0}
\psset{linecolor=dialinecolor}
\psclip{\pswedge[linestyle=none,fillstyle=none](11.626967,11.812564){1.009097}{42.825007}{177.008637}}
\psellipse(11.626967,11.812564)(0.713539,0.713539)
\endpsclip
\psset{linewidth=0.01} \psset{linestyle=dashed,dash=0.1 0.1}
\psset{linestyle=dashed,dash=0.1 0.1} \setlinecaps{0}
\newrgbcolor{dialinecolor}{0 0 0}
\psset{linecolor=dialinecolor}
\psclip{\pswedge[linestyle=none,fillstyle=none](11.153119,11.851846){0.331860}{182.600160}{257.667807}}
\psellipse(11.153119,11.851846)(0.234660,0.234660)
\endpsclip
\psset{linewidth=0.04} \psset{linestyle=solid}
\psset{linestyle=solid} \setlinecaps{0}
\newrgbcolor{dialinecolor}{0 0 0}
\psset{linecolor=dialinecolor}
\psline(12.110300,12.222700)(12.397400,11.751200)
\newrgbcolor{dialinecolor}{1 1 1}
\psset{linecolor=dialinecolor}
\pspolygon*(11.784500,12.246500)(11.784500,12.629840)(12.063110,12.629840)(12.063110,12.246500)
\psset{linewidth=0.04} \psset{linestyle=solid}
\psset{linestyle=solid} \setlinejoinmode{0}
\newrgbcolor{dialinecolor}{1 1 1}
\psset{linecolor=dialinecolor}
\pspolygon(11.784500,12.246500)(11.784500,12.629840)(12.063110,12.629840)(12.063110,12.246500)
\newrgbcolor{dialinecolor}{1 1 1}
\psset{linecolor=dialinecolor}
\pspolygon*(11.625900,12.544100)(11.625900,12.655540)(11.775920,12.655540)(11.775920,12.544100)
\psset{linewidth=0.04} \psset{linestyle=solid}
\psset{linestyle=solid} \setlinejoinmode{0}
\newrgbcolor{dialinecolor}{1 1 1}
\psset{linecolor=dialinecolor}
\pspolygon(11.625900,12.544100)(11.625900,12.655540)(11.775920,12.655540)(11.775920,12.544100)
\psset{linewidth=0.04} \psset{linestyle=solid}
\psset{linestyle=solid} \setlinecaps{0}
\newrgbcolor{dialinecolor}{0 0 0}
\psset{linecolor=dialinecolor}
\psclip{\pswedge[linestyle=none,fillstyle=none](11.396517,12.237180){0.672267}{304.431817}{74.030749}}
\psellipse(11.396517,12.237180)(0.475364,0.475364)
\endpsclip
\psset{linewidth=0.04} \psset{linestyle=solid}
\psset{linestyle=solid} \setlinecaps{0}
\newrgbcolor{dialinecolor}{0 0 0}
\psset{linecolor=dialinecolor}
\psclip{\pswedge[linestyle=none,fillstyle=none](11.348284,12.242246){0.998145}{301.923616}{79.401310}}
\psellipse(11.348284,12.242246)(0.705795,0.705795)
\endpsclip
\psset{linewidth=0.01} \psset{linestyle=dashed,dash=1 1}
\psset{linestyle=dashed,dash=0.1 0.1} \setlinecaps{0}
\newrgbcolor{dialinecolor}{0 0 0}
\psset{linecolor=dialinecolor}
\psline(3.849680,11.522900)(4.525580,11.542600)
\psset{linewidth=0.04} \psset{linestyle=solid}
\psset{linestyle=solid} \setlinecaps{0}
\newrgbcolor{dialinecolor}{0 0 0}
\psset{linecolor=dialinecolor}
\psline(3.016530,11.443200)(3.255300,11.442600)
\psset{linewidth=0.04} \psset{linestyle=solid}
\psset{linestyle=solid} \setlinecaps{0}
\newrgbcolor{dialinecolor}{0 0 0}
\psset{linecolor=dialinecolor}
\psline(3.012250,11.610400)(3.263870,11.609700)
\newrgbcolor{dialinecolor}{1 1 1}
\psset{linecolor=dialinecolor}
\pspolygon*(2.934970,11.475000)(2.934970,11.571158)(3.307681,11.571158)(3.307681,11.475000)
\psset{linewidth=0.04} \psset{linestyle=solid}
\psset{linestyle=solid} \setlinejoinmode{0}
\newrgbcolor{dialinecolor}{1 1 1}
\psset{linecolor=dialinecolor}
\pspolygon(2.934970,11.475000)(2.934970,11.571158)(3.307681,11.571158)(3.307681,11.475000)
\psset{linewidth=0.01} \psset{linestyle=dashed,dash=1 1}
\psset{linestyle=dashed,dash=0.1 0.1} \setlinecaps{0}
\newrgbcolor{dialinecolor}{0 0 0}
\psset{linecolor=dialinecolor}
\psline(2.864620,11.527695)(3.408490,11.522895)
\newrgbcolor{dialinecolor}{1 1 1}
\psset{linecolor=dialinecolor}
\pspolygon*(10.976400,12.575700)(10.976400,12.983460)(11.226400,12.983460)(11.226400,12.575700)
\psset{linewidth=0.04} \psset{linestyle=solid}
\psset{linestyle=solid} \setlinejoinmode{0}
\newrgbcolor{dialinecolor}{1 1 1}
\psset{linecolor=dialinecolor}
\pspolygon(10.976400,12.575700)(10.976400,12.983460)(11.226400,12.983460)(11.226400,12.575700)
\psset{linewidth=0.01} \psset{linestyle=dashed,dash=1 1}
\psset{linestyle=dashed,dash=0.1 0.1} \setlinecaps{0}
\newrgbcolor{dialinecolor}{0 0 0}
\psset{linecolor=dialinecolor}
\psclip{\pswedge[linestyle=none,fillstyle=none](11.346238,12.232401){0.858498}{289.252340}{152.170247}}
\psellipse(11.346238,12.232401)(0.607050,0.607050)
\endpsclip
\psset{linewidth=0.04} \psset{linestyle=solid}
\psset{linestyle=solid} \setlinecaps{0}
\newrgbcolor{dialinecolor}{0 0 0}
\psset{linecolor=dialinecolor}
\psclip{\pswedge[linestyle=none,fillstyle=none](11.565130,12.335837){0.962329}{116.275004}{157.512885}}
\psellipse(11.565130,12.335837)(0.680469,0.680469)
\endpsclip
\psset{linewidth=0.04} \psset{linestyle=solid}
\psset{linestyle=solid} \setlinecaps{0}
\newrgbcolor{dialinecolor}{0 0 0}
\psset{linecolor=dialinecolor}
\psclip{\pswedge[linestyle=none,fillstyle=none](11.042936,12.595817){0.360501}{31.234413}{114.704110}}
\psellipse(11.042936,12.595817)(0.254913,0.254913)
\endpsclip
\psset{linewidth=0.04} \psset{linestyle=solid}
\psset{linestyle=solid} \setlinecaps{0}
\newrgbcolor{dialinecolor}{0 0 0}
\psset{linecolor=dialinecolor}
\psline(11.529602,11.955953)(11.549900,11.944500)
\psset{linewidth=0.04} \psset{linestyle=solid} \setlinejoinmode{0}
\setlinecaps{0}
\newrgbcolor{dialinecolor}{0 0 0}
\psset{linecolor=dialinecolor}
\psline(11.474767,12.073008)(11.568551,11.933976)(11.401056,11.942368)
}\endpspicture \] Consider the situation of a composition $s \circ
s'$, where the orientation of an output circle of $s'$ is opposite to
that of the corresponding input circle of $s$. When applying the
definition of the composition from section \ref{sec1}, we need to
follow the orientation of the output circle and identify it with that
of the corresponding input circle. In the above picture, this is
achieved by flipping an input circle, which has introduced twists in
some of the chords. In general, while sewing input to output, one
needs to give chords an extra twist, if the orientations do not
match.

\begin{thm} \label{theorem2}
Let $A$ be a finite dimensional, unital, associative, and
commutative algebra endowed with a non-degenerate and invariant
inner product. Then, the normalized Hochschild cochain complex of
$A$ is an algebra over the PROP, $C_\ast \mathcal{S}$, of Sullivan
chord diagrams.
\end{thm}

\begin{proof}
The description of the map $\overline{\alpha}:C_*\mathcal S \to
\mathcal End_{\overline{ HC^*}(A;A)}$ is identical to that of the
previously established action $\alpha:C_*\mathcal S^c\to \mathcal
End_{\overline{HC^*}(A;A)}$ in the associative case. It remains to
show that $\overline{\alpha}$ is a map of operads. One can see
that $\overline{\alpha}$ resects composition for the same reasons
$\alpha$ did. In fact, commutativity of $A$ does not play a role
in this. Commutativity of $A$, however, plays an important role in
showing that $\overline{\alpha}$ respect the differentials. Recall
that for $s\in C_*\mathcal S$, the formula \eqref{dif(a(s))}
describes the differential $D(\alpha(s))$. There are two cases to
consider. In case there is at most one binding tensor factor of
$\Delta(c^i_j)=\sum_{ (c^i_j)}(c^i_j)'\otimes (c^i_j)''$, the
terms in \eqref{dif(a(s))} cancel each other. This is due to the
new feature of commutativity of $A$, as seen in observation (i) on
page \pageref{observat-i}. In case both tensor factors of
$\Delta(c^i_j)$ are bound, we obtain the terms which correspond to
$\alpha(\partial(s))$.
\end{proof}

Note that the chord associated to the orientation-reversal $\sim$
squares to the identity, $(\sim)^2= id$. Therefore,
$\overline{HC^*}(A;A)$ decomposes into the eigenspaces
$\overline{HC^*}(A;A)_+ \oplus \overline{HC^*}(A;A)_{-}$, where
\begin{eqnarray*}
\overline{ HC^*}(A;A)_+&=& span\{c_1\otimes\dots\otimes c_n\otimes
a+(-1)^{\epsilon} c_n\otimes\dots \otimes c_1\otimes a\}\\%
\overline{ HC^*}(A;A)_-&=& span\{c_1\otimes\dots\otimes c_n\otimes
a-(-1)^{\epsilon} c_n\otimes\dots \otimes c_1\otimes a\}%
\end{eqnarray*}
where $c_i\in A^*, a\in A$, and $\epsilon=\frac{n(n+1)}{2}$. The
operator $\Delta$ maps each eigenspace into the other, {\it i.e.}
$\Delta\left(\overline{ HC^*} (A;A)_\pm\right)\subset\overline{
HC^*}(A;A)_\mp$. In the notation of Example \ref{brace-smile},
$\sim$ anticommutes with $\Delta$, {\it i.e.}
$\sim\circ\Delta=-\Delta \circ\sim$. Moreover, we have $\sim\circ
\smile= \smile \circ \tau_2\circ( \sim\otimes \sim)$ and
$(\sim\otimes\sim)\circ\tau_2\circ\vee_0=\vee_0\circ\sim$. Since
$\sim$ commutes with the Hochschild boundary operator, we have,

\begin{cor}
Under the assumptions of Theorem \ref{theorem2}, the Hochschild
cohomology of $A$ is a Frobenius algebra, which is endowed with a
compatible $BV$ operator, $\Delta$, and an involution $\sim$. The
operator $\Delta$ maps each eigenspace of $\sim$ into the other, {\it
i.e.} $\Delta\left(HH^* (A;A)_\pm\right)\subset HH^*(A;A)_\mp$, where
$HH^*(A;A)_{\pm}$ are the $\pm 1$ eigenspaces of $\sim$. The map
$\sim$ is both an anti-algebra and an anti-coalgebra map. That is to
say $\widetilde{f\smile g}=\widetilde{ g}\smile \widetilde{f}$, and
$\vee_0\left(\widetilde{f}\right)= \sum_{(f)} \widetilde{f}''\otimes
\widetilde{f'}$, where $\vee_0(f)= \sum_{(f)} f'\otimes f''$.
\end{cor}

\end{document}